%% file: BorceuxClementino.tex
\documentclass[12pt]{article}

\usepackage{amssymb}
\usepackage{amsmath}
\input{diagram}

\newcommand{\Set}{\mathsf{Set}}

\newcommand{\ZERO}{\mathbf{0}}

\newcommand{\id}{\mathsf{id}}
\newcommand{\ONE}{\mathbf{1}}
\newcommand{\Cat}{\mathsf{Cat}}

\newcommand{\DisPar}{\mathsf{DisPar}}
\newcommand{\Stab}{\mathsf{Stab}}
\newcommand{\DisParE}{\DisPar(\calE)}
\newcommand{\StabE}{\Stab(\calE)}

\newcommand{\proof}{\medbreak\noindent%
{\it Proof}\hspace{1em plus .5em}}

\newcommand{\qed}{\mbox{}\hfill\mbox{}%
\mbox{$\square$}\medbreak}

\newcommand{\calC}{{\mathcal C}}
\newcommand{\calT}{{\mathcal T}}
\newcommand{\calF}{{\mathcal F}}

\newcommand{\calP}{{\mathcal P}}
\newcommand{\calA}{{\mathcal A}}
\newcommand{\calZ}{{\mathcal Z}}
\newcommand{\calS}{{\mathcal S}}
\newcommand{\calX}{{\mathcal X}}
\newcommand{\calY}{{\mathcal Y}}
\newcommand{\calE}{{\mathcal E}}
\newcommand{\calB}{{\mathcal B}}

\newcommand{\BbbA}{\mathbb A}
\newcommand{\BbbS}{\mathbb S}
\newcommand{\BbbT}{\mathbb T}

\newtheorem{theorem}{Theorem}[section]
\newtheorem{definition}[theorem]{Definition}
\newtheorem{proposition}[theorem]{Proposition}
\newtheorem{lemma}[theorem]{Lemma}
\newtheorem{corollary}[theorem]{Corollary}
\newtheorem{example}[theorem]{Example}

\newtheorem{examples}[theorem]{Examples}

\let\Label=\label
\let\Bibitem=\bibitem

\title{On coherent systems of subobjects\\with application to torsion theory}
\author{Francis Borceux and Maria Manuel Clementino\footnote{This work was partially supported by the {\it Centro de
Matem\'atica da Universidade de Coimbra} (CMUC) -- UIDB/00324/2020, funded by the Portuguese Government
through FCT/MCTES.}}

\date{To Marta Bunge, a friendly colleague and an inspiring mathematician}

\begin{document}

\maketitle

 \begin{abstract}
In a coherent category, the posets of subobjects have very strong properties. We emphasize the validity of these properties, in general categories, for well-behaved classes of subobjects. As an example of application, we investigate the problem of the various torsion theories which can be universally associated with a pretorsion one.\\

\noindent {\it Keywords}: coherent system of subobjects, (pre)torsion theory, torsion objects, torsion free objects, stable category, category of fractions.

\noindent {\it Math. Subject Classification (2020)}: 18E40, 18B50, 18A20, 18E35.
\end{abstract}

\section{Introduction}\Label{1}

A coherent category (see \cite{PTJ}) has in particular distributive lattices of subobjects whose operations are preserved by pullbacks. In this paper, we want to draw attention to the case of arbitrary categories in which there exist distributive lattices of particular subobjects -- we call them {\em distinguished} -- whose operations are preserved by pullbacks. We call this a {\em coherent system of subobjects}.  A well-known example is that of complemented subobjects in a lextensive category.  But there are many other examples of interest.

 In the case of small categories, we call a subcategory $\calS$  of $\calA$  {\em saturated} when an arrow of $\calA$ lies in $\calS$ as soon as its domain or its codomain is in $\calS$. The saturated subcategories constitute a distributive lattice whose operations are preserved by pullbacks. This result generalizes to the case of internal categories in an arbitrary coherent category $\calC$, yielding an interesting bunch of examples of coherent systems of subobjects. This particularizes further to preordered objects in a coherent category.

We apply these considerations to the study of the stable category associated with a pretorsion theory. A {\em pretorsion theory} (see \cite{FFG}) in a category $\calC$ consists of giving two classes of objects, respectively called the {\em torsion objects} and the {\em torsion free objects}, together with adequate axioms. The objects which are both torsion and torsion free are called {\em trivial}. When all trivial objects are zero objects, the pretorsion theory  is called a {\em torsion theory}. The problem of the {\em stable category} associated with a torsion theory is that of constructing naturally a torsion theory from a given pretorsion theory (see \cite{BCG3}).

Of course, to construct a torsion theory from a pretorsion one, one must identify all trivial objects to a zero object.  In \cite{FF, BCG3}, this is done with the extra requirement that, when $A=S\amalg T$ is the disjoint union of two complemented subobjects, that union is preserved; in particular, if $T$ is trivial,  $A$ is identified with $S$ in the stable category. In our paper, we consider the more flexible requirement that when $A=S\cup T$ is the union of two distinguished subobjects with $T$ trivial, then $T$ is identified with $\ZERO$ and $A$ is identified with $S$.

We prove a general theorem on the existence and the universality of the stable category, which applies in particular to the case of preordered objects in an exact coherent category (see \cite{FF, BCG1}), and to the case of internal categories in a Grothendieck topos (see \cite{BCGT, FC}). But the flexibility of our approach allows also choosing only trivial subobjects (i.e. only $\ZERO$ and the object itself) as distinguished subobjects. In that case, our theorem yields the torsion theory universally associated with the given pretorsion theory. And, under a very mild assumption, we prove that the stable category can then be obtained as the category of fractions which inverts all morphisms $\ZERO \ar A$, with $A$ a trivial object.

\section{Coherent systems of subobjects}\Label{2}

The structure that we want to promote in this paper is the following one\footnote{Since this cannot hurt in this paper, we choose to use freely the common abuse of language which does not distinguish subobjects and monomorphisms.}.

\begin{definition}\Label{dist}
Let $\calC$ be a category with a strict initial object $\ZERO$. By a {\em coherent system of subobjects} in $\calC$ is meant, for each object $C\in\calC$, the choice of a class of so-called {\em distinguished   subobjects of $C$}, in such a way that:
\begin{description}
\item[(CS1)] $\ZERO$ and $A$ are distinguished in $A$;
\item[(CS2)] the union of two distinguished subobjects exists and is distinguished;
\item[(CS3)] the pullback of a distinguished subobject along an arbitrary morphism exists and is distinguished;
\item[(CS4)] pulling back distinguished subobjects preserves their union;
\item[(CS5)] given two monomorphisms
$$
A\mono B \mono C
$$
if $A$ is distinguished in $B$ and $B$ is distinguished in $C$, then $A$ is distinguished in $C$.
\end{description}
A coherent system of subobjects is called {\em effective} when moreover (see \cite{BARR2})
\begin{description}
\item[(CS6)] the union of distinguished objects is effective, that is, given two distinguished subobjects $S$ and $T$ of $A$, the following pullback is also a pushout:
\begin{diagram}[80]
S\cap T | \emono | S ||
\smono | | \smono   ||
T | \emono | S\cup T ||
\end{diagram}
\end{description}
\end{definition}

Let us recall that, given a strict initial object $\ZERO$, every morphism with domain $\ZERO$ is a monomorphism. Indeed if there exists a morphism $X\ar\ZERO$, that morphism is necessarily unique since it must be an isomorphism, with the unique morphism $\ZERO\ar X$ as inverse.

\begin{proposition}\Label{propdis}
Let $\calC$ be a category with a strict initial object, provided with a coherent system of subobjects. The following properties hold:
\begin{description}
\item[(CS7)] the intersection of two distinguished subobjects exists and is distinguished;
\item[(CS8)] given three distinguished subobjects $R$, $S$, $T$ of $A$
$$
R\cap(S\cup T)=(R\cap S)\cup(R\cap T);
$$
\item[(CS9)] if $S$ is a distinguished subobject of $A$, $S$ is distinguished in every subobject $S\subseteq T \subseteq A$.
\end{description}
\end{proposition}

\proof
The case of the intersection follows from (CS3) and (CS5). The distributivity law is a special instance of (CS4). The last assertion is obtained via (CS3) when pulling back $S\subseteq A$ along $T\subseteq A$.
\qed

We shall also meet the following related notion, which is reminiscent of the notion of extremal epimorphism:

\begin{definition}\Label{epidis}
Let $\calC$ be a category with a strict initial object, provided with a coherent system of subobjects. An epimorphism $f\colon A \ar B$ is called {\em distinguished} when it does not factor through any proper  distinguished subobject of $B$.
\end{definition}

\begin{proposition}\Label{epiort}
Let $\calC$ be a category with a strict initial object, provided with a coherent system of subobjects. An epimorphism is distinguished if and only if it is left orthogonal to every distinguished monomorphism.
\end{proposition}

\proof
The classical proof for extremal versus strong epimorphism applies as such.
\qed

\begin{proposition}\Label{propepi}
Let $\calC$ be a category with a strict initial object, provided with a coherent system of subobjects.
\begin{enumerate}
\item Distinguished epimorphisms are stable under composition.
\item If $gf$ is a distinguished epimorphism, so is $g$.
\end{enumerate}
\end{proposition}

\proof
Once more the standard proof for extremal epimorphisms applies as such.
\qed

\begin{examples}\Label{exdist}
Here is a first list of rather immediate examples of coherent systems of subobjects; the last three are effective.
\begin{enumerate}
\item All subobjects in a coherent category.
\item The complemented subobjects in a lextensive category (see \cite{CLW}).
\item The open (respectively closed, clopen) subspaces in the category of topological spaces.
\item Only $\ZERO$ and the object itself in a category with a strict initial object. We call this the {\em indiscrete} coherent system of subobjects.
\hspace*{1cm}
\end{enumerate}
\end{examples}

\proof
In the lextensive case, given two complemented subobjects $S$, $T$ of $A$, their union is the disjoint coproduct
$$
(S\cap\complement T)\amalg(S\cap T)\amalg(\complement S\cap T).
$$
The rest is obvious.
\qed

Despite its triviality, the last of these examples will play a significant role in Section~\ref{7}. But we want to focus now on the following example, which suggested to us the notion of coherent system of subobjects.

\begin{example}\Label{preord}
A sub-preordered set $(S,\leq)$ of $(A,\leq)$ with the induced preorder is called {\em open} (see \cite{FF}) when
$$
a\leq b \mbox{ and } b\in S \implies a\in S
$$
and {\em closed} when
$$
a\leq b \mbox{ and } a\in S \implies b\in S.
$$
Clopen means, as usual, open and closed.
The open (respectively closed, clopen) sub-preordered sets in the sense of \cite{FF} constitute effective coherent systems of subobjects.
\end{example}

\proof
The proof of having effective coherent systems of subobjects is a special case of our next example. Indeed, a preordered set is a small category with at most one arrow between every two objects.
\qed

We generalize now the notions of open, closed and clopen subobjects to the case of subcategories. To avoid any ambiguity of terminology in the closed case, we prefer to use {\em left saturated}, {\em right saturated} and {\em saturated}  instead of open, closed and clopen.

\begin{example}\Label{cat}
A subcategory $\calS\subseteq\calC$ is {\em left saturated} when
$$
f\colon A \to B \mbox{ and } B\in \calS \implies f\in\calS
$$
and right saturated when
$$
f\colon A \to B \mbox{ and } A\in \calS \implies f\in\calS.
$$
$\calS$ is {\em saturated} in $\calC$ when it is both left and right saturated.
Left saturated, right saturated and saturated subcategories constitute effective coherent systems of subobjects in $\Cat$.
\end{example}

\proof
Only the case of the union requires a comment.
To construct the union of two subcategories  $\calS$ and $\calT$ of $\calC$, one considers first the graph constructed on the set theoretical union of both sets of objects, and the set theoretical union of both sets of arrows. One adds further all the possible composites in $\calC$ of chains of consecutive arrows in this last union. When both $\calS$ and $\calT$ are left saturated, the last arrow of a chain is thus in $\calS$ or in $\calT$, and therefore by left saturation, so do all the arrows of the chain and thus also  their composite. Analogously for right saturation, starting with the first arrow of the chain. Thus the set of arrows of $\calS\cup\calT$ is the set theoretical union of the sets of arrows of $\calS$ and $\calT$. Since $\Set$ is coherent with effective unions, this forces at once all the properties for having effective coherent systems of subobjects.
\qed

In particular, trivially:

\begin{proposition}\Label{full}
A left saturated, right saturated or saturated subcategory is full.
\qed
\end{proposition}

\begin{definition}\Label{intlex}
In a category $\calC$ with finite limits, a left saturated internal subcategory $\BbbS\subseteq \BbbA$ is one which satisfies the axiom
$$
\mbox{If }\Bigl[f\in A_1 \wedge d_1(f)\in  S_0\Bigr] \mbox{~then~} \Bigl[f\in S_1\Bigr]
$$
in the cartesian internal logic of $\calC$.
The right saturated case is obtained when using instead $d_0$.
\end{definition}

We shall use the notation
$$
A_2 \Ar m A_1 \Bktriadjar{d_0}{n}{d_1} A_0
$$
to indicate an internal category $\BbbA$. The internal category is an internal preordered object when the pair $(d_0,d_1)$ turns $A_1$ in a subobject of $A_0\times A_0$.

Given an internal subcategory $\BbbS$ of $\BbbA$, the left saturated notion  of Example~\ref{preord} (which is a cartesian notion; see \cite{PTJ}) translates as the existence of a (unique) factorization $\delta_1$ in the following pullback diagram
\begin{diagram}[80]
| | A_1\times_{A_0}S_0 | \ear | S_0 ||
| \Swdotar{\delta_1} | \smono | \mbox{p.b.} | \smono ||
S_1 | \emono | A_1 | \eaR{d_1} | A_0 ||
\end{diagram}
An analogous conclusion holds, using instead $d_0$ and a corresponding factorization $\delta_0$, in the right saturated case.

In the situation  of Definition \ref{intlex}, one could of course be tempted to define a saturated internal subcategory as one which is both left and right saturated. But in the case of a coherent category $\calC$ (see \cite{PTJ}), thus in particular in a topos, it sounds more sensible to define:

\begin{definition}\Label{satcoh}
In a coherent category $\calC$, a saturated internal subcategory $\BbbS\subseteq \BbbA$ is one which satisfies the axiom
$$
\mbox{If }\Bigl[f\in A_1 \wedge \bigl (d_0(f)\in S_0 \vee d_1(f)\in S_0\bigr)\Bigr]
\mbox{~then~} \Bigl[f\in S_1\Bigr]
$$
in the internal coherent logic of $\calC$.
\end{definition}

This can be translated as the existence of a (unique) factorization $\delta$  in the following diagram
\begin{diagram}[80]
| | A'_1 | \ear | \movevertexright{(S_0\times A_0)\cup(A_0\times S_0)}| |  ||
| \Swdotar{\delta} | \smono | \mbox{p.b.} | \smono | | ||
S_1 | \emono | A_1 | \eaR{(d_0,d_1)} | \movevertexright{A_0\times A_0} | | ||
\end{diagram}
For example in the topos of sheaves on a locale $L$, this means that if $f\in A_1(u)$ while $u\in L$ is covered by all the levels $v\leq u$ where the domain or the codomain of the restriction of $f$ is in $S_0(v)$, then $f\in S_1(u)$.

\begin{lemma}\Label{equivsat}
When $\calC$ is a coherent category, given an internal subcategory $\BbbS\subseteq\BbbA$, if $\BbbS$ is saturated in $\BbbA$, it is both left and right saturated. The converse holds when $\calC$ has effective unions.
\end{lemma}

\proof
With the notation of Definitions \ref{intlex} and  \ref{satcoh}, compose $(d_0,d_1)$ with both projections of the product. This yields first
$$
p_1^{-1}(S_0)=S_0\times A_0,~~~
p_2^{-1}(S_0)=A_0\times S_0.
$$
By axiom (CS4) this implies
$$
A'_1=(A_1\times_{A_0}S_1)\cup(S_1\times_{A_0}A_1).
$$
Since $S_1\mono A_1$ is a monomorphism, the existence of $\delta$ implies that of $\delta_0$ and $\delta_1$. The converse holds when $\calC$ has effective unions.
\qed

\begin{example}\Label{int}
Let $\calC$ be a coherent category.  The notions of left saturated, right saturated and saturated  internal subcategories yield coherent systems of subobjects in $\Cat(\calC)$. Restricting one's attention to internal preordered objects, one obtains corresponding coherent systems of subobjects in the category of preordered objects in $\calC$.
\end{example}

\proof
Given two saturated (respectively, left saturated, right saturated)  internal subcategories $\BbbS$, $\BbbT$ of $\BbbA$, let us  first observe that their union as internal subcategories admits $S_1\cup T_1$ in $\calC$ as object of arrows, and of course $S_0\cup T_0$ as object of objects. Indeed the pullback defining the object of composable pairs of morphisms in $S_1\cup T_1$
$$
(S_1\cup T_1)\times_{A_0}(S_1\cup T_1)
$$
can be split in four pieces by coherence of $\calC$. The two  pieces $S_1\times_{A_0}S_1$ and $T_1\times_{A_0}T_1$  yield composites lying respectively in $S_1$ and $T_1$ by the category axioms. The other two pieces $S_1\times_{A_0}T_1$ and $T_1\times_{A_0}S_1$ yield composites lying in $S_1\cup T_1$ in the saturated case, or already in $S_1$ or $T_1$ in the left or right saturated case.
The rest follows at once from the coherence of the category $\calC$.
\qed

It is immediate that, in $\Cat$, every saturated subcategory $\calS$ is complemented: its complement is the full subcategory generated by those objects which are not connected by a chain of arrows to any object of $\calS$. This is not the case for left or right saturated subcategories: here is an example of a right saturated subcategory which is not complemented.
$$
\{A \rightarrow B\} \mono \{A \rightarrow B \leftarrow C\}
$$

In an arbitrary category with finite limits, even in a topos, it is no longer the case that saturated subcategories are complemented. For example in the topos of sheaves on the Sierpinski space, consider the following preordered sheaf, where the restriction applies $a_i$ on $a$ and $b_i$ on $b$:
$$
\{a_1<b_1,a_2<b_2\}\ar\{a<b\}
$$
The two subsheaves
$$
\{a_i<b_i\}\ar\{a<b\}
$$
are saturated (or clopen) in the sense of Example \ref{int}, but are not complemented. Their union is of course the whole sheaf.
This shows that our notion of saturated (clopen) sub-preordered object differs from that in \cite{BCG1}, where complementarity of the subobject is forced. Of course, both notions coincide in the $\Set$ case.

\section{Variations on torsion theories}\Label{3}

Let us recall that an ideal $\calZ$ in a category $\calC$ is a class of arrows such that,
for every arrow $f\in\calZ$, one has $fu\in\calZ$ and
$vf\in\calZ$, for all arrows $u$, $v$ composable with $f$ (see \cite{Ehresmann}). When $\calC$ has a zero object, the zero morphisms constitute an ideal.

Given an ideal $\calZ$ in a category $\calC$, an arrow $k$ is the $\calZ$-kernel of an arrow $f$ when $fk\in\calZ$ and, if $fm\in\calZ$ for some arrow $m$, then $m$ factors uniquely through $k$. The uniqueness condition forces $k$ to be a monomorphism. When $\calZ$ is the ideal of zero morphisms, we recapture the usual notion of kernel. There is of course a dual notion of $\calZ$-cokernel. A pair of composable morphisms
$$
K\Ar k A \Ar q Q
$$
is a short $\calZ$-exact sequence when $k$ is the $\calZ$-kernel of $q$ and $q$ is the $\calZ$-cokernel of $k$.

 The following definition was introduced in \cite{FF} and then thoroughly investigated in \cite{FFG}:

\begin{definition}\Label{pretor}
A pretorsion theory in a category $\calC$ consists of a pair $(\calT,\calF)$ of classes of objects, both of them closed under isomorphisms, whose elements are called  the {\em torsion} and the {\em torsion-free} objects of the pretorsion theory, respectively. The objects in $\calT\cap\calF$ are called {\em trivial}, and the ideal $\calZ$ of {\em trivial morphisms} is that of those arrows factoring through a trivial object.\\
These data must satisfy the following two axioms:
\begin{description}
\item[(PT1)] every arrow $f\colon A\ar B$ with $A\in\calT$ and $B\in\calF$ is trivial;
\item[(PT2)] for every object $A\in\calC$, there exists a short $\calZ$-exact sequence
$$
\tau(A)\Ar{\varepsilon_A} A \Ar{\eta_A} \phi(A)
$$
with $\tau(A)\in\calT$ and $\phi(A)\in\calF$.
\end{description}
When the trivial objects are zero objects, the pretorsion theory is called a {\em torsion} theory.\footnote{The terminology {\em pretorsion} is thus somehow unfortunate, since the axioms are the same as for a torsion theory, but with respect to a more flexible choice of ideal.}\\
A {\em torsion functor} is a functor, between two categories provided with a pretorsion theory, which respects torsion objects, torsion free objects and the canonical $\calZ$-exact sequences of axiom (PT2).
\end{definition}

\begin{proposition}\Label{tauphi}
In the conditions of Definition \ref{pretor}, the objects $\tau(A)$ and $\phi(A)$ are defined uniquely up to an isomorphism. This extends in functors $\tau$ and $\phi$ which present respectively the full category of torsion objects as a coreflexive category of $\calC$, and the full subcategory of torsion free objects as a reflexive subcategory of $\calC$.
\end{proposition}

\proof
See \cite{FFG}.
\qed

It is trivial to observe that $\tau$ and $\phi$ extend as functors, which turn  respectively  the full subcategories of torsion/torsion free objects in a coreflective/reflective full subcategory of $\calC$ (see \cite{FFG}).

\begin{example}\Label{preordex}
One gets a pretorsion theory on the category of preordered sets when choosing (see \cite{FF})
\begin{itemize}
\item the equivalence relations as torsion objects;
\item the partial orders as torsion free objects.
\end{itemize}
These two notions are cartesian ones (see \cite{PTJ}), thus make sense in every category with finite limits. They yield a pretorsion theory in every Barr exact category (see \cite{BARR, BCG1}).
\end{example}

\proof
Given a preordered set $(A,R)$ and writing $R^\circ$ for the opposite relation, $\tau(A,R)=(A,R\cap R^\circ)$ while $\phi(A,R)$ is the quotient of $A$ by the equivalence relation $R\cap R^\circ$, this quotient being provided with the image-preorder of $R$.
\qed

\begin{example}\Label{catex}
One gets a pretorsion theory on the category of small categories when choosing (see \cite{BCGT})
\begin{itemize}
\item the groupoids as torsion objects (every arrow is an isomorphism);
\item the skeletal categories as torsion free objects (every isomorphism is an automorphism).
\end{itemize}
These two notions are cartesian ones (see \cite{PTJ}), thus make sense in every category with finite limits. They yield  a pretorsion theory in every Grothendieck topos (see \cite{FC}).
\end{example}

\proof
Given a small category $\calC$, $\tau(\calC)$ is the groupoid of isomorphisms of $\calC$ while $\phi(\calC)$ is the quotient of $\calC$ which identifies the domain and the codomain of every isomorphism (see \cite{BCGT, FC}).
\qed

More examples of pretorsion theories can be found in \cite{FFG, XA}.

The notion of {\em stable category} for a pretorsion theory has been introduced in \cite{FF} and further investigated in \cite{BCG1, BCG2, BCG3}. The question is, given a category $\calX$ provided with a pretorsion theory, to construct  ``in the best possible way'' a category $\calY$ provided with a torsion theory, and a morphism of pretorsion theories from $\calX$ to $\calY$. The first way to interpret ``in the best possible way'' is to look for the universal solution to the problem: this is what we do in Section~\ref{7}. The idea is of course to identify all the trivial objects to $\ZERO$ \ldots\ and check when this yields a solution to the problem.

But one can be interested in a more involved problem. For example, when an object $A$ can be written as a union $A=S\cup T$ of two subobjects, with $T$ a trivial object, one could want $A$ to be identified with $S$ in the stable category \ldots\  since $T$ will be identified with $\ZERO$. This yields of course a different universal problem, in which some compatibility with the union of subobjects is requested. In \cite{FF, BCG1, BCG3}, a universal solution of that type is produced in the special case where $S$ and $T$ are complemented subobjects and the union is a disjoint one.

In this paper, we shall give evidence that complemented subobjects do not play any canonical role: a stable category compatible with the union of distinguished subobjects can be constructed with respect to an arbitrary coherent system of subobjects, provided it is sufficiently compatible with the pretorsion theory.

\section{The compatibility conditions}\Label{4}

The various constructions of a stable category as in \cite{FF, BCG1, BCG2} underline clearly which properties are essential to get the expected result. In this section, we focus on these properties and show that they hold in our main cases of interest.  We use freely the notation of Section~\ref{3}, without recalling it.

\begin{definition}\Label{comp}
Let $\calC$ be a category with a strict initial object, provided with both a pretorsion theory and a coherent system of subobjects. These data are called {\em compatible} when the following conditions are satisfied.
\begin{description}
\item[(CC1)] The class of trivial objects is closed under distinguished subobjects.
\item[(CC2)] Given a morphism $f\colon S\cup T\ar B$, with $S$, $T$ two distinguished subobjects of $A$,  $f$ is trivial as soon as its restrictions on $S$ and $T$ are trivial.
\item[(CC3)] The functor  $\phi$ preserves distinguished subobjects.
\item[(CC4)] Given a pullback square
\begin{diagram}[80]
P |\Ear v | S||
\Smono u | | \smonO s ||
A | \eepI{\eta_A} | \phi(A) ||
\end{diagram}
with $s$ a distinguished monomorphism, $v$ is a distinguished epimorphism.
\item[(CC5)] Let $s$ be a distinguished monomorphism such that $\phi(s)$ is an isomorphism; then $s$ is an isomorphism.
\end{description}
\end{definition}

Let us infer some consequences of such a situation.

\begin{proposition}\Label{C1}
In the conditions of Definition~\ref{comp}, the initial object is trivial.
\end{proposition}

\proof
Trivial objects exist by axiom (PT2) for a pretorsion theory. We conclude by axiom (CS1) for a coherent system of subobjects, and the compatibility condition (CC1).
\qed

\begin{proposition}\Label{C7}
In the conditions of Definition~\ref{comp}, if $sf$ is trivial with $s\colon S\mono A$ a distinguished subobject, then $f$ is trivial.
\end{proposition}

\proof
The composite $sf$ factors as $gh$ through a trivial object $D$. Therefore $f$ factors through the distinguished subobject $g^{-1}(S)$ of $D$, which is trivial by condition (CC1).
\qed

\begin{proposition}\Label{exCC3}
In the conditions of Definition~\ref{comp}, the class of torsion free objects is closed under distinguished subobjects.
\end{proposition}

\proof
Consider the following diagram, where $A$ is torsion free and $S$ is distinguished in $A$.
\begin{diagram}[80]
\tau(S) | \Emono{\varepsilon_S} | S | \Eepi{\eta_S} | \phi(S) ||
\Sar{\tau_S} | | \Smono s | | \saR{\phi(s)} ||
\tau(A) | \emonO{\varepsilon_A} | A | \eeqL{\eta_A} | \phi(A) ||
\end{diagram}
The composite $s\varepsilon_S=\phi(s)\eta_S\varepsilon_S$ is trivial, thus $\varepsilon_S$ is trivial because $s$ is a distinguished subobject. Therefore $\id_S\varepsilon_S$ is trivial, proving that $\id_S$ factors through the $\calZ$-cokernel $\eta_S$. So the epimorphism $\eta_S$ is a section, thus an isomorphism.
\qed

\begin{proposition}\Label{epieta}
In the conditions of Definition~\ref{comp}, the canonical morphisms $\eta_A\colon A \epi \phi(A)$ are distinguished epimorphisms.
\end{proposition}

\proof
Simply choose $s=\id_{\phi(A)}$ in condition (CC4).
\qed

\begin{proposition}\Label{exCC4}
In the conditions of Definition~\ref{comp}, given the pullback diagram of condition (CC4), one has $S\cong\phi(P)$ and $v\cong\eta_P$.
\end{proposition}

\proof
Consider the following diagram, where $S$ is isomorphic to $\phi(S)$ by Proposition~\ref{exCC3}
\begin{diagram}[80]
| | \phi(P) ||
|  \Neepi{\eta_P} |  \Sar{\phi(v)} ||
P | \Eepi v |\crossarrows{S \cong\phi(S)}{\Escurvar{\movename(0,10){\phi(u)}}} |  ||
\Smono u | | \Smono s |   ||
A | \eepI{\eta_A} | \phi(A)   ||
\end{diagram}
By condition~(CC3), $\phi(s)$ and $\phi(u)$ are distinguished subobjects. By Proposition~\ref{propdis}, $\phi(v)$ is a distinguished subobject as well. By condition (CC4), $v$ is a distinguished epimorphism factoring through the distinguished subobject $\phi(v)$, thus $\phi(v)$ is an isomorphism. Therefore $\phi(P)\cong\phi(S)$ and $v\cong\eta_P$.
\qed

\begin{proposition}\Label{C5}
In the conditions of Definition~\ref{comp}, consider the following diagram
\begin{diagram}[80]
\tau(S)| \Emono{\varepsilon_{S}}|S|\Eepi{\eta_S}|\phi(S)||
\Sar{\tau(s)}| | \Smono{s}| | \saR{\phi(s)}||
\tau(A)| \Emono{\varepsilon_{A}}|A|\Eepi{\eta_{A}}|\phi(A)||
\end{diagram}
where the two rows are canonical $\calZ$-exact sequences and $s$ is a distinguished monomorphism. Then $\tau(s)$ and $\phi(s)$ are distinguished monomorphisms and both squares are pullbacks.
\end{proposition}

\proof
Let us begin with the right hand square. We know already, by condition (CC3), that  $\phi(s)$ is a distinguished monomorphism. We consider the pullback $(P,u,v)$ of $\phi(s)$ along $\eta_A$ and the corresponding factorization $w$.
\begin{diagram}[80]
S | | \Eepi[125]{\eta_S} | | \phi(S) ||
| \Semono w | | \Neepi{v} ||
\Sar[130]{s} | | P | | \smonO[130]{\phi(s)} ||
|  \Swmono u ||
A | | \eepI[125]{\eta_A} | | \phi(A) ||
\end{diagram}
Thus $u$ is a distinguished subobject by (CS3) and $w$ as well, by Proposition~\ref{propdis}. By Proposition \ref{exCC4}, $\phi(P)\cong\phi(S)$ and  $v\cong\eta_P$. But this implies that $\phi(w)$ is isomorphic to the identity on $\phi(S)$. By condition~(CC5), $w$ is an isomorphism and the outer square is a pullback.

For the left hand square, if $su=\varepsilon_Av$, composing with $\eta_A$ yields a trivial morphism $\eta_A\varepsilon_Av=\phi(s)\eta_Su$. Since $\phi(s)$ is a distinguished subobject, it follows from Proposition~\ref{C7} that $\eta_S u$  is trivial. Thus $u$ factors uniquely through $\varepsilon_S$. That factorization is also a factorization of $v$ through $\tau(s)$, because $\varepsilon_A$ is a monomorphism. Thus the left hand square is a pullback. And since $s$ is a distinguished subobject, so is $\tau(s)$ by axiom (CS3).
\qed

\begin{corollary}\Label{taudissub}
In the conditions of Definition~\ref{comp}, the functor $\tau$ preserves
distinguished subobjects.
\qed
\end{corollary}

\begin{corollary}\Label{tausub}
In the conditions of definition \ref{comp}, the class of torsion objects is closed under distinguished subobjects.
\end{corollary}

\proof
In the diagram of Proposition \ref{C5}, if $A$ is a torsion object, $\varepsilon_A$ is an isomorphism and thus by pullback, $\varepsilon_S$ as well.
\qed

\begin{corollary}\Label{taudissub}
In the conditions of Definition~\ref{comp}, the functor $\tau$ preserves
the union and the intersection of distinguished subobjects.
\end{corollary}

\proof
Considering again the diagram of Proposition~\ref{C5}, $\tau$ acts by pullbacks along $\varepsilon_A$ at the level of distinguished subobjects of $A$. One concludes by axiom (CS4).
\qed

\begin{corollary}\Label{phiuni}
In the conditions of Definition~\ref{comp}, the functor $\phi$ preserves the union and the intersection of distinguished subobjects.
\end{corollary}

\proof
Let $R$ and $S$ be two distinguished subobjects of $A$. By condition (CC3),  $\phi(S)$ and $\phi(R)$ are distinguished subobjects of $\phi_A$, thus also their union $\phi(S)\cup\phi(R)$, by axiom (CS2). By Proposition~\ref{C5}, the inverse images of $\phi(S)$ and $\phi(R)$ along $\eta_A$ are $S$ and $R$. Thus by axiom~(CS4), the inverse image of $\phi(A)\cup\phi(B)$ along $\eta_A$ is $S\cup R$. By Proposition ~\ref{exCC4}, $\phi(S\cup R)\cong\phi(S)\cup\phi(R)$. An analogous argument holds for the intersection.
\qed

\begin{proposition}\Label{C6}
In the conditions of Definition~\ref{comp}, consider the following diagram
\begin{diagram}[80]
X| \Ear{k}|S|\Ear{{p}}|Y||
\Smono{x}| | \Smono{s}| | \smonO{y}||
\tau(A)| \Emono{\varepsilon_{A}}|A|\Eepi{\eta_{A}}|\phi(A)||
\end{diagram}
where the bottom row is a canonical $\calZ$-exact sequence, the vertical morphisms are distinguished monomorphisms and the two squares are pullbacks. Then the upper row is a canonical $\calZ$-exact sequence as well.
\end{proposition}

\proof
By Proposition \ref{exCC4}, we have at once $p\cong\eta_S$.
By Proposition~\ref{C5}, the pullback of $s$ along $\varepsilon_A$ is $\tau(S)$. Thus up to isomorphism, $X\cong\tau(S)$ and $x\cong\varepsilon_S$.
\qed

\begin{example}\Label{comptriv}
Let $\calC$ be a category with a strict initial object. Any pretorsion theory on $\calC$ where the initial object is trivial is compatible with the indiscrete coherent system of subobjects (see \ref{exdist}).
\end{example}

\proof
Trivial because $\ZERO$ is strict initial.
\qed

\begin{example}\Label{complext}
If $\calC$ is a lextensive category, a pretorsion theory on $\calC$ and the coherent system of complemented subobjects (see \ref{exdist}) are compatible provided that the class of trivial objects is closed under complemented subobjects and binary coproducts.
\end{example}

\proof
The assumptions take care of conditions (CC1) and (CC2).  Lemma~4 in \cite{BCG3} implies our Propositions \ref{exCC4}, \ref{C5} and \ref{C6}, from which conditions (CC3), (CC4) and (CC5).
\qed

\begin{example}\Label{comppreord}
In the category of preordered objects in an exact coherent category, the (equivalence/partial order) pretorsion theory and the coherent system of saturated subobjects are compatible.
\end{example}

\proof
When useful, we use the more convenient notation $\BbbA=(A,R)$ for a preordered object, and $\BbbS=(S,R_S)$ for a subobject of $(A,R)$ provided with the induced preorder. We use freely the internal logic of the ambient coherent category $\calC$.

(CC1) holds because $\ZERO$ is strict initial.

(CC2) A morphism $f\colon A \ar B$ is trivial when it factors through $(B,\Delta_B)$. From which the result, since by saturation, the preorder relation of $S\cup T$ is the union of the preorder relations of $S$ and $T$.

(CC3) Consider $s\colon (S,R_S)\mono (A,R)$ saturated.  First, $\phi(s)$ is a monomorphism. Two elements $\overline x$, $\overline y$ of $\phi(S,R_S)$ identified by $\phi(s)$ correspond to two elements $x$, $y$ in $S\subseteq A$ such that the pair $(x,y)$ lies in $R\cap R^\circ$ (see Example~\ref{preordex}). But since $(S,R_S)$ is saturated in $(A,R)$, $(S,R^\circ_S)$ is saturated in $(A,R^\circ)$ and the pair $(x,y)$ lies in $R_S\cap R^\circ_S$. This implies $\overline x=\overline y$. Thus $\phi(s)$ is a monomorphism, which we shall write as an inclusion.

Next, choose $(\overline x,\overline y)\in\eta_A(R)$ in $\phi(A,R)$ and assume that $\overline x\in\phi(S,R_S)$ or $\overline y\in \phi(S,R_S)$.  This means the existence of elements $x$, $y$ in $A$,  mapped respectively to $\overline x$ and $\overline y$, with thus $(x,y)\in R$ and one of the two elements  $x$, $y$ in $S$. By saturation, both elements $x$, $y$ are in $S$ and thus both elements $\overline x$, $\overline y$ are in $\phi(S)$.

(CC4)
Consider the following pullback diagram, with $S$ saturated in $\phi(A,R)$. We know already that $P$ is saturated in $A$. There is no restriction in writing these distinguished monomorphisms as canonical inclusions.
\begin{diagram}[80]
| | T ||
| \Near{h} | \smonO t ||
\calP | \Ear g | S||
\Smono p | | \smonO s||
A | \eaR{\eta_A} | \phi(A)||
\end{diagram}
Suppose that $g$ factors as $th$ through a saturated subobject $T$. Given an element $x\in S$, there is an element $y\in A$ such that $\eta_A(y)=x$, thus $y$ lies in the pullback $P$ and $g(y)=x$. This means that $x=t\bigl(h(y)\bigr)$, so that the morphism of preordered objects $t$ is an isomorphism in the ambient category $\calC$. Since moreover $T$ is saturated in $S$, $t$ is an isomorphism of preordered objects, and so, $g$ is a distinguished epimorphism.

(CC5)
Let again $s\colon(S,R_S)\mono(A,R)$ be a saturated subobject, which we write as a canonical inclusion. If $\phi(s)$ is an isomorphism, for every  $a\in A$, there exists $a'\in S$ such that  $\eta_A(a)=\eta_S(a')$. This means that the pair $(a,a')$ is in $R\cap R^\circ$ with $a'\in S$. Since $S$ is saturated in $A$, $a\in S$ and $s$ is an isomorphism.
\qed

\begin{example}\Label{compcat}
In the category of internal categories in a Grothendieck topos $\calE$, the (groupoid/ skeletal) pretorsion theory and the coherent system of saturated internal subcategories are compatible.
\end{example}

\proof
Let us first write down the proof in $\Cat$.

(CC1) A trivial category is one with only automorphisms. A saturated subcategory has only automorphisms as well, because it is full.

]CC2] A functor $F\colon\calS\cup\calT\ar\calB$ is trivial when it factors through the subcategory of automorphisms of $\calB$. Unions of saturated subcategories are computed as in $\Set$ at both the level of objects and the level of arrows (see \ref{cat}, from which the result).

In \cite{BCGT}, it is proved that $\phi(\calA)$  is constructed from $\calA$  as the quotient which identifies isomorphic objects. An arrow $a\ar b$  in $\phi(\calA)$, as for every quotient in $\Cat$, is represented  by a finite chain of arrows in $\calA$, which become composable in the quotient. But in this specific case, it is proved in \cite{BCGT} that a unique  ``reduced'' representation exists: a chain which does not contain any identity arrow nor any consecutive pair of arrows which are composable in $\calA$, and where thus -- by definition of the quotient -- the codomain of each arrow of the chain is isomorphic to the domain of the next one, while the domain of the first arrow is isomorphic to $a$ and the codomain of the last arrow is isomorphic to $b$. The identities in $\phi(\calA)$ are represented by the empty sequences.

(CC3) Consider a saturated subcategory $S\colon \calS\mono \calA$. First, two objects of $\calS$ identified in $\phi(\calA)$ are isomorphic in $\calA$, thus also in $\calS$ by saturation. This proves that $\phi(S)$ is injective on the objects.  Next,  if two arrows in reduced form in $\phi(\calS)$ are identified in $\phi(\calA)$, they are equal by uniqueness of the reduced form. Thus $\phi(S)$ is a monomorphism. Moreover given an arrow of $\phi(\calA)$ expressed in reduced form, if one of the objects involved in the chain is in $\calS$, by saturation, so are all the objects and all the arrows of the chain. Thus $\phi(\calS)$ is saturated in $\phi(\calA)$.

(CC4) We refer to the same diagram as for proving condition (CC4) in  Example~\ref{comppreord}, using for clarity calligraphic letters and upper case letters to indicate respectively the corresponding categories and functors.
The subcategory $\calS$ is thus saturated in $\phi(\calA)$ and by (CS3), its pullback $\calP$ is saturated in $\calA$. Let us view $S$ and $P$ as inclusions of full subcategories. Since $\eta_\calA$ is surjective on objects, an object $x\in \calS$ is equal in $\phi(\calA)$ to an object of the form $\eta_\calA(y)$, with thus $y$ an object of $\calA$. But then $y$ lies in the pullback $\calP$ and $g(y)=x$. We get further $x=t\bigl(h(y)\bigr)$ and $t$ is surjective, thus bijective, on the objects. But $t$ is also full by saturation, thus is an isomorphism.

(CC5) Let $S\colon\calS\mono\calA$ be a saturated subcategory such that $\phi(S)$ is an isomorphism. Write $S$ as a canonical inclusion. For every object $a\in\calA$, there is thus an object $a'\in\calS$ such that $\phi(a)$ and $\phi(a')$ are isomorphic in $\phi(\calA)$. By saturation of $\calS$ in $\calA$, this forces $a\in\calS$. Thus $\calS$ and $\calA$ have the same objects and $\calS$ is full in $\calA$, as a saturated subcategory. Therefore $\calS=\calA$.

Let us now switch to the case of a Grothendieck topos. In \cite{FC}, it is observed that in a topos of presheaves, the (groupoid/skeletal) pretorsion theory has canonical $\calZ$-exact sequences computed pointwise as in $\Set$. But being a monomorphism, a saturated subcategory, a pullback, a groupoid, a skeletal category are notions expressed in terms of finite limits and unions. Since  limits and colimits (and thus unions) are computed pointwise in a topos of presheaves, we get the expected result in every topos of presheaves, just because it holds pointwise  in $\Set$.

Again in \cite{FC} it is observed that the canonical $\calZ$-exact sequences for the (groupoid/skeletal) pretorsion theory in a topos of sheaves are obtained by constructing these sequences in the corresponding topos of presheaves, and applying the associated sheaf functor. Since this last functor preserves finite limits and colimits (and thus unions), we conclude the proof in the case of every Grothendieck topos.
\qed

\section{The stable category }\Label{5}

To avoid repeating it each time, let us put:\\[5pt]
{\bf Blanket assumption for this whole section}\\
{\it Let $\calE$ be a category with a strict initial object,
provided with a pretorsion theory and a coherent system of distinguished subobjects which are compatible.}\vspace{5pt}

First, let us construct the category of distinguished partial morphisms of $\calE$, which we shall denote as $\DisParE$.
\begin{itemize}
\item the objects are those of $\calE$;
\item a morphism from $A$ to $B$ is a triple $(S_0,S_1,f)$ where
\begin{itemize}
\item $S_0$ and $S_1$ are distinguished subobjects of $A$;
\item $S_0\cup S_1=A$;
\item $f\colon S_1\ar B$ is a morphism in $\calE$;
\item $f$ restricted to $S_0\cap S_1$ is a trivial morphism.
\end{itemize}
\item the composition law is made via a pullback:
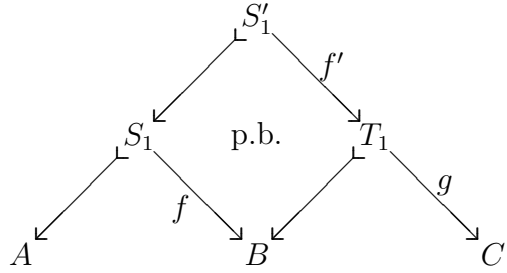
\begin{diagram}[56]
| | | | S'_1||
| | | \swmono | | \Sear{f'} ||
| | S_1 | | \mbox{p.b.} | | T_1 ||
| \swmono | | \seaR f | | \swmono | | \Sear g ||
A | | | | B | | | | C ||
\end{diagram}
$$
(T_0,T_1,g)\circ(S_0,S_1,f)=( S_0\cup S'_0,S'_1,gf')
$$
with $S'_0=f^{-1}(T_0)$ and $S'_1=f^{-1}(T_1)$.
\end{itemize}
Checking the category axioms is just routine computation.

There is an obvious functor $\iota\colon\calE\ar\DisParE$, which is the identity on objects and maps an arrow $f\colon A\ar B$ on $(\ZERO,A,f)$.

The partial morphisms should be thought as morphisms defined only on $S_1$, but that will be extended by zero on $S_0$ in the stable category,  in order to become defined on the whole of $A$.

Next we construct the stable category, denoted by $\StabE$. For this, we declare equivalent two parallel morphisms in $\DisParE$
$$
(S_0,S_1,f),(T_0,T_1,g)\colon A \biar B
$$
when there exists a so-called ``congruence'' diagram
\begin{diagram}[80]
\movevertexleft{U_0\cap S_1} | \emono | S_1 ||
| \nemono | \smono | \Sear f ||
U_1 | \emono | A | | B ||
| \semono | \nmono|  \neaR g ||
\movevertexleft{U_0\cap T_1} | \emono | T_1 ||
| | \makebox[0pt]{(Congruence diagram)} ||
\end{diagram}
with the properties:
\begin{itemize}
\item $U_0$ and $U_1$ are distinguished subobjects of $A$;
\item $U_0\cup U_1=A$;
\item $U_1$ is contained in both $S_1$ and $T_1$;
\item $f$ and $g$ coincide on $U_1$;
\item the restriction of $f$ to $U_0\cap S_1$ is trivial;
\item the restriction of $g$ to $U_0\cap T_1$ is trivial.
\end{itemize}
It is lengthy, but routine computation, to observe that this yields a congruence (see \cite{Ehresmann}) on $\DisParE$. The stable category $\StabE$ is the quotient of $\DisParE$ by that congruence. We write $\pi\colon\DisParE\ar\StabE$ for the quotient functor and $\sigma\colon \calE\ar\StabE$ for the composite $\pi\iota$. For objects and morphisms in $\calE$, we shall generally  write $A$ and $f$ in $\StabE$ instead of $\sigma(A)$ and $\sigma(f)$

Let us now prove the various properties of a stable category,  as announced at the end of Section~\ref{3}.

\begin{proposition}\Label{zero}
The category $\StabE$ admits $\ZERO$ as a zero object.
\end{proposition}

\proof
Writing $0_A\colon \ZERO\ar A$ for the unique morphism in $\calE$, the two partial distinguished morphisms
$$
(A,\ZERO,\id_{\ZERO})\colon A \ar \ZERO,~~~
(\ZERO,\ZERO,0_A)\colon \ZERO \ar A
$$
are the only possible ones in $\DisParE$, since $\ZERO$ is strict initial. It is immediate to observe that they are inverse isomorphisms in $\StabE$.
\qed

\begin{corollary}\Label{zeroar}
The zero morphism from $A$ to $B$ in $\StabE$ is represented by $(A,\ZERO,0_B)$.
\qed
\end{corollary}

\begin{lemma}\Label{zeropar}
A morphism $(S_0,S_1,f)\colon A \ar B$ in $\DisParE$ is mapped on a zero morphism by the quotient functor $\pi$ if and only if $f$ is trivial.
\end{lemma}

\proof
Consider a congruence diagram as above, expressing the possible  equivalence of $(S_0,S_1,f)$ and $(A,\ZERO,0_B)$. Again since $\ZERO$ is strict initial in $\calE$, we must have $U_1=\ZERO$  and thus $U_0=A$, which forces at once the result.
\qed

\begin{corollary}\Label{zeromap}
A morphism $f\colon A \ar B$ in $\calE$ is mapped on a zero morphism in $\StabE$ if and only if it is trivial.
\qed
\end{corollary}

\begin{proposition}\Label{zeroob}
An object $A\in\calE$ is mapped on a zero object in $\StabE$ if and only if it is trivial.
\end{proposition}

\proof
Expressing that the two partial morphisms in Proposition~\ref{zero} yield inverse isomorphisms in $\StabE$ reduces, by Lemma~\ref{zeropar}, to the identity on $A$ being trivial. But if the identity on $A$ factors through a trivial object $X$, $A$ is a retract of $X$ and therefore is trivial (see Corollary~2.8 in\cite{FFG}).
\qed

\begin{theorem}\Label{stabth}
In $\StabE$, choosing as torsion and torsion free objects the same objects as in $\calE$, one gets a torsion theory on $\StabE$ which turns the functor $\sigma\colon\calE\ar\StabE$ in a torsion functor.
\end{theorem}

\proof
Consider a morphism $A\ar B$ in $\StabE$, represented by $(S_0,S_1,f)$ in $\DisParE$. If $A$ is torsion and $B$ is torsion free, then $S_1$ is torsion by Corollary~\ref{tausub} and thus $f$ is trivial by axiom (PT1) in $\calE$ (see Definition~\ref{pretor}). By Lemma~\ref{zeropar}, axiom (PT1) holds in $\StabE$.

The core of the proof is about axiom (PT2). We fix an object $A\in\calE$ and must prove that the image by the functor $\sigma$ of its canonical $\calZ$-exact sequence in $\calE$
$$
\tau(A)\Mono{\varepsilon_A} A \Epi{\eta_A} \phi(A)
$$
is an exact sequence in $\StabE$.

First, the kernel part. Consider a morphism from $X$ to $A$ in $\StabE$ represented by a morphism $(S_0,S_1,x)$ in $\DisParE$ and whose composite with $\eta_A$ is zero. That composite is represented by $(S_0,S_1,\eta_Ax)$; by Lemma~\ref{zeropar}, saying that it is zero means that $\eta_Ax$ is trivial. Therefore $x$ factors uniquely in $\calE$ as $x=\varepsilon_Ay$. It follows at once that $(S_0,S_1,y)$ yields a factorization of $(S_0,S_1,x)$ through $\varepsilon_A$ in $\StabE$.

To prove the uniqueness, consider another factorization $(T_0,T_1,z)$. Write $U_0$, $U_1$ for the two distinguished subobjects in a congruence diagram exhibiting the equivalence between $(S_0,S_1,\varepsilon_Ay)$ and $(T_0,T_1,\varepsilon_Az)$. Since $\varepsilon_Ay$ and $\varepsilon_Az$ coincide on $U_1$, so do $y$ and $z$ because $\varepsilon_A$ is a monomorphism. The rest follows easily from Lemma~\ref{C7}.

To prove the cokernel part, let us choose in $\StabE$ a morphism from $A$ to $X$, represented by $(S_0,S_1,x)$ in $\DisParE$, and whose composite with $\varepsilon_A$ is zero. We consider then the situation
\begin{diagram}[80]
\tau(A) | \Emono{\varepsilon_A} | A | \Eepi{\eta_A} | \phi(A) ||
\Nmono{\tau(s)} | | \Nmono s | |  \nmonO{\phi(s)} ||
\tau(S_1) | \Emono{\varepsilon_{S_1}} | S_1 | \Eepi{\eta_{S_1}} | \phi(S_1) ||
| | | \seaR x | \sdotaR y ||
| | | | X ||
\end{diagram}
where, by Property~\ref{C5}, both squares are pullbacks and the vertical arrows are distinguished subobjects. In particular,
$$
(S_0,S_1,x)\circ \varepsilon_A
=\bigl(\varepsilon_A^{-1}(S_0), \tau(S_1), x\varepsilon_{S_1}\bigr).
$$
Since this composite is zero, $x\varepsilon_{S_1}$ is trivial, by Lemma~\ref{zeropar}. This implies the existence of a unique factorization $y$ of $x$ through $\eta_{S_1}$. Then $\bigl(\phi(S_0), \phi(S_1),y\bigr)$ yields the expected factorization of $(S_0,S_1,x)$ through $\eta_A$ in $\StabE$.

It remains to prove the uniqueness of the factorization. Let thus $(U_0,U_1,z)$ be another factorization of $(S_0,S_1,x)$ through $\eta_A$ in $\StabE$. Consider the following diagram
\begin{diagram}[80]
| | \movevertexleft{V'_1=\tau(V_1)} | \Emono{\varepsilon_{V_1}} | V_1 | \Eepi{\eta_{V_1}} | \movevertexright{\phi(V_1)=U_1} ||
| | \crossarrows{\Smono{\movename(0,5){\tau(v_1)}}}{\Enemono{\movename(-35,-10){w_1^V}}} | | \smonO{v_1} | | \Smono{\phi(v_1)} | \Sear z ||
W_1 | | \tau(A) | \Emono{\varepsilon_A} | A | \Eepi{\eta_A} | \phi(A) | | X ||
| | \crossarrows{\Nmono{\movename(0,-5){\tau(s_1)}}}{\esemonO{\movename(-35,7){w_1^S}}} | |  \nmonO{s_1} | | \Nmono{\phi(s_1)} |  \neaR y ||
| | \tau(S_1) | \emonO{\varepsilon_{S_1}} | S_1 | \eepI{\eta_{S_1}} | \phi(S_1) ||
\end{diagram}
where $V_1$ and $V'_1$ are obtained from $U_1$ by pullbacks. By Proposition~\ref{C6}, the upper line is a $\calZ$-exact sequence, allowing as in Proposition~\ref{C5} to rewrite it in terms of $V_1$, and the vertical morphisms in terms of $v_1$.
The equalities
\begin{multline*}
(V_0,V_1,z\eta_{V_1})
=(U_0,U_1,z)\circ\eta_A
=(S_0,S_1,x)
\\
=\bigl(\phi(S_0),\phi(S_1),y\bigr)\circ\eta_A
=\bigl(S_0,S_1,y\eta_{S_1})
\end{multline*}
mean the existence of a congruence diagram built from a distinguished covering $A=W_0\cup W_1$ of $A$.

We know by Corollary~\ref{phiuni} that $\phi(A)=\phi(W_0)\cup \phi(W_1)$ and we shall now prove that this covering allows constructing a congruence diagram proving the expected uniqueness of the factorization. Consider first the canonical $\calZ$-exact sequence of $W_1$.
\begin{diagram}[80]
\tau(W_1) | \Emono{\varepsilon_{W_1}} | W_1 | \Eepi{\eta_{W_1}} | \phi(W_1) ||
| | | \seaR{\eta_{V_1}\circ w_1^V} | \saR{\alpha_V} ||
| | | | \phi(V_1) ||
\end{diagram}
The composite $V_1 w_1^V \varepsilon_{W_1}$ is trivial by axiom (PT1), from which the factorization $\alpha_V$. The same argument can be repeated with $S_1$ and $w_1^S$, yielding $\alpha_S$ such that $\alpha_S\eta_{W_1}=\eta_{S_1}w_1^S$. Then
$$
z\alpha_V\eta_{W_1}
=z\eta_{V_1}w_1^V
=y\eta_{S_1}w_1^S
=y\alpha_S\eta_{W_1}.
$$
Since $\eta_{W_1}$ is an isomorphism, this implies $z\alpha_V=y\alpha_S$. This is the first condition for a congruence diagram based on $\phi(W_0)\cup \phi(W_1)$.

It remains to check the triviality of $z$ and $y$ when restricted to $\phi(W_0)$. For this, having in mind Corollaries~\ref{taudissub} and \ref{phiuni}, we consider the diagram
\begin{diagram}
\tau(W_0)\cap \tau(W_1) |
\Ear{\movename(0,3){\varepsilon_{W_0\cap W_1}}} |
W_0 \cap W_1 |
\Eepi{\movename(0,3){\eta_{W_0\cap W_1}}} |
\phi(W_0)\cap \phi(W_1) |
\Emono{\movename(0,3){\phi(i) }} |
\phi(W_1) ||
| | | \Sear{\beta} |
\saR{\delta} | |
\saR y ||
| | | | Z |
\eaR{\gamma} |
X ||
\end{diagram}
where $i$ is the inclusion of $W_0\cap W_1$ in $W_1$. From the congruence diagram built on the covering $A=W_0\cup W_1$, we know that $y\phi(i)\eta_{W_0\cap W_1}$ is trivial. Thus this morphism factors as $\gamma\beta$ through some trivial object $D$. But then $\beta$, and thus $\beta\varepsilon_{W_0\cap W_1}$ are trivial.  This implies that $\beta$ factors via a unique morphism $\delta$ through $\eta_{W_0\cap W_1}$. Since $\eta_{W_0\cap W_1}$ is an epimorphism, the various commutativities imply $\gamma\delta=y\phi(i)$. Thus $y\phi(i)$ factors through the trivial object $D$ and is trivial. The case of $z$ is perfectly analogous.
\qed

All the results of this section apply thus to Examples \ref{comptriv}, \ref{complext}, \ref{comppreord} and \ref{compcat}.

\section{The universal property}\Label{6}

Again we put:\\[5pt]
{\bf Blanket assumption for this whole section}\\
{\it Let $\calE$ be a category with a strict initial object,
provided with a pretorsion theory and a coherent system of distinguished subobjects which are compatible.}\vspace{5pt}

Let us begin with the motivating property that we emphasized in our discussion on stable categories at the end of Section~\ref{3}.

\begin{proposition}\Label{union}
Consider an object $A$ in $\calE$ which is the union $A=S\cup T$ of two distinguished subobjects. When $T$ is trivial, $A$ and $S$ become isomorphic in $\StabE$.
\end{proposition}

\proof
Write $s\colon S\mono A$ for the canonical inclusion. The two distinguished partial morphisms $(T,S,s)$ from $A$ to $S$ and $(\ZERO,S,\id_S)$ from $S$ to $A$ exhibit the isomorphism in $\StabE$.
\qed

This property somehow suggests the following definition.

\begin{definition}\Label{zeropush}
Let $\calX$ be a category with a zero object. A commutative square $ca=db$ is called a {\em zero-pushout}
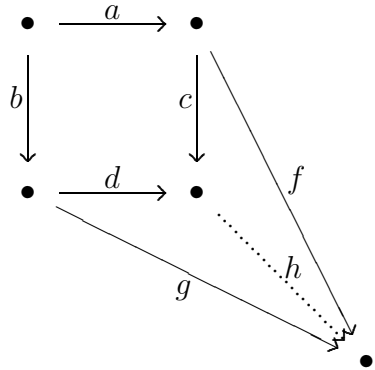
\begin{diagram}[80]
\bullet | \Ear a | \bullet ||
\Sar b | | \Sar c ||
\bullet | \Ear d | \bullet | \Ssear f ||
|  | \eseaR g | \Sedotar h ||
| | | | \bullet ||
\end{diagram}
when given $arrows$ $f$, $g$, with $fa=gb$, there exists a unique morphism $h$ such that $hc=f$ and $hd=g$, provided that at least one of the two morphisms $f$, $g$ is zero.
\end{definition}

\begin{proposition}\Label{sigmapush}
Given two distinguished subobjects $S$, $T$ of $A$ in $\calE$, the functor $\sigma\colon\calE\ar\StabE$ transforms the bicartesian square
\begin{diagram}[80]
S\cap T | \Emono{s'} | S ||
\Smono{t'} | | \smonO s ||
T | \emonO t | S\cup T ||
\end{diagram}
in a zero-pushout.
\end{proposition}

\proof
Let us consider the following diagram, where the plain arrows are in $\calE$ and the dotted arrows in $\StabE$.
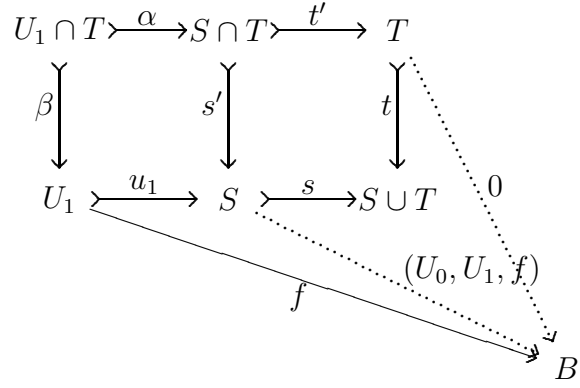
\begin{diagram}[80]
U_1\cap T | \Emono{\alpha} | S\cap T | \Emono{t'} | T ||
\Smono{\beta} | | \Smono{s'} | | \Smono{t} ||
U_1 | \Emono{u_1} | S | \Emono s | S\cup T | \Ssedotar 0 ||
| | | \eeseaR f | \Esedotar{(U_0,U_1,f)} ||
| | | | | | B ||
\end{diagram}
In $\StabE$, we assume thus $(U_0,U_1,f)\circ s'=0\circ t'$. By Lemma~\ref{zeropar}, $f\beta$ is trivial, thus $f$ is trivial on $U_1\cap T$. But by definition of the arrows in $\DisParE$, $f$ is also trivial on $U_1\cap U_0$. By condition~(CC2), $f$ is thus trivial on
$$
(U_1\cap  T)\cup (U_1\cap U_0) = U_1 \cap (T\cup U_0).
$$
Since $U_1\cup (T\cup U_0) = S\cup T$, this proves that we have a morphism in $\StabE$
$$
(T\cup U_0, U_1, f)\colon S\cup T \dotar B
$$
and it follows at once that this is the expected factorization.

It remains to prove the uniqueness of that factorization. Let $(V_0,V_1,g)$ be another factorization, with thus in particular $V_0\cup V_1=S\cup T$. On one hand the equality $(V_0,V_1,g)\circ s=(U_0,U_1,f)$ implies the existence of a distinguished covering $S=W_0\cup W_1$  and a corresponding congruence diagram
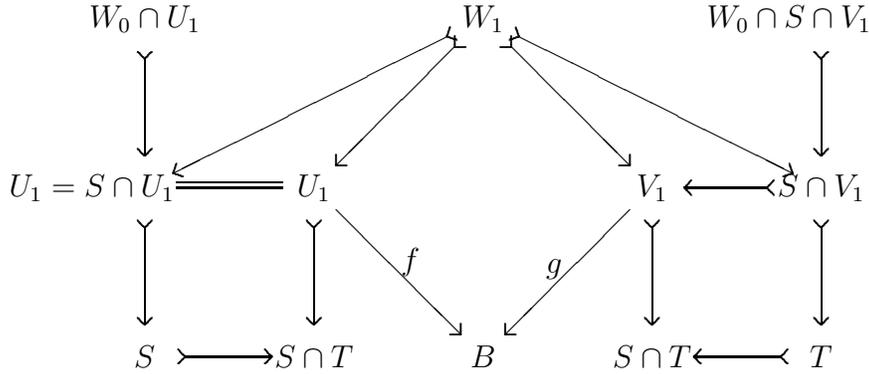
\begin{diagram}[80]
W_0\cap U_1| | | | W_1 | | | \movevertexright{W_0\cap S \cap V_1} ||
\smono| |  \wswmono | \swmono| |\semono |  \esemono| | \smono ||
\movevertexleft{U_1=S\cap U_1}| \eeql | U_1 | | | | V_1 | \wmono | S\cap V_1 ||
\smono | | \smono | \Sear f | | \Swar g | \smono | | \smono ||
S | \emono | S\cap T | | B | | S\cap T | \wmono | T ||
\end{diagram}
By this congruence diagram,  $f$ and $g$ coincide on $W_1$, while $f$ is trivial on $W_0\cap U_1$ and $g$ is trivial on $W_0\cap S \cap V_1$. On the other hand the equality $(V_0,V_1,g)\circ y=0$ means that $g$ is trivial on $T\cap V_1$. And we know already that $f$ is trivial on $T\cap U_1$.
Observe now that
$$
S\cup T= W_0\cup W_1\cup =W_1\cup (W_0\cup T)
$$
is a covering of $S\cup T$, while $f$ and $g$ coincide already on $W_1$. To conclude the proof, it remains to see that $f$ and $g$ are trivial when restricted respectively to $(W_0\cup T)\cap U_1$ and $(W_0\cup T)\cap  V_1$. We have
$$
(W_0\cup T)\cap U_1=(W_0\cap U_1) \cup (T\cap U_1)
$$
and we know already that $f$ is trivial on both pieces of this union; condition~(CC2) yields the conclusion for $f$. In the same way
$$
(W_0\cup T)\cap  V_1= (W_0\cap V_1) \cup (T \cap V_1)
$$
and we know already that $g$ is trivial on $T\cap V_1$ and $W_0\cap S \cap V_1$. This takes already care of the term $T\cap V_1$; for the other term, simply split it in two parts and apply again condition~(CC2)
$$
W_0\cap V_1
=(W_0\cap V_1)\cap (S\cup T)
=(W_0\cap V_1\cap S)\cup(W_0\cap V_1\cap T).
$$
The morphism  $g$ is already known to be trivial on each part.
\qed

\begin{theorem}\Label{univ}
The functor $\sigma\colon \calE\ar\StabE$ is universal among all the functors
$
F\colon \calE \ar \calX
$
where
\begin{enumerate}
\item $\calX$ is provided with a torsion theory;
\item $F$ is a torsion functor;
\item given two distinguished subobjects $S$, $T$ of $A$ in $\calE$, $F$ transforms the square
\begin{diagram}[80]
S\cap T | \emono | S ||
\smono | | \smono ||
T | \emono | S\cup T ||
\end{diagram}
in a zero-pushout.
\end{enumerate}
\end{theorem}

\proof
Call $G\colon\StabE\ar \calX$ the expected unique factorization of $F$ through $\sigma$. Since $\sigma$ is the identity on objects, we must put $G(A)=F(A)$ for every object $A\in\calE$. In the same way for every arrow $f\in\calE$, the factorization requirement imposes $G\bigl(\sigma(f)\bigr)=F(f)$.

Next consider a morphism from $A$ to $B$ in $\StabE$ represented by $(S_0,S_1,f)$ in $\DisParE$. We have in particular the following left hand diagram in $\calE$, with $S_0\cup S_1 = A$.
\begin{diagram}[80]
S_0\cap S_1 | \Emono{s'_1} | S_1 | | | |
F(S_0\cap S_1) | \Ear{F(s'_1)} | F(S_1)
||
\Smono{s'_0} | | \Smono{s_1} | | | |
\Sar{F(s'_0)} | | \Sar{F(s_1)}
||
S_0 | \emonO{s_0} | A | \Ssear f | | |
F(S_0) | \Ear{F(s_0)} | F(A) | \Ssear{F(f)}
||
| | | | | |
| | \eseaR 0 | \Sedotar{\varphi}
||
| | | | B  | |
| | | | F(B)
||
\end{diagram}
yielding the corresponding right hand diagram in $\calX$,
with the square a zero-pushout. By definition of a morphism in $\DisParE$, $fs'_1$ is trivial, thus $F(f)F(s'_1)=0$ because $F$ is a torsion functor. By the zero-pushout property, we get a unique factorization $\varphi$ which we choose as $G(S_0,S_1,f)$.

We must of course verify that this definition is independent of the choice of the distinguished partial morphism representing the morphism in $\StabE$. Using analogous notation, let the partial morphism $(T_0,T_1,g)$ be equivalent to $(U_0,U_1,f)$ in $\DisParE$ and call $\psi$ the corresponding factorization from $F(A)$ to $F(B)$. We have thus a congruence diagram with $U_0\cup U_1=A$
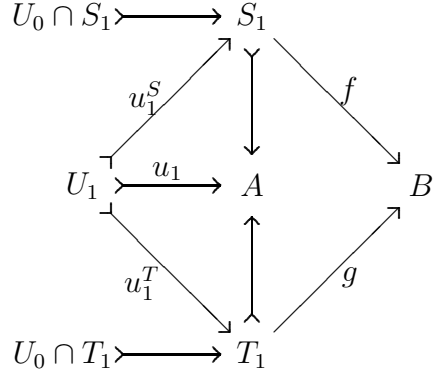
\begin{diagram}[80]
\movevertexleft{U_0\cap S_1} | \emono | S_1 ||
| \Nemono{u_1^S} | \smono | \Sear f ||
U_1 | \Emono{u_1} | A | | B ||
| \semonO{u_1^T} | \nmono|  \neaR g ||
\movevertexleft{U_0\cap T_1} | \emono | T_1 ||
\end{diagram}

Still with analogous notation, we consider further the following diagram
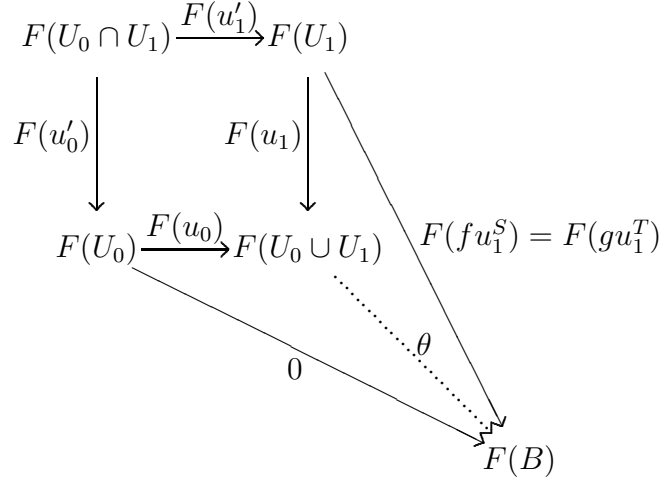
\begin{diagram}
F(U_0\cap U_1) | \Ear{F(u'_1)} | F(U_1) ||
\Sar{F(u'_0)} | | \Sar{F(u_1)} ||
F(U_0) | \Ear{F(u_0)} | F(U_0\cup U_1) | \Ssear{F(fu_1^S)=F(gu_1^T)}||
| | \eseaR 0 | \Sedotar{\theta}||
| | | | F(B) ||
\end{diagram}
yielding the unique factorization $\theta$ through the zero-pushout. But $U_0\cup U_1=A$ so that $\theta$ is a morphism from $F(A)$ to $F(B)$. Let us prove that both morphisms $\varphi=G(U_0,U_1,f)$ and $\psi=G(T_0,T_1,g)$ are equal to $\theta$. Of course it suffices to do the job for one of them: we do it for $\varphi$.

We have first
$$
\varphi F(u_1) =\varphi F(s_1)F(u_1^S)=F(f)F(u_1^S).
$$
To prove that $\varphi\circ F(u_0)=0$, observe that
$$
U_0=U_0\cap A = U_0\cap (S_0\cup S_1) = (U_0 \cap S_0)\cup (U_0\cap S_1).
$$
Using the zero-pushout in $\calX$ constructed from these last two distinguished subobjects in $\calE$, it suffices to prove that $\varphi\circ F(u_0)$ is zero on  $F(U_0\cap S_0)$ and $F( U_0\cap S_1)$. In the first case it is because $\varphi$ is already zero at the level $F(S_0)$ and in the second case, because $F(f)$ is zero at the level   $F(U_0\cap S_1)$.

Checking the functoriality of $G$  is just routine verification. Its uniqueness  is attested by that of $\varphi$ at the beginning of the proof. And by definition of the torsion theory on $\StabE$ and assumption on $F$, $G$ is a torsion functor.
\qed

All the results of this section apply again to Examples \ref{comptriv}, \ref{complext}, \ref{comppreord} and \ref{compcat}.

\section{The associated torsion theory}\Label{7}

In this section, we want to focus on the torsion theory universally associated with a pretorsion theory, without any further requirement of compatibility with unions, coproducts, or whatever.

\begin{theorem}\Label{basic}
Let $\calE$ be a category provided with a pretorsion theory and a strict initial object. Suppose that the initial object is trivial. There exists a category provided with a torsion theory, universally associated with those data.
\end{theorem}

\proof
The universal solution is the stable category corresponding to the choice of the indiscrete coherent system of subobobjects (see Example \ref{comptriv}).  Indeed, given two distinguished subobjects of $A$, the only possible bicartesion squares as in Theorem~\ref{univ} are
\begin{diagram}[80]
\ZERO | \eeql| \ZERO  |       |  A | \eeql | A |  |  \ZERO| \emono | A ||
\seql | | \seql |                  | \seql | | \seql |  |  \seql | | \seql ||
\ZERO  | \eeql | \ZERO |    | A | \eeql | A  |  |  \ZERO | \emono | A ||
\end{diagram}
and these are trivially sent to actual pushouts by every torsion functor.
\qed

\begin{corollary}\Label{descrip}
In the conditions of Theorem~\ref{basic}, the universal stable category $\StabE$ admits the following description:
\begin{itemize}
\item the objects are those of $\calE$;
\item the arrows $A\ar B$ in $\StabE$ are the non-trivial arrows from $A$ to $B$ in $\calE$, together with a formally added arrow $0_{AB}$.
\end{itemize}
The composition in $\StabE$ of two non-trivial morphisms of $\calE$ is their composition in $\calE$ when this composite is non-trivial, and the corresponding zero morphism  otherwise. Composing whatever arrow with a zero arrow is the corresponding zero arrow.
\end{corollary}

\proof
The morphisms from $A$ to $B$ in $\DisParE$ are
\begin{itemize}
\item $(\ZERO,A,f)$ for every morphism $f\colon A \ar B$;
\item $(A,\ZERO,0_B)$;
\item $(A,A,f)$ for every trivial morphism $f\colon A \ar B$.
\end{itemize}
By Lemma~\ref{zeropar}, all cases where $f$ is trivial are identified with $(A,\ZERO,0_B)$ in $\StabE$. And when $f$, $g$ are non-trivial, the only possible choice for $U_0$ in a congruence diagram is $U_0=\ZERO$, thus $U_1=A$, which yields $f=g$. The conclusion follows at once, when putting  $0_{AB}=\sigma(A,\ZERO,B)$
\qed

\begin{proposition}\Label{frac}
Consider a category $\calE$ provided with a torsion theory, a strict initial object and a terminal object. Suppose that the initial and the terminal object are both trivial.  The stable category $\StabE$ of Theorem \ref{basic} is then the category of fractions $\calE[\Sigma^{-1}]$, where $\Sigma$ is the class of morphisms $\ZERO\mono A$, with $A$ trivial. Equivalently, $\Sigma$ can be chosen as just the singleton $\{\ZERO\mono\ONE\}$.
\end{proposition}

\proof
When $A$ is trivial, $0_A\colon \ZERO \ar A$ becomes an isomorphism in $\StabE$, by Proposition~\ref{zeroob}. This is in particular the case for the unique morphism $\xi\colon \ZERO \ar \ONE$.

Choose now a functor $F\colon\calE\ar\calX$ which transforms the  morphism $\xi$ in an isomorphism.
\begin{diagram}[80]
\calE | \Ear{\sigma} | \movevertexright{\StabE} ||
| \seaR F | \sdotaR G ||
| | \calX ||
\end{diagram}
We must prove the existence of a unique factorization $G$. Since $\sigma$ is bijective on objects, we must put $G(A)=F(A)$ for every object $A$. And of course by functoriality, we must put $F(\xi^{-1})=F(\xi)^{-1}$.

Let us use the description of $\StabE$ as in Corollary \ref{descrip}, writing $1_A$ for the unique morphism from $A$ to $\ONE$ in $\calE$. If $f$ is a non trivial morphism in $\calE$, we must again put $G(f)=F(f)$. And when it goes about the zero morphism $0_{AB}$ in $\StabE$, we have the situation
\begin{diagram}[80]
A | \Edotar{0_{AB}} | B   | | | | F(A) | \Edotar{F(0_{AB})} | F(B) ||
\Sar{\sigma(1_A)} | | \naR{\sigma(0_B)} | | | |\Sar{F(1_A)} | | \naR{F(0_B)} ||
\ONE | \edotaR{\xi^{-1}} | \ZERO | | |  | F(\ONE) | \eaR{F(\xi)^{-1}}| F(\ZERO) ||
\end{diagram}
where the left hand square is commutative in $\StabE$. We are thus forced to define
$$
G(0_{AB})=F(0_B)\circ F(\xi)^{-1}\circ F(1_A).
$$
This shows in particular that $F(u)=F(0_B)\circ F(\xi)^{-1}\circ F(1_A)$ for every trivial morphism $u\colon A \ar B$.

Checking the functoriality of $G$ is then routine. For example consider the most involved case: two non trivial morphisms $f$, $g$ in $\calE$ whose composite is trivial. We have thus a factorization $fg=uv$ through a trivial object $D$, yielding the following situation in $\calX$
\begin{diagram}[80]
 F(A) | | | \Ear[190] {F(f)} | |  | F(B) ||
 | \Sear{F(1_A)} | | | | \Near{F(0_B)} ||
 \Sar[130]{F(u)} | | F(\ONE) | \Ear{F(\xi)^{-1}} | F(\ZERO) |  | \saR[130]{F(g)} ||
 | \Near{\movename(14,8){F(1_D)}} |  \wswaR{F(0_D)} | | | \Sear{\movename(-15,6){F(0_C)}}||
 F(D) | | | \eaR[190]{F(v)}  | | | F(C) ||
\end{diagram}
Since $f$ and $g$ are non trivial, we have
$$
G(\sigma(g))\circ G(\sigma(f))=F(g)\circ F(f) = F(gf) = F(vu) = F(v)\circ F(u).
$$
Since $D$ is trivial, so are $u$ and $v$ and thus we know already that
$$
F(u)=F(0_B)\circ F(\xi)^{-1}\circ F(1_A),~~~
F(v)=F(0_C)\circ F(\xi)^{-1}\circ F(1_B) .
$$
The composite of these two morphisms reduces to  $F(0_C)\circ F(\xi)^{-1}\circ F(1_A)$ since $1_B\circ 0_B =\id_{\ZERO}$. But this is precisely $G(0_{AC})$, that is $G(\sigma(gf))$ since $gf$ is trivial.
\qed

\vspace*{8mm}

\noindent {\sc Francis Borceux}, Universit\'{e} catholique de Louvain, Institut de Recherche em Math\'{e}matique et Physique, 1348 Louvain-la-Neuve, Belgium, francis.borceux@uclouvain.be\\

\noindent {\sc Maria Manuel Clementino}, Universidade de Coimbra, CMUC, Departamento de Matem\'{a}tica, 3001-501 Coimbra, Portugal,
 mmc@mat.uc.pt

\end{document}

%% file: diagram.tex
\catcode`\@=11

\newcount\printerresolution
\global\printerresolution=300

\newcount\headsize
\global\headsize=450

\font\dotfont = lcircle10 at 3pt

\newcount\epihead
\global\epihead=200

\newcount\defaultscale
\def\setdefaultscale#1{\global\defaultscale=#1}
\setdefaultscale{100}

\newcount\actualtextarrowspace
\newcount\actualtextarrowlength

\newcommand{\computetextparameters}%
{\global\actualtextarrowspace=\textarrowlength%
\global\advance\actualtextarrowspace by 3%
\global\actualtextarrowlength=\textarrowlength%
\global\multiply\actualtextarrowlength by 100}

\newcount\textarrowlength
\def\settextarrowlength#1{\global\textarrowlength=#1%
\computetextparameters}
\settextarrowlength{20}

\newcount\actualdisplayarrowspace
\newcount\actualdisplayarrowlength

\newcommand{\computedisplayparameters}%
{\global\actualdisplayarrowspace=\displayarrowlength%
\global\advance\actualdisplayarrowspace by 3%
\global\actualdisplayarrowlength=\displayarrowlength%
\global\multiply\actualdisplayarrowlength by 100}

\newcount\displayarrowlength
\def\setdisplayarrowlength#1{\global\displayarrowlength=#1%
\computedisplayparameters}
\setdisplayarrowlength{40}

\def\@ifnexttok#1#2#3{\let\@tempe #1\def\@tempa{#2}\def\@tempb{#3}%
\futurelet\@tempc\@ifntok}
\def\@ifntok{\ifx \@tempc \@tempe\let\@tempd\@tempa\else\let\@tempd\@tempb\fi%
\@tempd}

\def\DOT{{\dotfont q}}

\newcount\basicstep
\global\divide\basicstep by \printerresolution

\newcount\numberofsteps
\numberofsteps=\headsize
\global\divide\numberofsteps by \basicstep

\newcommand{\makehead}[3]{%
\begin{picture}(0,0)%
\multiput(0,0)(#1,#2){#3}{\DOT}%
\multiput(0,0)(-#2,#1){#3}{\DOT}%
\end{picture}}

\newcount\xstep
\newcount\ystep

\newsavebox{\northhead}
\savebox{\northhead}{%
\xstep=-\basicstep%
\multiply\xstep by 7071%
\divide\xstep by 10000%
\ystep=\xstep%
\makehead{\xstep}{\ystep}{\numberofsteps}}
\newcommand{\nhead}{\usebox{\northhead}}

\newsavebox{\easthead}
\savebox{\easthead}{%
\xstep=-\basicstep%
\multiply\xstep by 7071%
\divide\xstep by 10000%
\ystep=-\xstep%
\makehead{\xstep}{\ystep}{\numberofsteps}}
\newcommand{\ehead}{\usebox{\easthead}}

\newsavebox{\southhead}
\savebox{\southhead}{%
\xstep=\basicstep%
\multiply\xstep by 7071%
\divide\xstep by 10000%
\ystep=\xstep%
\makehead{\xstep}{\ystep}{\numberofsteps}}
\newcommand{\shead}{\usebox{\southhead}}

\newsavebox{\westhead}
\savebox{\westhead}{%
\xstep=\basicstep%
\multiply\xstep by 7071%
\divide\xstep by 10000%
\ystep=-\xstep%
\makehead{\xstep}{\ystep}{\numberofsteps}}
\newcommand{\whead}{\usebox{\westhead}}

\newsavebox{\northwesthead}
\savebox{\northwesthead}{%
\makehead{0}{-\basicstep}{\numberofsteps}}
\newcommand{\nwhead}{\usebox{\northwesthead}}

\newsavebox{\northeasthead}
\savebox{\northeasthead}{%
\makehead{-\basicstep}{0}{\numberofsteps}}
\newcommand{\nehead}{\usebox{\northeasthead}}

\newsavebox{\southwesthead}
\savebox{\southwesthead}{%
\makehead{\basicstep}{0}{\numberofsteps}}
\newcommand{\swhead}{\usebox{\southwesthead}}

\newsavebox{\southeasthead}
\savebox{\southeasthead}{%
\makehead{0}{\basicstep}{\numberofsteps}}
\newcommand{\sehead}{\usebox{\southeasthead}}

\newsavebox{\eastnortheasthead}
\savebox{\eastnortheasthead}{%
\xstep=-\basicstep%
\multiply\xstep by 9486%
\divide\xstep by 10000%
\ystep=\xstep%
\divide\ystep by -3%
\makehead{\xstep}{\ystep}{\numberofsteps}}
\newcommand{\enehead}{\usebox{\eastnortheasthead}}

\newsavebox{\northnortheasthead}
\savebox{\northnortheasthead}{%
\xstep=-\basicstep%
\multiply\xstep by 9486%
\divide\xstep by 10000%
\ystep=\xstep%
\divide\ystep by 3%
\makehead{\xstep}{\ystep}{\numberofsteps}}
\newcommand{\nnehead}{\usebox{\northnortheasthead}}

\newsavebox{\southsouthwesthead}
\savebox{\southsouthwesthead}{%
\xstep=\basicstep%
\multiply\xstep by 9486%
\divide\xstep by 10000%
\ystep=\xstep%
\divide\ystep by 3%
\makehead{\xstep}{\ystep}{\numberofsteps}}
\newcommand{\sswhead}{\usebox{\southsouthwesthead}}

\newsavebox{\westsouthwesthead}
\savebox{\westsouthwesthead}{%
\xstep=\basicstep%
\multiply\xstep by 9486%
\divide\xstep by 10000%
\ystep=\xstep%
\divide\ystep by -3%
\makehead{\xstep}{\ystep}{\numberofsteps}}
\newcommand{\wswhead}{\usebox{\westsouthwesthead}}

\newsavebox{\westnorthwesthead}
\savebox{\westnorthwesthead}{%
\xstep=\basicstep%
\multiply\xstep by 3162%
\divide\xstep by 10000%
\ystep=\xstep%
\multiply\ystep by -3%
\makehead{\xstep}{\ystep}{\numberofsteps}}
\newcommand{\wnwhead}{\usebox{\westnorthwesthead}}

\newsavebox{\eastsoutheasthead}
\savebox{\eastsoutheasthead}{%
\xstep=-\basicstep%
\multiply\xstep by 3162%
\divide\xstep by 10000%
\ystep=\xstep%
\multiply\ystep by -3%
\makehead{\xstep}{\ystep}{\numberofsteps}}
\newcommand{\esehead}{\usebox{\eastsoutheasthead}}

\newsavebox{\northnorthwesthead}
\savebox{\northnorthwesthead}{%
\xstep=-\basicstep%
\multiply\xstep by 3162%
\divide\xstep by 10000%
\ystep=\xstep%
\multiply\ystep by 3%
\makehead{\xstep}{\ystep}{\numberofsteps}}
\newcommand{\nnwhead}{\usebox{\northnorthwesthead}}

\newsavebox{\southsoutheasthead}
\savebox{\southsoutheasthead}{%
\xstep=\basicstep%
\multiply\xstep by 3162%
\divide\xstep by 10000%
\ystep=\xstep%
\multiply\ystep by 3%
\makehead{\xstep}{\ystep}{\numberofsteps}}
\newcommand{\ssehead}{\usebox{\southsoutheasthead}}

\newsavebox{\easteastnortheasthead}
\savebox{\easteastnortheasthead}{%
\xstep=-\basicstep%
\multiply\xstep by 8944%
\divide\xstep by 10000%
\ystep=\xstep%
\divide\ystep by -2%
\makehead{\xstep}{\ystep}{\numberofsteps}}
\newcommand{\eenehead}{\usebox{\easteastnortheasthead}}

\newsavebox{\northnorthnortheasthead}
\savebox{\northnorthnortheasthead}{%
\xstep=-\basicstep%
\multiply\xstep by 8944%
\divide\xstep by 10000%
\ystep=\xstep%
\divide\ystep by 2%
\makehead{\xstep}{\ystep}{\numberofsteps}}
\newcommand{\nnnehead}{\usebox{\northnorthnortheasthead}}

\newsavebox{\southsouthsouthwesthead}
\savebox{\southsouthsouthwesthead}{%
\xstep=\basicstep%
\multiply\xstep by 8944%
\divide\xstep by 10000%
\ystep=\xstep%
\divide\ystep by 2%
\makehead{\xstep}{\ystep}{\numberofsteps}}
\newcommand{\ssswhead}{\usebox{\southsouthsouthwesthead}}

\newsavebox{\westwestsouthwesthead}
\savebox{\westwestsouthwesthead}{%
\xstep=\basicstep%
\multiply\xstep by 8944%
\divide\xstep by 10000%
\ystep=\xstep%
\divide\ystep by -2%
\makehead{\xstep}{\ystep}{\numberofsteps}}
\newcommand{\wwswhead}{\usebox{\westwestsouthwesthead}}

\newsavebox{\westwestnorthwesthead}
\savebox{\westwestnorthwesthead}{%
\xstep=\basicstep%
\multiply\xstep by 4472%
\divide\xstep by 10000%
\ystep=\xstep%
\multiply\ystep by -2%
\makehead{\xstep}{\ystep}{\numberofsteps}}
\newcommand{\wwnwhead}{\usebox{\westwestnorthwesthead}}

\newsavebox{\easteastsoutheasthead}
\savebox{\easteastsoutheasthead}{%
\xstep=-\basicstep%
\multiply\xstep by 4472%
\divide\xstep by 10000%
\ystep=\xstep%
\multiply\ystep by -2%
\makehead{\xstep}{\ystep}{\numberofsteps}}
\newcommand{\eesehead}{\usebox{\easteastsoutheasthead}}

\newsavebox{\northnorthnorthwesthead}
\savebox{\northnorthnorthwesthead}{%
\xstep=-\basicstep%
\multiply\xstep by 4472%
\divide\xstep by 10000%
\ystep=\xstep%
\multiply\ystep by 2%
\makehead{\xstep}{\ystep}{\numberofsteps}}
\newcommand{\nnnwhead}{\usebox{\northnorthnorthwesthead}}

\newsavebox{\southsouthsoutheasthead}
\savebox{\southsouthsoutheasthead}{%
\xstep=\basicstep%
\multiply\xstep by 4472%
\divide\xstep by 10000%
\ystep=\xstep%
\multiply\ystep by 2%
\makehead{\xstep}{\ystep}{\numberofsteps}}
\newcommand{\sssehead}{\usebox{\southsouthsoutheasthead}}

\newsavebox{\northeasteastnortheasthead}
\savebox{\northeasteastnortheasthead}{%
\xstep=-\basicstep%
\multiply\xstep by 9806%
\divide\xstep by 10000%
\ystep=\xstep%
\divide\ystep by -5%
\makehead{\xstep}{\ystep}{\numberofsteps}}
\newcommand{\neenehead}{\usebox{\northeasteastnortheasthead}}

\newsavebox{\northeastnorthnortheasthead}
\savebox{\northeastnorthnortheasthead}{%
\xstep=-\basicstep%
\multiply\xstep by 9806%
\divide\xstep by 10000%
\ystep=\xstep%
\divide\ystep by 5%
\makehead{\xstep}{\ystep}{\numberofsteps}}
\newcommand{\nennehead}{\usebox{\northeastnorthnortheasthead}}

\newsavebox{\southwestsouthsouthwesthead}
\savebox{\southwestsouthsouthwesthead}{%
\xstep=\basicstep%
\multiply\xstep by 9806%
\divide\xstep by 10000%
\ystep=\xstep%
\divide\ystep by 5%
\makehead{\xstep}{\ystep}{\numberofsteps}}
\newcommand{\swsswhead}{\usebox{\southwestsouthsouthwesthead}}

\newsavebox{\southwestwestsouthwesthead}
\savebox{\southwestwestsouthwesthead}{%
\xstep=\basicstep%
\multiply\xstep by 9806%
\divide\xstep by 10000%
\ystep=\xstep%
\divide\ystep by -5%
\makehead{\xstep}{\ystep}{\numberofsteps}}
\newcommand{\swwswhead}{\usebox{\southwestwestsouthwesthead}}

\newsavebox{\northwestwestnorthwesthead}
\savebox{\northwestwestnorthwesthead}{%
\xstep=\basicstep%
\multiply\xstep by 1961%
\divide\xstep by 10000%
\ystep=\xstep%
\multiply\ystep by -5%
\makehead{\xstep}{\ystep}{\numberofsteps}}
\newcommand{\nwwnwhead}{\usebox{\northwestwestnorthwesthead}}

\newsavebox{\southeasteastsoutheasthead}
\savebox{\southeasteastsoutheasthead}{%
\xstep=-\basicstep%
\multiply\xstep by 1961%
\divide\xstep by 10000%
\ystep=\xstep%
\multiply\ystep by -5%
\makehead{\xstep}{\ystep}{\numberofsteps}}
\newcommand{\seesehead}{\usebox{\southeasteastsoutheasthead}}

\newsavebox{\northwestnorthnorthwesthead}
\savebox{\northwestnorthnorthwesthead}{%
\xstep=-\basicstep%
\multiply\xstep by 1961%
\divide\xstep by 10000%
\ystep=\xstep%
\multiply\ystep by 5%
\makehead{\xstep}{\ystep}{\numberofsteps}}
\newcommand{\nwnnwhead}{\usebox{\northwestnorthnorthwesthead}}

\newsavebox{\southeastsouthsoutheasthead}
\savebox{\southeastsouthsoutheasthead}{%
\xstep=\basicstep%
\multiply\xstep by 1961%
\divide\xstep by 10000%
\ystep=\xstep%
\multiply\ystep by 5%
\makehead{\xstep}{\ystep}{\numberofsteps}}
\newcommand{\sessehead}{\usebox{\southeastsouthsoutheasthead}}

\newsavebox{\isomorphismmark}
\newcommand{\isomark}[1]{\savebox{\isomorphismmark}{#1}}
\isomark{$\cong$}
\newcommand{\isosign}{\usebox{\isomorphismmark}}

\newif\ifuserdist
\newsavebox{\distributormark}
\newcommand{\distmark}[1]{\ifx#1\distcircle\userdistfalse\else%
\userdisttrue\savebox{\distributormark}{#1}\fi}
\newsavebox{\distributorcircle}
\savebox{\distributorcircle}{\begin{picture}(0,0)%
\put(0,0){\circle{4}}\end{picture}}
\newcommand{\distsign}{\ifuserdist\usebox{\distributormark}\else%
\usebox{\distributorcircle}\fi}
\userdistfalse

\newcount\monolength
\newcount\epilength
\newcount\bimolength
\newcount\secondmonolength
\newcount\secondepilength

\newcount\monotail%
\global\monotail=\headsize%
\global\multiply\monotail by 7070%
\global\divide\monotail by 10000%

\newcount\truemonotail%
\newcount\trueepihead%

\def\truetail{\truemonotail=\monotail%
\multiply\truemonotail by 100%
\divide\truemonotail by \SCALE}

\def\truehead{\trueepihead=\epihead%
\multiply\trueepihead by 100%
\divide\trueepihead by \SCALE}

\newcount\Monotail%
\global\Monotail=\headsize%
\global\divide\Monotail by 2%

\newcount\Epihead%
\global\Epihead=\epihead%
\global\multiply\Epihead by 7070%
\global\divide\Epihead by 10000%

\newcount\Truemonotail%
\newcount\Trueepihead%

\def\Truetail{\Truemonotail=\Monotail%
\multiply\Truemonotail by 100%
\divide\Truemonotail by \SCALE}%

\def\Truehead{\Trueepihead=\Epihead%
\multiply\Trueepihead by 100%
\divide\Trueepihead by \SCALE}

\newcount\MonoTail%
\global\MonoTail=\headsize%
\global\multiply\MonoTail by 6324%
\global\divide\MonoTail by 10000%

\newcount\EpiHead%
\global\EpiHead=\epihead%
\global\multiply\EpiHead by 8944%
\global\divide\EpiHead by 10000%

\newcount\TrueMonoTail%
\newcount\TrueEpiHead%

\def\TrueTail{\TrueMonoTail=\MonoTail%
\multiply\TrueMonoTail by 100%
\divide\TrueMonoTail by \SCALE}%

\def\TrueHead{\TrueEpiHead=\EpiHead%
\multiply\TrueEpiHead by 100%
\divide\TrueEpiHead by \SCALE}

\newcount\monotaiL%
\global\monotaiL=\headsize%
\global\multiply\monotaiL by 3162%
\global\divide\monotaiL by 10000%

\newcount\epiheaD%
\global\epiheaD=\epihead%
\global\multiply\epiheaD by 4472%
\global\divide\epiheaD by 10000%

\newcount\truemonotaiL%
\newcount\trueepiheaD%

\def\truetaiL{\truemonotaiL=\monotaiL%
\multiply\truemonotaiL by 100%
\divide\truemonotaiL by \SCALE}%

\def\trueheaD{\trueepiheaD=\epiheaD%
\multiply\trueepiheaD by 100%
\divide\trueepiheaD by \SCALE}

\newcount\SCALE%

\newcount\NUMBER%

\newcount\LINE%

\newcount\COLUMN%

\newcount\WIDTH%

\newcount\SOURCE%

\newcount\ARROW%

\newcount\TARGET%

\newcount\ARROWLENGTH%

\newcount\NUMBEROFDOTS%

\newcounter{x}%

\newcounter{y}%

\newcounter{z}%

\newcounter{horizontal}%

\newcounter{vertical}%

\newskip\itemlength%

\newskip\firstitem%

\newskip\seconditem%

\newcommand{\printarrow}{}%

\newcount\X%

\newcount\Y%

\newcount\Z%

\newcommand{\truex}[1]{%
\NUMBER=#1%
\multiply\NUMBER by 100%
\divide\NUMBER by \SCALE%
\setcounter{x}{\NUMBER}}%

\newcommand{\truey}[1]{%
\NUMBER=#1%
\multiply\NUMBER by 100%
\divide\NUMBER by \SCALE%
\setcounter{y}{\NUMBER}}%

\newcommand{\truez}[1]{%
\NUMBER=#1%
\multiply\NUMBER by 100%
\divide\NUMBER by \SCALE%
\setcounter{z}{\NUMBER}}%

\newcommand{\changecounters}[1]{%
\SOURCE=\ARROW%
\ARROW=\TARGET%
\settowidth{\itemlength}{#1}%
\ifdim \itemlength > 2800\unitlength%
\addtolength{\itemlength}{-2800\unitlength}%
\TARGET=\itemlength%
\divide\TARGET by 1310%
\multiply\TARGET by 100%
\divide\TARGET by \SCALE%
\else%
\TARGET=0%
\fi%
\ARROWLENGTH=5000%
\advance\ARROWLENGTH by -\SOURCE%
\advance\ARROWLENGTH by -\TARGET%
\divide\ARROWLENGTH by 100%
\advance\SOURCE by -\TARGET}%

\newcommand{\initialize}[1]{%
\LINE=0%
\COLUMN=0%
\WIDTH=0%
\ARROW=0%
\TARGET=0%
\changecounters{#1}%
\renewcommand{\printarrow}{#1}%
\begin{center}%
\vspace{2pt}%
\begin{picture}(0,0)}%

\newcommand{\DIAGV}[2]{%
\SCALE=#1%
\setlength{\unitlength}{655sp}%
\multiply\unitlength by \SCALE%
\divide\unitlength by 100%
\initialize{\mbox{$#2$}}}%

\newcommand{\n}[1]{%
\changecounters{\mbox{$#1$}}%
\put(\COLUMN,\LINE){\makebox(0,0){\printarrow}}%
\thinlines%
\renewcommand{\printarrow}{\mbox{$#1$}}%
\advance\COLUMN by 4000}%

\newcommand{\nn}[1]{%
\put(\COLUMN,\LINE){\makebox(0,0){\printarrow}}%
\thinlines%
\ifnum \WIDTH < \COLUMN%
\WIDTH=\COLUMN%
\else%
\fi%
\advance\LINE by -4000%
\COLUMN=0%
\ARROW=0%
\TARGET=0%
\changecounters{\mbox{$#1$}}%
\renewcommand{\printarrow}{\mbox{$#1$}}}%

\newcommand{\conclude}{%
\put(\COLUMN,\LINE){\makebox(0,0){\printarrow}}%
\thinlines%
\ifnum \WIDTH < \COLUMN%
\WIDTH=\COLUMN%
\else%
\fi%
\setcounter{horizontal}{\WIDTH}%
\setcounter{vertical}{-\LINE}%
\end{picture}}%

\newcommand{\diag}{%
\conclude%
\raisebox{0pt}[0pt][\value{vertical}\unitlength]{}%
\hspace*{\value{horizontal}\unitlength}%
\vspace{12pt}%
\end{center}%
\setlength{\unitlength}{1pt}%
}%

\newcommand{\diagv}[3]{%
\conclude%
\NUMBER=#1%
\rule{0pt}{\NUMBER pt}%
\hspace*{-#2pt}%
\raisebox{0pt}[0pt][\value{vertical}\unitlength]{}%
\hspace*{\value{horizontal}\unitlength}
\NUMBER=#3%
\advance\NUMBER by 12%
\vspace*{\NUMBER pt}%
\end{center}%
\setlength{\unitlength}{1pt}%
}%

\def\movename(#1,#2)#3{%
\hspace{#1pt}%
\raisebox{#2pt}[5pt][2pt]{\raisebox{#2pt}{$#3$}}%
\hspace{-#1pt}}%

\def\movearrow(#1,#2)#3{%
\makebox[0pt]{%
\hspace{#1pt}\hspace{#1pt}%
\raisebox{#2pt}[0pt][0pt]{\raisebox{#2pt}{$#3$}}}}%

\def\movevertex(#1,#2)#3{%
\mbox{\hspace{#1pt}%
\raisebox{#2pt}{\raisebox{#2pt}{$#3$}}%
\hspace{-#1pt}}}%

\newcommand{\movevertexleft}[1]{\makebox[2800\unitlength][r]{$#1$}}

\newcommand{\movevertexright}[1]{\makebox[2800\unitlength][l]{$#1$}}

\newcommand{\CrossLength}[2]{%
\settowidth{\firstitem}{#1}%
\settowidth{\seconditem}{#2}%
\ifdim\firstitem < \seconditem%
\itemlength=\seconditem%
\else%
\itemlength=\firstitem%
\fi%
\divide\itemlength by 2%
\hspace{\itemlength}}%

\def\basicDIAG#1|{\DIAGV{\defaultscale}{#1}\@ifnexttok|{\finishline}{\basicn}}
\def\basicDIAGV[#1]#2|{\DIAGV{#1}{#2}\@ifnexttok|{\finishline}{\basicn}}
\def\basicn#1|{\n{#1}\@ifnexttok|{\finishline}{\basicn}}
\def\basicnn#1|{\nn{#1}\@ifnexttok|{\finishline}{\basicn}}
\def\finishline#1{\@ifnextchar\end{\diag}%
{\@ifnextchar\spacing{\relax}{\basicnn}}}
\def\spacing(#1,#2,#3){\diagv{#1}{#2}{#3}}

\newif\ifcaption%

\newenvironment{diagram}{%
\iffloatdiag\relax\else
\global\def\diagramcaption##1{%
\global\captiontrue%
\global\def\@diagcaption{##1}}%
\global\def\@diagcaption{}\fi%
\@ifnextchar[{\basicDIAGV}{\basicDIAG}}%
{\iffloatdiag\relax\else%
\ifcaption
\begin{center}\mbox{}\@diagcaption\end{center}%
\else\relax\fi\fi\global\captionfalse}

\global\captionfalse

\gdef\@diaglabel{Diagram}
\gdef\diagramlabel#1{\gdef\@diaglabel{#1}}

\newcounter{Diagram}
\def\theDiagram{\@arabic\c@Diagram}
\def\fps@Diagram{tpbh}
\def\ftype@Diagram{1}
\def\ext@Diagram{lof}
\def\fnum@Diagram{\@diaglabel\ \theDiagram}
\def\Diagram{\@float{Diagram}}
\let\endDiagram\end@float
\@namedef{Diagram*}{\@dblfloat{Diagram}}
\@namedef{endDiagram*}{\end@dblfloat}

\newif\iffloatdiag

{\end{diagram}\caption{\@diagcaption}\end{Diagram}%
\global\floatdiagfalse}

{\end{diagram}\caption{\@diagcaption}\end{Diagram}%
\global\floatdiagfalse%
\typeout{*****}%
\typeout{WARNING: Floating diagram forced to be on the page}%
\typeout{*****}
}

{\end{diagram}\caption{\@diagcaption}\end{Diagram}%
\global\floatdiagfalse%
}

\global\floatdiagfalse

\def\dlabel#1{\immediate\write\@auxout{\string\newlabel{#1}%
{{\theChapter.\@arabic\c@Diagram}{\thepage}}}}

\newcommand{\TUP}[1]{\raisebox{0pt}[0pt][3pt]{}#1}

\newcommand{\TDOWN}[1]{\raisebox{0pt}[6pt][0pt]{}#1}

\newcommand{\tlowername}[2]%
{$\stackrel{\makebox[1pt]{#1}}%
{\begin{picture}(0,0)%
\put(0,0){\makebox(0,6)[t]{\makebox[1pt]{$\scriptstyle#2$}}}%
\end{picture}}$}%

\newcommand{\tcase}[1]{%
\setlength{\unitlength}{0.01pt}%
\makebox[\actualtextarrowspace pt]%
{\raisebox{2.5pt}{#1{\actualtextarrowlength}}}%
\setlength{\unitlength}{1pt}}%

\newcommand{\Tcase}[2]{%
\setlength{\unitlength}{0.01pt}%
\makebox[\actualtextarrowspace pt]%
{\raisebox{2.5pt}{$\stackrel{\scriptstyle #2}{#1{\actualtextarrowlength}}$}}%
\setlength{\unitlength}{1pt}}%

\newcommand{\tbicase}[1]{%
\setlength{\unitlength}{0.01pt}%
\makebox[\actualtextarrowspace pt]%
{\raisebox{1pt}{#1{\actualtextarrowlength}}}%
\setlength{\unitlength}{1pt}}%

\newcommand{\Tbicase}[3]{%
\setlength{\unitlength}{0.01pt}%
\makebox[\actualtextarrowspace pt]%
{\raisebox{-1pt}%
{$\stackrel{\scriptstyle #2}%
{\mbox{\tlowername{#1{\actualtextarrowlength}}%
{\scriptstyle #3}}}$}}%
\setlength{\unitlength}{1pt}}%

\newcommand{\ttricase}[1]{%
\setlength{\unitlength}{0.01pt}%
\makebox[\actualtextarrowspace pt]%
{\raisebox{4pt}{#1{\actualtextarrowlength}%
{\rule{0pt}{6pt}}{\rule{0pt}{6pt}}{\rule{0pt}{6pt}}}}%
\setlength{\unitlength}{1pt}}%

\newcommand{\Ttricase}[4]{%
\setlength{\unitlength}{0.01pt}%
\makebox[\actualtextarrowspace pt]%
{\raisebox{4pt}%
{#1{\actualtextarrowlength}{\rule{0pt}{6pt}\scriptstyle #2}%
{\rule{0pt}{6pt}\scriptstyle #3}{\rule{0pt}{6pt}\scriptstyle #4}}}%
\setlength{\unitlength}{1pt}}%

\newcommand{\tmulticase}[1]{%
\setlength{\unitlength}{0.01pt}%
\makebox[\actualtextarrowspace pt]%
{\raisebox{1pt}{#1{\actualtextarrowlength}}}%
\setlength{\unitlength}{1pt}}%

\newcommand{\DUP}[1]{\raisebox{0pt}[0pt][4pt]{}#1}

\newcommand{\DDOWN}[1]{\raisebox{0pt}[9pt][0pt]{}#1}

\newcommand{\dlowername}[2]%
{$\stackrel{\makebox[1pt]{#1}}%
{\begin{picture}(0,0)%
\put(0,0){\makebox(0,6)[t]{\makebox[1pt]{$\textstyle#2$}}}%
\end{picture}}$}%

\newcommand{\dcase}[1]{%
\setlength{\unitlength}{0.01pt}%
\makebox[\actualdisplayarrowspace pt]%
{\raisebox{2.5pt}{#1{\actualdisplayarrowlength}}}%
\setlength{\unitlength}{1pt}}%

\newcommand{\Dcase}[2]{%
\setlength{\unitlength}{0.01pt}%
\makebox[\actualdisplayarrowspace pt]%
{\raisebox{2.5pt}{$\stackrel{\textstyle #2}{#1{\actualdisplayarrowlength}}$}}%
\setlength{\unitlength}{1pt}}%

\newcommand{\dbicase}[1]{%
\setlength{\unitlength}{0.01pt}%
\makebox[\actualdisplayarrowspace pt]%
{\raisebox{1pt}{#1{\actualdisplayarrowlength}}}%
\setlength{\unitlength}{1pt}}%

\newcommand{\Dbicase}[3]{%
\setlength{\unitlength}{0.01pt}%
\makebox[\actualdisplayarrowspace pt]%
{\raisebox{-1pt}%
{$\stackrel{\textstyle #2}%
{\mbox{\tlowername{#1{\actualdisplayarrowlength}}%
{\textstyle #3}}}$}}%
\setlength{\unitlength}{1pt}}%

\newcommand{\dtricase}[1]{%
\setlength{\unitlength}{0.01pt}%
\makebox[\actualdisplayarrowspace pt]%
{\raisebox{4pt}{#1{\actualdisplayarrowlength}%
{\rule{0pt}{6pt}}{\rule{0pt}{6pt}}{\rule{0pt}{6pt}}}}%
\setlength{\unitlength}{1pt}}%

\newcommand{\Dtricase}[4]{%
\setlength{\unitlength}{0.01pt}%
\makebox[\actualdisplayarrowspace pt]%
{\raisebox{4pt}%
{#1{\actualdisplayarrowlength}{\rule{0pt}{6pt}\textstyle #2}%
{\rule{0pt}{6pt}\textstyle #3}{\rule{0pt}{6pt}\textstyle #4}}}%
\setlength{\unitlength}{1pt}}%

\newcommand{\dmulticase}[1]{%
\setlength{\unitlength}{0.01pt}%
\makebox[\actualdisplayarrowspace pt]%
{\raisebox{1pt}{#1{\actualdisplayarrowlength}}}%
\setlength{\unitlength}{1pt}}%

\newcommand{\AR}[1]%
{\begin{picture}(#1,0)%
\put(0,0){\line(1,0){#1}}%
\put(#1,0){\ehead}%
\end{picture}}%

\newcommand{\DIST}[1]%
{\begin{picture}(#1,0)%
\put(0,0){\line(1,0){#1}}%
\put(#1,0){\ehead}%
\NUMBER=#1%
\divide\NUMBER by 2%
\put(\NUMBER,0){\distsign}%
\end{picture}}%

\newcommand{\DOTAR}[1]%
{\NUMBEROFDOTS=#1%
\divide\NUMBEROFDOTS by 300%
\advance\NUMBEROFDOTS by 1%
\begin{picture}(#1,0)%
\multiput(0,0)(300,0){\NUMBEROFDOTS}{\circle*{100}}%
\put(#1,0){\ehead}%
\end{picture}}%

\newcommand{\MONO}[1]%
{\monolength=#1%
\advance\monolength by -\monotail%
\begin{picture}(#1,0)%
\put(\monotail,0){\line(1,0){\monolength}}%
\put(#1,0){\ehead}%
\put(\monotail,0){\ehead}%
\end{picture}}%

\newcommand{\EPI}[1]%
{\epilength=#1%
\advance\epilength by -\epihead%
\begin{picture}(#1,0)(-#1,0)%
\put(-#1,0){\line(1,0){\epilength}}%
\put(-\epihead,0){\ehead}%
\put(0,0){\ehead}%
\end{picture}}%

\newcommand{\BIMO}[1]%
{\monolength=#1%
\advance\monolength by -\monotail%
\epilength=\monolength%
\advance\epilength by -\epihead%
\begin{picture}(#1,0)(-#1,0)%
\put(-\monolength,0){\line(1,0){\epilength}}%
\put(-\monolength,0){\ehead}%
\put(-\epihead,0){\ehead}%
\put(0,0){\ehead}%
\end{picture}}%

\newcommand{\BIAR}[1]%
{\begin{picture}(#1,700)%
\put(0,0){\line(1,0){#1}}%
\put(#1,0){\ehead}%
\put(0,700){\line(1,0){#1}}%
\put(#1,700){\ehead}%
\end{picture}}%

\newcommand{\BIDIST}[1]%
{\begin{picture}(#1,700)%
\put(0,0){\line(1,0){#1}}%
\put(#1,0){\ehead}%
\put(0,700){\line(1,0){#1}}%
\put(#1,700){\ehead}%
\NUMBER=#1%
\divide\NUMBER by 2%
\put(\NUMBER,0){\distsign}%
\put(\NUMBER,700){\distsign}%
\end{picture}}%

\newcommand{\EQL}[1]%
{\begin{picture}(#1,0)%
\put(0,100){\line(1,0){#1}}%
\put(0,-100){\line(1,0){#1}}%
\end{picture}}%

\newcommand{\LLINE}[1]%
{\begin{picture}(#1,0)%
\put(0,0){\line(1,0){#1}}%
\end{picture}}%

\newcommand{\ADJAR}[1]%
{\begin{picture}(#1,700)%
\put(0,0){\line(1,0){#1}}%
\put(#1,0){\ehead}%
\put(#1,700){\line(-1,0){#1}}%
\put(0,700){\whead}
\end{picture}}%

\newcommand{\ADJDIST}[1]%
{\begin{picture}(#1,700)%
\put(0,0){\line(1,0){#1}}%
\put(#1,0){\ehead}%
\put(#1,700){\line(-1,0){#1}}%
\put(0,700){\whead}
\NUMBER=#1%
\divide\NUMBER by 2%
\put(\NUMBER,0){\distsign}%
\put(\NUMBER,700){\distsign}%
\end{picture}}%

\newcommand{\TRIAR}[4]%
{%
\begin{tabular}{c}%
\mbox{$#2$}\\ \AR{#1}\\%
\mbox{$#3$}\\ \AR{#1}\\%
\mbox{$#4$}\\ \AR{#1}%
\end{tabular}%
}%

\newcommand{\TRIADJAR}[4]%
{%
\begin{tabular}{c}%
\mbox{$#2$}\\ \BKAR{#1}\\%
\mbox{$#3$}\\ \AR{#1}\\%
\mbox{$#4$}\\ \BKAR{#1}%
\end{tabular}%
}%

\newcommand{\QUADRIAR}[1]%
{%
\begin{tabular}{c}%
\AR{#1}\\%
\rule{0pt}{3pt}\\ \AR{#1}\\%
\rule{0pt}{3pt}\\ \AR{#1}\\%
\rule{0pt}{3pt}\\ \AR{#1}%
\end{tabular}%
}%

\newcommand{\QUADRIADJAR}[1]%
{%
\begin{tabular}{c}%
\BKAR{#1}\\%
\rule{0pt}{3pt}\\ \AR{#1}\\%
\rule{0pt}{3pt}\\ \BKAR{#1}\\%
\rule{0pt}{3pt}\\ \AR{#1}%
\end{tabular}%
}%

\newcommand{\QUINTIAR}[1]%
{%
\begin{tabular}{c}%
\AR{#1}\\%
\rule{0pt}{2pt}\\ \AR{#1}\\%
\rule{0pt}{2pt}\\ \AR{#1}\\%
\rule{0pt}{2pt}\\ \AR{#1}\\%
\rule{0pt}{2pt}\\ \AR{#1}%
\end{tabular%
}%

\newcommand{\QUINTIADJAR}[1]%
{%
\begin{tabular}{c}%
\AR{#1}\\%
\rule{0pt}{2pt}\\ \BKAR{#1}\\%
\rule{0pt}{2pt}\\ \AR{#1}\\%
\rule{0pt}{2pt}\\ \BKAR{#1}\\%
\rule{0pt}{2pt}\\ \AR{#1}%
\end{tabular}}%
}%

\newcommand{\ar}{\ifinner\tcase{\AR}\else\dcase{\AR}\fi}%

\newcommand{\Ar}[1]{\ifinner\Tcase{\AR}{#1}\else\Dcase{\AR}{#1}\fi}%

\newcommand{\dist}{\ifinner\tcase{\DIST}\else\dcase{\DIST}\fi}%

\newcommand{\Dist}[1]{\ifinner\Tcase{\DIST}{\TUP{#1}}%
\else\Dcase{\DIST}{\TUP{#1}}\fi}%

\newcommand{\dotar}{\ifinner\tcase{\DOTAR}\else\dcase{\DOTAR}\fi}%

\newcommand{\Dotar}[1]{\ifinner\Tcase{\DOTAR}{#1}%
\else\Dcase{\DOTAR}{#1}\fi}%

\newcommand{\mono}{\ifinner\tcase{\MONO}\else\dcase{\MONO}\fi}%

\newcommand{\Mono}[1]{\ifinner\Tcase{\MONO}{#1}\else\Dcase{\MONO}{#1}\fi}%

\newcommand{\epi}{\ifinner\tcase{\EPI}\else\dcase{\EPI}\fi}%

\newcommand{\Epi}[1]{\ifinner\Tcase{\EPI}{#1}\else\Dcase{\EPI}{#1}\fi}%

\newcommand{\bimo}{\ifinner\tcase{\BIMO}\else\dcase{\BIMO}\fi}%

\newcommand{\Bimo}[1]{\ifinner\Tcase{\BIMO}{#1}%
\else\Dcase{\BIMO}{#1}\fi}%

\newcommand{\iso}{\ifinner\Tcase{\AR}{\isosign}\else\Dcase{\AR}{\isosign}\fi}%

\newcommand{\Iso}[1]{\ifinner\Tcase{\AR}{\isosign{#1}}%
\else\Dcase{\AR}{\isosign{#1}}\fi}%

\newcommand{\biar}{\ifinner\tbicase{\BIAR}\else\dbicase{\BIAR}\fi}%

\newcommand{\Biar}[2]{\ifinner\Tbicase{\BIAR}{#1}{#2}%
\else\Dbicase{\BIAR}{#1}{#2}\fi}%

\newcommand{\bidist}{\ifinner\tbicase{\BIDIST}\else\dbicase{\BIDIST}\fi}%

\newcommand{\Bidist}[2]{\ifinner\Tbicase{\BIDIST}{\TUP{#1}}{\TDOWN{#2}}%
\else\Dbicase{\BIDIST}{\DUP{#1}}{\DDOWN{#2}}\fi}%

\newcommand{\eql}{\ifinner\tcase{\EQL}\else\dcase{\EQL}\fi}%

\newcommand{\Eql}[1]{\ifinner\Tcase{\EQL}{\TUP{#1}}%
\else\Dcase{\EQL}{\DUP{#1}}\fi}%

\newcommand{\lline}{\ifinner\tcase{\LLINE}\else\dcase{\LLINE}\fi}%

\newcommand{\Line}[1]{\ifinner\Tcase{\LLINE}{\TUP{#1}}%
\else\Dcase{\LLINE}{\DUP{#1}}\fi}%

\newcommand{\adjar}{\ifinner\tbicase{\ADJAR}\else\dbicase{\ADJAR}\fi}%

\newcommand{\Adjar}[2]{\ifinner\Tbicase{\ADJAR}{#1}{#2}%
\else\Dbicase{\ADJAR}{#1}{#2}\fi}%

\newcommand{\adjdist}{\ifinner\tbicase{\ADJDIST}\else\dbicase{\ADJDIST}\fi}%

\newcommand{\Adjdist}[2]{\ifinner\Tbicase{\ADJDIST}{\TUP{#1}}{\TDOWN{#2}}%
\else\Dbicase{\ADJDIST}{\DUP{#1}}{\DDOWN{#2}}\fi}%

\newcommand{\triar}{\ifinner\ttricase{\TRIAR}\else\dtricase{\TRIAR}\fi}%

\newcommand{\Triar}[3]{\ifinner\Ttricase{\TRIAR}{#1}{#2}{#3}%
\else\Dtricase{\TRIAR}{#1}{#2}{#3}\fi}%

\newcommand{\triadjar}{\ifinner\ttricase{\TRIADJAR}%
\else\dtricase{\TRIADJAR}\fi}%

\newcommand{\Triadjar}[3]{\ifinner\Ttricase{\TRIADJAR}{#1}{#2}{#3}%
\else\Dtricase{\TRIADJAR}{#1}{#2}{#3}\fi}%

\newcommand{\quadriar}{\ifinner\tquadricase{\QUADRIAR}%
\else\dmulticase{\QUADRIAR}\fi}%

\newcommand{\quadriadjar}{\ifinner\tmulticase{\QUADRIADJAR}%
\else\dmulticase{\QUADRIADJAR}\fi}%

\newcommand{\quintiar}{\ifinner\tmulticase{\QUINTIAR}%
\else\dmulticase{\QUINTIAR}\fi}%

\newcommand{\quintiadjar}{\ifinner\tmulticase{\QUINTIADJAR}%
\else\dmulticase{\QUINTIADJAR}\fi}%

\newcommand{\BKAR}[1]%
{\begin{picture}(#1,0)%
\put(#1,0){\line(-1,0){#1}}%
\put(0,0){\whead}%
\end{picture}}%

\newcommand{\BKDIST}[1]%
{\begin{picture}(#1,0)%
\put(#1,0){\line(-1,0){#1}}%
\put(0,0){\whead}%
\NUMBER=#1%
\divide\NUMBER by 2%
\put(\NUMBER,0){\distsign}%
\end{picture}}%

\newcommand{\BKDOTAR}[1]%
{\NUMBEROFDOTS=#1%
\divide\NUMBEROFDOTS by 300%
\advance\NUMBEROFDOTS by 1%
\begin{picture}(#1,0)%
\multiput(#1,0)(-300,0){\NUMBEROFDOTS}{\circle*{100}}%
\put(0,0){\whead}%
\end{picture}}%

\newcommand{\BKMONO}[1]%
{\monolength=#1%
\advance\monolength by -\monotail%
\begin{picture}(#1,0)(-#1,0)%
\put(-\monotail,0){\line(-1,0){\monolength}}%
\put(-\monotail,0){\whead}%
\put(-#1,0){\whead}%
\end{picture}}%

\newcommand{\BKEPI}[1]%
{\epilength=#1%
\advance\epilength by -\epihead%
\begin{picture}(#1,0)%
\put(#1,0){\line(-1,0){\epilength}}%
\put(\epihead,0){\whead}%
\put(0,0){\whead}%
\end{picture}}%

\newcommand{\BKBIMO}[1]%
{\monolength=#1%
\advance\monolength by -\monotail%
\epilength=\monolength%
\advance\epilength by -\epihead%
\begin{picture}(#1,0)%
\put(\monolength,0){\line(-1,0){\epilength}}%
\put(\monolength,0){\whead}%
\put(\epihead,0){\whead}%
\put(0,0){\whead}%
\end{picture}}%

\newcommand{\BKBIAR}[1]%
{\begin{picture}(#1,700)%
\put(#1,0){\line(-1,0){#1}}%
\put(0,0){\whead}%
\put(#1,700){\line(-1,0){#1}}%
\put(0,700){\whead}%
\end{picture}}%

\newcommand{\BKBIDIST}[1]%
{\begin{picture}(#1,700)%
\put(#1,0){\line(-1,0){#1}}%
\put(0,0){\whead}%
\put(#1,700){\line(-1,0){#1}}%
\put(0,700){\whead}%
\NUMBER=#1%
\divide\NUMBER by 2%
\put(\NUMBER,0){\distsign}%
\put(\NUMBER,700){\distsign}%
\end{picture}}%

\newcommand{\BKADJAR}[1]%
{\begin{picture}(#1,700)%
\put(0,700){\line(1,0){#1}}%
\put(#1,700){\ehead}%
\put(#1,0){\line(-1,0){#1}}%
\put(0,0){\whead}%
\end{picture}}%

\newcommand{\BKADJDIST}[1]%
{\begin{picture}(#1,700)%
\put(0,700){\line(1,0){#1}}%
\put(#1,700){\ehead}%
\put(#1,0){\line(-1,0){#1}}%
\put(0,0){\whead}%
\NUMBER=#1%
\divide\NUMBER by 2%
\put(\NUMBER,0){\distsign}%
\put(\NUMBER,700){\distsign}%
\end{picture}}%

\newcommand{\BKTRIAR}[4]%
{%
\begin{tabular}{c}%
\mbox{$#2$}\\ \BKAR{#1}\\%
\mbox{$#3$}\\ \BKAR{#1}\\%
\mbox{$#4$}\\ \BKAR{#1}%
\end{tabular}%
}%

\newcommand{\BKTRIADJAR}[4]%
{%
\begin{tabular}{c}%
\mbox{$#2$}\\ \AR{#1}\\%
\mbox{$#3$}\\ \BKAR{#1}\\%
\mbox{$#4$}\\ \AR{#1}%
\end{tabular}%
}%

\newcommand{\BKQUADRIAR}[1]%
{%
\begin{tabular}{c}%
\BKAR{#1}\\%
\rule{0pt}{3pt}\\ \BKAR{#1}\\%
\rule{0pt}{3pt}\\ \BKAR{#1}\\%
\rule{0pt}{3pt}\\ \BKAR{#1}%
\end{tabular}%
}%

\newcommand{\BKQUADRIADJAR}[1]%
{%
\begin{tabular}{c}%
\AR{#1}\\%
\rule{0pt}{3pt}\\ \BKAR{#1}\\%
\rule{0pt}{3pt}\\ \AR{#1}\\%
\rule{0pt}{3pt}\\ \BKAR{#1}%
\end{tabular}%
}%

\newcommand{\BKQUINTIAR}[1]%
{%
\begin{tabular}{c}%
\BKAR{#1}\\%
\rule{0pt}{2pt}\\ \BKAR{#1}\\%
\rule{0pt}{2pt}\\ \BKAR{#1}\\%
\rule{0pt}{2pt}\\ \BKAR{#1}\\%
\rule{0pt}{2pt}\\ \BKAR{#1}%
\end{tabular}%
}%

\newcommand{\BKQUINTIADJAR}[1]%
{%
\begin{tabular}{c}%
\BKAR{#1}\\%
\rule{0pt}{2pt}\\ \AR{#1}\\%
\rule{0pt}{2pt}\\ \BKAR{#1}\\%
\rule{0pt}{2pt}\\ \AR{#1}\\%
\rule{0pt}{2pt}\\ \BKAR{#1}%
\end{tabular}%
}%

\newcommand{\bkar}{\ifinner\tcase{\BKAR}\else\dcase{\BKAR}\fi}%

\newcommand{\Bkar}[1]{\ifinner\Tcase{\BKAR}{#1}\else\Dcase{\BKAR}{#1}\fi}%

\newcommand{\bkdist}{\ifinner\tcase{\BKDIST}\else\dcase{\BKDIST}\fi}%

\newcommand{\Bkdist}[1]{\ifinner\Tcase{\BKDIST}{\TUP{#1}}%
\else\Dcase{\BKDIST}{\TUP{#1}}\fi}%

\newcommand{\bkdotar}{\ifinner\tcase{\BKDOTAR}\else\dcase{\BKDOTAR}\fi}%

\newcommand{\Bkdotar}[1]{\ifinner\Tcase{\BKDOTAR}{#1}%
\else\Dcase{\BKDOTAR}{#1}\fi}%

\newcommand{\bkmono}{\ifinner\tcase{\BKMONO}\else\dcase{\BKMONO}\fi}%

\newcommand{\Bkmono}[1]{\ifinner\Tcase{\BKMONO}{#1}%
\else\Dcase{\BKMONO}{#1}\fi}%

\newcommand{\bkepi}{\ifinner\tcase{\BKEPI}\else\dcase{\BKEPI}\fi}%

\newcommand{\Bkepi}[1]{\ifinner\Tcase{\BKEPI}{#1}%
\else\Dcase{\BKEPI}{#1}\fi}%

\newcommand{\bkbimo}{\ifinner\tcase{\BKBIMO}\else\dcase{\BKBIMO}\fi}%

\newcommand{\Bkbimo}[1]{\ifinner\Tcase{\BKBIMO}{\hspace{9pt}#1}%
\else\Dcase{\BKBIMO}{\hspace{9pt}#1}\fi}%

\newcommand{\bkiso}{\ifinner\Tcase{\BKAR}{\isosign}%
\else\Dcase{\BKAR}{\isosign}\fi}%

\newcommand{\Bkiso}[1]{\ifinner\Tcase{\BKAR}{\isosign{#1}}%
\else\Dcase{\BKAR}{\isosign{#1}}\fi}%

\newcommand{\bkbiar}{\ifinner\tbicase{\BKBIAR}\else\dbicase{\BKBIAR}\fi}%

\newcommand{\Bkbiar}[2]{\ifinner\Tbicase{\BKBIAR}{#1}{#2}%
\else\Dbicase{\BKBIAR}{#1}{#2}\fi}%

\newcommand{\bkbidist}{\ifinner\tbicase{\BKBIDIST}%
\else\dbicase{\BKBIDIST}\fi}%

\newcommand{\Bkbidist}[2]{\ifinner\Tbicase{\BKBIDIST}{\TUP{#1}}{\TDOWN{#2}}%
\else\Tbicase{\BKBIDIST}{\DUP{#1}}{\DDOWN{#2}}\fi}%

\newcommand{\bkadjar}{\ifinner\tbicase{\BKADJAR}%
\else\dbicase{\BKADJAR}\fi}%

\newcommand{\Bkadjar}[2]{\ifinner\Tbicase{\BKADJAR}{#1}{#2}%
\else\Dbicase{\BKADJAR}{#1}{#2}\fi}%

\newcommand{\bkadjdist}{\ifinner\tbicase{\BKADJDIST}%
\else\dbicase{\BKADJDIST}\fi}%

\newcommand{\Bkadjdist}[2]{\ifinner\Tbicase{\BKADJDIST}{\TUP{#1}}{\TDOWN{#2}}%
\else\Dbicase{\BKADJDIST}{\TUP{#1}}{\TDOWN{#2}}\fi}%

\newcommand{\bktriar}{\ifinner\ttricase{\BKTRIAR}\else\dtricase{\BKTRIAR}\fi}%

\newcommand{\Bktriar}[3]{\ifinner\Ttricase{\BKTRIAR}{#1}{#2}{#3}%
\else\Dtricase{\BKTRIAR}{#1}{#2}{#3}\fi}%

\newcommand{\bktriadjar}{\ifinner\ttricase{\BKTRIADJAR}%
\else\dtricase{\BKTRIADJAR}\fi}%

\newcommand{\Bktriadjar}[3]{\ifinner\Ttricase{\BKTRIADJAR}{#1}{#2}{#3}%
\else\Dtricase{\BKTRIADJAR}{#1}{#2}{#3}\fi}%

\newcommand{\bkquadriar}{\ifinner\tquadricase{\BKQUADRIAR}%
\else\dmulticase{\BKQUADRIAR}\fi}%

\newcommand{\bkquadriadjar}{\ifinner\tmulticase{\BKQUADRIADJAR}%
\else\dmulticase{\BKQUADRIADJAR}\fi}%

\newcommand{\bkquintiar}{\ifinner\tmulticase{\BKQUINTIAR}%
\else\dmulticase{\BKQUINTIAR}\fi}%

\newcommand{\bkquintiadjar}{\ifinner\tmulticase{\BKQUINTIADJAR}%
\else\dmulticase{\BKQUINTIADJAR}\fi}%

\def\maketopbox#1{\vbox to 0\unitlength%
{\vfill\hbox to 0\unitlength{\hss#1\hss}\vss}}

\def\makebottombox#1{\vbox to 0\unitlength%
{\vss\hbox to 0\unitlength{\hss#1\hss}\vfill}}

\newcommand{\hcase}[2]%
{\makebox[0pt]%
{\raisebox{0pt}[0pt][0pt]%
{\begin{picture}(0,0)%
\put(0,0){\makebox(0,0){#1{#2}}}%
\end{picture}%
}}}%

\newcommand{\Hcase}[3]%
{\makebox[0pt]%
{\raisebox{0pt}[0pt][0pt]%
{\begin{picture}(0,0)%
\put(0,0){\makebox(0,0){#1{#3}}}%
\truex{200}%
\put(0,\value{x}){\makebottombox{$#2$}}%
\end{picture}%
}}}%

\newcommand{\hcasE}[3]%
{\makebox[0pt]%
{\raisebox{0pt}[0pt][0pt]%
{\begin{picture}(0,0)%
\put(0,0){\makebox(0,0){#1{#3}}}%
\truex{200}%
\put(0,-\value{x}){\maketopbox{$#2$}}%
\end{picture}%
}}}%

\newcommand{\Hisocase}[4]%
{\makebox[0pt]
{\raisebox{0pt}[0pt][0pt]%
{\begin{picture}(0,0)%
\put(0,0){\makebox(0,0){#1{#4}}}%
\truex{200}%
\put(0,\value{x}){\makebottombox{$#2$}}%
\put(0,-\value{x}){\maketopbox{$#3$}}%
\end{picture}%
}}}%

\newcommand{\hbicase}[2]%
{\makebox[0pt]%
{\raisebox{0pt}[0pt][0pt]%
{\begin{picture}(0,0)%
\put(0,0){\makebox(0,0){#1{#2}}}%
\end{picture}%
}}}%

\newcommand{\Hbicase}[4]%
{\makebox[0pt]
{\raisebox{0pt}[0pt][0pt]%
{\begin{picture}(0,0)%
\put(0,0){\makebox(0,0){#1{#4}}}%
\truex{700}%
\put(0,\value{x}){\makebottombox{$#2$}}%
\put(0,-\value{x}){\maketopbox{$#3$}}%
\end{picture}%
}}}%

\newcommand{\htricase}[2]%
{\makebox[0pt]%
{\raisebox{-12pt}[0pt][0pt]{#1{#2}{}{}{}{}}}}%

\newcommand{\Htricase}[6]%
{\makebox[0pt]%
{\raisebox{-12pt}[0pt][0pt]{#1{#6}{#2}{#3}{#4}{#5}}}}%

\newcommand{\hmulticase}[2]%
{\makebox[0pt]{\raisebox{-12pt}[0pt][0pt]{#1{#2}}}}%

\newcommand{\EAR}[1]%
{\begin{picture}(#1,0)%
\put(0,0){\line(1,0){#1}}%
\put(#1,0){\ehead}%
\end{picture}}%

\newcommand{\EDIST}[1]%
{\begin{picture}(#1,0)%
\put(0,0){\line(1,0){#1}}%
\put(#1,0){\ehead}%
\truex{400}
\NUMBER=#1%
\divide\NUMBER by 2%
\put(\NUMBER,0){\distsign}
\end{picture}}%

\newcommand{\EDOTAR}[1]%
{\truex{100}\truey{300}%
\NUMBEROFDOTS=#1%
\divide\NUMBEROFDOTS by \value{y}%
\advance\NUMBEROFDOTS by 1%
\begin{picture}(#1,0)%
\multiput(0,0)(\value{y},0){\NUMBEROFDOTS}%
{\circle*{\value{x}}}%
\put(#1,0){\ehead}%
\end{picture}}%

\newcommand{\EMONO}[1]%
{\truetail
\monolength=#1%
\advance\monolength by -\truemonotail%
\begin{picture}(#1,0)%
\put(\truemonotail,0){\line(1,0){\monolength}}%
\put(#1,0){\ehead}%
\put(\truemonotail,0){\ehead}%
\end{picture}}%

\newcommand{\EEPI}[1]%
{\truehead%
\epilength=#1%
\advance\epilength by -\trueepihead%
\begin{picture}(#1,0)(-#1,0)%
\put(-#1,0){\line(1,0){\epilength}}%
\put(-\trueepihead,0){\ehead}%
\put(0,0){\ehead}%
\end{picture}}%

\newcommand{\EBIMO}[1]%
{\truehead\truetail%
\monolength=#1%
\advance\monolength by -\truemonotail%
\epilength=\monolength%
\advance\epilength by -\trueepihead%
\begin{picture}(#1,0)(-#1,0)%
\put(-\monolength,0){\line(1,0){\epilength}}%
\put(-\monolength,0){\ehead}%
\put(-\trueepihead,0){\ehead}%
\put(0,0){\ehead}%
\end{picture}}%

\newcommand{\EBIAR}[1]%
{\truex{700}%
\begin{picture}(#1,\value{x})%
\put(0,0){\line(1,0){#1}}%
\put(#1,0){\ehead}%
\put(0,\value{x}){\line(1,0){#1}}%
\put(#1,\value{x}){\ehead}%
\end{picture}}%

\newcommand{\EBIDIST}[1]%
{\truex{700}%
\begin{picture}(#1,\value{x})%
\put(0,0){\line(1,0){#1}}%
\put(#1,0){\ehead}%
\put(0,\value{x}){\line(1,0){#1}}%
\put(#1,\value{x}){\ehead}%
\NUMBER=#1%
\divide\NUMBER by 2%
\put(\NUMBER,0){\distsign}
\put(\NUMBER,\value{x}){\distsign}%
\end{picture}}%

\newcommand{\EEQL}[1]%
{\begin{picture}(#1,0)%
\truex{200}%
\put(0,\value{x}){\line(1,0){#1}}%
\put(0,0){\line(1,0){#1}}%
\end{picture}}%

\newcommand{\ELINE}[1]%
{\begin{picture}(#1,0)%
\truex{100}%
\put(0,\value{x}){\line(1,0){#1}}%
\end{picture}}%

\newcommand{\EADJAR}[1]%
{\truex{700}%
\begin{picture}(#1,\value{x})%
\put(0,0){\line(1,0){#1}}%
\put(#1,0){\ehead}%
\put(#1,\value{x}){\line(-1,0){#1}}%
\put(0,\value{x}){\whead}%
\end{picture}}%

\newcommand{\EADJDIST}[1]%
{\truex{700}%
\begin{picture}(#1,\value{x})%
\put(0,0){\line(1,0){#1}}%
\put(#1,0){\ehead}%
\put(#1,\value{x}){\line(-1,0){#1}}%
\put(0,\value{x}){\whead}%
\NUMBER=#1%
\divide\NUMBER by 2%
\put(\NUMBER,0){\distsign}
\put(\NUMBER,\value{x}){\distsign}%
\end{picture}}%

\newcommand{\ETRIAR}[5]%
{\truex{1200}%
\X=\value{x}%
\multiply\X by 2%
\begin{picture}(#1,\X)%
\put(0,0){\line(1,0){#1}}%
\put(#1,0){\ehead}%
\put(0,\value{x}){\line(1,0){#1}}%
\put(#1,\value{x}){\ehead}%
\put(0,\X){\line(1,0){#1}}%
\put(#1,\X){\ehead}%
\X=#1%
\divide\X by 2%
\truex{-600}%
\put(\X,\value{x}){\makebox(0,0){$#5$}}%
\truex{600}%
\put(\X,\value{x}){\makebox(0,0){$#4$}}%
\truex{1800}%
\put(\X,\value{x}){\makebox(0,0){$#3$}}%
\truex{3000}%
\put(\X,\value{x}){\makebox(0,0){$#2$}}%
\end{picture}}%

\newcommand{\ETRIADJAR}[5]%
{\truex{1200}%
\X=\value{x}%
\multiply\X by 2%
\begin{picture}(#1,\X)%
\put(#1,0){\line(-1,0){#1}}%
\put(0,0){\whead}%
\put(0,\value{x}){\line(1,0){#1}}%
\put(#1,\value{x}){\ehead}%
\put(#1,\X){\line(-1,0){#1}}%
\put(0,\X){\whead}%
\X=#1%
\divide\X by 2%
\truex{-600}%
\put(\X,\value{x}){\makebox(0,0){$#5$}}%
\truex{600}%
\put(\X,\value{x}){\makebox(0,0){$#4$}}%
\truex{1800}%
\put(\X,\value{x}){\makebox(0,0){$#3$}}%
\truex{3000}%
\put(\X,\value{x}){\makebox(0,0){$#2$}}%
\end{picture}}%

\newcommand{\EQUADRIAR}[1]%
{\truex{800}%
\truey{1600}%
\X=\value{x}%
\multiply\X by 3%
\begin{picture}(#1,\X)%
\put(0,0){\line(1,0){#1}}%
\put(#1,0){\ehead}%
\put(0,\value{x}){\line(1,0){#1}}%
\put(#1,\value{x}){\ehead}%
\put(0,\value{y}){\line(1,0){#1}}%
\put(#1,\value{y}){\ehead}%
\put(0,\X){\line(1,0){#1}}%
\put(#1,\X){\ehead}%
\end{picture}}%

\newcommand{\EQUADRIADJAR}[1]%
{\truex{800}%
\truey{1600}%
\X=\value{x}%
\multiply\X by 3%
\begin{picture}(#1,\X)%
\put(0,0){\line(1,0){#1}}%
\put(#1,0){\ehead}%
\put(#1,\value{x}){\line(-1,0){#1}}%
\put(0,\value{x}){\whead}%
\put(0,\value{y}){\line(1,0){#1}}%
\put(#1,\value{y}){\ehead}%
\put(#1,\X){\line(-1,0){#1}}%
\put(0,\X){\whead}%
\end{picture}}%

\newcommand{\EQUINTIAR}[1]%
{\truex{600}%
\truey{1200}%
\truez{1800}%
\X=\value{x}%
\multiply\X by 4%
\begin{picture}(#1,\X)%
\put(0,0){\line(1,0){#1}}%
\put(#1,0){\ehead}%
\put(0,\value{x}){\line(1,0){#1}}%
\put(#1,\value{x}){\ehead}%
\put(0,\value{y}){\line(1,0){#1}}%
\put(#1,\value{y}){\ehead}%
\put(0,\value{z}){\line(1,0){#1}}%
\put(#1,\value{z}){\ehead}%
\put(0,\X){\line(1,0){#1}}%
\put(#1,\X){\ehead}%
\end{picture}}%

\newcommand{\EQUINTIADJAR}[1]%
{\truex{600}%
\truey{1200}%
\truez{1800}%
\X=\value{x}%
\multiply\X by 4%
\begin{picture}(#1,\X)%
\put(0,0){\line(1,0){#1}}%
\put(#1,0){\ehead}%
\put(#1,\value{x}){\line(-1,0){#1}}%
\put(0,\value{x}){\whead}%
\put(0,\value{y}){\line(1,0){#1}}%
\put(#1,\value{y}){\ehead}%
\put(#1,\value{z}){\line(-1,0){#1}}%
\put(0,\value{z}){\whead}%
\put(0,\X){\line(1,0){#1}}%
\put(#1,\X){\ehead}%
\end{picture}}%

\def\basicear[#1]{%
\Z=#1%
\multiply \Z by 100%
\hcase{\EAR}{\Z}}%

\newcommand{\ear}{\@ifnextchar[{\basicear}%
{\hspace{\SOURCE\unitlength}\basicear[\ARROWLENGTH]}}%

\def\basicEar[#1]#2{%
\Z=#1%
\multiply \Z by 100%
\Hcase{\EAR}{#2}{\Z}}%

\newcommand{\Ear}{\@ifnextchar[{\basicEar}%
{\hspace{\SOURCE\unitlength}\basicEar[\ARROWLENGTH]}}%

\def\basiceaR[#1]#2{%
\Z=#1%
\multiply \Z by 100%
\hcasE{\EAR}{#2}{\Z}}%

\newcommand{\eaR}{\@ifnextchar[{\basiceaR}%
{\hspace{\SOURCE\unitlength}\basiceaR[\ARROWLENGTH]}}%

\def\basicedist[#1]{%
\Z=#1%
\multiply \Z by 100%
\hcase{\EDIST}{\Z}}%

\newcommand{\edist}{\@ifnextchar[{\basicedist}%
{\hspace{\SOURCE\unitlength}\basicedist[\ARROWLENGTH]}}%

\def\basicEdist[#1]#2{%
\Z=#1%
\multiply \Z by 100%
\Hcase{\EDIST}{\DUP{#2}}{\Z}}%

\newcommand{\Edist}{\@ifnextchar[{\basicEdist}%
{\hspace{\SOURCE\unitlength}\basicEdist[\ARROWLENGTH]}}%

\def\basicedisT[#1]#2{%
\Z=#1%
\multiply \Z by 100%
\hcasE{\EDIST}{\DDOWN{#2}}{\Z}}%

\newcommand{\edisT}{\@ifnextchar[{\basicedisT}%
{\hspace{\SOURCE\unitlength}\basicedisT[\ARROWLENGTH]}}%

\def\basicedotar[#1]{%
\Z=#1%
\multiply \Z by 100%
\hcase{\EDOTAR}{\Z}}%

\newcommand{\edotar}{\@ifnextchar[{\basicedotar}%
{\hspace{\SOURCE\unitlength}\basicedotar[\ARROWLENGTH]}}%

\def\basicEdotar[#1]#2{%
\Z=#1%
\multiply \Z by 100%
\Hcase{\EDOTAR}{#2}{\Z}}%

\newcommand{\Edotar}{\@ifnextchar[{\basicEdotar}%
{\hspace{\SOURCE\unitlength}\basicEdotar[\ARROWLENGTH]}}%

\def\basicedotaR[#1]#2{%
\Z=#1%
\multiply \Z by 100%
\hcasE{\EDOTAR}{#2}{\Z}}%

\newcommand{\edotaR}{\@ifnextchar[{\basicedotaR}%
{\hspace{\SOURCE\unitlength}\basicedotaR[\ARROWLENGTH]}}%

\def\basicemono[#1]{%
\Z=#1%
\multiply \Z by 100%
\hcase{\EMONO}{\Z}}%

\newcommand{\emono}{\@ifnextchar[{\basicemono}%
{\hspace{\SOURCE\unitlength}\basicemono[\ARROWLENGTH]}}%

\def\basicEmono[#1]#2{%
\Z=#1%
\multiply \Z by 100%
\Hcase{\EMONO}{#2}{\Z}}%

\newcommand{\Emono}{\@ifnextchar[{\basicEmono}%
{\hspace{\SOURCE\unitlength}\basicEmono[\ARROWLENGTH]}}%

\def\basicemonO[#1]#2{%
\Z=#1%
\multiply \Z by 100%
\hcasE{\EMONO}{#2}{\Z}}%

\newcommand{\emonO}{\@ifnextchar[{\basicemonO}%
{\hspace{\SOURCE\unitlength}\basicemonO[\ARROWLENGTH]}}%

\def\basiceepi[#1]{%
\Z=#1%
\multiply \Z by 100%
\hcase{\EEPI}{\Z}}%

\newcommand{\eepi}{\@ifnextchar[{\basiceepi}%
{\hspace{\SOURCE\unitlength}\basiceepi[\ARROWLENGTH]}}%

\def\basicEepi[#1]#2{%
\Z=#1%
\multiply \Z by 100%
\Hcase{\EEPI}{#2}{\Z}}%

\newcommand{\Eepi}{\@ifnextchar[{\basicEepi}%
{\hspace{\SOURCE\unitlength}\basicEepi[\ARROWLENGTH]}}%

\def\basiceepI[#1]#2{%
\Z=#1%
\multiply \Z by 100%
\hcasE{\EEPI}{#2}{\Z}}%

\newcommand{\eepI}{\@ifnextchar[{\basiceepI}%
{\hspace{\SOURCE\unitlength}\basiceepI[\ARROWLENGTH]}}%

\def\basicebimo[#1]{%
\Z=#1%
\multiply \Z by 100%
\hcase{\EBIMO}{\Z}}%

\newcommand{\ebimo}{\@ifnextchar[{\basicebimo}%
{\hspace{\SOURCE\unitlength}\basicebimo[\ARROWLENGTH]}}%

\def\basicEbimo[#1]#2{%
\Z=#1%
\multiply \Z by 100%
\Hcase{\EBIMO}{#2}{\Z}}%

\newcommand{\Ebimo}{\@ifnextchar[{\basicEbimo}%
{\hspace{\SOURCE\unitlength}\basicEbimo[\ARROWLENGTH]}}%

\def\basicebimO[#1]#2{%
\Z=#1%
\multiply \Z by 100%
\hcasE{\EBIMO}{#2}{\Z}}%

\newcommand{\ebimO}{\@ifnextchar[{\basicebimO}%
{\hspace{\SOURCE\unitlength}\basicebimO[\ARROWLENGTH]}}%

\def\basiceiso[#1]{%
\Z=#1%
\multiply \Z by 100%
\Hisocase{\EAR}{\isosign}{}{\Z}}%

\newcommand{\eiso}{\@ifnextchar[{\basiceiso}%
{\hspace{\SOURCE\unitlength}\basiceiso[\ARROWLENGTH]}}%

\def\basicEiso[#1]#2{%
\Z=#1%
\multiply \Z by 100%
\Hisocase{\EAR}{#2}{\isosign}{\Z}}%

\newcommand{\Eiso}{\@ifnextchar[{\basicEiso}%
{\hspace{\SOURCE\unitlength}\basicEiso[\ARROWLENGTH]}}%

\def\basiceisO[#1]#2{%
\Z=#1%
\multiply \Z by 100%
\Hisocase{\EAR}{\isosign}{#2}{\Z}}%

\newcommand{\eisO}{\@ifnextchar[{\basiceisO}%
{\hspace{\SOURCE\unitlength}\basiceisO[\ARROWLENGTH]}}%

\def\basiceeql[#1]{%
\Z=#1%
\multiply \Z by 100%
\hcase{\EEQL}{\Z}}%

\newcommand{\eeql}{\@ifnextchar[{\basiceeql}%
{\hspace{\SOURCE\unitlength}\basiceeql[\ARROWLENGTH]}}%

\def\basicEeql[#1]#2{%
\Z=#1%
\multiply \Z by 100%
\Hcase{\EEQL}{\DUP{#2}}{\Z}}%

\newcommand{\Eeql}{\@ifnextchar[{\basicEeql}%
{\hspace{\SOURCE\unitlength}\basicEeql[\ARROWLENGTH]}}%

\def\basiceeqL[#1]#2{%
\Z=#1%
\multiply \Z by 100%
\hcasE{\EEQL}{#2}{\Z}}%

\newcommand{\eeqL}{\@ifnextchar[{\basiceeqL}%
{\hspace{\SOURCE\unitlength}\basiceeqL[\ARROWLENGTH]}}%

\def\basiceline[#1]{%
\Z=#1%
\multiply \Z by 100%
\hcase{\ELINE}{\Z}}%

\newcommand{\eline}{\@ifnextchar[{\basiceline}%
{\hspace{\SOURCE\unitlength}\basiceline[\ARROWLENGTH]}}%

\def\basicEline[#1]#2{%
\Z=#1%
\multiply \Z by 100%
\Hcase{\ELINE}{\DUP{#2}}{\Z}}%

\newcommand{\Eline}{\@ifnextchar[{\basicEline}%
{\hspace{\SOURCE\unitlength}\basicEline[\ARROWLENGTH]}}%

\def\basicelinE[#1]#2{%
\Z=#1%
\multiply \Z by 100%
\hcasE{\ELINE}{#2}{\Z}}%

\newcommand{\elinE}{\@ifnextchar[{\basicelinE}%
{\hspace{\SOURCE\unitlength}\basicelinE[\ARROWLENGTH]}}%

\def\basicebiar[#1]{%
\Z=#1%
\multiply \Z by 100%
\hbicase{\EBIAR}{\Z}}%

\newcommand{\ebiar}{\@ifnextchar[{\basicebiar}%
{\hspace{\SOURCE\unitlength}\basicebiar[\ARROWLENGTH]}}%

\def\basicEbiar[#1]#2#3{%
\Z=#1%
\multiply \Z by 100%
\Hbicase{\EBIAR}{#2}{#3}{\Z}}%

\newcommand{\Ebiar}{\@ifnextchar[{\basicEbiar}%
{\hspace{\SOURCE\unitlength}\basicEbiar[\ARROWLENGTH]}}%

\def\basicebidist[#1]{%
\Z=#1%
\multiply \Z by 100%
\hbicase{\EBIDIST}{\Z}}%

\newcommand{\ebidist}{\@ifnextchar[{\basicebidist}%
{\hspace{\SOURCE\unitlength}\basicebidist[\ARROWLENGTH]}}%

\def\basicEbidist[#1]#2#3{%
\Z=#1%
\multiply \Z by 100%
\Hbicase{\EBIDIST}{\DUP{#2}}{\DDOWN{#3}}{\Z}}%

\newcommand{\Ebidist}{\@ifnextchar[{\basicEbidist}%
{\hspace{\SOURCE\unitlength}\basicEbidist[\ARROWLENGTH]}}%

\def\basiceadjar[#1]{%
\Z=#1%
\multiply \Z by 100%
\hbicase{\EADJAR}{\Z}}%

\newcommand{\eadjar}{\@ifnextchar[{\basiceadjar}%
{\hspace{\SOURCE\unitlength}\basiceadjar[\ARROWLENGTH]}}%

\def\basicEadjar[#1]#2#3{%
\Z=#1%
\multiply \Z by 100%
\Hbicase{\EADJAR}{#2}{#3}{\Z}}%

\newcommand{\Eadjar}{\@ifnextchar[{\basicEadjar}%
{\hspace{\SOURCE\unitlength}\basicEadjar[\ARROWLENGTH]}}%

\def\basiceadjdist[#1]{%
\Z=#1%
\multiply \Z by 100%
\hbicase{\EADJDIST}{\Z}}%

\newcommand{\eadjdist}{\@ifnextchar[{\basiceadjdist}%
{\hspace{\SOURCE\unitlength}\basiceadjdist[\ARROWLENGTH]}}%

\def\basicEadjdist[#1]#2#3{%
\Z=#1%
\multiply \Z by 100%
\Hbicase{\EADJDIST}{\DUP{#2}}{\DDOWN{#3}}{\Z}}%

\newcommand{\Eadjdist}{\@ifnextchar[{\basicEadjdist}%
{\hspace{\SOURCE\unitlength}\basicEadjdist[\ARROWLENGTH]}}%

\def\basicetriar[#1]{%
\Z=#1%
\multiply \Z by 100%
\htricase{\ETRIAR}{\Z}}%

\newcommand{\etriar}{\@ifnextchar[{\basicetriar}%
{\hspace{\SOURCE\unitlength}\basicetriar[\ARROWLENGTH]}}%

\def\basicEtriar[#1]#2#3#4{%
\Z=#1%
\multiply \Z by 100%
\Htricase{\ETRIAR}{#2}{#3}{#4}{}{\Z}}%

\newcommand{\Etriar}{\@ifnextchar[{\basicEtriar}%
{\hspace{\SOURCE\unitlength}\basicEtriar[\ARROWLENGTH]}}%

\def\basicetriaR[#1]#2#3#4{%
\Z=#1%
\multiply \Z by 100%
\Htricase{\ETRIAR}{}{#2}{#3}{#4}{\Z}}%

\newcommand{\etriaR}{\@ifnextchar[{\basicetriaR}%
{\hspace{\SOURCE\unitlength}\basicetriaR[\ARROWLENGTH]}}%

\def\basicetriadjar[#1]{%
\Z=#1%
\multiply \Z by 100%
\htricase{\ETRIADJAR}{\Z}}%

\newcommand{\etriadjar}{\@ifnextchar[{\basicetriadjar}%
{\hspace{\SOURCE\unitlength}\basicetriadjar[\ARROWLENGTH]}}%

\def\basicEtriadjar[#1]#2#3#4{%
\Z=#1%
\multiply \Z by 100%
\Htricase{\ETRIADJAR}{#2}{#3}{#4}{}{\Z}}%

\newcommand{\Etriadjar}{\@ifnextchar[{\basicEtriadjar}%
{\hspace{\SOURCE\unitlength}\basicEtriadjar[\ARROWLENGTH]}}%

\def\basicetriadjaR[#1]#2#3#4{%
\Z=#1%
\multiply \Z by 100%
\Htricase{\ETRIADJAR}{}{#2}{#3}{#4}{\Z}}%

\newcommand{\etriadjaR}{\@ifnextchar[{\basicetriadjaR}%
{\hspace{\SOURCE\unitlength}\basicetriadjaR[\ARROWLENGTH]}}%

\def\basicequadriar[#1]{%
\Z=#1%
\multiply \Z by 100%
\hmulticase{\EQUADRIAR}{\Z}}%

\newcommand{\equadriar}{\@ifnextchar[{\basicequadriar}%
{\hspace{\SOURCE\unitlength}\basicequadriar[\ARROWLENGTH]}}%

\def\basicequadriadjar[#1]{%
\Z=#1%
\multiply \Z by 100%
\hmulticase{\EQUADRIADJAR}{\Z}}%

\newcommand{\equadriadjar}{\@ifnextchar[{\basicequadriadjar}%
{\hspace{\SOURCE\unitlength}\basicequadriadjar[\ARROWLENGTH]}}%

\def\basicequintiar[#1]{%
\Z=#1%
\multiply \Z by 100%
\hmulticase{\EQUINTIAR}{\Z}}%

\newcommand{\equintiar}{\@ifnextchar[{\basicequintiar}%
{\hspace{\SOURCE\unitlength}\basicequintiar[\ARROWLENGTH]}}%

\def\basicequintiadjar[#1]{%
\Z=#1%
\multiply \Z by 100%
\hmulticase{\EQUINTIADJAR}{\Z}}%

\newcommand{\equintiadjar}{\@ifnextchar[{\basicequintiadjar}%
{\hspace{\SOURCE\unitlength}\basicequintiadjar[\ARROWLENGTH]}}%

\newcommand{\WAR}[1]%
{\begin{picture}(#1,0)%
\put(#1,0){\line(-1,0){#1}}%
\put(0,0){\whead}%
\end{picture}}%

\newcommand{\WDIST}[1]%
{\begin{picture}(#1,0)%
\put(#1,0){\line(-1,0){#1}}%
\put(0,0){\whead}%
\NUMBER=#1%
\divide\NUMBER by 2%
\put(\NUMBER,0){\distsign}%
\end{picture}}%

\newcommand{\WDOTAR}[1]%
{\truex{100}\truey{300}%
\NUMBEROFDOTS=#1%
\divide\NUMBEROFDOTS by \value{y}%
\advance\NUMBEROFDOTS by 1%
\begin{picture}(#1,0)%
\multiput(#1,0)(-\value{y},0){\NUMBEROFDOTS}%
{\circle*{\value{x}}}%
\put(0,0){\whead}%
\end{picture}}%

\newcommand{\WMONO}[1]%
{\truetail%
\monolength=#1%
\advance\monolength by -\truemonotail%
\begin{picture}(#1,0)(-#1,0)%
\put(-\truemonotail,0){\line(-1,0){\monolength}}%
\put(-\truemonotail,0){\whead}%
\put(-#1,0){\whead}%
\end{picture}}%

\newcommand{\WEPI}[1]%
{\truehead%
\epilength=#1%
\advance\epilength by -\trueepihead%
\begin{picture}(#1,0)%
\put(#1,0){\line(-1,0){\epilength}}%
\put(\trueepihead,0){\whead}%
\put(0,0){\whead}%
\end{picture}}%

\newcommand{\WBIMO}[1]%
{\truehead\truetail%
\monolength=#1
\advance\monolength by -\truemonotail%
\epilength=\monolength%
\advance\epilength by -\trueepihead%
\begin{picture}(#1,0)%
\put(\monolength,0){\line(-1,0){\epilength}}%
\put(\monolength,0){\whead}%
\put(\trueepihead,0){\whead}%
\put(0,0){\whead}%
\end{picture}}%

\newcommand{\WBIAR}[1]%
{\truex{700}%
\begin{picture}(#1,\value{x})%
\put(#1,0){\line(-1,0){#1}}%
\put(0,0){\whead}%
\put(#1,\value{x}){\line(-1,0){#1}}%
\put(0,\value{x}){\whead}%
\end{picture}}%

\newcommand{\WBIDIST}[1]%
{\truex{700}%
\begin{picture}(#1,\value{x})%
\put(#1,0){\line(-1,0){#1}}%
\put(0,0){\whead}%
\put(#1,\value{x}){\line(-1,0){#1}}%
\put(0,\value{x}){\whead}%
\NUMBER=#1%
\divide\NUMBER by 2%
\put(\NUMBER,0){\distsign}%
\put(\NUMBER,\value{x}){\distsign}%
\end{picture}}%

\newcommand{\WADJAR}[1]%
{\truex{700}%
\begin{picture}(#1,\value{x})%
\put(0,\value{x}){\line(1,0){#1}}%
\put(#1,\value{x}){\ehead}%
\put(#1,0){\line(-1,0){#1}}%
\put(0,0){\whead}%
\end{picture}}%

\newcommand{\WADJDIST}[1]%
{\truex{700}%
\begin{picture}(#1,\value{x})%
\put(0,\value{x}){\line(1,0){#1}}%
\put(#1,\value{x}){\ehead}%
\put(#1,0){\line(-1,0){#1}}%
\put(0,0){\whead}%
\NUMBER=#1%
\divide\NUMBER by 2%
\put(\NUMBER,0){\distsign}%
\put(\NUMBER,\value{x}){\distsign}%
\end{picture}}%

\newcommand{\WTRIAR}[5]%
{\truex{1200}%
\X=\value{x}%
\multiply\X by 2%
\begin{picture}(#1,\X)%
\put(#1,0){\line(-1,0){#1}}%
\put(0,0){\whead}%
\put(#1,\value{x}){\line(-1,0){#1}}%
\put(0,\value{x}){\whead}%
\put(#1,\X){\line(-1,0){#1}}%
\put(0,\X){\whead}%
\X=#1%
\divide\X by 2%
\truex{-600}%
\put(\X,\value{x}){\makebox(0,0){$#5$}}%
\truex{600}%
\put(\X,\value{x}){\makebox(0,0){$#4$}}%
\truex{1800}%
\put(\X,\value{x}){\makebox(0,0){$#3$}}%
\truex{3000}%
\put(\X,\value{x}){\makebox(0,0){$#2$}}%
\end{picture}}%

\newcommand{\WTRIADJAR}[5]%
{\truex{1200}%
\X=\value{x}%
\multiply\X by 2%
\begin{picture}(#1,\X)%
\put(0,0){\line(1,0){#1}}%
\put(#1,0){\ehead}%
\put(#1,\value{x}){\line(-1,0){#1}}%
\put(0,\value{x}){\whead}%
\put(0,\X){\line(1,0){#1}}%
\put(#1,\X){\ehead}%
\X=#1%
\divide\X by 2%
\truex{-600}%
\put(\X,\value{x}){\makebox(0,0){$#5$}}%
\truex{600}%
\put(\X,\value{x}){\makebox(0,0){$#4$}}%
\truex{1800}%
\put(\X,\value{x}){\makebox(0,0){$#3$}}%
\truex{3000}%
\put(\X,\value{x}){\makebox(0,0){$#2$}}%
\end{picture}}%

\newcommand{\WQUADRIAR}[1]%
{\truex{800}%
\truey{1600}%
\X=\value{x}%
\multiply\X by 3%
\begin{picture}(#1,\X)%
\put(#1,0){\line(-1,0){#1}}%
\put(0,0){\whead}%
\put(#1,\value{x}){\line(-1,0){#1}}%
\put(0,\value{x}){\whead}%
\put(#1,\value{y}){\line(-1,0){#1}}%
\put(0,\value{y}){\whead}%
\put(#1,\X){\line(-1,0){#1}}%
\put(0,\X){\whead}%
\end{picture}}%

\newcommand{\WQUADRIADJAR}[1]%
{\truex{800}%
\truey{1600}%
\X=\value{x}%
\multiply\X by 3%
\begin{picture}(#1,\X)%
\put(#1,0){\line(-1,0){#1}}%
\put(0,0){\whead}%
\put(0,\value{x}){\line(1,0){#1}}%
\put(#1,\value{x}){\ehead}%
\put(#1,\value{y}){\line(-1,0){#1}}%
\put(0,\value{y}){\whead}%
\put(0,\X){\line(1,0){#1}}%
\put(#1,\X){\ehead}%
\end{picture}}%

\newcommand{\WQUINTIAR}[1]%
{\truex{600}%
\truey{1200}%
\truez{1800}%
\X=\value{x}%
\multiply\X by 4%
\begin{picture}(#1,\X)%
\put(#1,0){\line(-1,0){#1}}%
\put(0,0){\whead}%
\put(#1,\value{x}){\line(-1,0){#1}}%
\put(0,\value{x}){\whead}%
\put(#1,\value{y}){\line(-1,0){#1}}%
\put(0,\value{y}){\whead}%
\put(#1,\value{z}){\line(-1,0){#1}}%
\put(0,\value{z}){\whead}%
\put(#1,\X){\line(-1,0){#1}}%
\put(0,\X){\whead}%
\end{picture}}%

\newcommand{\WQUINTIADJAR}[1]%
{\truex{600}%
\truey{1200}%
\truez{1800}%
\X=\value{x}%
\multiply\X by 4%
\begin{picture}(#1,\X)%
\put(#1,0){\line(-1,0){#1}}%
\put(0,0){\whead}%
\put(0,\value{x}){\line(1,0){#1}}%
\put(#1,\value{x}){\ehead}%
\put(#1,\value{y}){\line(-1,0){#1}}%
\put(0,\value{y}){\whead}%
\put(0,\value{z}){\line(1,0){#1}}%
\put(#1,\value{z}){\ehead}%
\put(#1,\X){\line(-1,0){#1}}%
\put(0,\X){\whead}%
\end{picture}}%

\def\basicwar[#1]{%
\Z=#1%
\multiply \Z by 100%
\hcase{\WAR}{\Z}}%

\newcommand{\war}{\@ifnextchar[{\basicwar}%
{\hspace{\SOURCE\unitlength}\basicwar[\ARROWLENGTH]}}%

\def\basicWar[#1]#2{%
\Z=#1%
\multiply \Z by 100%
\Hcase{\WAR}{#2}{\Z}}%

\newcommand{\War}{\@ifnextchar[{\basicWar}%
{\hspace{\SOURCE\unitlength}\basicWar[\ARROWLENGTH]}}%

\def\basicwaR[#1]#2{%
\Z=#1%
\multiply \Z by 100%
\hcasE{\WAR}{#2}{\Z}}%

\newcommand{\waR}{\@ifnextchar[{\basicwaR}%
{\hspace{\SOURCE\unitlength}\basicwaR[\ARROWLENGTH]}}%

\def\basicwdist[#1]{%
\Z=#1%
\multiply \Z by 100%
\hcase{\WDIST}{\Z}}%

\newcommand{\wdist}{\@ifnextchar[{\basicwdist}%
{\hspace{\SOURCE\unitlength}\basicwdist[\ARROWLENGTH]}}%

\def\basicWdist[#1]#2{%
\Z=#1%
\multiply \Z by 100%
\Hcase{\WDIST}{\DUP{#2}}{\Z}}%

\newcommand{\Wdist}{\@ifnextchar[{\basicWdist}%
{\hspace{\SOURCE\unitlength}\basicWdist[\ARROWLENGTH]}}%

\def\basicwdisT[#1]#2{%
\Z=#1%
\multiply \Z by 100%
\hcasE{\WDIST}{\DDOWN{#2}}{\Z}}%

\newcommand{\wdisT}{\@ifnextchar[{\basicwdisT}%
{\hspace{\SOURCE\unitlength}\basicwdisT[\ARROWLENGTH]}}%

\def\basicwdotar[#1]{%
\Z=#1%
\multiply \Z by 100%
\hcase{\WDOTAR}{\Z}}%

\newcommand{\wdotar}{\@ifnextchar[{\basicwdotar}%
{\hspace{\SOURCE\unitlength}\basicwdotar[\ARROWLENGTH]}}%

\def\basicWdotar[#1]#2{%
\Z=#1%
\multiply \Z by 100%
\Hcase{\WDOTAR}{#2}{\Z}}%

\newcommand{\Wdotar}{\@ifnextchar[{\basicWdotar}%
{\hspace{\SOURCE\unitlength}\basicWdotar[\ARROWLENGTH]}}%

\def\basicwdotaR[#1]#2{%
\Z=#1%
\multiply \Z by 100%
\hcasE{\WDOTAR}{#2}{\Z}}%

\newcommand{\wdotaR}{\@ifnextchar[{\basicwdotaR}%
{\hspace{\SOURCE\unitlength}\basicwdotaR[\ARROWLENGTH]}}%

\def\basicwmono[#1]{%
\Z=#1%
\multiply \Z by 100%
\hcase{\WMONO}{\Z}}%

\newcommand{\wmono}{\@ifnextchar[{\basicwmono}%
{\hspace{\SOURCE\unitlength}\basicwmono[\ARROWLENGTH]}}%

\def\basicWmono[#1]#2{%
\Z=#1%
\multiply \Z by 100%
\Hcase{\WMONO}{#2}{\Z}}%

\newcommand{\Wmono}{\@ifnextchar[{\basicWmono}%
{\hspace{\SOURCE\unitlength}\basicWmono[\ARROWLENGTH]}}%

\def\basicwmonO[#1]#2{%
\Z=#1%
\multiply \Z by 100%
\hcasE{\WMONO}{#2}{\Z}}%

\newcommand{\wmonO}{\@ifnextchar[{\basicwmonO}%
{\hspace{\SOURCE\unitlength}\basicwmonO[\ARROWLENGTH]}}%

\def\basicwepi[#1]{%
\Z=#1%
\multiply \Z by 100%
\hcase{\WEPI}{\Z}}%

\newcommand{\wepi}{\@ifnextchar[{\basicwepi}%
{\hspace{\SOURCE\unitlength}\basicwepi[\ARROWLENGTH]}}%

\def\basicWepi[#1]#2{%
\Z=#1%
\multiply \Z by 100%
\Hcase{\WEPI}{#2}{\Z}}%

\newcommand{\Wepi}{\@ifnextchar[{\basicWepi}%
{\hspace{\SOURCE\unitlength}\basicWepi[\ARROWLENGTH]}}%

\def\basicwepI[#1]#2{%
\Z=#1%
\multiply \Z by 100%
\hcasE{\WEPI}{#2}{\Z}}%

\newcommand{\wepI}{\@ifnextchar[{\basicwepI}%
{\hspace{\SOURCE\unitlength}\basicwepI[\ARROWLENGTH]}}%

\def\basicwbimo[#1]{%
\Z=#1%
\multiply \Z by 100%
\hcase{\WBIMO}{\Z}}%

\newcommand{\wbimo}{\@ifnextchar[{\basicwbimo}%
{\hspace{\SOURCE\unitlength}\basicwbimo[\ARROWLENGTH]}}%

\def\basicWbimo[#1]#2{%
\Z=#1%
\multiply \Z by 100%
\Hcase{\WBIMO}{#2}{\Z}}%

\newcommand{\Wbimo}{\@ifnextchar[{\basicWbimo}%
{\hspace{\SOURCE\unitlength}\basicWbimo[\ARROWLENGTH]}}%

\def\basicwbimO[#1]#2{%
\Z=#1%
\multiply \Z by 100%
\hcasE{\WBIMO}{#2}{\Z}}%

\newcommand{\wbimO}{\@ifnextchar[{\basicwbimO}%
{\hspace{\SOURCE\unitlength}\basicwbimO[\ARROWLENGTH]}}%

\def\basicwiso[#1]{%
\Z=#1%
\multiply \Z by 100%
\Hisocase{\WAR}{\isosign}{}{\Z}}%

\newcommand{\wiso}{\@ifnextchar[{\basicwiso}%
{\hspace{\SOURCE\unitlength}\basicwiso[\ARROWLENGTH]}}%

\def\basicWiso[#1]#2{%
\Z=#1%
\multiply \Z by 100%
\Hisocase{\WAR}{#2}{\isosign}{\Z}}%

\newcommand{\Wiso}{\@ifnextchar[{\basicWiso}%
{\hspace{\SOURCE\unitlength}\basicWiso[\ARROWLENGTH]}}%

\def\basicwisO[#1]#2{%
\Z=#1%
\multiply \Z by 100%
\Hisocase{\WAR}{\isosign}{#2}{\Z}}%

\newcommand{\wisO}{\@ifnextchar[{\basicwisO}%
{\hspace{\SOURCE\unitlength}\basicwisO[\ARROWLENGTH]}}%

\def\basicwbiar[#1]{%
\Z=#1%
\multiply \Z by 100%
\hbicase{\WBIAR}{\Z}}%

\newcommand{\wbiar}{\@ifnextchar[{\basicwbiar}%
{\hspace{\SOURCE\unitlength}\basicwbiar[\ARROWLENGTH]}}%

\def\basicWbiar[#1]#2#3{%
\Z=#1%
\multiply \Z by 100%
\Hbicase{\WBIAR}{#2}{#3}{\Z}}%

\newcommand{\Wbiar}{\@ifnextchar[{\basicWbiar}%
{\hspace{\SOURCE\unitlength}\basicWbiar[\ARROWLENGTH]}}%

\def\basicwbidist[#1]{%
\Z=#1%
\multiply \Z by 100%
\hbicase{\WBIDIST}{\Z}}%

\newcommand{\wbidist}{\@ifnextchar[{\basicwbidist}%
{\hspace{\SOURCE\unitlength}\basicwbidist[\ARROWLENGTH]}}%

\def\basicWbidist[#1]#2#3{%
\Z=#1%
\multiply \Z by 100%
\Hbicase{\WBIDIST}{\DUP{#2}}{\DDOWN{#3}}{\Z}}%

\newcommand{\Wbidist}{\@ifnextchar[{\basicWbidist}%
{\hspace{\SOURCE\unitlength}\basicWbidist[\ARROWLENGTH]}}%

\def\basicwadjar[#1]{%
\Z=#1%
\multiply \Z by 100%
\hbicase{\WADJAR}{\Z}}%

\newcommand{\wadjar}{\@ifnextchar[{\basicwadjar}%
{\hspace{\SOURCE\unitlength}\basicwadjar[\ARROWLENGTH]}}%

\def\basicWadjar[#1]#2#3{%
\Z=#1%
\multiply \Z by 100%
\Hbicase{\WADJAR}{#2}{#3}{\Z}}%

\newcommand{\Wadjar}{\@ifnextchar[{\basicWadjar}%
{\hspace{\SOURCE\unitlength}\basicWadjar[\ARROWLENGTH]}}%

\def\basicwadjdist[#1]{%
\Z=#1%
\multiply \Z by 100%
\hbicase{\WADJDIST}{\Z}}%

\newcommand{\wadjdist}{\@ifnextchar[{\basicwadjdist}%
{\hspace{\SOURCE\unitlength}\basicwadjdist[\ARROWLENGTH]}}%

\def\basicWadjdist[#1]#2#3{%
\Z=#1%
\multiply \Z by 100%
\Hbicase{\WADJDIST}{\DUP{#2}}{\DDOWN{#3}}{\Z}}%

\newcommand{\Wadjdist}{\@ifnextchar[{\basicWadjdist}%
{\hspace{\SOURCE\unitlength}\basicWadjdist[\ARROWLENGTH]}}%

\def\basicwtriar[#1]{%
\Z=#1%
\multiply \Z by 100%
\htricase{\WTRIAR}{\Z}}%

\newcommand{\wtriar}{\@ifnextchar[{\basicwtriar}%
{\hspace{\SOURCE\unitlength}\basicwtriar[\ARROWLENGTH]}}%

\def\basicWtriar[#1]#2#3#4{%
\Z=#1%
\multiply \Z by 100%
\Htricase{\WTRIAR}{#2}{#3}{#4}{}{\Z}}%

\newcommand{\Wtriar}{\@ifnextchar[{\basicWtriar}%
{\hspace{\SOURCE\unitlength}\basicWtriar[\ARROWLENGTH]}}%

\def\basicwtriaR[#1]#2#3#4{%
\Z=#1%
\multiply \Z by 100%
\Htricase{\WTRIAR}{}{#2}{#3}{#4}{\Z}}%

\newcommand{\wtriaR}{\@ifnextchar[{\basicwtriaR}%
{\hspace{\SOURCE\unitlength}\basicwtriaR[\ARROWLENGTH]}}%

\def\basicwtriadjar[#1]{%
\Z=#1%
\multiply \Z by 100%
\htricase{\WTRIADJAR}{\Z}}%

\newcommand{\wtriadjar}{\@ifnextchar[{\basicwtriadjar}%
{\hspace{\SOURCE\unitlength}\basicwtriadjar[\ARROWLENGTH]}}%

\def\basicWtriadjar[#1]#2#3#4{%
\Z=#1%
\multiply \Z by 100%
\Htricase{\WTRIADJAR}{#2}{#3}{#4}{}{\Z}}%

\newcommand{\Wtriadjar}{\@ifnextchar[{\basicWtriadjar}%
{\hspace{\SOURCE\unitlength}\basicWtriadjar[\ARROWLENGTH]}}%

\def\basicwtriadjaR[#1]#2#3#4{%
\Z=#1%
\multiply \Z by 100%
\Htricase{\WTRIADJAR}{}{#2}{#3}{#4}{\Z}}%

\newcommand{\wtriadjaR}{\@ifnextchar[{\basicwtriadjaR}%
{\hspace{\SOURCE\unitlength}\basicwtriadjaR[\ARROWLENGTH]}}%

\def\basicwquadriar[#1]{%
\Z=#1%
\multiply \Z by 100%
\hmulticase{\WQUADRIAR}{\Z}}%

\newcommand{\wquadriar}{\@ifnextchar[{\basicwquadriar}%
{\hspace{\SOURCE\unitlength}\basicwquadriar[\ARROWLENGTH]}}%

\def\basicwquadriadjar[#1]{%
\Z=#1%
\multiply \Z by 100%
\hmulticase{\WQUADRIADJAR}{\Z}}%

\newcommand{\wquadriadjar}{\@ifnextchar[{\basicwquadriadjar}%
{\hspace{\SOURCE\unitlength}\basicwquadriadjar[\ARROWLENGTH]}}%

\def\basicwquintiar[#1]{%
\Z=#1%
\multiply \Z by 100%
\hmulticase{\WQUINTIAR}{\Z}}%

\newcommand{\wquintiar}{\@ifnextchar[{\basicwquintiar}%
{\hspace{\SOURCE\unitlength}\basicwquintiar[\ARROWLENGTH]}}%

\def\basicwquintiadjar[#1]{%
\Z=#1%
\multiply \Z by 100%
\hmulticase{\WQUINTIADJAR}{\Z}}%

\newcommand{\wquintiadjar}{\@ifnextchar[{\basicwquintiadjar}%
{\hspace{\SOURCE\unitlength}\basicwquintiadjar[\ARROWLENGTH]}}%

\newcommand{\vcase}[2]{#1{#2}}%

\newcommand{\Vcase}[3]{\makebox[0pt]%
{\makebox[0pt][r]{\raisebox{0pt}[0pt][0pt]{${#2}\hspace{2pt}$}}}#1{#3}}%

\newcommand{\vcasE}[3]{\makebox[0pt]%
{#1{#3}\makebox[0pt][l]{\raisebox{0pt}[0pt][0pt]{\hspace{2pt}$#2$}}}}%

\newcommand{\Visocase}[4]{\makebox[0pt]%
{\makebox[0pt][r]{\raisebox{0pt}[0pt][0pt]{$#2$\hspace{2pt}}}#1{#4}%
\makebox[0pt][l]{\raisebox{0pt}[0pt][0pt]{\hspace{2pt}$#3$}}}}%

\newcommand{\vbicase}[2]{\makebox[0pt]{{#1{#2}}}}%

\newcommand{\Vbicase}[4]{\makebox[0pt]%
{\makebox[0pt][r]{\raisebox{0pt}[0pt][0pt]{$#2$\hspace{5.5pt}}}#1{#4}%
\makebox[0pt][l]{\raisebox{0pt}[0pt][0pt]{\hspace{6.5pt}$#3$}}}}%

\newcommand{\vtricase}[2]%
{\makebox[0pt]{#1{}{}{}{}{#2}}}%

\newcommand{\Vtricase}[6]%
{\makebox[0pt]{#1{#2}{#3}{#4}{#5}{#6}}}%

\newcommand{\vmulticase}[2]{\makebox[0pt]{#1{#2}}}%

\newcommand{\SAR}[1]%
{\begin{picture}(0,0)%
\put(0,0){\makebox(0,0)%
{\begin{picture}(0,#1)%
\put(0,#1){\line(0,-1){#1}}%
\put(0,0){\shead}%
\end{picture}}}\end{picture}}%

\newcommand{\SDIST}[1]%
{\begin{picture}(0,0)%
\put(0,0){\makebox(0,0)%
{\begin{picture}(0,#1)%
\put(0,#1){\line(0,-1){#1}}%
\put(0,0){\shead}%
\end{picture}}}%
\put(0,0){\distsign}%
\end{picture}}%

\newcommand{\SDOTAR}[1]%
{\truex{100}\truey{300}%
\NUMBEROFDOTS=#1%
\divide\NUMBEROFDOTS by \value{y}%
\advance\NUMBEROFDOTS by 1%
\begin{picture}(0,0)%
\put(0,0){\makebox(0,0)%
{\begin{picture}(0,#1)%
\multiput(0,#1)(0,-\value{y}){\NUMBEROFDOTS}%
{\circle*{\value{x}}}%
\put(0,0){\shead}%
\end{picture}}}\end{picture}}%

\newcommand{\SMONO}[1]%
{\truetail%
\monolength=#1%
\advance\monolength by -\truemonotail%
\begin{picture}(0,0)%
\put(0,0){\makebox(0,0)%
{\begin{picture}(0,#1)%
\put(0,\monolength){\line(0,-1){\monolength}}%
\put(0,\monolength){\shead}%
\put(0,0){\shead}%
\end{picture}}}\end{picture}}%

\newcommand{\SEPI}[1]%
{\truehead%
\epilength=#1%
\advance\epilength by -\trueepihead%
\begin{picture}(0,0)%
\put(0,0){\makebox(0,0)%
{\begin{picture}(0,#1)%
\put(0,#1){\line(0,-1){\epilength}}%
\put(0,\trueepihead){\shead}%
\put(0,0){\shead}%
\end{picture}}}\end{picture}}%

\newcommand{\SBIMO}[1]%
{\truehead\truetail%
\monolength=#1%
\advance\monolength by -\truemonotail%
\epilength=\monolength%
\advance\epilength by -\trueepihead%
\begin{picture}(0,0)%
\put(0,0){\makebox(0,0)%
{\begin{picture}(0,#1)%
\put(0,\monolength){\line(0,-1){\epilength}}%
\put(0,\monolength){\shead}%
\put(0,\trueepihead){\shead}%
\put(0,0){\shead}%
\end{picture}}}\end{picture}}%

\newcommand{\SBIAR}[1]%
{\begin{picture}(0,0)%
\truex{350}%
\put(0,0){\makebox(0,0)%
{\begin{picture}(0,#1)%
\put(-\value{x},#1){\line(0,-1){#1}}%
\put(-\value{x},0){\shead}%
\put(\value{x},#1){\line(0,-1){#1}}%
\put(\value{x},0){\shead}%
\end{picture}}}\end{picture}}%

\newcommand{\SBIDIST}[1]%
{\begin{picture}(0,0)%
\truex{350}%
\put(0,0){\makebox(0,0)%
{\begin{picture}(0,#1)%
\put(-\value{x},#1){\line(0,-1){#1}}%
\put(-\value{x},0){\shead}%
\put(\value{x},#1){\line(0,-1){#1}}%
\put(\value{x},0){\shead}%
\end{picture}}}%
\put(-\value{x},0){\distsign}%
\put(\value{x},0){\distsign}%
\end{picture}}%

\newcommand{\SEQL}[1]%
{\begin{picture}(0,0)%
\truex{100}%
\put(0,0){\makebox(0,0)%
{\begin{picture}(0,#1)\put(-\value{x},#1){\line(0,-1){#1}}%
\put(\value{x},#1){\line(0,-1){#1}}%
\end{picture}}}\end{picture}}%

\newcommand{\SLINE}[1]%
{\begin{picture}(0,0)%
\put(0,0){\makebox(0,0)%
{\begin{picture}(0,#1)\put(0,#1){\line(0,-1){#1}}%
\end{picture}}}\end{picture}}%

\newcommand{\SADJAR}[1]{\begin{picture}(0,0)%
\truex{350}%
\put(0,0){\makebox(0,0)%
{\begin{picture}(0,#1)%
\put(-\value{x},#1){\line(0,-1){#1}}%
\put(-\value{x},0){\shead}%
\put(\value{x},0){\line(0,1){#1}}%
\put(\value{x},#1){\nhead}%
\end{picture}}}\end{picture}}%

\newcommand{\SADJDIST}[1]{\begin{picture}(0,0)%
\truex{350}%
\put(0,0){\makebox(0,0)%
{\begin{picture}(0,#1)%
\put(-\value{x},#1){\line(0,-1){#1}}%
\put(-\value{x},0){\shead}%
\put(\value{x},0){\line(0,1){#1}}%
\put(\value{x},#1){\nhead}%
\end{picture}}}%
\put(-\value{x},0){\distsign}%
\put(\value{x},0){\distsign}%
\end{picture}}%

\newcommand{\STRIAR}[5]%
{\truex{1200}%
\truey{1800}%
\truez{600}%
\begin{picture}(0,0)%
\put(0,0){\makebox(0,0)%
{\begin{picture}(0,#5)%
\X=#5\divide\X by 2%
\put(-\value{x},#5){\line(0,-1){#5}}%
\put(-\value{x},0){\shead}%
\put(0,#5){\line(0,-1){#5}}%
\put(0,0){\shead}%
\put(\value{x},#5){\line(0,-1){#5}}%
\put(\value{x},0){\shead}%
\put(-\value{y},\X){\makebox(0,0){$#1$}}%
\put(\value{y},\X){\makebox(0,0){$#4$}}%
\put(-\value{z},\X){\makebox(0,0){$#2$}}%
\put(\value{z},\X){\makebox(0,0){$#3$}}%
\end{picture}}}\end{picture}}%

\newcommand{\STRIADJAR}[5]%
{\truex{1200}%
\truey{1800}%
\truez{600}%
\begin{picture}(0,0)%
\put(0,0){\makebox(0,0)%
{\begin{picture}(0,#5)%
\X=#5\divide\X by 2%
\put(-\value{x},0){\line(0,1){#5}}%
\put(-\value{x},#5){\nhead}%
\put(0,#5){\line(0,-1){#5}}%
\put(0,0){\shead}%
\put(\value{x},0){\line(0,1){#5}}%
\put(\value{x},#5){\nhead}%
\put(-\value{y},\X){\makebox(0,0){$#1$}}%
\put(\value{y},\X){\makebox(0,0){$#4$}}%
\put(-\value{z},\X){\makebox(0,0){$#2$}}%
\put(\value{z},\X){\makebox(0,0){$#3$}}%
\end{picture}}}\end{picture}}%

\newcommand{\SQUADRIAR}[1]%
{\truex{400}%
\truey{1200}%
\begin{picture}(0,0)%
\put(0,0){\makebox(0,0)%
{\begin{picture}(0,#1)%
\put(-\value{y},#1){\line(0,-1){#1}}%
\put(-\value{y},0){\shead}%
\put(-\value{x},#1){\line(0,-1){#1}}%
\put(-\value{x},0){\shead}%
\put(\value{x},#1){\line(0,-1){#1}}%
\put(\value{x},0){\shead}%
\put(\value{y},#1){\line(0,-1){#1}}%
\put(\value{y},0){\shead}%
\end{picture}}}\end{picture}}%

\newcommand{\SQUADRIADJAR}[1]%
{\truex{400}%
\truey{1200}%
\begin{picture}(0,0)%
\put(0,0){\makebox(0,0)%
{\begin{picture}(0,#1)%
\put(-\value{y},#1){\line(0,-1){#1}}%
\put(-\value{y},0){\shead}%
\put(-\value{x},0){\line(0,1){#1}}%
\put(-\value{x},#1){\nhead}%
\put(\value{x},#1){\line(0,-1){#1}}%
\put(\value{x},0){\shead}%
\put(\value{y},0){\line(0,1){#1}}%
\put(\value{y},#1){\nhead}%
\end{picture}}}\end{picture}}%

\newcommand{\SQUINTIAR}[1]%
{\truex{600}%
\truey{1200}%
\begin{picture}(0,0)%
\put(0,0){\makebox(0,0)%
{\begin{picture}(0,#1)%
\put(-\value{y},#1){\line(0,-1){#1}}%
\put(-\value{y},0){\shead}%
\put(-\value{x},#1){\line(0,-1){#1}}%
\put(-\value{x},0){\shead}%
\put(0,#1){\line(0,-1){#1}}%
\put(0,0){\shead}%
\put(\value{x},#1){\line(0,-1){#1}}%
\put(\value{x},0){\shead}%
\put(\value{y},#1){\line(0,-1){#1}}%
\put(\value{y},0){\shead}%
\end{picture}}}\end{picture}}%

\newcommand{\SQUINTIADJAR}[1]%
{\truex{600}%
\truey{1200}%
\begin{picture}(0,0)%
\put(0,0){\makebox(0,0)%
{\begin{picture}(0,#1)%
\put(-\value{y},#1){\line(0,-1){#1}}%
\put(-\value{y},0){\shead}%
\put(-\value{x},0){\line(0,1){#1}}%
\put(-\value{x},#1){\nhead}%
\put(0,#1){\line(0,-1){#1}}%
\put(0,0){\shead}%
\put(\value{x},0){\line(0,1){#1}}%
\put(\value{x},#1){\nhead}%
\put(\value{y},#1){\line(0,-1){#1}}%
\put(\value{y},0){\shead}%
\end{picture}}}\end{picture}}%

\def\basicsar[#1]{\vcase{\SAR}{#100}}%

\newcommand{\sar}{\@ifnextchar[{\basicsar}{\basicsar[50]}}%

\def\basicSar[#1]#2{\Vcase{\SAR}{#2}{#100}}%

\newcommand{\Sar}{\@ifnextchar[{\basicSar}{\basicSar[50]}}%

\def\basicsaR[#1]#2{\vcasE{\SAR}{#2}{#100}}%

\newcommand{\saR}{\@ifnextchar[{\basicsaR}{\basicsaR[50]}}%

\def\basicsdist[#1]{\vcase{\SDIST}{#100}}%

\newcommand{\sdist}{\@ifnextchar[{\basicsdist}{\basicsdist[50]}}%

\def\basicSdist[#1]#2{\Vcase{\SDIST}{#2\hspace*{2pt}}{#100}}%

\newcommand{\Sdist}{\@ifnextchar[{\basicSdist}{\basicSdist[50]}}%

\def\basicsdisT[#1]#2{\vcasE{\SDIST}{\hspace*{2pt}#2}{#100}}%

\newcommand{\sdisT}{\@ifnextchar[{\basicsdisT}{\basicsdisT[50]}}%

\def\basicsdotar[#1]{\vcase{\SDOTAR}{#100}}%

\newcommand{\sdotar}{\@ifnextchar[{\basicsdotar}{\basicsdotar[50]}}%

\def\basicSdotar[#1]#2{\Vcase{\SDOTAR}{#2}{#100}}%

\newcommand{\Sdotar}{\@ifnextchar[{\basicSdotar}{\basicSdotar[50]}}%

\def\basicsdotaR[#1]#2{\vcasE{\SDOTAR}{#2}{#100}}%

\newcommand{\sdotaR}{\@ifnextchar[{\basicsdotaR}{\basicsdotaR[50]}}%

\def\basicsmono[#1]{\vcase{\SMONO}{#100}}%

\newcommand{\smono}{\@ifnextchar[{\basicsmono}{\basicsmono[50]}}%

\def\basicSmono[#1]#2{\Vcase{\SMONO}{#2}{#100}}%

\newcommand{\Smono}{\@ifnextchar[{\basicSmono}{\basicSmono[50]}}%

\def\basicsmonO[#1]#2{\vcasE{\SMONO}{#2}{#100}}%

\newcommand{\smonO}{\@ifnextchar[{\basicsmonO}{\basicsmonO[50]}}%

\def\basicsepi[#1]{\vcase{\SEPI}{#100}}%

\newcommand{\sepi}{\@ifnextchar[{\basicsepi}{\basicsepi[50]}}%

\def\basicSepi[#1]#2{\Vcase{\SEPI}{#2}{#100}}%

\newcommand{\Sepi}{\@ifnextchar[{\basicSepi}{\basicSepi[50]}}%

\def\basicsepI[#1]#2{\vcasE{\SEPI}{#2}{#100}}%

\newcommand{\sepI}{\@ifnextchar[{\basicsepI}{\basicsepI[50]}}%

\def\basicsbimo[#1]{\vcase{\SBIMO}{#100}}%

\newcommand{\sbimo}{\@ifnextchar[{\basicsbimo}{\basicsbimo[50]}}%

\def\basicSbimo[#1]#2{\Vcase{\SBIMO}{#2}{#100}}%

\newcommand{\Sbimo}{\@ifnextchar[{\basicSbimo}{\basicSbimo[50]}}%

\def\basicsbimO[#1]#2{\vcasE{\SBIMO}{#2}{#100}}%

\newcommand{\sbimO}{\@ifnextchar[{\basicsbimO}{\basicsbimO[50]}}%

\def\basicsiso[#1]{\Visocase{\SAR}{\isosign}{}{#100}}%

\newcommand{\siso}{\@ifnextchar[{\basicsiso}{\basicsiso[50]}}%

\def\basicSiso[#1]#2{\Visocase{\SAR}{#2}{\isosign}{#100}}%

\newcommand{\Siso}{\@ifnextchar[{\basicSiso}{\basicSiso[50]}}%

\def\basicsisO[#1]#2{\Visocase{\SAR}{\isosign}{#2}{#100}}%

\newcommand{\sisO}{\@ifnextchar[{\basicsisO}{\basicsisO[50]}}%

\def\basicseql[#1]{\vcase{\SEQL}{#100}}%

\newcommand{\seql}{\@ifnextchar[{\basicseql}{\basicseql[50]}}%

\def\basicSeql[#1]#2{\Vcase{\SEQL}{#2\hspace*{2pt}}{#100}}%

\newcommand{\Seql}{\@ifnextchar[{\basicSeql}{\basicSeql[50]}}%

\def\basicseqL[#1]#2{\vcasE{\SEQL}{\hspace*{2pt}#2}{#100}}%

\newcommand{\seqL}{\@ifnextchar[{\basicseqL}{\basicseqL[50]}}%

\def\basicsline[#1]{\vcase{\SLINE}{#100}}%

\newcommand{\sline}{\@ifnextchar[{\basicsline}{\basicsline[50]}}%

\def\basicSline[#1]#2{\Vcase{\SLINE}{#2\hspace*{2pt}}{#100}}%

\newcommand{\Sline}{\@ifnextchar[{\basicSline}{\basicSline[50]}}%

\def\basicslinE[#1]#2{\vcasE{\SLINE}{\hspace*{2pt}#2}{#100}}%

\newcommand{\slinE}{\@ifnextchar[{\basicslinE}{\basicslinE[50]}}%

\def\basicsbiar[#1]{\vbicase{\SBIAR}{#100}}%

\newcommand{\sbiar}{\@ifnextchar[{\basicsbiar}{\basicsbiar[50]}}%

\def\basicSbiar[#1]#2#3{\Vbicase{\SBIAR}{#2}{#3}{#100}}%

\newcommand{\Sbiar}{\@ifnextchar[{\basicSbiar}{\basicSbiar[50]}}%

\def\basicsbidist[#1]{\vbicase{\SBIDIST}{#100}}%

\newcommand{\sbidist}{\@ifnextchar[{\basicsbidist}{\basicsbidist[50]}}%

\def\basicSbidist[#1]#2#3%
{\Vbicase{\SBIDIST}{#2\hspace*{2pt}}{\hspace*{2pt}#3}{#100}}%

\newcommand{\Sbidist}{\@ifnextchar[{\basicSbidist}{\basicSbidist[50]}}%

\def\basicsadjar[#1]{\vbicase{\SADJAR}{#100}}%

\newcommand{\sadjar}{\@ifnextchar[{\basicsadjar}{\basicsadjar[50]}}%

\def\basicSadjar[#1]#2#3{\Vbicase{\SADJAR}{#2}{#3}{#100}}%

\newcommand{\Sadjar}{\@ifnextchar[{\basicSadjar}{\basicSadjar[50]}}%

\def\basicsadjdist[#1]{\vbicase{\SADJDIST}{#100}}%

\newcommand{\sadjdist}{\@ifnextchar[{\basicsadjdist}{\basicsadjdist[50]}}%

\def\basicSadjdist[#1]#2#3%
{\Vbicase{\SADJDIST}{#2\hspace*{2pt}}{\hspace*{2pt}#3}{#100}}%

\newcommand{\Sadjdist}{\@ifnextchar[{\basicSadjdist}{\basicSadjdist[50]}}%

\def\basicstriar[#1]{\vtricase{\STRIAR}{#100}}%

\newcommand{\striar}{\@ifnextchar[{\basicstriar}{\basicstriar[50]}}%

\def\basicStriar[#1]#2#3#4{\Vtricase{\STRIAR}{#2}{#3}{#4}{}{#100}}%

\newcommand{\Striar}{\@ifnextchar[{\basicStriar}{\basicStriar[50]}}%

\def\basicstriaR[#1]#2#3#4{\Vtricase{\STRIAR}{}{#2}{#3}{#4}{#100}}%

\newcommand{\striaR}{\@ifnextchar[{\basicstriaR}{\basicstriaR[50]}}%

\def\basicstriadjar[#1]{\vtricase{\STRIADJAR}{#100}}%

\newcommand{\striadjar}{\@ifnextchar[{\basicstriadjar}{\basicstriadjar[50]}}%

\def\basicStriadjar[#1]#2#3#4{\Vtricase{\STRIADJAR}{#2}{#3}{#4}{}{#100}}%

\newcommand{\Striadjar}{\@ifnextchar[{\basicStriadjar}{\basicStriadjar[50]}}%

\def\basicstriadjaR[#1]#2#3#4{\Vtricase{\STRIADJAR}{}{#2}{#3}{#4}{#100}}%

\newcommand{\striadjaR}{\@ifnextchar[{\basicstriadjaR}{\basicstriadjaR[50]}}%

\def\basicsquadriar[#1]{%
\Z=#1%
\multiply \Z by 100%
\vmulticase{\SQUADRIAR}{\Z}}%

\newcommand{\squadriar}{\@ifnextchar[{\basicsquadriar}%
{\hspace{\SOURCE\unitlength}\basicsquadriar[\ARROWLENGTH]}}%

\def\basicsquadriadjar[#1]{%
\Z=#1%
\multiply \Z by 100%
\vmulticase{\SQUADRIADJAR}{\Z}}%

\newcommand{\squadriadjar}{\@ifnextchar[{\basicsquadriadjar}%
{\hspace{\SOURCE\unitlength}\basicsquadriadjar[\ARROWLENGTH]}}%

\def\basicsquintiar[#1]{%
\Z=#1%
\multiply \Z by 100%
\vmulticase{\SQUINTIAR}{\Z}}%

\newcommand{\squintiar}{\@ifnextchar[{\basicsquintiar}%
{\hspace{\SOURCE\unitlength}\basicsquintiar[\ARROWLENGTH]}}%

\def\basicsquintiadjar[#1]{%
\Z=#1%
\multiply \Z by 100%
\vmulticase{\SQUINTIADJAR}{\Z}}%

\newcommand{\squintiadjar}{\@ifnextchar[{\basicsquintiadjar}%
{\hspace{\SOURCE\unitlength}\basicsquintiadjar[\ARROWLENGTH]}}%

\newcommand{\NAR}[1]%
{\begin{picture}(0,0)%
\put(0,0){\makebox(0,0)%
{\begin{picture}(0,#1)%
\put(0,0){\line(0,1){#1}}%
\put(0,#1){\nhead}%
\end{picture}}}\end{picture}}%

\newcommand{\NDIST}[1]%
{\begin{picture}(0,0)%
\put(0,0){\makebox(0,0)%
{\begin{picture}(0,#1)%
\put(0,0){\line(0,1){#1}}%
\put(0,#1){\nhead}%
\end{picture}}}
\put(0,0){\distsign}%
\end{picture}}%

\newcommand{\NDOTAR}[1]%
{\truex{100}\truey{300}%
\NUMBEROFDOTS=#1%
\divide\NUMBEROFDOTS by \value{y}%
\advance\NUMBEROFDOTS by 1%
\begin{picture}(0,0)%
\put(0,0){\makebox(0,0)%
{\begin{picture}(0,#1)%
\multiput(0,0)(0,\value{y}){\NUMBEROFDOTS}%
{\circle*{\value{x}}}%
\put(0,#1){\nhead}%
\end{picture}}}\end{picture}}%

\newcommand{\NMONO}[1]%
{\truetail%
\monolength=#1%
\advance\monolength by -\truemonotail%
\begin{picture}(0,0)%
\put(0,0){\makebox(0,0)%
{\begin{picture}(0,#1)%
\put(0,\truemonotail){\line(0,1){\monolength}}%
\put(0,#1){\nhead}%
\put(0,\truemonotail){\nhead}%
\end{picture}}}\end{picture}}%

\newcommand{\NEPI}[1]%
{\truehead%
\epilength=#1%
\advance\epilength by -\trueepihead%
\begin{picture}(0,0)%
\put(0,0){\makebox(0,0)%
{\begin{picture}(0,#1)%
\put(0,0){\line(0,1){\epilength}}%
\put(0,#1){\nhead}%
\put(0,\epilength){\nhead}%
\end{picture}}}\end{picture}}%

\newcommand{\NBIMO}[1]%
{\truehead\truetail%
\epilength=#1%
\advance\epilength by -\trueepihead%
\monolength=\epilength%
\advance\monolength by -\truemonotail%
\begin{picture}(0,0)%
\put(0,0){\makebox(0,0)%
{\begin{picture}(0,#1)%
\put(0,\truemonotail){\line(0,1){\monolength}}%
\put(0,#1){\nhead}%
\put(0,\truemonotail){\nhead}%
\put(0,\epilength){\nhead}%
\end{picture}}}\end{picture}}%

\newcommand{\NBIAR}[1]%
{\begin{picture}(0,0)%
\truex{350}%
\put(0,0){\makebox(0,0)%
{\begin{picture}(0,#1)%
\put(-\value{x},0){\line(0,1){#1}}%
\put(-\value{x},#1){\nhead}%
\put(\value{x},0){\line(0,1){#1}}%
\put(\value{x},#1){\nhead}%
\end{picture}}}\end{picture}}%

\newcommand{\NBIDIST}[1]%
{\begin{picture}(0,0)%
\truex{350}%
\put(0,0){\makebox(0,0)%
{\begin{picture}(0,#1)%
\put(-\value{x},0){\line(0,1){#1}}%
\put(-\value{x},#1){\nhead}%
\put(\value{x},0){\line(0,1){#1}}%
\put(\value{x},#1){\nhead}%
\end{picture}}}
\put(-\value{x},0){\distsign}%
\put(\value{x},0){\distsign}%
\end{picture}}%

\newcommand{\NADJAR}[1]{\begin{picture}(0,0)%
\truex{350}%
\put(0,0){\makebox(0,0)%
{\begin{picture}(0,#1)%
\put(\value{x},#1){\line(0,-1){#1}}%
\put(\value{x},0){\shead}%
\put(-\value{x},0){\line(0,1){#1}}%
\put(-\value{x},#1){\nhead}%
\end{picture}}}\end{picture}}%

\newcommand{\NADJDIST}[1]{\begin{picture}(0,0)%
\truex{350}%
\put(0,0){\makebox(0,0)%
{\begin{picture}(0,#1)%
\put(\value{x},#1){\line(0,-1){#1}}%
\put(\value{x},0){\shead}%
\put(-\value{x},0){\line(0,1){#1}}%
\put(-\value{x},#1){\nhead}%
\end{picture}}}
\put(-\value{x},0){\distsign}%
\put(\value{x},0){\distsign}%
\end{picture}}%

\newcommand{\NTRIAR}[5]%
{\truex{1200}%
\truey{1800}%
\truez{600}%
\begin{picture}(0,0)%
\put(0,0){\makebox(0,0)%
{\begin{picture}(0,#5)%
\X=#5\divide\X by 2%
\put(-\value{x},0){\line(0,1){#5}}%
\put(-\value{x},#5){\nhead}%
\put(0,0){\line(0,1){#5}}%
\put(0,#5){\nhead}%
\put(\value{x},0){\line(0,1){#5}}%
\put(\value{x},#5){\nhead}%
\put(-\value{y},\X){\makebox(0,0){$#1$}}%
\put(\value{y},\X){\makebox(0,0){$#4$}}%
\put(-\value{z},\X){\makebox(0,0){$#2$}}%
\put(\value{z},\X){\makebox(0,0){$#3$}}%
\end{picture}}}\end{picture}}%

\newcommand{\NTRIADJAR}[5]%
{\truex{1200}%
\truey{1800}%
\truez{600}%
\begin{picture}(0,0)%
\put(0,0){\makebox(0,0)%
{\begin{picture}(0,#5)%
\X=#5\divide\X by 2%
\put(-\value{x},#5){\line(0,-1){#5}}%
\put(-\value{x},0){\shead}%
\put(0,0){\line(0,1){#5}}%
\put(0,#5){\nhead}%
\put(\value{x},#5){\line(0,-1){#5}}%
\put(\value{x},0){\shead}%
\put(-\value{y},\X){\makebox(0,0){$#1$}}%
\put(\value{y},\X){\makebox(0,0){$#4$}}%
\put(-\value{z},\X){\makebox(0,0){$#2$}}%
\put(\value{z},\X){\makebox(0,0){$#3$}}%
\end{picture}}}\end{picture}}%

\newcommand{\NQUADRIAR}[1]%
{\truex{400}%
\truey{1200}%
\begin{picture}(0,0)%
\put(0,0){\makebox(0,0)%
{\begin{picture}(0,#1)%
\put(-\value{y},0){\line(0,1){#1}}%
\put(-\value{y},#1){\nhead}%
\put(-\value{x},0){\line(0,1){#1}}%
\put(-\value{x},#1){\nhead}%
\put(\value{x},0){\line(0,1){#1}}%
\put(\value{x},#1){\nhead}%
\put(\value{y},0){\line(0,1){#1}}%
\put(\value{y},#1){\nhead}%
\end{picture}}}\end{picture}}%

\newcommand{\NQUADRIADJAR}[1]%
{\truex{400}%
\truey{1200}%
\begin{picture}(0,0)%
\put(0,0){\makebox(0,0)%
{\begin{picture}(0,#1)%
\put(-\value{y},0){\line(0,1){#1}}%
\put(-\value{y},#1){\nhead}%
\put(-\value{x},#1){\line(0,-1){#1}}%
\put(-\value{x},0){\shead}%
\put(\value{x},0){\line(0,1){#1}}%
\put(\value{x},#1){\nhead}%
\put(\value{y},#1){\line(0,-1){#1}}%
\put(\value{y},0){\shead}%
\end{picture}}}\end{picture}}%

\newcommand{\NQUINTIAR}[1]%
{\truex{600}%
\truey{1200}%
\begin{picture}(0,0)%
\put(0,0){\makebox(0,0)%
{\begin{picture}(0,#1)%
\put(-\value{y},0){\line(0,1){#1}}%
\put(-\value{y},#1){\nhead}%
\put(-\value{x},0){\line(0,1){#1}}%
\put(-\value{x},#1){\nhead}%
\put(0,0){\line(0,1){#1}}%
\put(0,#1){\nhead}%
\put(\value{x},0){\line(0,1){#1}}%
\put(\value{x},#1){\nhead}%
\put(\value{y},0){\line(0,1){#1}}%
\put(\value{y},#1){\nhead}%
\end{picture}}}\end{picture}}%

\newcommand{\NQUINTIADJAR}[1]%
{\truex{600}%
\truey{1200}%
\begin{picture}(0,0)%
\put(0,0){\makebox(0,0)%
{\begin{picture}(0,#1)%
\put(-\value{y},0){\line(0,1){#1}}%
\put(-\value{y},#1){\nhead}%
\put(-\value{x},#1){\line(0,-1){#1}}%
\put(-\value{x},0){\shead}%
\put(0,0){\line(0,1){#1}}%
\put(0,#1){\nhead}%
\put(\value{x},#1){\line(0,-1){#1}}%
\put(\value{x},0){\shead}%
\put(\value{y},0){\line(0,1){#1}}%
\put(\value{y},#1){\nhead}%
\end{picture}}}\end{picture}}%

\def\basicnar[#1]{\vcase{\NAR}{#100}}%

\newcommand{\nar}{\@ifnextchar[{\basicnar}{\basicnar[50]}}%

\def\basicNar[#1]#2{\Vcase{\NAR}{#2}{#100}}%

\newcommand{\Nar}{\@ifnextchar[{\basicNar}{\basicNar[50]}}%

\def\basicnaR[#1]#2{\vcasE{\NAR}{#2}{#100}}%

\newcommand{\naR}{\@ifnextchar[{\basicnaR}{\basicnaR[50]}}%

\def\basicndist[#1]{\vcase{\NDIST}{#100}}%

\newcommand{\ndist}{\@ifnextchar[{\basicndist}{\basicndist[50]}}%

\def\basicNdist[#1]#2{\Vcase{\NDIST}{#2\hspace*{2pt}}{#100}}%

\newcommand{\Ndist}{\@ifnextchar[{\basicNdist}{\basicNdist[50]}}%

\def\basicndisT[#1]#2{\vcasE{\NDIST}{\hspace*{2pt}#2}{#100}}%

\newcommand{\ndisT}{\@ifnextchar[{\basicndisT}{\basicndisT[50]}}%

\def\basicndotar[#1]{\vcase{\NDOTAR}{#100}}%

\newcommand{\ndotar}{\@ifnextchar[{\basicndotar}{\basicndotar[50]}}%

\def\basicNdotar[#1]#2{\Vcase{\NDOTAR}{#2}{#100}}%

\newcommand{\Ndotar}{\@ifnextchar[{\basicNdotar}{\basicNdotar[50]}}%

\def\basicndotaR[#1]#2{\vcasE{\NDOTAR}{#2}{#100}}%

\newcommand{\ndotaR}{\@ifnextchar[{\basicndotaR}{\basicndotaR[50]}}%

\def\basicnmono[#1]{\vcase{\NMONO}{#100}}%

\newcommand{\nmono}{\@ifnextchar[{\basicnmono}%
{\basicnmono[50]}}%

\def\basicNmono[#1]#2{\Vcase{\NMONO}{#2}{#100}}%

\newcommand{\Nmono}{\@ifnextchar[{\basicNmono}{\basicNmono[50]}}%

\def\basicnmonO[#1]#2{\vcasE{\NMONO}{#2}{#100}}%

\newcommand{\nmonO}{\@ifnextchar[{\basicnmonO}{\basicnmonO[50]}}%

\def\basicnepi[#1]{\vcase{\NEPI}{#100}}%

\newcommand{\nepi}{\@ifnextchar[{\basicnepi}{\basicnepi[50]}}%

\def\basicNepi[#1]#2{\Vcase{\NEPI}{#2}{#100}}%

\newcommand{\Nepi}{\@ifnextchar[{\basicNepi}{\basicNepi[50]}}%

\def\basicnepI[#1]#2{\vcasE{\NEPI}{#2}{#100}}%

\newcommand{\nepI}{\@ifnextchar[{\basicnepI}{\basicnepI[50]}}%

\def\basicnbimo[#1]{\vcase{\NBIMO}{#100}}%

\newcommand{\nbimo}{\@ifnextchar[{\basicnbimo}{\basicnbimo[50]}}%

\def\basicNbimo[#1]#2{\Vcase{\NBIMO}{#2}{#100}}%

\newcommand{\Nbimo}{\@ifnextchar[{\basicNbimo}{\basicNbimo[50]}}%

\def\basicnbimO[#1]#2{\vcasE{\NBIMO}{#2}{#100}}%

\newcommand{\nbimO}{\@ifnextchar[{\basicnbimO}{\basicnbimO[50]}}%

\def\basicniso[#1]{\Visocase{\NAR}{\isosign}{}{#100}}%

\newcommand{\niso}{\@ifnextchar[{\basicniso}{\basicniso[50]}}%

\def\basicNiso[#1]#2{\Visocase{\NAR}{#2}{\isosign}{#100}}%

\newcommand{\Niso}{\@ifnextchar[{\basicNiso}{\basicNiso[50]}}%

\def\basicnisO[#1]#2{\Visocase{\NAR}{\isosign}{#2}{#100}}%

\newcommand{\nisO}{\@ifnextchar[{\basicnisO}{\basicnisO[50]}}%

\def\basicnbiar[#1]{\vbicase{\NBIAR}{#100}}%

\newcommand{\nbiar}{\@ifnextchar[{\basicnbiar}{\basicnbiar[50]}}%

\def\basicNbiar[#1]#2#3{\Vbicase{\NBIAR}{#2}{#3}{#100}}%

\newcommand{\Nbiar}{\@ifnextchar[{\basicNbiar}{\basicNbiar[50]}}%

\def\basicnbidist[#1]{\vbicase{\NBIDIST}{#100}}%

\newcommand{\nbidist}{\@ifnextchar[{\basicnbidist}{\basicnbidist[50]}}%

\def\basicNbidist[#1]#2#3%
{\Vbicase{\NBIDIST}{#2\hspace*{2pt}}{\hspace*{2pt}#3}{#100}}%

\newcommand{\Nbidist}{\@ifnextchar[{\basicNbidist}{\basicNbidist[50]}}%

\def\basicnadjar[#1]{\vbicase{\NADJAR}{#100}}%

\newcommand{\nadjar}{\@ifnextchar[{\basicnadjar}{\basicnadjar[50]}}%

\def\basicNadjar[#1]#2#3{\Vbicase{\NADJAR}{#2}{#3}{#100}}%

\newcommand{\Nadjar}{\@ifnextchar[{\basicNadjar}{\basicNadjar[50]}}%

\def\basicnadjdist[#1]{\vbicase{\NADJDIST}{#100}}%

\newcommand{\nadjdist}{\@ifnextchar[{\basicnadjdist}{\basicnadjdist[50]}}%

\def\basicNadjdist[#1]#2#3%
{\Vbicase{\NADJDIST}{#2\hspace*{2pt}}{\hspace*{2pt}#3}{#100}}%

\newcommand{\Nadjdist}{\@ifnextchar[{\basicNadjdist}{\basicNadjdist[50]}}%

\def\basicntriar[#1]{\vtricase{\NTRIAR}{#100}}%

\newcommand{\ntriar}{\@ifnextchar[{\basicntriar}{\basicntriar[50]}}%

\def\basicNtriar[#1]#2#3#4{\Vtricase{\NTRIAR}{#2}{#3}{#4}{}{#100}}%

\newcommand{\Ntriar}{\@ifnextchar[{\basicNtriar}{\basicNtriar[50]}}%

\def\basicntriaR[#1]#2#3#4{\Vtricase{\NTRIAR}{}{#2}{#3}{#4}{#100}}%

\newcommand{\ntriaR}{\@ifnextchar[{\basicntriaR}{\basicntriaR[50]}}%

\def\basicntriadjar[#1]{\vtricase{\NTRIADJAR}{#100}}%

\newcommand{\ntriadjar}{\@ifnextchar[{\basicntriadjar}{\basicntriadjar[50]}}%

\def\basicNtriadjar[#1]#2#3#4{\Vtricase{\NTRIADJAR}{#2}{#3}{#4}{}{#100}}%

\newcommand{\Ntriadjar}{\@ifnextchar[{\basicNtriadjar}{\basicNtriadjar[50]}}%

\def\basicntriadjaR[#1]#2#3#4{\Vtricase{\NTRIADJAR}{}{#2}{#3}{#4}{#100}}%

\newcommand{\ntriadjaR}{\@ifnextchar[{\basicntriadjaR}{\basicntriadjaR[50]}}%

\def\basicnquadriar[#1]{%
\Z=#1%
\multiply \Z by 100%
\vmulticase{\NQUADRIAR}{\Z}}%

\newcommand{\nquadriar}{\@ifnextchar[{\basicnquadriar}%
{\hspace{\SOURCE\unitlength}\basicnquadriar[\ARROWLENGTH]}}%

\def\basicnquadriadjar[#1]{%
\Z=#1%
\multiply \Z by 100%
\vmulticase{\NQUADRIADJAR}{\Z}}%

\newcommand{\nquadriadjar}{\@ifnextchar[{\basicnquadriadjar}%
{\hspace{\SOURCE\unitlength}\basicnquadriadjar[\ARROWLENGTH]}}%

\def\basicnquintiar[#1]{%
\Z=#1%
\multiply \Z by 100%
\vmulticase{\NQUINTIAR}{\Z}}%

\newcommand{\nquintiar}{\@ifnextchar[{\basicnquintiar}%
{\hspace{\SOURCE\unitlength}\basicnquintiar[\ARROWLENGTH]}}%

\def\basicnquintiadjar[#1]{%
\Z=#1%
\multiply \Z by 100%
\vmulticase{\NQUINTIADJAR}{\Z}}%

\newcommand{\nquintiadjar}{\@ifnextchar[{\basicnquintiadjar}%
{\hspace{\SOURCE\unitlength}\basicnquintiadjar[\ARROWLENGTH]}}%

\newcommand{\fdcase}[4]{\begin{picture}(0,0)%
\put(0,0){#1{#4}}%
\truex{200}\truey{600}\truez{600}%
\put(-\value{x},-\value{x}){\makebox(0,\value{z})[r]{${#2}$}}%
\put(\value{x},-\value{y}){\makebox(0,\value{z})[l]{${#3}$}}%
\end{picture}}%

\newcommand{\fdbicase}[4]{\begin{picture}(0,0)%
\put(0,0){#1{#4}}%
\truex{600}\truey{150}\truez{600}%
\put(-\value{x},\value{y}){\makebox(0,\value{z})[r]{${#2}$}}%
\truex{300}\truey{1050}\truez{600}%
\put(\value{x},-\value{y}){\makebox(0,\value{z})[l]{${#3}$}}%
\end{picture}}%

\newcommand{\NEAR}[1]{%
\Y=#1%
\divide\Y by 2%
\begin{picture}(0,0)%
\put(-\Y,-\Y){\line(1,1){#1}}%
\put(\Y,\Y){\nehead}%
\end{picture}}%

\newcommand{\NEDIST}[1]{%
\Y=#1%
\divide\Y by 2%
\begin{picture}(0,0)%
\put(-\Y,-\Y){\line(1,1){#1}}%
\put(\Y,\Y){\nehead}%
\put(0,0){\distsign}%
\end{picture}}%

\newcommand{\NEDOTAR}[1]%
{\truex{100}\truey{212}%
\Y=#1%
\divide\Y by 2%
\NUMBEROFDOTS=#1%
\divide\NUMBEROFDOTS by \value{y}%
\advance\NUMBEROFDOTS by 1%
\begin{picture}(0,0)%
\multiput(-\Y,-\Y)(\value{y},\value{y}){\NUMBEROFDOTS}%
{\circle*{\value{x}}}%
\put(\Y,\Y){\nehead}%
\end{picture}}%

\newcommand{\NEMONO}[1]{%
\Y=#1%
\divide \Y by 2%
\Truetail%
\bimolength=#1%
\advance\bimolength by -\Truemonotail%
\monolength=\bimolength%
\advance\monolength by -\Y%
\begin{picture}(0,0)%
\put(-\monolength,-\monolength){\line(1,1){\bimolength}}%
\put(-\monolength,-\monolength){\nehead}%
\put(\Y,\Y){\nehead}%
\end{picture}}%

\newcommand{\NEEPI}[1]{%
\Y=#1%
\divide\Y by 2%
\Truehead%
\bimolength=#1%
\advance\bimolength by -\Trueepihead%
\epilength=\bimolength%
\advance\epilength by -\Y%
\begin{picture}(0,0)%
\put(-\Y,-\Y){\line(1,1){\bimolength}}%
\put(\epilength,\epilength){\nehead}%
\put(\Y,\Y){\nehead}%
\end{picture}}%

\newcommand{\NEBIMO}[1]{%
\Y=#1%
\divide\Y by 2%
\Truetail\Truehead%
\bimolength=#1%
\advance\bimolength by -\Truemonotail%
\monolength=\bimolength%
\advance\monolength by -\Y%
\advance\bimolength by -\Trueepihead%
\epilength=\bimolength%
\advance\epilength by -\monolength%
\begin{picture}(0,0)%
\put(-\monolength,-\monolength){\line(1,1){\bimolength}}%
\put(-\monolength,-\monolength){\nehead}%
\put(\epilength,\epilength){\nehead}%
\put(\Y,\Y){\nehead}%
\end{picture}}%

\newcommand{\NEBIAR}[1]{%
\Y=#1%
\divide\Y by 2%
\begin{picture}(0,0)%
\put(-\Y,-\Y){\begin{picture}(0,0)%
\truex{247}%
\put(-\value{x},\value{x}){\line(1,1){#1}}%
\put(\value{x},-\value{x}){\line(1,1){#1}}%
\monolength=#1%
\advance\monolength by -\value{x}%
\epilength=#1%
\advance\epilength by \value{x}%
\put(\monolength,\epilength){\nehead}%
\put(\epilength,\monolength){\nehead}%
\end{picture}}\end{picture}}%

\newcommand{\NEBIDIST}[1]{%
\Y=#1%
\divide\Y by 2%
\begin{picture}(0,0)%
\put(-\Y,-\Y){\begin{picture}(0,0)%
\truex{247}%
\monolength=#1%
\advance\monolength by -\value{x}%
\epilength=#1%
\advance\epilength by \value{x}%
\put(\value{x},-\value{x}){\line(1,1){#1}}%
\put(\epilength,\monolength){\nehead}%
\end{picture}}%
\put(-\Y,-\Y){\begin{picture}(0,0)%
\truex{247}%
\monolength=#1%
\advance\monolength by \value{x}%
\epilength=#1%
\advance\epilength by -\value{x}%
\put(-\value{x},\value{x}){\line(1,1){#1}}%
\put(\epilength,\monolength){\nehead}%
\end{picture}}%
\put(-\value{x},\value{x}){\distsign}%
\put(\value{x},-\value{x}){\distsign}%
\end{picture}}%

\newcommand{\NEEQL}[1]{%
\Y=#1%
\divide\Y by 2%
\begin{picture}(0,0)%
\put(-\Y,-\Y){\begin{picture}(0,0)%
\truex{70}%
\put(-\value{x},\value{x}){\line(1,1){#1}}%
\put(\value{x},-\value{x}){\line(1,1){#1}}%
\end{picture}}\end{picture}}%

\newcommand{\NELINE}[1]{%
\Y=#1%
\divide\Y by 2%
\begin{picture}(0,0)%
\put(-\Y,-\Y){\begin{picture}(0,0)%
\put(0,0){\line(1,1){#1}}%
\end{picture}}\end{picture}}%

\newcommand{\NEADJAR}[1]{%
\Y=#1%
\divide\Y by 2%
\begin{picture}(0,0)%
\put(-\Y,-\Y){\begin{picture}(0,0)%
\truex{247}%
\monolength=#1%
\advance\monolength by -\value{x}%
\epilength=#1%
\advance\epilength by \value{x}%
\put(\value{x},-\value{x}){\line(1,1){#1}}%
\put(\epilength,\monolength){\nehead}%
\end{picture}}%
\put(\Y,\Y){\begin{picture}(0,0)%
\truex{247}%
\monolength=#1%
\advance\monolength by -\value{x}%
\epilength=#1%
\advance\epilength by \value{x}%
\put(-\value{x},\value{x}){\line(-1,-1){#1}}%
\put(-\epilength,-\monolength){\swhead}%
\end{picture}}\end{picture}}%

\newcommand{\NEADJDIST}[1]{%
\Y=#1%
\divide\Y by 2%
\begin{picture}(0,0)%
\put(-\Y,-\Y){\begin{picture}(0,0)%
\truex{247}%
\monolength=#1%
\advance\monolength by -\value{x}%
\epilength=#1%
\advance\epilength by \value{x}%
\put(\value{x},-\value{x}){\line(1,1){#1}}%
\put(\epilength,\monolength){\nehead}%
\end{picture}}%
\put(\Y,\Y){\begin{picture}(0,0)%
\truex{247}%
\monolength=#1%
\advance\monolength by -\value{x}%
\epilength=#1%
\advance\epilength by \value{x}%
\put(-\value{x},\value{x}){\line(-1,-1){#1}}%
\put(-\epilength,-\monolength){\swhead}%
\end{picture}}%
\put(-\value{x},\value{x}){\distsign}%
\put(\value{x},-\value{x}){\distsign}%
\end{picture}}%

\def\basicnear[#1]{\fdcase{\NEAR}{}{}{#100}}%

\newcommand{\near}{\@ifnextchar[{\basicnear}{\basicnear[59]}}%

\def\basicNear[#1]#2{\fdcase{\NEAR}{#2}{}{#100}}%

\newcommand{\Near}{\@ifnextchar[{\basicNear}{\basicNear[59]}}%

\def\basicneaR[#1]#2{\fdcase{\NEAR}{}{#2}{#100}}%

\newcommand{\neaR}{\@ifnextchar[{\basicneaR}{\basicneaR[59]}}%

\def\basicnedist[#1]{\fdcase{\NEDIST}{}{}{#100}}%

\newcommand{\nedist}{\@ifnextchar[{\basicnedist}{\basicnedist[59]}}%

\def\basicNedist[#1]#2{\fdcase{\NEDIST}{#2}{}{#100}}%

\newcommand{\Nedist}{\@ifnextchar[{\basicNedist}{\basicNedist[59]}}%

\def\basicnedisT[#1]#2{\fdcase{\NEDIST}{}{#2}{#100}}%

\newcommand{\nedisT}{\@ifnextchar[{\basicnedisT}{\basicnedisT[59]}}%

\def\basicnedotar[#1]{\fdcase{\NEDOTAR}{}{}{#100}}%

\newcommand{\nedotar}{\@ifnextchar[{\basicnedotar}{\basicnedotar[59]}}%

\def\basicNedotar[#1]#2{\fdcase{\NEDOTAR}{#2}{}{#100}}%

\newcommand{\Nedotar}{\@ifnextchar[{\basicNedotar}{\basicNedotar[59]}}%

\def\basicnedotaR[#1]#2{\fdcase{\NEDOTAR}{}{#2}{#100}}%

\newcommand{\nedotaR}{\@ifnextchar[{\basicnedotaR}{\basicnedotaR[59]}}%

\def\basicnemono[#1]{\fdcase{\NEMONO}{}{}{#100}}%

\newcommand{\nemono}{\@ifnextchar[{\basicnemono}{\basicnemono[59]}}%

\def\basicNemono[#1]#2{\fdcase{\NEMONO}{#2}{}{#100}}%

\newcommand{\Nemono}{\@ifnextchar[{\basicNemono}{\basicNemono[59]}}%

\def\basicnemonO[#1]#2{\fdcase{\NEMONO}{}{#2}{#100}}%

\newcommand{\nemonO}{\@ifnextchar[{\basicnemonO}{\basicnemonO[59]}}%

\def\basicneepi[#1]{\fdcase{\NEEPI}{}{}{#100}}%

\newcommand{\neepi}{\@ifnextchar[{\basicneepi}{\basicneepi[59]}}%

\def\basicNeepi[#1]#2{\fdcase{\NEEPI}{#2}{}{#100}}%

\newcommand{\Neepi}{\@ifnextchar[{\basicNeepi}{\basicNeepi[59]}}%

\def\basicneepI[#1]#2{\fdcase{\NEEPI}{}{#2}{#100}}%

\newcommand{\neepI}{\@ifnextchar[{\basicneepI}{\basicneepI[59]}}%

\def\basicnebimo[#1]{\fdcase{\NEBIMO}{}{}{#100}}%

\newcommand{\nebimo}{\@ifnextchar[{\basicnebimo}{\basicnebimo[59]}}%

\def\basicNebimo[#1]#2{\fdcase{\NEBIMO}{#2}{}{#100}}%

\newcommand{\Nebimo}{\@ifnextchar[{\basicNebimo}{\basicNebimo[59]}}%

\def\basicnebimO[#1]#2{\fdcase{\NEBIMO}{}{#2}{#100}}%

\newcommand{\nebimO}{\@ifnextchar[{\basicnebimO}{\basicnebimO[59]}}%

\def\basicneiso[#1]{\fdcase{\NEAR}{\hspace{-2pt}\isosign}{}{#100}}%

\newcommand{\neiso}{\@ifnextchar[{\basicneiso}{\basicneiso[59]}}%

\def\basicNeiso[#1]#2{\fdcase{\NEAR}{#2}{\isosign}{#100}}%

\newcommand{\Neiso}{\@ifnextchar[{\basicNeiso}{\basicNeiso[59]}}%

\def\basicneisO[#1]#2{\fdcase{\NEAR}{\hspace{-2pt}\isosign}{#2}{#100}}%

\newcommand{\neisO}{\@ifnextchar[{\basicneisO}{\basicneisO[59]}}%

\def\basicneeql[#1]{\fdcase{\NEEQL}{}{}{#100}}%

\newcommand{\neeql}{\@ifnextchar[{\basicneeql}{\basicneeql[59]}}%

\def\basicNeeql[#1]#2{\fdcase{\NEEQL}{#2}{}{#100}}%

\newcommand{\Neeql}{\@ifnextchar[{\basicNeeql}{\basicNeeql[59]}}%

\def\basicneeqL[#1]#2{\fdcase{\NEEQL}{}{#2}{#100}}%

\newcommand{\neeqL}{\@ifnextchar[{\basicneeqL}{\basicneeqL[59]}}%

\def\basicneline[#1]{\fdcase{\NELINE}{}{}{#100}}%

\newcommand{\neline}{\@ifnextchar[{\basicneline}{\basicneline[59]}}%

\def\basicNeline[#1]#2{\fdcase{\NELINE}{#2}{}{#100}}%

\newcommand{\Neline}{\@ifnextchar[{\basicNeline}{\basicNeline[59]}}%

\def\basicnelinE[#1]#2{\fdcase{\NELINE}{}{#2}{#100}}%

\newcommand{\nelinE}{\@ifnextchar[{\basicnelinE}{\basicneeqL[59]}}%

\def\basicnebiar[#1]{\fdbicase{\NEBIAR}{}{}{#100}}%

\newcommand{\nebiar}{\@ifnextchar[{\basicnebiar}{\basicnebiar[59]}}%

\def\basicNebiar[#1]#2#3{\fdbicase{\NEBIAR}{#2}{#3}{#100}}%

\newcommand{\Nebiar}{\@ifnextchar[{\basicNebiar}{\basicNebiar[59]}}%

\def\basicneadjar[#1]{\fdbicase{\NEADJAR}{}{}{#100}}%

\newcommand{\neadjar}{\@ifnextchar[{\basicneadjar}{\basicneadjar[59]}}%

\def\basicNeadjar[#1]#2#3{\fdbicase{\NEADJAR}{#2}{#3}{#100}}%

\newcommand{\Neadjar}{\@ifnextchar[{\basicNeadjar}{\basicNeadjar[59]}}%

\def\basicnebidist[#1]{\fdbicase{\NEBIDIST}{}{}{#100}}%

\newcommand{\nebidist}{\@ifnextchar[{\basicnebidist}{\basicnebidist[59]}}%

\def\basicNebidist[#1]#2#3{\fdbicase{\NEBIDIST}{#2}{#3}{#100}}%

\newcommand{\Nebidist}{\@ifnextchar[{\basicNebidist}{\basicNebidist[59]}}%

\def\basicneadjdist[#1]{\fdbicase{\NEADJDIST}{}{}{#100}}%

\newcommand{\neadjdist}{\@ifnextchar[{\basicneadjdist}{\basicneadjdist[59]}}%

\def\basicNeadjdist[#1]#2#3{\fdbicase{\NEADJDIST}{#2}{#3}{#100}}%

\newcommand{\Neadjdist}{\@ifnextchar[{\basicNeadjdist}{\basicNeadjdist[59]}}%

\newcommand{\SWAR}[1]{%
\Y=#1%
\divide\Y by 2%
\begin{picture}(0,0)%
\put(\Y,\Y){\line(-1,-1){#1}}%
\put(-\Y,-\Y){\swhead}%
\end{picture}}%

\newcommand{\SWDIST}[1]{%
\Y=#1%
\divide\Y by 2%
\begin{picture}(0,0)%
\put(\Y,\Y){\line(-1,-1){#1}}%
\put(-\Y,-\Y){\swhead}%
\put(0,0){\distsign}%
\end{picture}}%

\newcommand{\SWDOTAR}[1]%
{\truex{100}\truey{212}%
\Y=#1%
\divide\Y by 2%
\NUMBEROFDOTS=#1%
\divide\NUMBEROFDOTS by \value{y}%
\advance\NUMBEROFDOTS by 1%
\begin{picture}(0,0)%
\multiput(\Y,\Y)(-\value{y},-\value{y}){\NUMBEROFDOTS}%
{\circle*{\value{x}}}%
\put(-\Y,-\Y){\swhead}%
\end{picture}}%

\newcommand{\SWMONO}[1]{%
\Y=#1%
\divide \Y by 2%
\Truetail%
\bimolength=#1%
\advance\bimolength by -\Truemonotail%
\monolength=\bimolength%
\advance\monolength by -\Y%
\begin{picture}(0,0)%
\put(\monolength,\monolength){\line(-1,-1){\bimolength}}%
\put(\monolength,\monolength){\swhead}%
\put(-\Y,-\Y){\swhead}%
\end{picture}}%

\newcommand{\SWEPI}[1]{%
\Y=#1%
\divide\Y by 2%
\Truehead%
\bimolength=#1%
\advance\bimolength by -\Trueepihead%
\epilength=\bimolength%
\advance\epilength by -\Y%
\begin{picture}(0,0)%
\put(\Y,\Y){\line(-1,-1){\bimolength}}%
\put(-\epilength,-\epilength){\swhead}%
\put(-\Y,-\Y){\swhead}%
\end{picture}}%

\newcommand{\SWBIMO}[1]{%
\Y=#1%
\divide\Y by 2%
\Truetail\Truehead%
\bimolength=#1%
\advance\bimolength by -\Truemonotail%
\monolength=\bimolength%
\advance\monolength by -\Y%
\advance\bimolength by -\Trueepihead%
\epilength=\bimolength%
\advance\epilength by -\monolength%
\begin{picture}(0,0)%
\put(\monolength,\monolength){\line(-1,-1){\bimolength}}%
\put(\monolength,\monolength){\swhead}%
\put(-\epilength,-\epilength){\swhead}%
\put(-\Y,-\Y){\swhead}%
\end{picture}}%

\newcommand{\SWBIAR}[1]{%
\Y=#1%
\divide\Y by 2%
\begin{picture}(0,0)%
\put(\Y,\Y){\begin{picture}(0,0)%
\truex{247}%
\put(\value{x},-\value{x}){\line(-1,-1){#1}}%
\put(-\value{x},\value{x}){\line(-1,-1){#1}}%
\monolength=#1%
\advance\monolength by -\value{x}%
\epilength=#1%
\advance\epilength by \value{x}%
\put(-\monolength,-\epilength){\swhead}%
\put(-\epilength,-\monolength){\swhead}%
\end{picture}}\end{picture}}%

\newcommand{\SWBIDIST}[1]{%
\Y=#1%
\divide\Y by 2%
\begin{picture}(0,0)%
\put(\Y,\Y){\begin{picture}(0,0)%
\truex{247}%
\monolength=#1%
\advance\monolength by -\value{x}%
\epilength=#1%
\advance\epilength by \value{x}%
\put(-\value{x},\value{x}){\line(-1,-1){#1}}%
\put(-\epilength,-\monolength){\swhead}%
\end{picture}}%
\put(\Y,\Y){\begin{picture}(0,0)%
\truex{247}%
\monolength=#1%
\advance\monolength by \value{x}%
\epilength=#1%
\advance\epilength by -\value{x}%
\put(\value{x},-\value{x}){\line(-1,-1){#1}}%
\put(-\epilength,-\monolength){\swhead}%
\end{picture}}%
\put(\value{x},-\value{x}){\distsign}%
\put(-\value{x},\value{x}){\distsign}%
\end{picture}}%

\newcommand{\SWADJAR}[1]{%
\Y=#1%
\divide\Y by 2%
\begin{picture}(0,0)%
\put(\Y,\Y){\begin{picture}(0,0)%
\truex{247}%
\monolength=#1%
\advance\monolength by -\value{x}%
\epilength=#1%
\advance\epilength by \value{x}%
\put(\value{x},-\value{x}){\line(-1,-1){#1}}%
\put(-\monolength,-\epilength){\swhead}%
\end{picture}}%
\put(-\Y,-\Y){\begin{picture}(0,0)%
\truex{247}%
\monolength=#1%
\advance\monolength by -\value{x}%
\epilength=#1%
\advance\epilength by \value{x}%
\put(-\value{x},\value{x}){\line(1,1){#1}}%
\put(\monolength,\epilength){\nehead}%
\end{picture}}\end{picture}}%

\newcommand{\SWADJDIST}[1]{%
\Y=#1%
\divide\Y by 2%
\begin{picture}(0,0)%
\put(\Y,\Y){\begin{picture}(0,0)%
\truex{247}%
\monolength=#1%
\advance\monolength by -\value{x}%
\epilength=#1%
\advance\epilength by \value{x}%
\put(\value{x},-\value{x}){\line(-1,-1){#1}}%
\put(-\monolength,-\epilength){\swhead}%
\end{picture}}%
\put(-\Y,-\Y){\begin{picture}(0,0)%
\truex{247}%
\monolength=#1%
\advance\monolength by -\value{x}%
\epilength=#1%
\advance\epilength by \value{x}%
\put(-\value{x},\value{x}){\line(1,1){#1}}%
\put(\monolength,\epilength){\nehead}%
\end{picture}}%
\put(-\value{x},\value{x}){\distsign}%
\put(\value{x},-\value{x}){\distsign}%
\end{picture}}%

\def\basicswar[#1]{\fdcase{\SWAR}{}{}{#100}}%

\newcommand{\swar}{\@ifnextchar[{\basicswar}{\basicswar[59]}}%

\def\basicSwar[#1]#2{\fdcase{\SWAR}{#2}{}{#100}}%

\newcommand{\Swar}{\@ifnextchar[{\basicSwar}{\basicSwar[59]}}%

\def\basicswaR[#1]#2{\fdcase{\SWAR}{}{#2}{#100}}%

\newcommand{\swaR}{\@ifnextchar[{\basicswaR}{\basicswaR[59]}}%

\def\basicswdist[#1]{\fdcase{\SWDIST}{}{}{#100}}%

\newcommand{\swdist}{\@ifnextchar[{\basicswdist}{\basicswdist[59]}}%

\def\basicSwdist[#1]#2{\fdcase{\SWDIST}{#2}{}{#100}}%

\newcommand{\Swdist}{\@ifnextchar[{\basicSwdist}{\basicSwdist[59]}}%

\def\basicswdisT[#1]#2{\fdcase{\SWDIST}{}{#2}{#100}}%

\newcommand{\swdisT}{\@ifnextchar[{\basicswdisT}{\basicswdisT[59]}}%

\def\basicswdotar[#1]{\fdcase{\SWDOTAR}{}{}{#100}}%

\newcommand{\swdotar}{\@ifnextchar[{\basicswdotar}{\basicswdotar[59]}}%

\def\basicSwdotar[#1]#2{\fdcase{\SWDOTAR}{#2}{}{#100}}%

\newcommand{\Swdotar}{\@ifnextchar[{\basicSwdotar}{\basicSwdotar[59]}}%

\def\basicswdotaR[#1]#2{\fdcase{\SWDOTAR}{}{#2}{#100}}%

\newcommand{\swdotaR}{\@ifnextchar[{\basicswdotaR}{\basicswdotaR[59]}}%

\def\basicswmono[#1]{\fdcase{\SWMONO}{}{}{#100}}%

\newcommand{\swmono}{\@ifnextchar[{\basicswmono}{\basicswmono[59]}}%

\def\basicSwmono[#1]#2{\fdcase{\SWMONO}{#2}{}{#100}}%

\newcommand{\Swmono}{\@ifnextchar[{\basicSwmono}{\basicSwmono[59]}}%

\def\basicswmonO[#1]#2{\fdcase{\SWMONO}{}{#2}{#100}}%

\newcommand{\swmonO}{\@ifnextchar[{\basicswmonO}{\basicswmonO[59]}}%

\def\basicswepi[#1]{\fdcase{\SWEPI}{}{}{#100}}%

\newcommand{\swepi}{\@ifnextchar[{\basicswepi}{\basicswepi[59]}}%

\def\basicSwepi[#1]#2{\fdcase{\SWEPI}{#2}{}{#100}}%

\newcommand{\Swepi}{\@ifnextchar[{\basicSwepi}{\basicSwepi[59]}}%

\def\basicswepI[#1]#2{\fdcase{\SWEPI}{}{#2}{#100}}%

\newcommand{\swepI}{\@ifnextchar[{\basicswepI}{\basicswepI[59]}}%

\def\basicswbimo[#1]{\fdcase{\SWBIMO}{}{}{#100}}%

\newcommand{\swbimo}{\@ifnextchar[{\basicswbimo}{\basicswbimo[59]}}%

\def\basicSwbimo[#1]#2{\fdcase{\SWBIMO}{#2}{}{#100}}%

\newcommand{\Swbimo}{\@ifnextchar[{\basicSwbimo}{\basicSwbimo[59]}}%

\def\basicswbimO[#1]#2{\fdcase{\SWBIMO}{}{#2}{#100}}%

\newcommand{\swbimO}{\@ifnextchar[{\basicswbimO}{\basicswbimO[59]}}%

\def\basicswiso[#1]{\fdcase{\SWAR}{\hspace{-2pt}\isosign}{}{#100}}%

\newcommand{\swiso}{\@ifnextchar[{\basicswiso}{\basicswiso[59]}}%

\def\basicSwiso[#1]#2{\fdcase{\SWAR}{#2}{\isosign}{#100}}%

\newcommand{\Swiso}{\@ifnextchar[{\basicSwiso}{\basicSwiso[59]}}%

\def\basicswisO[#1]#2{\fdcase{\SWAR}{\hspace{-2pt}\isosign}{#2}{#100}}%

\newcommand{\swisO}{\@ifnextchar[{\basicswisO}{\basicswisO[59]}}%

\def\basicswbiar[#1]{\fdbicase{\SWBIAR}{}{}{#100}}%

\newcommand{\swbiar}{\@ifnextchar[{\basicswbiar}{\basicswbiar[59]}}%

\def\basicSwbiar[#1]#2#3{\fdbicase{\SWBIAR}{#2}{#3}{#100}}%

\newcommand{\Swbiar}{\@ifnextchar[{\basicSwbiar}{\basicSwbiar[59]}}%

\def\basicswadjar[#1]{\fdbicase{\SWADJAR}{}{}{#100}}%

\newcommand{\swadjar}{\@ifnextchar[{\basicswadjar}{\basicswadjar[59]}}%

\def\basicSwadjar[#1]#2#3{\fdbicase{\SWADJAR}{#2}{#3}{#100}}%

\newcommand{\Swadjar}{\@ifnextchar[{\basicSwadjar}{\basicSwadjar[59]}}%

\def\basicswbidist[#1]{\fdbicase{\SWBIDIST}{}{}{#100}}%

\newcommand{\swbidist}{\@ifnextchar[{\basicswbidist}{\basicswbidist[59]}}%

\def\basicSwbidist[#1]#2#3{\fdbicase{\SWBIDIST}{#2}{#3}{#100}}%

\newcommand{\Swbidist}{\@ifnextchar[{\basicSwbidist}{\basicSwbidist[59]}}%

\def\basicswadjdist[#1]{\fdbicase{\SWADJDIST}{}{}{#100}}%

\newcommand{\swadjdist}{\@ifnextchar[{\basicswadjdist}{\basicswadjdist[59]}}%

\def\basicSwadjdist[#1]#2#3{\fdbicase{\SWADJDIST}{#2}{#3}{#100}}%

\newcommand{\Swadjdist}{\@ifnextchar[{\basicSwadjdist}{\basicSwadjdist[59]}}%

\newcommand{\sdcase}[4]{\begin{picture}(0,0)%
\put(0,0){#1{#4}}%
\truex{100}\truez{600}%
\put(\value{x},\value{x}){\makebox(0,\value{z})[l]{${#2}$}}%
\truex{300}\truey{800}%
\put(-\value{x},-\value{y}){\makebox(0,\value{z})[r]{${#3}$}}%
\end{picture}}%

\newcommand{\sdbicase}[4]{\begin{picture}(0,0)%
\put(0,0){#1{#4}}%
\truex{350}\truey{600}\truez{950}%
\put(\value{x},\value{x}){\makebox(0,\value{y})[l]{${#2}$}}%
\truex{450}\truey{600}\truez{1050}%
\put(-\value{x},-\value{z}){\makebox(0,\value{y})[r]{${#3}$}}%
\end{picture}}%

\newcommand{\SEAR}[1]{%
\Y=#1%
\divide\Y by 2%
\begin{picture}(0,0)%
\put(-\Y,\Y){\line(1,-1){#1}}%
\put(\Y,-\Y){\sehead}%
\end{picture}}%

\newcommand{\SEDIST}[1]{%
\Y=#1%
\divide\Y by 2%
\begin{picture}(0,0)%
\put(-\Y,\Y){\line(1,-1){#1}}%
\put(\Y,-\Y){\sehead}%
\put(0,0){\distsign}%
\end{picture}}%

\newcommand{\SEDOTAR}[1]%
{\truex{100}\truey{212}%
\Y=#1%
\divide\Y by 2%
\NUMBEROFDOTS=#1%
\divide\NUMBEROFDOTS by \value{y}%
\advance\NUMBEROFDOTS by 1%
\begin{picture}(0,0)%
\multiput(-\Y,\Y)(\value{y},-\value{y}){\NUMBEROFDOTS}%
{\circle*{\value{x}}}%
\put(\Y,-\Y){\sehead}%
\end{picture}}%

\newcommand{\SEMONO}[1]{%
\Y=#1%
\divide \Y by 2%
\Truetail%
\bimolength=#1%
\advance\bimolength by -\Truemonotail%
\monolength=\bimolength%
\advance\monolength by -\Y%
\begin{picture}(0,0)%
\put(-\monolength,\monolength){\line(1,-1){\bimolength}}%
\put(-\monolength,\monolength){\sehead}%
\put(\Y,-\Y){\sehead}%
\end{picture}}%

\newcommand{\SEEPI}[1]{%
\Y=#1%
\divide\Y by 2%
\Truehead%
\bimolength=#1%
\advance\bimolength by -\Trueepihead%
\epilength=\bimolength%
\advance\epilength by -\Y%
\begin{picture}(0,0)%
\put(-\Y,\Y){\line(1,-1){\bimolength}}%
\put(\epilength,-\epilength){\sehead}%
\put(\Y,-\Y){\sehead}%
\end{picture}}%

\newcommand{\SEBIMO}[1]{%
\Y=#1%
\divide\Y by 2%
\Truetail\Truehead%
\bimolength=#1%
\advance\bimolength by -\Truemonotail%
\monolength=\bimolength%
\advance\monolength by -\Y%
\advance\bimolength by -\Trueepihead%
\epilength=\bimolength%
\advance\epilength by -\monolength%
\begin{picture}(0,0)%
\put(-\monolength,\monolength){\line(1,-1){\bimolength}}%
\put(-\monolength,\monolength){\sehead}%
\put(\epilength,-\epilength){\sehead}%
\put(\Y,-\Y){\sehead}%
\end{picture}}%

\newcommand{\SEBIAR}[1]{%
\Y=#1%
\divide\Y by 2%
\begin{picture}(0,0)%
\put(-\Y,\Y){\begin{picture}(0,0)%
\truex{247}%
\put(-\value{x},-\value{x}){\line(1,-1){#1}}%
\put(\value{x},\value{x}){\line(1,-1){#1}}%
\monolength=#1%
\advance\monolength by -\value{x}%
\epilength=#1%
\advance\epilength by \value{x}%
\put(\monolength,-\epilength){\sehead}%
\put(\epilength,-\monolength){\sehead}%
\end{picture}}\end{picture}}%

\newcommand{\SEBIDIST}[1]{%
\Y=#1%
\divide\Y by 2%
\begin{picture}(0,0)%
\put(-\Y,\Y){\begin{picture}(0,0)%
\truex{247}%
\monolength=#1%
\advance\monolength by -\value{x}%
\epilength=#1%
\advance\epilength by \value{x}%
\put(\value{x},\value{x}){\line(1,-1){#1}}%
\put(\epilength,-\monolength){\sehead}%
\end{picture}}%
\put(-\Y,\Y){\begin{picture}(0,0)%
\truex{247}%
\monolength=#1%
\advance\monolength by \value{x}%
\epilength=#1%
\advance\epilength by -\value{x}%
\put(-\value{x},-\value{x}){\line(1,-1){#1}}%
\put(\epilength,-\monolength){\sehead}%
\end{picture}}%
\put(-\value{x},-\value{x}){\distsign}%
\put(\value{x},\value{x}){\distsign}%
\end{picture}}%

\newcommand{\SEEQL}[1]{%
\Y=#1%
\divide\Y by 2%
\begin{picture}(0,0)%
\put(-\Y,\Y){\begin{picture}(0,0)%
\truex{70}%
\put(-\value{x},-\value{x}){\line(1,-1){#1}}%
\put(\value{x},\value{x}){\line(1,-1){#1}}%
\end{picture}}\end{picture}}%

\newcommand{\SELINE}[1]{%
\Y=#1%
\divide\Y by 2%
\begin{picture}(0,0)%
\put(-\Y,\Y){\begin{picture}(0,0)%
\put(0,0){\line(1,-1){#1}}%
\end{picture}}\end{picture}}%

\newcommand{\SEADJAR}[1]{%
\Y=#1%
\divide\Y by 2%
\begin{picture}(0,0)%
\put(-\Y,\Y){\begin{picture}(0,0)%
\truex{247}%
\monolength=#1%
\advance\monolength by -\value{x}%
\epilength=#1%
\advance\epilength by \value{x}%
\put(-\value{x},-\value{x}){\line(1,-1){#1}}%
\put(\monolength,-\epilength){\sehead}%
\end{picture}}%
\put(\Y,-\Y){\begin{picture}(0,0)%
\truex{247}%
\monolength=#1%
\advance\monolength by -\value{x}%
\epilength=#1%
\advance\epilength by \value{x}%
\put(\value{x},\value{x}){\line(-1,1){#1}}%
\put(-\monolength,\epilength){\nwhead}%
\end{picture}}\end{picture}}%

\newcommand{\SEADJDIST}[1]{%
\Y=#1%
\divide\Y by 2%
\begin{picture}(0,0)%
\put(-\Y,\Y){\begin{picture}(0,0)%
\truex{247}%
\monolength=#1%
\advance\monolength by -\value{x}%
\epilength=#1%
\advance\epilength by \value{x}%
\put(-\value{x},-\value{x}){\line(1,-1){#1}}%
\put(\monolength,-\epilength){\sehead}%
\end{picture}}%
\put(\Y,-\Y){\begin{picture}(0,0)%
\truex{247}%
\monolength=#1%
\advance\monolength by -\value{x}%
\epilength=#1%
\advance\epilength by \value{x}%
\put(\value{x},\value{x}){\line(-1,1){#1}}%
\put(-\monolength,\epilength){\nwhead}%
\end{picture}}%
\put(-\value{x},-\value{x}){\distsign}%
\put(\value{x},\value{x}){\distsign}%
\end{picture}}%

\def\basicsear[#1]{\sdcase{\SEAR}{}{}{#100}}%

\newcommand{\sear}{\@ifnextchar[{\basicsear}{\basicsear[59]}}%

\def\basicSear[#1]#2{\sdcase{\SEAR}{#2}{}{#100}}%

\newcommand{\Sear}{\@ifnextchar[{\basicSear}{\basicSear[59]}}%

\def\basicseaR[#1]#2{\sdcase{\SEAR}{}{#2}{#100}}%

\newcommand{\seaR}{\@ifnextchar[{\basicseaR}{\basicseaR[59]}}%

\def\basicsedist[#1]{\sdcase{\SEDIST}{}{}{#100}}%

\newcommand{\sedist}{\@ifnextchar[{\basicsedist}{\basicsedist[59]}}%

\def\basicSedist[#1]#2{\sdcase{\SEDIST}{#2}{}{#100}}%

\newcommand{\Sedist}{\@ifnextchar[{\basicSedist}{\basicSedist[59]}}%

\def\basicsedisT[#1]#2{\sdcase{\SEDIST}{}{#2}{#100}}%

\newcommand{\sedisT}{\@ifnextchar[{\basicsedisT}{\basicsedisT[59]}}%

\def\basicsedotar[#1]{\sdcase{\SEDOTAR}{}{}{#100}}%

\newcommand{\sedotar}{\@ifnextchar[{\basicsedotar}{\basicsedotar[59]}}%

\def\basicSedotar[#1]#2{\sdcase{\SEDOTAR}{#2}{}{#100}}%

\newcommand{\Sedotar}{\@ifnextchar[{\basicSedotar}{\basicSedotar[59]}}%

\def\basicsedotaR[#1]#2{\sdcase{\SEDOTAR}{}{#2}{#100}}%

\newcommand{\sedotaR}{\@ifnextchar[{\basicsedotaR}{\basicsedotaR[59]}}%

\def\basicsemono[#1]{\sdcase{\SEMONO}{}{}{#100}}%

\newcommand{\semono}{\@ifnextchar[{\basicsemono}{\basicsemono[59]}}%

\def\basicSemono[#1]#2{\sdcase{\SEMONO}{#2}{}{#100}}%

\newcommand{\Semono}{\@ifnextchar[{\basicSemono}{\basicSemono[59]}}%

\def\basicsemonO[#1]#2{\sdcase{\SEMONO}{}{#2}{#100}}%

\newcommand{\semonO}{\@ifnextchar[{\basicsemonO}{\basicsemonO[59]}}%

\def\basicseepi[#1]{\sdcase{\SEEPI}{}{}{#100}}%

\newcommand{\seepi}{\@ifnextchar[{\basicseepi}{\basicseepi[59]}}%

\def\basicSeepi[#1]#2{\sdcase{\SEEPI}{#2}{}{#100}}%

\newcommand{\Seepi}{\@ifnextchar[{\basicSeepi}{\basicSeepi[59]}}%

\def\basicseepI[#1]#2{\sdcase{\SEEPI}{}{#2}{#100}}%

\newcommand{\seepI}{\@ifnextchar[{\basicseepI}{\basicseepI[59]}}%

\def\basicsebimo[#1]{\sdcase{\SEBIMO}{}{}{#100}}%

\newcommand{\sebimo}{\@ifnextchar[{\basicsebimo}{\basicsebimo[59]}}%

\def\basicSebimo[#1]#2{\sdcase{\SEBIMO}{#2}{}{#100}}%

\newcommand{\Sebimo}{\@ifnextchar[{\basicSebimo}{\basicSebimo[59]}}%

\def\basicsebimO[#1]#2{\sdcase{\SEBIMO}{}{#2}{#100}}%

\newcommand{\sebimO}{\@ifnextchar[{\basicsebimO}{\basicsebimO[59]}}%

\def\basicseiso[#1]{\sdcase{\SEAR}{\hspace{-2pt}\isosign}{}{#100}}%

\newcommand{\seiso}{\@ifnextchar[{\basicseiso}{\basicseiso[59]}}%

\def\basicSeiso[#1]#2{\sdcase{\SEAR}{#2}{\isosign}{#100}}%

\newcommand{\Seiso}{\@ifnextchar[{\basicSeiso}{\basicSeiso[59]}}%

\def\basicseisO[#1]#2{\sdcase{\SEAR}{\hspace{-2pt}\isosign}{#2}{#100}}%

\newcommand{\seisO}{\@ifnextchar[{\basicseisO}{\basicseisO[59]}}%

\def\basicseeql[#1]{\sdcase{\SEEQL}{}{}{#100}}%

\newcommand{\seeql}{\@ifnextchar[{\basicseeql}{\basicseeql[59]}}%

\def\basicSeeql[#1]#2{\sdcase{\SEEQL}{#2}{}{#100}}%

\newcommand{\Seeql}{\@ifnextchar[{\basicSeeql}{\basicSeeql[59]}}%

\def\basicseeqL[#1]#2{\sdcase{\SEEQL}{}{#2}{#100}}%

\newcommand{\seeqL}{\@ifnextchar[{\basicseeqL}{\basicseeqL[59]}}%

\def\basicseline[#1]{\sdcase{\SELINE}{}{}{#100}}%

\newcommand{\seline}{\@ifnextchar[{\basicseline}{\basicseline[59]}}%

\def\basicSeline[#1]#2{\sdcase{\SELINE}{#2}{}{#100}}%

\newcommand{\Seline}{\@ifnextchar[{\basicSeline}{\basicSeline[59]}}%

\def\basicselinE[#1]#2{\sdcase{\SELINE}{}{#2}{#100}}%

\newcommand{\selinE}{\@ifnextchar[{\basicselinE}{\basicselinE[59]}}%

\def\basicsebiar[#1]{\sdbicase{\SEBIAR}{}{}{#100}}%

\newcommand{\sebiar}{\@ifnextchar[{\basicsebiar}{\basicsebiar[59]}}%

\def\basicSebiar[#1]#2#3{\sdbicase{\SEBIAR}{#2}{#3}{#100}}%

\newcommand{\Sebiar}{\@ifnextchar[{\basicSebiar}{\basicSebiar[59]}}%

\def\basicseadjar[#1]{\sdbicase{\SEADJAR}{}{}{#100}}%

\newcommand{\seadjar}{\@ifnextchar[{\basicseadjar}{\basicseadjar[59]}}%

\def\basicSeadjar[#1]#2#3{\sdbicase{\SEADJAR}{#2}{#3}{#100}}%

\newcommand{\Seadjar}{\@ifnextchar[{\basicSeadjar}{\basicSeadjar[59]}}%

\def\basicsebidist[#1]{\sdbicase{\SEBIDIST}{}{}{#100}}%

\newcommand{\sebidist}{\@ifnextchar[{\basicsebidist}{\basicsebidist[59]}}%

\def\basicSebidist[#1]#2#3{\sdbicase{\SEBIDIST}{#2}{#3}{#100}}%

\newcommand{\Sebidist}{\@ifnextchar[{\basicSebidist}{\basicSebidist[59]}}%

\def\basicseadjdist[#1]{\sdbicase{\SEADJDIST}{}{}{#100}}%

\newcommand{\seadjdist}{\@ifnextchar[{\basicseadjdist}{\basicseadjdist[59]}}%

\def\basicSeadjdist[#1]#2#3{\sdbicase{\SEADJDIST}{#2}{#3}{#100}}%

\newcommand{\Seadjdist}{\@ifnextchar[{\basicSeadjdist}{\basicSeadjdist[59]}}%

\newcommand{\NWAR}[1]{%
\Y=#1%
\divide\Y by 2%
\begin{picture}(0,0)%
\put(\Y,-\Y){\line(-1,1){#1}}%
\put(-\Y,\Y){\nwhead}%
\end{picture}}%

\newcommand{\NWDIST}[1]{%
\Y=#1%
\divide\Y by 2%
\begin{picture}(0,0)%
\put(\Y,-\Y){\line(-1,1){#1}}%
\put(-\Y,\Y){\nwhead}%
\put(0,0){\distsign}%
\end{picture}}%

\newcommand{\NWDOTAR}[1]%
{\truex{100}\truey{212}%
\Y=#1%
\divide\Y by 2%
\NUMBEROFDOTS=#1%
\divide\NUMBEROFDOTS by \value{y}%
\advance\NUMBEROFDOTS by 1%
\begin{picture}(0,0)%
\multiput(\Y,-\Y)(-\value{y},\value{y}){\NUMBEROFDOTS}%
{\circle*{\value{x}}}%
\put(-\Y,\Y){\nwhead}%
\end{picture}}%

\newcommand{\NWMONO}[1]{%
\Y=#1%
\divide \Y by 2%
\Truetail%
\bimolength=#1%
\advance\bimolength by -\Truemonotail%
\monolength=\bimolength%
\advance\monolength by -\Y%
\begin{picture}(0,0)%
\put(\monolength,-\monolength){\line(-1,1){\bimolength}}%
\put(\monolength,-\monolength){\nwhead}%
\put(-\Y,\Y){\nwhead}%
\end{picture}}%

\newcommand{\NWEPI}[1]{%
\Y=#1%
\divide\Y by 2%
\Truehead%
\bimolength=#1%
\advance\bimolength by -\Trueepihead%
\epilength=\bimolength%
\advance\epilength by -\Y%
\begin{picture}(0,0)%
\put(\Y,-\Y){\line(-1,1){\bimolength}}%
\put(-\epilength,\epilength){\nwhead}%
\put(-\Y,\Y){\nwhead}%
\end{picture}}%

\newcommand{\NWBIMO}[1]{%
\Y=#1%
\divide\Y by 2%
\Truetail\Truehead%
\bimolength=#1%
\advance\bimolength by -\Truemonotail%
\monolength=\bimolength%
\advance\monolength by -\Y%
\advance\bimolength by -\Trueepihead%
\epilength=\bimolength%
\advance\epilength by -\monolength%
\begin{picture}(0,0)%
\put(\monolength,-\monolength){\line(-1,1){\bimolength}}%
\put(\monolength,-\monolength){\nwhead}%
\put(-\epilength,\epilength){\nwhead}%
\put(-\Y,\Y){\nwhead}%
\end{picture}}%

\newcommand{\NWBIAR}[1]{%
\Y=#1%
\divide\Y by 2%
\begin{picture}(0,0)%
\put(\Y,-\Y){\begin{picture}(0,0)%
\truex{247}%
\put(-\value{x},-\value{x}){\line(-1,1){#1}}%
\put(\value{x},\value{x}){\line(-1,1){#1}}%
\monolength=#1%
\advance\monolength by -\value{x}%
\epilength=#1%
\advance\epilength by \value{x}%
\put(-\monolength,\epilength){\nwhead}%
\put(-\epilength,\monolength){\nwhead}%
\end{picture}}\end{picture}}%

\newcommand{\NWBIDIST}[1]{%
\Y=#1%
\divide\Y by 2%
\begin{picture}(0,0)%
\put(\Y,-\Y){\begin{picture}(0,0)%
\truex{247}%
\monolength=#1%
\advance\monolength by -\value{x}%
\epilength=#1%
\advance\epilength by \value{x}%
\put(-\value{x},-\value{x}){\line(-1,1){#1}}%
\put(-\epilength,\monolength){\nwhead}%
\end{picture}}%
\put(\Y,-\Y){\begin{picture}(0,0)%
\truex{247}%
\monolength=#1%
\advance\monolength by \value{x}%
\epilength=#1%
\advance\epilength by -\value{x}%
\put(\value{x},\value{x}){\line(-1,1){#1}}%
\put(-\epilength,\monolength){\nwhead}%
\end{picture}}%
\put(-\value{x},-\value{x}){\distsign}%
\put(\value{x},\value{x}){\distsign}%
\end{picture}}%

\newcommand{\NWADJAR}[1]{%
\Y=#1%
\divide\Y by 2%
\begin{picture}(0,0)%
\put(\Y,-\Y){\begin{picture}(0,0)%
\truex{247}%
\monolength=#1%
\advance\monolength by -\value{x}%
\epilength=#1%
\advance\epilength by \value{x}%
\put(-\value{x},-\value{x}){\line(-1,1){#1}}%
\put(-\epilength,\monolength){\nwhead}%
\end{picture}}%
\put(-\Y,\Y){\begin{picture}(0,0)%
\truex{247}%
\monolength=#1%
\advance\monolength by -\value{x}%
\epilength=#1%
\advance\epilength by \value{x}%
\put(\value{x},\value{x}){\line(1,-1){#1}}%
\put(\epilength,-\monolength){\sehead}%
\end{picture}}\end{picture}}%

\newcommand{\NWADJDIST}[1]{%
\Y=#1%
\divide\Y by 2%
\begin{picture}(0,0)%
\put(\Y,-\Y){\begin{picture}(0,0)%
\truex{247}%
\monolength=#1%
\advance\monolength by -\value{x}%
\epilength=#1%
\advance\epilength by \value{x}%
\put(-\value{x},-\value{x}){\line(-1,1){#1}}%
\put(-\epilength,\monolength){\nwhead}%
\end{picture}}%
\put(-\Y,\Y){\begin{picture}(0,0)%
\truex{247}%
\monolength=#1%
\advance\monolength by -\value{x}%
\epilength=#1%
\advance\epilength by \value{x}%
\put(\value{x},\value{x}){\line(1,-1){#1}}%
\put(\epilength,-\monolength){\sehead}%
\end{picture}}%
\put(-\value{x},-\value{x}){\distsign}%
\put(\value{x},\value{x}){\distsign}%
\end{picture}}%

\def\basicnwar[#1]{\sdcase{\NWAR}{}{}{#100}}%

\newcommand{\nwar}{\@ifnextchar[{\basicnwar}{\basicnwar[59]}}%

\def\basicNwar[#1]#2{\sdcase{\NWAR}{#2}{}{#100}}%

\newcommand{\Nwar}{\@ifnextchar[{\basicNwar}{\basicNwar[59]}}%

\def\basicnwaR[#1]#2{\sdcase{\NWAR}{}{#2}{#100}}%

\newcommand{\nwaR}{\@ifnextchar[{\basicnwaR}{\basicnwaR[59]}}%

\def\basicnwdist[#1]{\sdcase{\NWDIST}{}{}{#100}}%

\newcommand{\nwdist}{\@ifnextchar[{\basicnwdist}{\basicnwdist[59]}}%

\def\basicNwdist[#1]#2{\sdcase{\NWDIST}{#2}{}{#100}}%

\newcommand{\Nwdist}{\@ifnextchar[{\basicNwdist}{\basicNwdist[59]}}%

\def\basicnwdisT[#1]#2{\sdcase{\NWDIST}{}{#2}{#100}}%

\newcommand{\nwdisT}{\@ifnextchar[{\basicnwdisT}{\basicnwdisT[59]}}%

\def\basicnwdotar[#1]{\sdcase{\NWDOTAR}{}{}{#100}}%

\newcommand{\nwdotar}{\@ifnextchar[{\basicnwdotar}{\basicnwdotar[59]}}%

\def\basicNwdotar[#1]#2{\sdcase{\NWDOTAR}{#2}{}{#100}}%

\newcommand{\Nwdotar}{\@ifnextchar[{\basicNwdotar}{\basicNwdotar[59]}}%

\def\basicnwdotaR[#1]#2{\sdcase{\NWDOTAR}{}{#2}{#100}}%

\newcommand{\nwdotaR}{\@ifnextchar[{\basicnwdotaR}{\basicnwdotaR[59]}}%

\def\basicnwmono[#1]{\sdcase{\NWMONO}{}{}{#100}}%

\newcommand{\nwmono}{\@ifnextchar[{\basicnwmono}{\basicnwmono[59]}}%

\def\basicNwmono[#1]#2{\sdcase{\NWMONO}{#2}{}{#100}}%

\newcommand{\Nwmono}{\@ifnextchar[{\basicNwmono}{\basicNwmono[59]}}%

\def\basicnwmonO[#1]#2{\sdcase{\NWMONO}{}{#2}{#100}}%

\newcommand{\nwmonO}{\@ifnextchar[{\basicnwmonO}{\basicnwmonO[59]}}%

\def\basicnwepi[#1]{\sdcase{\NWEPI}{}{}{#100}}%

\newcommand{\nwepi}{\@ifnextchar[{\basicnwepi}{\basicnwepi[59]}}%

\def\basicNwepi[#1]#2{\sdcase{\NWEPI}{#2}{}{#100}}%

\newcommand{\Nwepi}{\@ifnextchar[{\basicNwepi}{\basicNwepi[59]}}%

\def\basicnwepI[#1]#2{\sdcase{\NWEPI}{}{#2}{#100}}%

\newcommand{\nwepI}{\@ifnextchar[{\basicnwepI}{\basicnwepI[59]}}%

\def\basicnwbimo[#1]{\sdcase{\NWBIMO}{}{}{#100}}%

\newcommand{\nwbimo}{\@ifnextchar[{\basicnwbimo}{\basicnwbimo[59]}}%

\def\basicNwbimo[#1]#2{\sdcase{\NWBIMO}{#2}{}{#100}}%

\newcommand{\Nwbimo}{\@ifnextchar[{\basicNwbimo}{\basicNwbimo[59]}}%

\def\basicnwbimO[#1]#2{\sdcase{\NWBIMO}{}{#2}{#100}}%

\newcommand{\nwbimO}{\@ifnextchar[{\basicnwbimO}{\basicnwbimO[59]}}%

\def\basicnwiso[#1]{\sdcase{\NWAR}{\hspace{-2pt}\isosign}{}{#100}}%

\newcommand{\nwiso}{\@ifnextchar[{\basicnwiso}{\basicnwiso[59]}}%

\def\basicNwiso[#1]#2{\sdcase{\NWAR}{#2}{\isosign}{#100}}%

\newcommand{\Nwiso}{\@ifnextchar[{\basicNwiso}{\basicNwiso[59]}}%

\def\basicnwisO[#1]#2{\sdcase{\NWAR}{\hspace{-2pt}\isosign}{#2}{#100}}%

\newcommand{\nwisO}{\@ifnextchar[{\basicnwisO}{\basicnwisO[59]}}%

\def\basicnwbiar[#1]{\sdbicase{\NWBIAR}{}{}{#100}}%

\newcommand{\nwbiar}{\@ifnextchar[{\basicnwbiar}{\basicnwbiar[59]}}%

\def\basicNwbiar[#1]#2#3{\sdbicase{\NWBIAR}{#2}{#3}{#100}}%

\newcommand{\Nwbiar}{\@ifnextchar[{\basicNwbiar}{\basicNwbiar[59]}}%

\def\basicnwadjar[#1]{\sdbicase{\NWADJAR}{}{}{#100}}%

\newcommand{\nwadjar}{\@ifnextchar[{\basicnwadjar}{\basicnwadjar[59]}}%

\def\basicNwadjar[#1]#2#3{\sdbicase{\NWADJAR}{#2}{#3}{#100}}%

\newcommand{\Nwadjar}{\@ifnextchar[{\basicNwadjar}{\basicNwadjar[59]}}%

\def\basicnwbidist[#1]{\sdbicase{\NWBIDIST}{}{}{#100}}%

\newcommand{\nwbidist}{\@ifnextchar[{\basicnwbidist}{\basicnwbidist[59]}}%

\def\basicNwbidist[#1]#2#3{\sdbicase{\NWBIDIST}{#2}{#3}{#100}}%

\newcommand{\Nwbidist}{\@ifnextchar[{\basicNwbidist}{\basicNwbidist[59]}}%

\def\basicnwadjdist[#1]{\sdbicase{\NWADJDIST}{}{}{#100}}%

\newcommand{\nwadjdist}{\@ifnextchar[{\basicnwadjdist}{\basicnwadjdist[59]}}%

\def\basicNwadjdist[#1]#2#3{\sdbicase{\NWADJDIST}{#2}{#3}{#100}}%

\newcommand{\Nwadjdist}{\@ifnextchar[{\basicNwadjdist}{\basicNwadjdist[59]}}%

\newcommand{\ENEAR}[3]{%
\Y=#3%
\divide\Y by 2%
\Z=\Y%
\divide\Z by 2%
\begin{picture}(0,0)%
\put(-\Y,-\Z){\line(2,1){#3}}%
\put(\Y,\Z){\enehead}%
\truex{200}\truey{800}\truez{600}%
\put(-\value{x},\value{x}){\makebox(0,\value{z})[r]{${#1}$}}%
\put(\value{x},-\value{y}){\makebox(0,\value{z})[l]{${#2}$}}%
\end{picture}}%

\newcommand{\ENEDIST}[3]{%
\Y=#3%
\divide\Y by 2%
\Z=\Y%
\divide\Z by 2%
\begin{picture}(0,0)%
\put(-\Y,-\Z){\line(2,1){#3}}%
\put(\Y,\Z){\enehead}%
\put(0,0){\distsign}%
\truex{200}\truey{800}\truez{600}%
\put(-\value{x},\value{x}){\makebox(0,\value{z})[r]{${#1}$}}%
\put(\value{x},-\value{y}){\makebox(0,\value{z})[l]{${#2}$}}%
\end{picture}}%

\newcommand{\ENEDOTAR}[3]{%
\truex{100}\truey{268}\truez{134}%
\Y=#3%
\divide\Y by 2%
\Z=\Y%
\divide\Z by 2%
\NUMBEROFDOTS=#3%
\divide\NUMBEROFDOTS by \value{y}%
\advance\NUMBEROFDOTS by 1%
\begin{picture}(0,0)%
\multiput(-\Y,-\Z)(\value{y},\value{z}){\NUMBEROFDOTS}%
{\circle*{\value{x}}}%
\put(\Y,\Z){\enehead}%
\truex{200}\truey{800}\truez{600}%
\put(-\value{x},\value{x}){\makebox(0,\value{z})[r]{${#1}$}}%
\put(\value{x},-\value{y}){\makebox(0,\value{z})[l]{${#2}$}}%
\end{picture}}%

\newcommand{\ENEMONO}[3]{%
\Y=#3%
\divide\Y by 2%
\Z=\Y%
\divide\Z by 2%
\TrueTail%
\bimolength=#3%
\advance\bimolength by -\TrueMonoTail%
\monolength=\bimolength%
\advance\monolength by -\Y%
\secondmonolength=\monolength%
\divide\secondmonolength by 2%
\begin{picture}(0,0)%
\put(-\monolength,-\secondmonolength){\line(2,1){\bimolength}}%
\put(-\monolength,-\secondmonolength){\enehead}%
\put(\Y,\Z){\enehead}%
\truex{200}\truey{800}\truez{600}%
\put(-\value{x},\value{x}){\makebox(0,\value{z})[r]{${#1}$}}%
\put(\value{x},-\value{y}){\makebox(0,\value{z})[l]{${#2}$}}%
\end{picture}}%

\newcommand{\ENEEPI}[3]{%
\Y=#3%
\divide\Y by 2%
\Z=\Y%
\divide\Z by 2%
\TrueHead%
\bimolength=#3%
\advance\bimolength by -\TrueEpiHead%
\epilength=\bimolength%
\advance\epilength by -\Y%
\secondepilength=\epilength%
\divide\secondepilength by 2%
\begin{picture}(0,0)%
\put(-\Y,-\Z){\line(2,1){\bimolength}}%
\put(\epilength,\secondepilength){\enehead}%
\put(\Y,\Z){\enehead}%
\truex{200}\truey{800}\truez{600}%
\put(-\value{x},\value{x}){\makebox(0,\value{z})[r]{${#1}$}}%
\put(\value{x},-\value{y}){\makebox(0,\value{z})[l]{${#2}$}}%
\end{picture}}%

\newcommand{\ENEBIMO}[3]{%
\Y=#3%
\divide\Y by 2%
\Z=\Y%
\divide\Z by 2%
\TrueTail\TrueHead%
\bimolength=#3%
\advance\bimolength by -\TrueMonoTail%
\monolength=\bimolength%
\advance\monolength by -\Y%
\advance\bimolength by -\TrueEpiHead%
\epilength=\bimolength%
\advance\epilength by -\monolength%
\secondmonolength=\monolength%
\divide\secondmonolength by 2%
\secondepilength=\epilength%
\divide\secondepilength by 2%
\begin{picture}(0,0)%
\put(-\monolength,-\secondmonolength){\line(2,1){\bimolength}}%
\put(-\monolength,-\secondmonolength){\enehead}%
\put(\epilength,\secondepilength){\enehead}%
\put(\Y,\Z){\enehead}%
\truex{200}\truey{800}\truez{600}%
\put(-\value{x},\value{x}){\makebox(0,\value{z})[r]{${#1}$}}%
\put(\value{x},-\value{y}){\makebox(0,\value{z})[l]{${#2}$}}%
\end{picture}}%

\newcommand{\ENEEQL}[3]{%
\Y=#3%
\divide\Y by 2%
\Z=\Y%
\divide\Z by 2%
\begin{picture}(0,0)%
\put(-\Y,-\Z){\begin{picture}(0,0)%
\truex{44}\truey{89}%
\put(-\value{x},\value{y}){\line(2,1){#3}}%
\put(\value{x},-\value{y}){\line(2,1){#3}}%
\end{picture}}%
\truex{200}\truey{800}\truez{600}%
\put(-\value{x},\value{x}){\makebox(0,\value{z})[r]{${#1}$}}%
\put(\value{x},-\value{y}){\makebox(0,\value{z})[l]{${#2}$}}%
\end{picture}}%

\newcommand{\ENEBIAR}[3]{%
\Y=#3%
\divide\Y by 2%
\Z=\Y%
\divide\Z by 2%
\begin{picture}(0,0)%
\put(-\Y,-\Z){\begin{picture}(0,0)%
\truex{156}\truey{313}%
\put(-\value{x},\value{y}){\line(2,1){#3}}%
\put(\value{x},-\value{y}){\line(2,1){#3}}%
\monolength=#3%
\advance\monolength by -\value{x}%
\epilength=#3%
\advance\epilength by \value{x}%
\secondmonolength=\Y%
\advance\secondmonolength by -\value{y}%
\secondepilength=\Y%
\advance\secondepilength by \value{y}%
\put(\monolength,\secondepilength){\enehead}%
\put(\epilength,\secondmonolength){\enehead}%
\end{picture}}
\truex{300}\truey{1000}\truez{600}%
\put(-\value{x},\value{x}){\makebox(0,\value{z})[r]{${#1}$}}%
\put(\value{x},-\value{y}){\makebox(0,\value{z})[l]{${#2}$}}%
\end{picture}}%

\newcommand{\ENEBIDIST}[3]{%
\Y=#3%
\divide\Y by 2%
\Z=\Y%
\divide\Z by 2%
\begin{picture}(0,0)%
\truex{156}\truey{313}\truez{400}%
\put(-\Y,-\Z){\begin{picture}(0,0)%
\put(-\value{x},\value{y}){\line(2,1){#3}}%
\put(\value{x},-\value{y}){\line(2,1){#3}}%
\monolength=#3%
\advance\monolength by -\value{x}%
\epilength=#3%
\advance\epilength by \value{x}%
\secondmonolength=\Y%
\advance\secondmonolength by -\value{y}%
\secondepilength=\Y%
\advance\secondepilength by \value{y}%
\put(\monolength,\secondepilength){\enehead}%
\put(\epilength,\secondmonolength){\enehead}%
\end{picture}}
\put(-\value{x},\value{y}){\distsign}%
\put(\value{x},-\value{y}){\distsign}%
\truex{300}\truey{1000}\truez{600}%
\put(-\value{x},\value{x}){\makebox(0,\value{z})[r]{${#1}$}}%
\put(\value{x},-\value{y}){\makebox(0,\value{z})[l]{${#2}$}}%
\end{picture}}%

\newcommand{\ENEADJAR}[3]{%
\Y=#3%
\divide\Y by 2%
\Z=\Y%
\divide\Z by 2%
\begin{picture}(0,0)%
\put(-\Y,-\Z){\begin{picture}(0,0)%
\truex{156}\truey{313}%
\monolength=#3%
\advance\monolength by -\value{x}%
\epilength=#3%
\advance\epilength by \value{x}%
\secondmonolength=\Y%
\advance\secondmonolength by -\value{y}%
\secondepilength=\Y%
\advance\secondepilength by \value{y}%
\put(\value{x},-\value{y}){\line(2,1){#3}}%
\put(\epilength,\secondmonolength){\enehead}%
\put(\monolength,\secondepilength){\line(-2,-1){#3}}%
\put(-\value{x},\value{y}){\wswhead}%
\end{picture}}
\truex{300}\truey{1000}\truez{600}%
\put(-\value{x},\value{x}){\makebox(0,\value{z})[r]{${#1}$}}%
\put(\value{x},-\value{y}){\makebox(0,\value{z})[l]{${#2}$}}%
\end{picture}}%

\newcommand{\ENEADJDIST}[3]{%
\Y=#3%
\divide\Y by 2%
\Z=\Y%
\divide\Z by 2%
\begin{picture}(0,0)%
\truex{156}\truey{313}\truez{400}%
\put(-\Y,-\Z){\begin{picture}(0,0)%
\monolength=#3%
\advance\monolength by -\value{x}%
\epilength=#3%
\advance\epilength by \value{x}%
\secondmonolength=\Y%
\advance\secondmonolength by -\value{y}%
\secondepilength=\Y%
\advance\secondepilength by \value{y}%
\put(\value{x},-\value{y}){\line(2,1){#3}}%
\put(\epilength,\secondmonolength){\enehead}%
\put(\monolength,\secondepilength){\line(-2,-1){#3}}%
\put(-\value{x},\value{y}){\wswhead}%
\end{picture}}
\put(-\value{x},\value{y}){\distsign}%
\put(\value{x},-\value{y}){\distsign}%
\truex{300}\truey{1000}\truez{600}%
\put(-\value{x},\value{x}){\makebox(0,\value{z})[r]{${#1}$}}%
\put(\value{x},-\value{y}){\makebox(0,\value{z})[l]{${#2}$}}%
\end{picture}}%

\def\basicenear[#1]{\ENEAR{}{}{#100}}%

\newcommand{\enear}{\@ifnextchar[{\basicenear}{\basicenear[133]}}%

\def\basicEnear[#1]#2{\ENEAR{#2}{}{#100}}%

\newcommand{\Enear}{\@ifnextchar[{\basicEnear}{\basicEnear[133]}}%

\def\basiceneaR[#1]#2{\ENEAR{}{#2}{#100}}%

\newcommand{\eneaR}{\@ifnextchar[{\basiceneaR}{\basiceneaR[133]}}%

\def\basicenedist[#1]{\ENEDIST{}{}{#100}}%

\newcommand{\enedist}{\@ifnextchar[{\basicenedist}{\basicenedist[133]}}%

\def\basicEnedist[#1]#2{\ENEDIST{#2}{}{#100}}%

\newcommand{\Enedist}{\@ifnextchar[{\basicEnedist}{\basicEnedist[133]}}%

\def\basicenedisT[#1]#2{\ENEDIST{}{#2}{#100}}%

\newcommand{\enedisT}{\@ifnextchar[{\basicenedisT}{\basicenedisT[133]}}%

\def\basicenedotar[#1]{\ENEDOTAR{}{}{#100}}%

\newcommand{\enedotar}{\@ifnextchar[{\basicenedotar}{\basicenedotar[133]}}%

\def\basicEnedotar[#1]#2{\ENEDOTAR{#2}{}{#100}}%

\newcommand{\Enedotar}{\@ifnextchar[{\basicEnedotar}{\basicEnedotar[133]}}%

\def\basicenedotaR[#1]#2{\ENEDOTAR{}{#2}{#100}}%

\newcommand{\enedotaR}{\@ifnextchar[{\basicenedotaR}{\basicenedotaR[133]}}%

\def\basicenemono[#1]{\ENEMONO{}{}{#100}}%

\newcommand{\enemono}{\@ifnextchar[{\basicenemono}{\basicenemono[133]}}%

\def\basicEnemono[#1]#2{\ENEMONO{#2}{}{#100}}%

\newcommand{\Enemono}{\@ifnextchar[{\basicEnemono}{\basicEnemono[133]}}%

\def\basicenemonO[#1]#2{\ENEMONO{}{#2}{#100}}%

\newcommand{\enemonO}{\@ifnextchar[{\basicenemonO}{\basicenemonO[133]}}%

\def\basiceneepi[#1]{\ENEEPI{}{}{#100}}%

\newcommand{\eneepi}{\@ifnextchar[{\basiceneepi}{\basiceneepi[133]}}%

\def\basicEneepi[#1]#2{\ENEEPI{#2}{}{#100}}%

\newcommand{\Eneepi}{\@ifnextchar[{\basicEneepi}{\basicEneepi[133]}}%

\def\basiceneepI[#1]#2{\ENEEPI{}{#2}{#100}}%

\newcommand{\eneepI}{\@ifnextchar[{\basiceneepI}{\basiceneepI[133]}}%

\def\basicenebimo[#1]{\ENEBIMO{}{}{#100}}%

\newcommand{\enebimo}{\@ifnextchar[{\basicenebimo}{\basicenebimo[133]}}%

\def\basicEnebimo[#1]#2{\ENEBIMO{#2}{}{#100}}%

\newcommand{\Enebimo}{\@ifnextchar[{\basicEnebimo}{\basicEnebimo[133]}}%

\def\basicenebimO[#1]#2{\ENEBIMO{}{#2}{#100}}%

\newcommand{\enebimO}{\@ifnextchar[{\basicenebimO}{\basicenebimO[133]}}%

\def\basiceneiso[#1]{\ENEAR{\isosign}{}{#100}}%

\newcommand{\eneiso}{\@ifnextchar[{\basiceneiso}{\basiceneiso[133]}}%

\def\basicEneiso[#1]#2{\ENEAR{#2}{\isosign}{#100}}%

\newcommand{\Eneiso}{\@ifnextchar[{\basicEneiso}{\basicEneiso[133]}}%

\def\basiceneisO[#1]#2{\ENEAR{\isosign}{#2}{#100}}%

\newcommand{\eneisO}{\@ifnextchar[{\basiceneisO}{\basiceneisO[133]}}%

\def\basiceneeql[#1]{\ENEEQL{}{}{#100}}%

\newcommand{\eneeql}{\@ifnextchar[{\basiceneeql}{\basiceneeql[133]}}%

\def\basicEneeql[#1]#2{\ENEEQL{#2}{}{#100}}%

\newcommand{\Eneeql}{\@ifnextchar[{\basicEneeql}{\basicEneeql[133]}}%

\def\basiceneeqL[#1]#2{\ENEEQL{}{#2}{#100}}%

\newcommand{\eneeqL}{\@ifnextchar[{\basiceneeqL}{\basiceneeqL[133]}}%

\def\basicenebiar[#1]{\ENEBIAR{}{}{#100}}%

\newcommand{\enebiar}{\@ifnextchar[{\basicenebiar}{\basicenebiar[133]}}%

\def\basicEnebiar[#1]#2#3{\ENEBIAR{#2}{#3}{#100}}%

\newcommand{\Enebiar}{\@ifnextchar[{\basicEnebiar}{\basicEnebiar[133]}}%

\def\basicenebidist[#1]{\ENEBIDIST{}{}{#100}}%

\newcommand{\enebidist}{\@ifnextchar[{\basicenebidist}{\basicenebidist[133]}}%

\def\basicEnebidist[#1]#2#3{\ENEBIDIST{#2}{#3}{#100}}%

\newcommand{\Enebidist}{\@ifnextchar[{\basicEnebidist}{\basicEnebidist[133]}}%

\def\basiceneadjar[#1]{\ENEADJAR{}{}{#100}}%

\newcommand{\eneadjar}{\@ifnextchar[{\basiceneadjar}{\basiceneadjar[133]}}%

\def\basicEneadjar[#1]#2#3{\ENEADJAR{#2}{#3}{#100}}%

\newcommand{\Eneadjar}{\@ifnextchar[{\basicEneadjar}{\basicEneadjar[133]}}%

\def\basiceneadjdist[#1]{\ENEADJDIST{}{}{#100}}%

\newcommand{\eneadjdist}{\@ifnextchar[{\basiceneadjdist}{\basiceneadjdist[13
3]}}%

\def\basicEneadjdist[#1]#2#3{\ENEADJDIST{#2}{#3}{#100}}%

\newcommand{\Eneadjdist}{\@ifnextchar[{\basicEneadjdist}{\basicEneadjdist[13
3]}}%

\newcommand{\ESEAR}[3]{%
\Y=#3%
\divide\Y by 2%
\Z=\Y%
\divide\Z by 2%
\begin{picture}(0,0)%
\put(-\Y,\Z){\line(2,-1){#3}}%
\put(\Y,-\Z){\esehead}%
\truex{200}\truey{800}\truez{600}%
\put(\value{x},\value{x}){\makebox(0,\value{z})[l]{${#1}$}}%
\put(-\value{x},-\value{y}){\makebox(0,\value{z})[r]{${#2}$}}%
\end{picture}}%

\newcommand{\ESEDIST}[3]{%
\Y=#3%
\divide\Y by 2%
\Z=\Y%
\divide\Z by 2%
\begin{picture}(0,0)%
\put(-\Y,\Z){\line(2,-1){#3}}%
\put(\Y,-\Z){\esehead}%
\put(0,0){\distsign}%
\truex{200}\truey{800}\truez{600}%
\put(\value{x},\value{x}){\makebox(0,\value{z})[l]{${#1}$}}%
\put(-\value{x},-\value{y}){\makebox(0,\value{z})[r]{${#2}$}}%
\end{picture}}%

\newcommand{\ESEDOTAR}[3]{%
\truex{100}\truey{268}\truez{134}%
\Y=#3%
\divide\Y by 2%
\Z=\Y%
\divide\Z by 2%
\NUMBEROFDOTS=#3%
\divide\NUMBEROFDOTS by \value{y}%
\advance\NUMBEROFDOTS by 1%
\begin{picture}(0,0)%
\multiput(-\Y,\Z)(\value{y},-\value{z}){\NUMBEROFDOTS}%
{\circle*{\value{x}}}%
\put(\Y,-\Z){\esehead}%
\truex{200}\truey{800}\truez{600}%
\put(\value{x},\value{x}){\makebox(0,\value{z})[l]{${#1}$}}%
\put(-\value{x},-\value{y}){\makebox(0,\value{z})[r]{${#2}$}}%
\end{picture}}%

\newcommand{\ESEMONO}[3]{%
\Y=#3%
\divide\Y by 2%
\Z=\Y%
\divide\Z by 2%
\TrueTail%
\bimolength=#3%
\advance\bimolength by -\TrueMonoTail%
\monolength=\bimolength%
\advance\monolength by -\Y%
\secondmonolength=\monolength%
\divide\secondmonolength by 2%
\begin{picture}(0,0)%
\put(-\monolength,\secondmonolength){\line(2,-1){\bimolength}}%
\put(-\monolength,\secondmonolength){\esehead}%
\put(\Y,-\Z){\esehead}%
\truex{200}\truey{800}\truez{600}%
\put(\value{x},\value{x}){\makebox(0,\value{z})[l]{${#1}$}}%
\put(-\value{x},-\value{y}){\makebox(0,\value{z})[r]{${#2}$}}%
\end{picture}}%

\newcommand{\ESEEPI}[3]{%
\Y=#3%
\divide\Y by 2%
\Z=\Y%
\divide\Z by 2%
\TrueHead%
\bimolength=#3%
\advance\bimolength by -\TrueEpiHead%
\epilength=\bimolength%
\advance\epilength by -\Y%
\secondepilength=\epilength%
\divide\secondepilength by 2%
\begin{picture}(0,0)%
\put(-\Y,\Z){\line(2,-1){\bimolength}}%
\put(\epilength,-\secondepilength){\esehead}%
\put(\Y,-\Z){\esehead}%
\truex{200}\truey{800}\truez{600}%
\put(\value{x},\value{x}){\makebox(0,\value{z})[l]{${#1}$}}%
\put(-\value{x},-\value{y}){\makebox(0,\value{z})[r]{${#2}$}}%
\end{picture}}%

\newcommand{\ESEBIMO}[3]{%
\Y=#3%
\divide\Y by 2%
\Z=\Y%
\divide\Z by 2%
\TrueTail\TrueHead%
\bimolength=#3%
\advance\bimolength by -\TrueMonoTail%
\monolength=\bimolength%
\advance\monolength by -\Y%
\advance\bimolength by -\TrueEpiHead%
\epilength=\bimolength%
\advance\epilength by -\monolength%
\secondmonolength=\monolength%
\divide\secondmonolength by 2%
\secondepilength=\epilength%
\divide\secondepilength by 2%
\begin{picture}(0,0)%
\put(-\monolength,\secondmonolength){\line(2,-1){\bimolength}}%
\put(-\monolength,\secondmonolength){\esehead}%
\put(\epilength,-\secondepilength){\esehead}%
\put(\Y,-\Z){\esehead}%
\truex{200}\truey{800}\truez{600}%
\put(\value{x},\value{x}){\makebox(0,\value{z})[l]{${#1}$}}%
\put(-\value{x},-\value{y}){\makebox(0,\value{z})[r]{${#2}$}}%
\end{picture}}%

\newcommand{\ESEEQL}[3]{%
\Y=#3%
\divide\Y by 2%
\Z=\Y%
\divide\Z by 2%
\begin{picture}(0,0)%
\put(-\Y,\Z){\begin{picture}(0,0)%
\truex{44}\truey{89}%
\put(-\value{x},-\value{y}){\line(2,-1){#3}}%
\put(\value{x},\value{y}){\line(2,-1){#3}}%
\end{picture}}%
\truex{200}\truey{800}\truez{600}%
\put(\value{x},\value{x}){\makebox(0,\value{z})[l]{${#1}$}}%
\put(-\value{x},-\value{y}){\makebox(0,\value{z})[r]{${#2}$}}%
\end{picture}}%

\newcommand{\ESEBIAR}[3]{%
\Y=#3%
\divide\Y by 2%
\Z=\Y%
\divide\Z by 2%
\begin{picture}(0,0)%
\put(-\Y,\Z){\begin{picture}(0,0)%
\truex{156}\truey{313}%
\put(-\value{x},-\value{y}){\line(2,-1){#3}}%
\put(\value{x},\value{y}){\line(2,-1){#3}}%
\monolength=#3%
\advance\monolength by -\value{x}%
\epilength=#3%
\advance\epilength by \value{x}%
\secondmonolength=\Y%
\advance\secondmonolength by -\value{y}%
\secondepilength=\Y%
\advance\secondepilength by \value{y}%
\put(\monolength,-\secondepilength){\esehead}%
\put(\epilength,-\secondmonolength){\esehead}%
\end{picture}}
\truex{400}\truey{1000}\truez{600}%
\put(\value{x},\value{x}){\makebox(0,\value{z})[l]{${#1}$}}%
\put(-\value{x},-\value{y}){\makebox(0,\value{z})[r]{${#2}$}}%
\end{picture}}%

\newcommand{\ESEBIDIST}[3]{%
\Y=#3%
\divide\Y by 2%
\Z=\Y%
\divide\Z by 2%
\begin{picture}(0,0)%
\truex{156}\truey{313}\truez{400}%
\put(-\Y,\Z){\begin{picture}(0,0)%
\put(-\value{x},-\value{y}){\line(2,-1){#3}}%
\put(\value{x},\value{y}){\line(2,-1){#3}}%
\monolength=#3%
\advance\monolength by -\value{x}%
\epilength=#3%
\advance\epilength by \value{x}%
\secondmonolength=\Y%
\advance\secondmonolength by -\value{y}%
\secondepilength=\Y%
\advance\secondepilength by \value{y}%
\put(\monolength,-\secondepilength){\esehead}%
\put(\epilength,-\secondmonolength){\esehead}%
\end{picture}}
\put(\value{x},\value{y}){\distsign}%
\put(-\value{x},-\value{y}){\distsign}%
\truex{400}\truey{1000}\truez{600}%
\put(\value{x},\value{x}){\makebox(0,\value{z})[l]{${#1}$}}%
\put(-\value{x},-\value{y}){\makebox(0,\value{z})[r]{${#2}$}}%
\end{picture}}%

\newcommand{\ESEADJAR}[3]{%
\Y=#3%
\divide\Y by 2%
\Z=\Y%
\divide\Z by 2%
\begin{picture}(0,0)%
\put(-\Y,\Z){\begin{picture}(0,0)%
\truex{156}\truey{313}%
\monolength=#3%
\advance\monolength by -\value{x}%
\epilength=#3%
\advance\epilength by \value{x}%
\secondmonolength=\Y%
\advance\secondmonolength by -\value{y}%
\secondepilength=\Y%
\advance\secondepilength by \value{y}%
\put(-\value{x},-\value{y}){\line(2,-1){#3}}%
\put(\monolength,-\secondepilength){\esehead}%
\put(\epilength,-\secondmonolength){\line(-2,1){#3}}%
\put(\value{x},\value{y}){\wnwhead}%
\end{picture}}
\truex{400}\truey{1000}\truez{600}%
\put(\value{x},\value{x}){\makebox(0,\value{z})[l]{${#1}$}}%
\put(-\value{x},-\value{y}){\makebox(0,\value{z})[r]{${#2}$}}%
\end{picture}}%

\newcommand{\ESEADJDIST}[3]{%
\Y=#3%
\divide\Y by 2%
\Z=\Y%
\divide\Z by 2%
\begin{picture}(0,0)%
\truex{156}\truey{313}\truez{400}%
\put(-\Y,\Z){\begin{picture}(0,0)%
\monolength=#3%
\advance\monolength by -\value{x}%
\epilength=#3%
\advance\epilength by \value{x}%
\secondmonolength=\Y%
\advance\secondmonolength by -\value{y}%
\secondepilength=\Y%
\advance\secondepilength by \value{y}%
\put(-\value{x},-\value{y}){\line(2,-1){#3}}%
\put(\monolength,-\secondepilength){\esehead}%
\put(\epilength,-\secondmonolength){\line(-2,1){#3}}%
\put(\value{x},\value{y}){\wnwhead}%
\end{picture}}
\put(\value{x},\value{y}){\distsign}%
\put(-\value{x},-\value{y}){\distsign}%
\truex{400}\truey{1000}\truez{600}%
\put(\value{x},\value{x}){\makebox(0,\value{z})[l]{${#1}$}}%
\put(-\value{x},-\value{y}){\makebox(0,\value{z})[r]{${#2}$}}%
\end{picture}}%

\def\basicesear[#1]{\ESEAR{}{}{#100}}%

\newcommand{\esear}{\@ifnextchar[{\basicesear}{\basicesear[133]}}%

\def\basicEsear[#1]#2{\ESEAR{#2}{}{#100}}%

\newcommand{\Esear}{\@ifnextchar[{\basicEsear}{\basicEsear[133]}}%

\def\basiceseaR[#1]#2{\ESEAR{}{#2}{#100}}%

\newcommand{\eseaR}{\@ifnextchar[{\basiceseaR}{\basiceseaR[133]}}%

\def\basicesedist[#1]{\ESEDIST{}{}{#100}}%

\newcommand{\esedist}{\@ifnextchar[{\basicesedist}{\basicesedist[133]}}%

\def\basicEsedist[#1]#2{\ESEDIST{#2}{}{#100}}%

\newcommand{\Esedist}{\@ifnextchar[{\basicEsedist}{\basicEsedist[133]}}%

\def\basicesedisT[#1]#2{\ESEDIST{}{#2}{#100}}%

\newcommand{\esedisT}{\@ifnextchar[{\basicesedisT}{\basicesedisT[133]}}%

\def\basicesedotar[#1]{\ESEDOTAR{}{}{#100}}%

\newcommand{\esedotar}{\@ifnextchar[{\basicesedotar}{\basicesedotar[133]}}%

\def\basicEsedotar[#1]#2{\ESEDOTAR{#2}{}{#100}}%

\newcommand{\Esedotar}{\@ifnextchar[{\basicEsedotar}{\basicEsedotar[133]}}%

\def\basicesedotaR[#1]#2{\ESEDOTAR{}{#2}{#100}}%

\newcommand{\esedotaR}{\@ifnextchar[{\basicesedotaR}{\basicesedotaR[133]}}%

\def\basicesemono[#1]{\ESEMONO{}{}{#100}}%

\newcommand{\esemono}{\@ifnextchar[{\basicesemono}{\basicesemono[133]}}%

\def\basicEsemono[#1]#2{\ESEMONO{#2}{}{#100}}%

\newcommand{\Esemono}{\@ifnextchar[{\basicEsemono}{\basicEsemono[133]}}%

\def\basicesemonO[#1]#2{\ESEMONO{}{#2}{#100}}%

\newcommand{\esemonO}{\@ifnextchar[{\basicesemonO}{\basicesemonO[133]}}%

\def\basiceseepi[#1]{\ESEEPI{}{}{#100}}%

\newcommand{\eseepi}{\@ifnextchar[{\basiceseepi}{\basiceseepi[133]}}%

\def\basicEseepi[#1]#2{\ESEEPI{#2}{}{#100}}%

\newcommand{\Eseepi}{\@ifnextchar[{\basicEseepi}{\basicEseepi[133]}}%

\def\basiceseepI[#1]#2{\ESEEPI{}{#2}{#100}}%

\newcommand{\eseepI}{\@ifnextchar[{\basiceseepI}{\basiceseepI[133]}}%

\def\basicesebimo[#1]{\ESEBIMO{}{}{#100}}%

\newcommand{\esebimo}{\@ifnextchar[{\basicesebimo}{\basicesebimo[133]}}%

\def\basicEsebimo[#1]#2{\ESEBIMO{#2}{}{#100}}%

\newcommand{\Esebimo}{\@ifnextchar[{\basicEsebimo}{\basicEsebimo[133]}}%

\def\basicesebimO[#1]#2{\ESEBIMO{}{#2}{#100}}%

\newcommand{\esebimO}{\@ifnextchar[{\basicesebimO}{\basicesebimO[133]}}%

\def\basiceseiso[#1]{\ESEAR{\isosign}{}{#100}}%

\newcommand{\eseiso}{\@ifnextchar[{\basiceseiso}{\basiceseiso[133]}}%

\def\basicEseiso[#1]#2{\ESEAR{#2}{\isosign}{#100}}%

\newcommand{\Eseiso}{\@ifnextchar[{\basicEseiso}{\basicEseiso[133]}}%

\def\basiceseisO[#1]#2{\ESEAR{\isosign}{#2}{#100}}%

\newcommand{\eseisO}{\@ifnextchar[{\basiceseisO}{\basiceseisO[133]}}%

\def\basiceseeql[#1]{\ESEEQL{}{}{#100}}%

\newcommand{\eseeql}{\@ifnextchar[{\basiceseeql}{\basiceseeql[133]}}%

\def\basicEseeql[#1]#2{\ESEEQL{#2}{}{#100}}%

\newcommand{\Eseeql}{\@ifnextchar[{\basicEseeql}{\basicEseeql[133]}}%

\def\basiceseeqL[#1]#2{\ESEEQL{}{#2}{#100}}%

\newcommand{\eseeqL}{\@ifnextchar[{\basiceseeqL}{\basiceseeqL[133]}}%

\def\basicesebiar[#1]{\ESEBIAR{}{}{#100}}%

\newcommand{\esebiar}{\@ifnextchar[{\basicesebiar}{\basicesebiar[133]}}%

\def\basicEsebiar[#1]#2#3{\ESEBIAR{#2}{#3}{#100}}%

\newcommand{\Esebiar}{\@ifnextchar[{\basicEsebiar}{\basicEsebiar[133]}}%

\def\basicesebidist[#1]{\ESEBIDIST{}{}{#100}}%

\newcommand{\esebidist}{\@ifnextchar[{\basicesebidist}{\basicesebidist[133]}}%

\def\basicEsebidist[#1]#2#3{\ESEBIDIST{#2}{#3}{#100}}%

\newcommand{\Esebidist}{\@ifnextchar[{\basicEsebidist}{\basicEsebidist[133]}}%

\def\basiceseadjar[#1]{\ESEADJAR{}{}{#100}}%

\newcommand{\eseadjar}{\@ifnextchar[{\basiceseadjar}{\basiceseadjar[133]}}%

\def\basicEseadjar[#1]#2#3{\ESEADJAR{#2}{#3}{#100}}%

\newcommand{\Eseadjar}{\@ifnextchar[{\basicEseadjar}{\basicEseadjar[133]}}%

\def\basiceseadjdist[#1]{\ESEADJDIST{}{}{#100}}%

\newcommand{\eseadjdist}{\@ifnextchar[{\basiceseadjdist}{\basiceseadjdist[13
3]}}%

\def\basicEseadjdist[#1]#2#3{\ESEADJDIST{#2}{#3}{#100}}%

\newcommand{\Eseadjdist}{\@ifnextchar[{\basicEseadjdist}{\basicEseadjdist[13
3]}}%

\newcommand{\WSWAR}[3]{%
\Y=#3%
\divide\Y by 2%
\Z=\Y%
\divide\Z by 2%
\begin{picture}(0,0)%
\put(\Y,\Z){\line(-2,-1){#3}}%
\put(-\Y,-\Z){\wswhead}%
\truex{200}\truey{800}\truez{600}%
\put(-\value{x},\value{x}){\makebox(0,\value{z})[r]{${#1}$}}%
\put(\value{x},-\value{y}){\makebox(0,\value{z})[l]{${#2}$}}%
\end{picture}}%

\newcommand{\WSWDIST}[3]{%
\Y=#3%
\divide\Y by 2%
\Z=\Y%
\divide\Z by 2%
\begin{picture}(0,0)%
\put(\Y,\Z){\line(-2,-1){#3}}%
\put(-\Y,-\Z){\wswhead}%
\truex{400}%
\put(0,0){\distsign}%
\truex{200}\truey{800}\truez{600}%
\put(-\value{x},\value{x}){\makebox(0,\value{z})[r]{${#1}$}}%
\put(\value{x},-\value{y}){\makebox(0,\value{z})[l]{${#2}$}}%
\end{picture}}%

\newcommand{\WSWDOTAR}[3]{%
\truex{100}\truey{268}\truez{134}%
\Y=#3%
\divide\Y by 2%
\Z=\Y%
\divide\Z by 2%
\NUMBEROFDOTS=#3%
\divide\NUMBEROFDOTS by \value{y}%
\advance\NUMBEROFDOTS by 1%
\begin{picture}(0,0)%
\multiput(\Y,\Z)(-\value{y},-\value{z}){\NUMBEROFDOTS}%
{\circle*{\value{x}}}%
\put(-\Y,-\Z){\wswhead}%
\truex{200}\truey{800}\truez{600}%
\put(-\value{x},\value{x}){\makebox(0,\value{z})[r]{${#1}$}}%
\put(\value{x},-\value{y}){\makebox(0,\value{z})[l]{${#2}$}}%
\end{picture}}%

\newcommand{\WSWMONO}[3]{%
\Y=#3%
\divide\Y by 2%
\Z=\Y%
\divide\Z by 2%
\TrueTail%
\bimolength=#3%
\advance\bimolength by -\TrueMonoTail%
\monolength=\bimolength%
\advance\monolength by -\Y%
\secondmonolength=\monolength%
\divide\secondmonolength by 2%
\begin{picture}(0,0)%
\put(\monolength,\secondmonolength){\line(-2,-1){\bimolength}}%
\put(\monolength,\secondmonolength){\wswhead}%
\put(-\Y,-\Z){\wswhead}%
\truex{200}\truey{800}\truez{600}%
\put(-\value{x},\value{x}){\makebox(0,\value{z})[r]{${#1}$}}%
\put(\value{x},-\value{y}){\makebox(0,\value{z})[l]{${#2}$}}%
\end{picture}}%

\newcommand{\WSWEPI}[3]{%
\Y=#3%
\divide\Y by 2%
\Z=\Y%
\divide\Z by 2%
\TrueHead%
\bimolength=#3%
\advance\bimolength by -\TrueEpiHead%
\epilength=\bimolength%
\advance\epilength by -\Y%
\secondepilength=\epilength%
\divide\secondepilength by 2%
\begin{picture}(0,0)%
\put(\Y,\Z){\line(-2,-1){\bimolength}}%
\put(-\epilength,-\secondepilength){\wswhead}%
\put(-\Y,-\Z){\wswhead}%
\truex{200}\truey{800}\truez{600}%
\put(-\value{x},\value{x}){\makebox(0,\value{z})[r]{${#1}$}}%
\put(\value{x},-\value{y}){\makebox(0,\value{z})[l]{${#2}$}}%
\end{picture}}%

\newcommand{\WSWBIMO}[3]{%
\Y=#3%
\divide\Y by 2%
\Z=\Y%
\divide\Z by 2%
\TrueTail\TrueHead%
\bimolength=#3%
\advance\bimolength by -\TrueMonoTail%
\monolength=\bimolength%
\advance\monolength by -\Y%
\advance\bimolength by -\TrueEpiHead%
\epilength=\bimolength%
\advance\epilength by -\monolength%
\secondmonolength=\monolength%
\divide\secondmonolength by 2%
\secondepilength=\epilength%
\divide\secondepilength by 2%
\begin{picture}(0,0)%
\put(\monolength,\secondmonolength){\line(-2,-1){\bimolength}}%
\put(\monolength,\secondmonolength){\wswhead}%
\put(-\epilength,-\secondepilength){\wswhead}%
\put(-\Y,-\Z){\wswhead}%
\truex{200}\truey{800}\truez{600}%
\put(-\value{x},\value{x}){\makebox(0,\value{z})[r]{${#1}$}}%
\put(\value{x},-\value{y}){\makebox(0,\value{z})[l]{${#2}$}}%
\end{picture}}%

\newcommand{\WSWBIAR}[3]{%
\Y=#3%
\divide\Y by 2%
\Z=\Y%
\divide\Z by 2%
\begin{picture}(0,0)%
\put(\Y,\Z){\begin{picture}(0,0)%
\truex{156}\truey{313}%
\put(-\value{x},\value{y}){\line(-2,-1){#3}}%
\put(\value{x},-\value{y}){\line(-2,-1){#3}}%
\monolength=#3%
\advance\monolength by -\value{x}%
\epilength=#3%
\advance\epilength by \value{x}%
\secondmonolength=\Y%
\advance\secondmonolength by -\value{y}%
\secondepilength=\Y%
\advance\secondepilength by \value{y}%
\put(-\monolength,-\secondepilength){\wswhead}%
\put(-\epilength,-\secondmonolength){\wswhead}%
\end{picture}}
\truex{300}\truey{1000}\truez{600}%
\put(-\value{x},\value{x}){\makebox(0,\value{z})[r]{${#1}$}}%
\put(\value{x},-\value{y}){\makebox(0,\value{z})[l]{${#2}$}}%
\end{picture}}%

\newcommand{\WSWBIDIST}[3]{%
\Y=#3%
\divide\Y by 2%
\Z=\Y%
\divide\Z by 2%
\begin{picture}(0,0)%
\truex{156}\truey{313}\truez{400}%
\put(\Y,\Z){\begin{picture}(0,0)%
\put(-\value{x},\value{y}){\line(-2,-1){#3}}%
\put(\value{x},-\value{y}){\line(-2,-1){#3}}%
\monolength=#3%
\advance\monolength by -\value{x}%
\epilength=#3%
\advance\epilength by \value{x}%
\secondmonolength=\Y%
\advance\secondmonolength by -\value{y}%
\secondepilength=\Y%
\advance\secondepilength by \value{y}%
\put(-\monolength,-\secondepilength){\wswhead}%
\put(-\epilength,-\secondmonolength){\wswhead}%
\end{picture}}
\put(-\value{x},\value{y}){\distsign}%
\put(\value{x},-\value{y}){\distsign}%
\truex{300}\truey{1000}\truez{600}%
\put(-\value{x},\value{x}){\makebox(0,\value{z})[r]{${#1}$}}%
\put(\value{x},-\value{y}){\makebox(0,\value{z})[l]{${#2}$}}%
\end{picture}}%

\newcommand{\WSWADJAR}[3]{%
\Y=#3%
\divide\Y by 2%
\Z=\Y%
\divide\Z by 2%
\begin{picture}(0,0)%
\put(\Y,\Z){\begin{picture}(0,0)%
\truex{156}\truey{313}%
\monolength=#3%
\advance\monolength by -\value{x}%
\epilength=#3%
\advance\epilength by \value{x}%
\secondmonolength=\Y%
\advance\secondmonolength by -\value{y}%
\secondepilength=\Y%
\advance\secondepilength by \value{y}%
\put(\value{x},-\value{y}){\line(-2,-1){#3}}%
\put(-\monolength,-\secondepilength){\wswhead}%
\put(-\epilength,-\secondmonolength){\line(2,1){#3}}%
\put(-\value{x},\value{y}){\enehead}%
\end{picture}}
\truex{300}\truey{1000}\truez{600}%
\put(-\value{x},\value{x}){\makebox(0,\value{z})[r]{${#1}$}}%
\put(\value{x},-\value{y}){\makebox(0,\value{z})[l]{${#2}$}}%
\end{picture}}%

\newcommand{\WSWADJDIST}[3]{%
\Y=#3%
\divide\Y by 2%
\Z=\Y%
\divide\Z by 2%
\begin{picture}(0,0)%
\truex{156}\truey{313}\truez{400}%
\put(\Y,\Z){\begin{picture}(0,0)%
\monolength=#3%
\advance\monolength by -\value{x}%
\epilength=#3%
\advance\epilength by \value{x}%
\secondmonolength=\Y%
\advance\secondmonolength by -\value{y}%
\secondepilength=\Y%
\advance\secondepilength by \value{y}%
\put(\value{x},-\value{y}){\line(-2,-1){#3}}%
\put(-\monolength,-\secondepilength){\wswhead}%
\put(-\epilength,-\secondmonolength){\line(2,1){#3}}%
\put(-\value{x},\value{y}){\enehead}%
\end{picture}}
\put(-\value{x},\value{y}){\distsign}%
\put(\value{x},-\value{y}){\distsign}%
\truex{300}\truey{1000}\truez{600}%
\put(-\value{x},\value{x}){\makebox(0,\value{z})[r]{${#1}$}}%
\put(\value{x},-\value{y}){\makebox(0,\value{z})[l]{${#2}$}}%
\end{picture}}%

\def\basicwswar[#1]{\WSWAR{}{}{#100}}%

\newcommand{\wswar}{\@ifnextchar[{\basicwswar}{\basicwswar[133]}}%

\def\basicWswar[#1]#2{\WSWAR{#2}{}{#100}}%

\newcommand{\Wswar}{\@ifnextchar[{\basicWswar}{\basicWswar[133]}}%

\def\basicwswaR[#1]#2{\WSWAR{}{#2}{#100}}%

\newcommand{\wswaR}{\@ifnextchar[{\basicwswaR}{\basicwswaR[133]}}%

\def\basicwswdist[#1]{\WSWDIST{}{}{#100}}%

\newcommand{\wswdist}{\@ifnextchar[{\basicwswdist}{\basicwswdist[133]}}%

\def\basicWswdist[#1]#2{\WSWDIST{#2}{}{#100}}%

\newcommand{\Wswdist}{\@ifnextchar[{\basicWswdist}{\basicWswdist[133]}}%

\def\basicwswdisT[#1]#2{\WSWDIST{}{#2}{#100}}%

\newcommand{\wswdisT}{\@ifnextchar[{\basicwswdisT}{\basicwswdisT[133]}}%

\def\basicwswdotar[#1]{\WSWDOTAR{}{}{#100}}%

\newcommand{\wswdotar}{\@ifnextchar[{\basicwswdotar}{\basicwswdotar[133]}}%

\def\basicWswdotar[#1]#2{\WSWDOTAR{#2}{}{#100}}%

\newcommand{\Wswdotar}{\@ifnextchar[{\basicWswdotar}{\basicWswdotar[133]}}%

\def\basicwswdotaR[#1]#2{\WSWDOTAR{}{#2}{#100}}%

\newcommand{\wswdotaR}{\@ifnextchar[{\basicwswdotaR}{\basicwswdotaR[133]}}%

\def\basicwswmono[#1]{\WSWMONO{}{}{#100}}%

\newcommand{\wswmono}{\@ifnextchar[{\basicwswmono}{\basicwswmono[133]}}%

\def\basicWswmono[#1]#2{\WSWMONO{#2}{}{#100}}%

\newcommand{\Wswmono}{\@ifnextchar[{\basicWswmono}{\basicWswmono[133]}}%

\def\basicwswmonO[#1]#2{\WSWMONO{}{#2}{#100}}%

\newcommand{\wswmonO}{\@ifnextchar[{\basicwswmonO}{\basicwswmonO[133]}}%

\def\basicwswepi[#1]{\WSWEPI{}{}{#100}}%

\newcommand{\wswepi}{\@ifnextchar[{\basicwswepi}{\basicwswepi[133]}}%

\def\basicWswepi[#1]#2{\WSWEPI{#2}{}{#100}}%

\newcommand{\Wswepi}{\@ifnextchar[{\basicWswepi}{\basicWswepi[133]}}%

\def\basicwswepI[#1]#2{\WSWEPI{}{#2}{#100}}%

\newcommand{\wswepI}{\@ifnextchar[{\basicwswepI}{\basicwswepI[133]}}%

\def\basicwswbimo[#1]{\WSWBIMO{}{}{#100}}%

\newcommand{\wswbimo}{\@ifnextchar[{\basicwswbimo}{\basicwswbimo[133]}}%

\def\basicWswbimo[#1]#2{\WSWBIMO{#2}{}{#100}}%

\newcommand{\Wswbimo}{\@ifnextchar[{\basicWswbimo}{\basicWswbimo[133]}}%

\def\basicwswbimO[#1]#2{\WSWBIMO{}{#2}{#100}}%

\newcommand{\wswbimO}{\@ifnextchar[{\basicwswbimO}{\basicwswbimO[133]}}%

\def\basicwswiso[#1]{\WSWAR{\isosign}{}{#100}}%

\newcommand{\wswiso}{\@ifnextchar[{\basicwswiso}{\basicwswiso[133]}}%

\def\basicWswiso[#1]#2{\WSWAR{#2}{\isosign}{#100}}%

\newcommand{\Wswiso}{\@ifnextchar[{\basicWswiso}{\basicWswiso[133]}}%

\def\basicwswisO[#1]#2{\WSWAR{\isosign}{#2}{#100}}%

\newcommand{\wswisO}{\@ifnextchar[{\basicwswisO}{\basicwswisO[133]}}%

\def\basicwswbiar[#1]{\WSWBIAR{}{}{#100}}%

\newcommand{\wswbiar}{\@ifnextchar[{\basicwswbiar}{\basicwswbiar[133]}}%

\def\basicWswbiar[#1]#2#3{\WSWBIAR{#2}{#3}{#100}}%

\newcommand{\Wswbiar}{\@ifnextchar[{\basicWswbiar}{\basicWswbiar[133]}}%

\def\basicwswbidist[#1]{\WSWBIDIST{}{}{#100}}%

\newcommand{\wswbidist}{\@ifnextchar[{\basicwswbidist}{\basicwswbidist[133]}}%

\def\basicWswbidist[#1]#2#3{\WSWBIDIST{#2}{#3}{#100}}%

\newcommand{\Wswbidist}{\@ifnextchar[{\basicWswbidist}{\basicWswbidist[133]}}%

\def\basicwswadjar[#1]{\WSWADJAR{}{}{#100}}%

\newcommand{\wswadjar}{\@ifnextchar[{\basicwswadjar}{\basicwswadjar[133]}}%

\def\basicWswadjar[#1]#2#3{\WSWADJAR{#2}{#3}{#100}}%

\newcommand{\Wswadjar}{\@ifnextchar[{\basicWswadjar}{\basicWswadjar[133]}}%

\def\basicwswadjdist[#1]{\WSWADJDIST{}{}{#100}}%

\newcommand{\wswadjdist}{\@ifnextchar[{\basicwswadjdist}{\basicwswadjdist[13
3]}}%

\def\basicWswadjdist[#1]#2#3{\WSWADJDIST{#2}{#3}{#100}}%

\newcommand{\Wswadjdist}{\@ifnextchar[{\basicWswadjdist}{\basicWswadjdist[13
3]}}%

\newcommand{\WNWAR}[3]{%
\Y=#3%
\divide\Y by 2%
\Z=\Y%
\divide\Z by 2%
\begin{picture}(0,0)%
\put(\Y,-\Z){\line(-2,1){#3}}%
\put(-\Y,\Z){\wnwhead}%
\truex{200}\truey{800}\truez{600}%
\put(\value{x},\value{x}){\makebox(0,\value{z})[l]{${#1}$}}%
\put(-\value{x},-\value{y}){\makebox(0,\value{z})[r]{${#2}$}}%
\end{picture}}%

\newcommand{\WNWDIST}[3]{%
\Y=#3%
\divide\Y by 2%
\Z=\Y%
\divide\Z by 2%
\begin{picture}(0,0)%
\put(\Y,-\Z){\line(-2,1){#3}}%
\put(-\Y,\Z){\wnwhead}%
\truex{400}%
\put(0,0){\distsign}%
\truex{200}\truey{800}\truez{600}%
\put(\value{x},\value{x}){\makebox(0,\value{z})[l]{${#1}$}}%
\put(-\value{x},-\value{y}){\makebox(0,\value{z})[r]{${#2}$}}%
\end{picture}}%

\newcommand{\WNWDOTAR}[3]{%
\truex{100}\truey{268}\truez{134}%
\Y=#3%
\divide\Y by 2%
\Z=\Y%
\divide\Z by 2%
\NUMBEROFDOTS=#3%
\divide\NUMBEROFDOTS by \value{y}%
\advance\NUMBEROFDOTS by 1%
\begin{picture}(0,0)%
\multiput(\Y,-\Z)(-\value{y},\value{z}){\NUMBEROFDOTS}%
{\circle*{\value{x}}}%
\put(-\Y,\Z){\wnwhead}%
\truex{200}\truey{800}\truez{600}%
\put(\value{x},\value{x}){\makebox(0,\value{z})[l]{${#1}$}}%
\put(-\value{x},-\value{y}){\makebox(0,\value{z})[r]{${#2}$}}%
\end{picture}}%

\newcommand{\WNWMONO}[3]{%
\Y=#3%
\divide\Y by 2%
\Z=\Y%
\divide\Z by 2%
\TrueTail%
\bimolength=#3%
\advance\bimolength by -\TrueMonoTail%
\monolength=\bimolength%
\advance\monolength by -\Y%
\secondmonolength=\monolength%
\divide\secondmonolength by 2%
\begin{picture}(0,0)%
\put(\monolength,-\secondmonolength){\line(-2,1){\bimolength}}%
\put(\monolength,-\secondmonolength){\wnwhead}%
\put(-\Y,\Z){\wnwhead}%
\truex{200}\truey{800}\truez{600}%
\put(\value{x},\value{x}){\makebox(0,\value{z})[l]{${#1}$}}%
\put(-\value{x},-\value{y}){\makebox(0,\value{z})[r]{${#2}$}}%
\end{picture}}%

\newcommand{\WNWEPI}[3]{%
\Y=#3%
\divide\Y by 2%
\Z=\Y%
\divide\Z by 2%
\TrueHead%
\bimolength=#3%
\advance\bimolength by -\TrueEpiHead%
\epilength=\bimolength%
\advance\epilength by -\Y%
\secondepilength=\epilength%
\divide\secondepilength by 2%
\begin{picture}(0,0)%
\put(\Y,-\Z){\line(-2,1){\bimolength}}%
\put(-\epilength,\secondepilength){\wnwhead}%
\put(-\Y,\Z){\wnwhead}%
\truex{200}\truey{800}\truez{600}%
\put(\value{x},\value{x}){\makebox(0,\value{z})[l]{${#1}$}}%
\put(-\value{x},-\value{y}){\makebox(0,\value{z})[r]{${#2}$}}%
\end{picture}}%

\newcommand{\WNWBIMO}[3]{%
\Y=#3%
\divide\Y by 2%
\Z=\Y%
\divide\Z by 2%
\TrueTail\TrueHead%
\bimolength=#3%
\advance\bimolength by -\TrueMonoTail%
\monolength=\bimolength%
\advance\monolength by -\Y%
\advance\bimolength by -\TrueEpiHead%
\epilength=\bimolength%
\advance\epilength by -\monolength%
\secondmonolength=\monolength%
\divide\secondmonolength by 2%
\secondepilength=\epilength%
\divide\secondepilength by 2%
\begin{picture}(0,0)%
\put(\monolength,-\secondmonolength){\line(-2,1){\bimolength}}%
\put(\monolength,-\secondmonolength){\wnwhead}%
\put(-\epilength,\secondepilength){\wnwhead}%
\put(-\Y,\Z){\wnwhead}%
\truex{200}\truey{800}\truez{600}%
\put(\value{x},\value{x}){\makebox(0,\value{z})[l]{${#1}$}}%
\put(-\value{x},-\value{y}){\makebox(0,\value{z})[r]{${#2}$}}%
\end{picture}}%

\newcommand{\WNWBIAR}[3]{%
\Y=#3%
\divide\Y by 2%
\Z=\Y%
\divide\Z by 2%
\begin{picture}(0,0)%
\put(\Y,-\Z){\begin{picture}(0,0)%
\truex{156}\truey{313}%
\put(-\value{x},-\value{y}){\line(-2,1){#3}}%
\put(\value{x},\value{y}){\line(-2,1){#3}}%
\monolength=#3%
\advance\monolength by -\value{x}%
\epilength=#3%
\advance\epilength by \value{x}%
\secondmonolength=\Y%
\advance\secondmonolength by -\value{y}%
\secondepilength=\Y%
\advance\secondepilength by \value{y}%
\put(-\monolength,\secondepilength){\wnwhead}%
\put(-\epilength,\secondmonolength){\wnwhead}%
\end{picture}}
\truex{400}\truey{1000}\truez{600}%
\put(\value{x},\value{x}){\makebox(0,\value{z})[l]{${#1}$}}%
\put(-\value{x},-\value{y}){\makebox(0,\value{z})[r]{${#2}$}}%
\end{picture}}%

\newcommand{\WNWBIDIST}[3]{%
\Y=#3%
\divide\Y by 2%
\Z=\Y%
\divide\Z by 2%
\begin{picture}(0,0)%
\truex{156}\truey{313}\truez{400}%
\put(\Y,-\Z){\begin{picture}(0,0)%
\put(-\value{x},-\value{y}){\line(-2,1){#3}}%
\put(\value{x},\value{y}){\line(-2,1){#3}}%
\monolength=#3%
\advance\monolength by -\value{x}%
\epilength=#3%
\advance\epilength by \value{x}%
\secondmonolength=\Y%
\advance\secondmonolength by -\value{y}%
\secondepilength=\Y%
\advance\secondepilength by \value{y}%
\put(-\monolength,\secondepilength){\wnwhead}%
\put(-\epilength,\secondmonolength){\wnwhead}%
\end{picture}}
\put(\value{x},\value{y}){\distsign}%
\put(-\value{x},-\value{y}){\distsign}%
\truex{400}\truey{1000}\truez{600}%
\put(\value{x},\value{x}){\makebox(0,\value{z})[l]{${#1}$}}%
\put(-\value{x},-\value{y}){\makebox(0,\value{z})[r]{${#2}$}}%
\end{picture}}%

\newcommand{\WNWADJAR}[3]{%
\Y=#3%
\divide\Y by 2%
\Z=\Y%
\divide\Z by 2%
\begin{picture}(0,0)%
\put(\Y,-\Z){\begin{picture}(0,0)%
\truex{156}\truey{313}%
\monolength=#3%
\advance\monolength by -\value{x}%
\epilength=#3%
\advance\epilength by \value{x}%
\secondmonolength=\Y%
\advance\secondmonolength by -\value{y}%
\secondepilength=\Y%
\advance\secondepilength by \value{y}%
\put(-\value{x},-\value{y}){\line(-2,1){#3}}%
\put(-\epilength,\secondmonolength){\wnwhead}%
\put(-\monolength,\secondepilength){\line(2,-1){#3}}%
\put(\value{x},\value{y}){\esehead}%
\end{picture}}
\truex{400}\truey{1000}\truez{600}%
\put(\value{x},\value{x}){\makebox(0,\value{z})[l]{${#1}$}}%
\put(-\value{x},-\value{y}){\makebox(0,\value{z})[r]{${#2}$}}%
\end{picture}}%

\newcommand{\WNWADJDIST}[3]{%
\Y=#3%
\divide\Y by 2%
\Z=\Y%
\divide\Z by 2%
\begin{picture}(0,0)%
\truex{156}\truey{313}\truez{400}%
\put(\Y,-\Z){\begin{picture}(0,0)%
\monolength=#3%
\advance\monolength by -\value{x}%
\epilength=#3%
\advance\epilength by \value{x}%
\secondmonolength=\Y%
\advance\secondmonolength by -\value{y}%
\secondepilength=\Y%
\advance\secondepilength by \value{y}%
\put(-\value{x},-\value{y}){\line(-2,1){#3}}%
\put(-\epilength,\secondmonolength){\wnwhead}%
\put(-\monolength,\secondepilength){\line(2,-1){#3}}%
\put(\value{x},\value{y}){\esehead}%
\end{picture}}
\put(\value{x},\value{y}){\distsign}%
\put(-\value{x},-\value{y}){\distsign}%
\truex{400}\truey{1000}\truez{600}%
\put(\value{x},\value{x}){\makebox(0,\value{z})[l]{${#1}$}}%
\put(-\value{x},-\value{y}){\makebox(0,\value{z})[r]{${#2}$}}%
\end{picture}}%

\def\basicwnwar[#1]{\WNWAR{}{}{#100}}%

\newcommand{\wnwar}{\@ifnextchar[{\basicwnwar}{\basicwnwar[133]}}%

\def\basicWnwar[#1]#2{\WNWAR{#2}{}{#100}}%

\newcommand{\Wnwar}{\@ifnextchar[{\basicWnwar}{\basicWnwar[133]}}%

\def\basicwnwaR[#1]#2{\WNWAR{}{#2}{#100}}%

\newcommand{\wnwaR}{\@ifnextchar[{\basicwnwaR}{\basicwnwaR[133]}}%

\def\basicwnwdist[#1]{\WNWDIST{}{}{#100}}%

\newcommand{\wnwdist}{\@ifnextchar[{\basicwnwdist}{\basicwnwdist[133]}}%

\def\basicWnwdist[#1]#2{\WNWDIST{#2}{}{#100}}%

\newcommand{\Wnwdist}{\@ifnextchar[{\basicWnwdist}{\basicWnwdist[133]}}%

\def\basicwnwdisT[#1]#2{\WNWDIST{}{#2}{#100}}%

\newcommand{\wnwdisT}{\@ifnextchar[{\basicwnwdisT}{\basicwnwdisT[133]}}%

\def\basicwnwdotar[#1]{\WNWDOTAR{}{}{#100}}%

\newcommand{\wnwdotar}{\@ifnextchar[{\basicwnwdotar}{\basicwnwdotar[133]}}%

\def\basicWnwdotar[#1]#2{\WNWDOTAR{#2}{}{#100}}%

\newcommand{\Wnwdotar}{\@ifnextchar[{\basicWnwdotar}{\basicWnwdotar[133]}}%

\def\basicwnwdotaR[#1]#2{\WNWDOTAR{}{#2}{#100}}%

\newcommand{\wnwdotaR}{\@ifnextchar[{\basicwnwdotaR}{\basicwnwdotaR[133]}}%

\def\basicwnwmono[#1]{\WNWMONO{}{}{#100}}%

\newcommand{\wnwmono}{\@ifnextchar[{\basicwnwmono}{\basicwnwmono[133]}}%

\def\basicWnwmono[#1]#2{\WNWMONO{#2}{}{#100}}%

\newcommand{\Wnwmono}{\@ifnextchar[{\basicWnwmono}{\basicWnwmono[133]}}%

\def\basicwnwmonO[#1]#2{\WNWMONO{}{#2}{#100}}%

\newcommand{\wnwmonO}{\@ifnextchar[{\basicwnwmonO}{\basicwnwmonO[133]}}%

\def\basicwnwepi[#1]{\WNWEPI{}{}{#100}}%

\newcommand{\wnwepi}{\@ifnextchar[{\basicwnwepi}{\basicwnwepi[133]}}%

\def\basicWnwepi[#1]#2{\WNWEPI{#2}{}{#100}}%

\newcommand{\Wnwepi}{\@ifnextchar[{\basicWnwepi}{\basicWnwepi[133]}}%

\def\basicwnwepI[#1]#2{\WNWEPI{}{#2}{#100}}%

\newcommand{\wnwepI}{\@ifnextchar[{\basicwnwepI}{\basicwnwepI[133]}}%

\def\basicwnwbimo[#1]{\WNWBIMO{}{}{#100}}%

\newcommand{\wnwbimo}{\@ifnextchar[{\basicwnwbimo}{\basicwnwbimo[133]}}%

\def\basicWnwbimo[#1]#2{\WNWBIMO{#2}{}{#100}}%

\newcommand{\Wnwbimo}{\@ifnextchar[{\basicWnwbimo}{\basicWnwbimo[133]}}%

\def\basicwnwbimO[#1]#2{\WNWBIMO{}{#2}{#100}}%

\newcommand{\wnwbimO}{\@ifnextchar[{\basicwnwbimO}{\basicwnwbimO[133]}}%

\def\basicwnwiso[#1]{\WNWAR{\isosign}{}{#100}}%

\newcommand{\wnwiso}{\@ifnextchar[{\basicwnwiso}{\basicwnwiso[133]}}%

\def\basicWnwiso[#1]#2{\WNWAR{#2}{\isosign}{#100}}%

\newcommand{\Wnwiso}{\@ifnextchar[{\basicWnwiso}{\basicWnwiso[133]}}%

\def\basicwnwisO[#1]#2{\WNWAR{\isosign}{#2}{#100}}%

\newcommand{\wnwisO}{\@ifnextchar[{\basicwnwisO}{\basicwnwisO[133]}}%

\def\basicwnwbiar[#1]{\WNWBIAR{}{}{#100}}%

\newcommand{\wnwbiar}{\@ifnextchar[{\basicwnwbiar}{\basicwnwbiar[133]}}%

\def\basicWnwbiar[#1]#2#3{\WNWBIAR{#2}{#3}{#100}}%

\newcommand{\Wnwbiar}{\@ifnextchar[{\basicWnwbiar}{\basicWnwbiar[133]}}%

\def\basicwnwbidist[#1]{\WNWBIDIST{}{}{#100}}%

\newcommand{\wnwbidist}{\@ifnextchar[{\basicwnwbidist}{\basicwnwbidist[133]}}%

\def\basicWnwbidist[#1]#2#3{\WNWBIDIST{#2}{#3}{#100}}%

\newcommand{\Wnwbidist}{\@ifnextchar[{\basicWnwbidist}{\basicWnwbidist[133]}}%

\def\basicwnwadjar[#1]{\WNWADJAR{}{}{#100}}%

\newcommand{\wnwadjar}{\@ifnextchar[{\basicwnwadjar}{\basicwnwadjar[133]}}%

\def\basicWnwadjar[#1]#2#3{\WNWADJAR{#2}{#3}{#100}}%

\newcommand{\Wnwadjar}{\@ifnextchar[{\basicWnwadjar}{\basicWnwadjar[133]}}%

\def\basicwnwadjdist[#1]{\WNWADJDIST{}{}{#100}}%

\newcommand{\wnwadjdist}{\@ifnextchar[{\basicwnwadjdist}{\basicwnwadjdist[13
3]}}%

\def\basicWnwadjdist[#1]#2#3{\WNWADJDIST{#2}{#3}{#100}}%

\newcommand{\Wnwadjdist}{\@ifnextchar[{\basicWnwadjdist}{\basicWnwadjdist[13
3]}}%

\newcommand{\NNEAR}[3]{%
\Z=#3%
\divide\Z by 2%
\begin{picture}(0,0)%
\put(-\Z,-#3){\line(1,2){#3}}%
\put(\Z,#3){\nnehead}%
\truex{200}\truey{800}\truez{600}%
\put(-\value{x},\value{x}){\makebox(0,\value{z})[r]{${#1}$}}%
\put(\value{x},-\value{y}){\makebox(0,\value{z})[l]{${#2}$}}%
\end{picture}}%

\newcommand{\NNEDIST}[3]{%
\Z=#3%
\divide\Z by 2%
\begin{picture}(0,0)%
\put(-\Z,-#3){\line(1,2){#3}}%
\put(\Z,#3){\nnehead}%
\truex{400}%
\put(0,0){\distsign}%
\truex{200}\truey{800}\truez{600}%
\put(-\value{x},\value{x}){\makebox(0,\value{z})[r]{${#1}$}}%
\put(\value{x},-\value{y}){\makebox(0,\value{z})[l]{${#2}$}}%
\end{picture}}%

\newcommand{\NNEDOTAR}[3]{%
\truex{100}\truey{268}\truez{134}%
\Z=#3%
\divide\Z by 2%
\NUMBEROFDOTS=#3%
\divide\NUMBEROFDOTS by \value{z}%
\advance\NUMBEROFDOTS by 1%
\begin{picture}(0,0)%
\multiput(-\Z,-#3)(\value{z},\value{y}){\NUMBEROFDOTS}%
{\circle*{\value{x}}}%
\put(\Z,#3){\nnehead}%
\truex{200}\truey{800}\truez{600}%
\put(-\value{x},\value{x}){\makebox(0,\value{z})[r]{${#1}$}}%
\put(\value{x},-\value{y}){\makebox(0,\value{z})[l]{${#2}$}}%
\end{picture}}%

\newcommand{\NNEMONO}[3]{%
\Z=#3%
\divide\Z by 2%
\truetaiL%
\bimolength=#3%
\advance\bimolength by -\truemonotaiL%
\monolength=\bimolength%
\advance\monolength by -\Z%
\secondmonolength=\monolength%
\multiply\secondmonolength by 2%
\begin{picture}(0,0)%
\put(-\monolength,-\secondmonolength){\line(1,2){\bimolength}}%
\put(-\monolength,-\secondmonolength){\nnehead}%
\put(\Z,#3){\nnehead}%
\truex{200}\truey{800}\truez{600}%
\put(-\value{x},\value{x}){\makebox(0,\value{z})[r]{${#1}$}}%
\put(\value{x},-\value{y}){\makebox(0,\value{z})[l]{${#2}$}}%
\end{picture}}%

\newcommand{\NNEEPI}[3]{%
\Z=#3%
\divide\Z by 2%
\trueheaD%
\bimolength=#3%
\advance\bimolength by -\trueepiheaD%
\epilength=\bimolength%
\advance\epilength by -\Z%
\secondepilength=\epilength%
\multiply\secondepilength by 2%
\begin{picture}(0,0)%
\put(-\Z,-#3){\line(1,2){\bimolength}}%
\put(\epilength,\secondepilength){\nnehead}%
\put(\Z,#3){\nnehead}%
\truex{200}\truey{800}\truez{600}%
\put(-\value{x},\value{x}){\makebox(0,\value{z})[r]{${#1}$}}%
\put(\value{x},-\value{y}){\makebox(0,\value{z})[l]{${#2}$}}%
\end{picture}}%

\newcommand{\NNEBIMO}[3]{%
\Z=#3%
\divide\Z by 2%
\truetaiL\trueheaD%
\bimolength=#3%
\advance\bimolength by -\truemonotaiL%
\monolength=\bimolength%
\advance\monolength by -\Z%
\advance\bimolength by -\trueepiheaD%
\epilength=\bimolength%
\advance\epilength by -\monolength%
\secondmonolength=\monolength%
\multiply\secondmonolength by 2%
\secondepilength=\epilength%
\multiply\secondepilength by 2%
\begin{picture}(0,0)%
\put(-\monolength,-\secondmonolength){\line(1,2){\bimolength}}%
\put(-\monolength,-\secondmonolength){\nnehead}%
\put(\epilength,\secondepilength){\nnehead}%
\put(\Z,#3){\nnehead}%
\truex{200}\truey{800}\truez{600}%
\put(-\value{x},\value{x}){\makebox(0,\value{z})[r]{${#1}$}}%
\put(\value{x},-\value{y}){\makebox(0,\value{z})[l]{${#2}$}}%
\end{picture}}%

\newcommand{\NNEEQL}[3]{%
\Z=#3%
\divide\Z by 2%
\begin{picture}(0,0)%
\put(-\Z,-#3){\begin{picture}(0,0)%
\truex{44}\truey{89}%
\put(-\value{y},\value{x}){\line(1,2){#3}}%
\put(\value{y},-\value{x}){\line(1,2){#3}}%
\end{picture}}%
\truex{200}\truey{800}\truez{600}%
\put(-\value{x},\value{x}){\makebox(0,\value{z})[r]{${#1}$}}%
\put(\value{x},-\value{y}){\makebox(0,\value{z})[l]{${#2}$}}%
\end{picture}}%

\newcommand{\NNEBIAR}[3]{%
\Y=#3%
\divide\Y by 2%
\Z=#3%
\multiply \Z by 2%
\begin{picture}(0,0)%
\put(-\Y,-#3){\begin{picture}(0,0)%
\truex{313}\truey{156}%
\put(-\value{x},\value{y}){\line(1,2){#3}}%
\put(\value{x},-\value{y}){\line(1,2){#3}}%
\monolength=#3%
\advance\monolength by -\value{x}%
\epilength=#3%
\advance\epilength by \value{x}%
\secondmonolength=\Z%
\advance\secondmonolength by -\value{y}%
\secondepilength=\Z%
\advance\secondepilength by \value{y}%
\put(\monolength,\secondepilength){\nnehead}%
\put(\epilength,\secondmonolength){\nnehead}%
\end{picture}}
\truex{300}\truey{1000}\truez{600}%
\put(-\value{x},\value{x}){\makebox(0,\value{z})[r]{${#1}$}}%
\put(\value{x},-\value{y}){\makebox(0,\value{z})[l]{${#2}$}}%
\end{picture}}%

\newcommand{\NNEBIDIST}[3]{%
\Y=#3%
\divide\Y by 2%
\Z=#3%
\multiply \Z by 2%
\begin{picture}(0,0)%
\truex{313}\truey{156}\truez{400}%
\put(-\Y,-#3){\begin{picture}(0,0)%
\put(-\value{x},\value{y}){\line(1,2){#3}}%
\put(\value{x},-\value{y}){\line(1,2){#3}}%
\monolength=#3%
\advance\monolength by -\value{x}%
\epilength=#3%
\advance\epilength by \value{x}%
\secondmonolength=\Z%
\advance\secondmonolength by -\value{y}%
\secondepilength=\Z%
\advance\secondepilength by \value{y}%
\put(\monolength,\secondepilength){\nnehead}%
\put(\epilength,\secondmonolength){\nnehead}%
\end{picture}}
\put(-\value{x},\value{y}){\distsign}%
\put(\value{x},-\value{y}){\distsign}%
\truex{300}\truey{1000}\truez{600}%
\put(-\value{x},\value{x}){\makebox(0,\value{z})[r]{${#1}$}}%
\put(\value{x},-\value{y}){\makebox(0,\value{z})[l]{${#2}$}}%
\end{picture}}%

\newcommand{\NNEADJAR}[3]{%
\Y=#3%
\divide\Y by 2%
\Z=#3%
\multiply \Z by 2%
\begin{picture}(0,0)%
\put(-\Y,-#3){\begin{picture}(0,0)%
\truex{313}\truey{156}%
\monolength=#3%
\advance\monolength by -\value{x}%
\epilength=#3%
\advance\epilength by \value{x}%
\secondmonolength=\Z%
\advance\secondmonolength by -\value{y}%
\secondepilength=\Z%
\advance\secondepilength by \value{y}%
\put(\value{x},-\value{y}){\line(1,2){#3}}%
\put(\epilength,\secondmonolength){\nnehead}%
\put(\monolength,\secondepilength){\line(-1,-2){#3}}%
\put(-\value{x},\value{y}){\sswhead}%
\end{picture}}
\truex{300}\truey{1000}\truez{600}%
\put(-\value{x},\value{x}){\makebox(0,\value{z})[r]{${#1}$}}%
\put(\value{x},-\value{y}){\makebox(0,\value{z})[l]{${#2}$}}%
\end{picture}}%

\newcommand{\NNEADJDIST}[3]{%
\Y=#3%
\divide\Y by 2%
\Z=#3%
\multiply \Z by 2%
\begin{picture}(0,0)%
\truex{313}\truey{156}\truez{400}%
\put(-\Y,-#3){\begin{picture}(0,0)%
\monolength=#3%
\advance\monolength by -\value{x}%
\epilength=#3%
\advance\epilength by \value{x}%
\secondmonolength=\Z%
\advance\secondmonolength by -\value{y}%
\secondepilength=\Z%
\advance\secondepilength by \value{y}%
\put(\value{x},-\value{y}){\line(1,2){#3}}%
\put(\epilength,\secondmonolength){\nnehead}%
\put(\monolength,\secondepilength){\line(-1,-2){#3}}%
\put(-\value{x},\value{y}){\sswhead}%
\end{picture}}
\put(-\value{x},\value{y}){\distsign}%
\put(\value{x},-\value{y}){\distsign}%
\truex{300}\truey{1000}\truez{600}%
\put(-\value{x},\value{x}){\makebox(0,\value{z})[r]{${#1}$}}%
\put(\value{x},-\value{y}){\makebox(0,\value{z})[l]{${#2}$}}%
\end{picture}}%

\def\basicnnear[#1]{\NNEAR{}{}{#100}}%

\newcommand{\nnear}{\@ifnextchar[{\basicnnear}{\basicnnear[67]}}%

\def\basicNnear[#1]#2{\NNEAR{#2}{}{#100}}%

\newcommand{\Nnear}{\@ifnextchar[{\basicNnear}{\basicNnear[67]}}%

\def\basicnneaR[#1]#2{\NNEAR{}{#2}{#100}}%

\newcommand{\nneaR}{\@ifnextchar[{\basicnneaR}{\basicnneaR[67]}}%

\def\basicnnedist[#1]{\NNEDIST{}{}{#100}}%

\newcommand{\nnedist}{\@ifnextchar[{\basicnnedist}{\basicnnedist[67]}}%

\def\basicNnedist[#1]#2{\NNEDIST{#2}{}{#100}}%

\newcommand{\Nnedist}{\@ifnextchar[{\basicNnedist}{\basicNnedist[67]}}%

\def\basicnnedisT[#1]#2{\NNEDIST{}{#2}{#100}}%

\newcommand{\nnedisT}{\@ifnextchar[{\basicnnedisT}{\basicnnedisT[67]}}%

\def\basicnnedotar[#1]{\NNEDOTAR{}{}{#100}}%

\newcommand{\nnedotar}{\@ifnextchar[{\basicnnedotar}{\basicnnedotar[67]}}%

\def\basicNnedotar[#1]#2{\NNEDOTAR{#2}{}{#100}}%

\newcommand{\Nnedotar}{\@ifnextchar[{\basicNnedotar}{\basicNnedotar[67]}}%

\def\basicnnedotaR[#1]#2{\NNEDOTAR{}{#2}{#100}}%

\newcommand{\nnedotaR}{\@ifnextchar[{\basicnnedotaR}{\basicnnedotaR[67]}}%

\def\basicnnemono[#1]{\NNEMONO{}{}{#100}}%

\newcommand{\nnemono}{\@ifnextchar[{\basicnnemono}{\basicnnemono[67]}}%

\def\basicNnemono[#1]#2{\NNEMONO{#2}{}{#100}}%

\newcommand{\Nnemono}{\@ifnextchar[{\basicNnemono}{\basicNnemono[67]}}%

\def\basicnnemonO[#1]#2{\NNEMONO{}{#2}{#100}}%

\newcommand{\nnemonO}{\@ifnextchar[{\basicnnemonO}{\basicnnemonO[67]}}%

\def\basicnneepi[#1]{\NNEEPI{}{}{#100}}%

\newcommand{\nneepi}{\@ifnextchar[{\basicnneepi}{\basicnneepi[67]}}%

\def\basicNneepi[#1]#2{\NNEEPI{#2}{}{#100}}%

\newcommand{\Nneepi}{\@ifnextchar[{\basicNneepi}{\basicNneepi[67]}}%

\def\basicnneepI[#1]#2{\NNEEPI{}{#2}{#100}}%

\newcommand{\nneepI}{\@ifnextchar[{\basicnneepI}{\basicnneepI[67]}}%

\def\basicnnebimo[#1]{\NNEBIMO{}{}{#100}}%

\newcommand{\nnebimo}{\@ifnextchar[{\basicnnebimo}{\basicnnebimo[67]}}%

\def\basicNnebimo[#1]#2{\NNEBIMO{#2}{}{#100}}%

\newcommand{\Nnebimo}{\@ifnextchar[{\basicNnebimo}{\basicNnebimo[67]}}%

\def\basicnnebimO[#1]#2{\NNEBIMO{}{#2}{#100}}%

\newcommand{\nnebimO}{\@ifnextchar[{\basicnnebimO}{\basicnnebimO[67]}}%

\def\basicnneiso[#1]{\NNEAR{\isosign}{}{#100}}%

\newcommand{\nneiso}{\@ifnextchar[{\basicnneiso}{\basicnneiso[67]}}%

\def\basicNneiso[#1]#2{\NNEAR{#2}{\isosign}{#100}}%

\newcommand{\Nneiso}{\@ifnextchar[{\basicNneiso}{\basicNneiso[67]}}%

\def\basicnneisO[#1]#2{\NNEAR{\isosign}{#2}{#100}}%

\newcommand{\nneisO}{\@ifnextchar[{\basicnneisO}{\basicnneisO[67]}}%

\def\basicnneeql[#1]{\NNEEQL{}{}{#100}}%

\newcommand{\nneeql}{\@ifnextchar[{\basicnneeql}{\basicnneeql[67]}}%

\def\basicNneeql[#1]#2{\NNEEQL{#2}{}{#100}}%

\newcommand{\Nneeql}{\@ifnextchar[{\basicNneeql}{\basicNneeql[67]}}%

\def\basicnneeqL[#1]#2{\NNEEQL{}{#2}{#100}}%

\newcommand{\nneeqL}{\@ifnextchar[{\basicnneeqL}{\basicnneeqL[67]}}%

\def\basicnnebiar[#1]{\NNEBIAR{}{}{#100}}%

\newcommand{\nnebiar}{\@ifnextchar[{\basicnnebiar}{\basicnnebiar[67]}}%

\def\basicNnebiar[#1]#2#3{\NNEBIAR{#2}{#3}{#100}}%

\newcommand{\Nnebiar}{\@ifnextchar[{\basicNnebiar}{\basicNnebiar[67]}}%

\def\basicnnebidist[#1]{\NNEBIDIST{}{}{#100}}%

\newcommand{\nnebidist}{\@ifnextchar[{\basicnnebidist}{\basicnnebidist[67]}}%

\def\basicNnebidist[#1]#2#3{\NNEBIDIST{#2}{#3}{#100}}%

\newcommand{\Nnebidist}{\@ifnextchar[{\basicNnebidist}{\basicNnebidist[67]}}%

\def\basicnneadjar[#1]{\NNEADJAR{}{}{#100}}%

\newcommand{\nneadjar}{\@ifnextchar[{\basicnneadjar}{\basicnneadjar[67]}}%

\def\basicNneadjar[#1]#2#3{\NNEADJAR{#2}{#3}{#100}}%

\newcommand{\Nneadjar}{\@ifnextchar[{\basicNneadjar}{\basicNneadjar[67]}}%

\def\basicnneadjdist[#1]{\NNEADJDIST{}{}{#100}}%

\newcommand{\nneadjdist}{\@ifnextchar[{\basicnneadjdist}{\basicnneadjdist[67]}}%

\def\basicNneadjdist[#1]#2#3{\NNEADJDIST{#2}{#3}{#100}}%

\newcommand{\Nneadjdist}{\@ifnextchar[{\basicNneadjdist}{\basicNneadjdist[67]}}%

\newcommand{\SSEAR}[3]{%
\Z=#3%
\divide\Z by 2%
\begin{picture}(0,0)%
\put(-\Z,#3){\line(1,-2){#3}}%
\put(\Z,-#3){\ssehead}%
\truex{200}\truey{800}\truez{600}%
\put(\value{x},\value{x}){\makebox(0,\value{z})[l]{${#1}$}}%
\put(-\value{x},-\value{y}){\makebox(0,\value{z})[r]{${#2}$}}%
\end{picture}}%

\newcommand{\SSEDIST}[3]{%
\Z=#3%
\divide\Z by 2%
\begin{picture}(0,0)%
\put(-\Z,#3){\line(1,-2){#3}}%
\put(\Z,-#3){\ssehead}%
\truex{400}%
\put(0,0){\distsign}%
\truex{200}\truey{800}\truez{600}%
\put(\value{x},\value{x}){\makebox(0,\value{z})[l]{${#1}$}}%
\put(-\value{x},-\value{y}){\makebox(0,\value{z})[r]{${#2}$}}%
\end{picture}}%

\newcommand{\SSEDOTAR}[3]{%
\truex{100}\truey{268}\truez{134}%
\Z=#3%
\divide\Z by 2%
\NUMBEROFDOTS=#3%
\divide\NUMBEROFDOTS by \value{z}%
\advance\NUMBEROFDOTS by 1%
\begin{picture}(0,0)%
\multiput(-\Z,#3)(\value{z},-\value{y}){\NUMBEROFDOTS}%
{\circle*{\value{x}}}%
\put(\Z,-#3){\ssehead}%
\truex{200}\truey{800}\truez{600}%
\put(\value{x},\value{x}){\makebox(0,\value{z})[l]{${#1}$}}%
\put(-\value{x},-\value{y}){\makebox(0,\value{z})[r]{${#2}$}}%
\end{picture}}%

\newcommand{\SSEMONO}[3]{%
\Z=#3%
\divide\Z by 2%
\truetaiL%
\bimolength=#3%
\advance\bimolength by -\truemonotaiL%
\monolength=\bimolength%
\advance\monolength by -\Z%
\secondmonolength=\monolength%
\multiply\secondmonolength by 2%
\begin{picture}(0,0)%
\put(-\monolength,\secondmonolength){\line(1,-2){\bimolength}}%
\put(-\monolength,\secondmonolength){\ssehead}%
\put(\Z,-#3){\ssehead}%
\truex{200}\truey{800}\truez{600}%
\put(\value{x},\value{x}){\makebox(0,\value{z})[l]{${#1}$}}%
\put(-\value{x},-\value{y}){\makebox(0,\value{z})[r]{${#2}$}}%
\end{picture}}%

\newcommand{\SSEEPI}[3]{%
\Z=#3%
\divide\Z by 2%
\trueheaD%
\bimolength=#3%
\advance\bimolength by -\trueepiheaD%
\epilength=\bimolength%
\advance\epilength by -\Z%
\secondepilength=\epilength%
\multiply\secondepilength by 2%
\begin{picture}(0,0)%
\put(-\Z,#3){\line(1,-2){\bimolength}}%
\put(\epilength,-\secondepilength){\ssehead}%
\put(\Z,-#3){\ssehead}%
\truex{200}\truey{800}\truez{600}%
\put(\value{x},\value{x}){\makebox(0,\value{z})[l]{${#1}$}}%
\put(-\value{x},-\value{y}){\makebox(0,\value{z})[r]{${#2}$}}%
\end{picture}}%

\newcommand{\SSEBIMO}[3]{%
\Z=#3%
\divide\Z by 2%
\truetaiL\trueheaD%
\bimolength=#3%
\advance\bimolength by -\truemonotaiL%
\monolength=\bimolength%
\advance\monolength by -\Z%
\advance\bimolength by -\trueepiheaD%
\epilength=\bimolength%
\advance\epilength by -\monolength%
\secondmonolength=\monolength%
\multiply\secondmonolength by 2%
\secondepilength=\epilength%
\multiply\secondepilength by 2%
\begin{picture}(0,0)%
\put(-\monolength,\secondmonolength){\line(1,-2){\bimolength}}%
\put(-\monolength,\secondmonolength){\ssehead}%
\put(\epilength,-\secondepilength){\ssehead}%
\put(\Z,-#3){\ssehead}%
\truex{200}\truey{800}\truez{600}%
\put(\value{x},\value{x}){\makebox(0,\value{z})[l]{${#1}$}}%
\put(-\value{x},-\value{y}){\makebox(0,\value{z})[r]{${#2}$}}%
\end{picture}}%

\newcommand{\SSEEQL}[3]{%
\Z=#3%
\divide\Z by 2%
\begin{picture}(0,0)%
\put(-\Z,#3){\begin{picture}(0,0)%
\truex{44}\truey{89}%
\put(-\value{y},-\value{x}){\line(1,-2){#3}}%
\put(\value{y},\value{x}){\line(1,-2){#3}}%
\end{picture}}%
\truex{200}\truey{800}\truez{600}%
\put(\value{x},\value{x}){\makebox(0,\value{z})[l]{${#1}$}}%
\put(-\value{x},-\value{y}){\makebox(0,\value{z})[r]{${#2}$}}%
\end{picture}}%

\newcommand{\SSEBIAR}[3]{%
\Y=#3%
\divide\Y by 2%
\Z=#3%
\multiply \Z by 2%
\begin{picture}(0,0)%
\put(-\Y,#3){\begin{picture}(0,0)%
\truex{313}\truey{156}%
\put(-\value{x},-\value{y}){\line(1,-2){#3}}%
\put(\value{x},\value{y}){\line(1,-2){#3}}%
\monolength=#3%
\advance\monolength by -\value{x}%
\epilength=#3%
\advance\epilength by \value{x}%
\secondmonolength=\Z%
\advance\secondmonolength by -\value{y}%
\secondepilength=\Z%
\advance\secondepilength by \value{y}%
\put(\monolength,-\secondepilength){\ssehead}%
\put(\epilength,-\secondmonolength){\ssehead}%
\end{picture}}
\truex{400}\truey{1000}\truez{600}%
\put(\value{x},\value{x}){\makebox(0,\value{z})[l]{${#1}$}}%
\put(-\value{x},-\value{y}){\makebox(0,\value{z})[r]{${#2}$}}%
\end{picture}}%

\newcommand{\SSEBIDIST}[3]{%
\Y=#3%
\divide\Y by 2%
\Z=#3%
\multiply \Z by 2%
\begin{picture}(0,0)%
\truex{313}\truey{156}\truez{400}%
\put(-\Y,#3){\begin{picture}(0,0)%
\put(-\value{x},-\value{y}){\line(1,-2){#3}}%
\put(\value{x},\value{y}){\line(1,-2){#3}}%
\monolength=#3%
\advance\monolength by -\value{x}%
\epilength=#3%
\advance\epilength by \value{x}%
\secondmonolength=\Z%
\advance\secondmonolength by -\value{y}%
\secondepilength=\Z%
\advance\secondepilength by \value{y}%
\put(\monolength,-\secondepilength){\ssehead}%
\put(\epilength,-\secondmonolength){\ssehead}%
\end{picture}}
\put(-\value{x},-\value{y}){\distsign}%
\put(\value{x},\value{y}){\distsign}%
\truex{500}\truey{1000}\truez{600}%
\put(\value{x},\value{x}){\makebox(0,\value{z})[l]{${#1}$}}%
\put(-\value{x},-\value{y}){\makebox(0,\value{z})[r]{${#2}$}}%
\end{picture}}%

\newcommand{\SSEADJAR}[3]{%
\Y=#3%
\divide\Y by 2%
\Z=#3%
\multiply \Z by 2%
\begin{picture}(0,0)%
\put(-\Y,#3){\begin{picture}(0,0)%
\truex{313}\truey{156}%
\monolength=#3%
\advance\monolength by -\value{x}%
\epilength=#3%
\advance\epilength by \value{x}%
\secondmonolength=\Z%
\advance\secondmonolength by -\value{y}%
\secondepilength=\Z%
\advance\secondepilength by \value{y}%
\put(-\value{x},-\value{y}){\line(1,-2){#3}}%
\put(\monolength,-\secondepilength){\ssehead}%
\put(\epilength,-\secondmonolength){\line(-1,2){#3}}%
\put(\value{x},\value{y}){\nnwhead}%
\end{picture}}
\truex{400}\truey{1000}\truez{600}%
\put(\value{x},\value{x}){\makebox(0,\value{z})[l]{${#1}$}}%
\put(-\value{x},-\value{y}){\makebox(0,\value{z})[r]{${#2}$}}%
\end{picture}}%

\newcommand{\SSEADJDIST}[3]{%
\Y=#3%
\divide\Y by 2%
\Z=#3%
\multiply \Z by 2%
\begin{picture}(0,0)%
\truex{313}\truey{156}\truez{400}%
\put(-\Y,#3){\begin{picture}(0,0)%
\monolength=#3%
\advance\monolength by -\value{x}%
\epilength=#3%
\advance\epilength by \value{x}%
\secondmonolength=\Z%
\advance\secondmonolength by -\value{y}%
\secondepilength=\Z%
\advance\secondepilength by \value{y}%
\put(-\value{x},-\value{y}){\line(1,-2){#3}}%
\put(\monolength,-\secondepilength){\ssehead}%
\put(\epilength,-\secondmonolength){\line(-1,2){#3}}%
\put(\value{x},\value{y}){\nnwhead}%
\end{picture}}
\put(\value{x},\value{y}){\distsign}%
\put(-\value{x},-\value{y}){\distsign}%
\truex{500}\truey{1000}\truez{600}%
\put(\value{x},\value{x}){\makebox(0,\value{z})[l]{${#1}$}}%
\put(-\value{x},-\value{y}){\makebox(0,\value{z})[r]{${#2}$}}%
\end{picture}}%

\def\basicssear[#1]{\SSEAR{}{}{#100}}%

\newcommand{\ssear}{\@ifnextchar[{\basicssear}{\basicssear[67]}}%

\def\basicSsear[#1]#2{\SSEAR{#2}{}{#100}}%

\newcommand{\Ssear}{\@ifnextchar[{\basicSsear}{\basicSsear[67]}}%

\def\basicsseaR[#1]#2{\SSEAR{}{#2}{#100}}%

\newcommand{\sseaR}{\@ifnextchar[{\basicsseaR}{\basicsseaR[67]}}%

\def\basicssedist[#1]{\SSEDIST{}{}{#100}}%

\newcommand{\ssedist}{\@ifnextchar[{\basicssedist}{\basicssedist[67]}}%

\def\basicSsedist[#1]#2{\SSEDIST{#2}{}{#100}}%

\newcommand{\Ssedist}{\@ifnextchar[{\basicSsedist}{\basicSsedist[67]}}%

\def\basicssedisT[#1]#2{\SSEDIST{}{#2}{#100}}%

\newcommand{\ssedisT}{\@ifnextchar[{\basicssedisT}{\basicssedisT[67]}}%

\def\basicssedotar[#1]{\SSEDOTAR{}{}{#100}}%

\newcommand{\ssedotar}{\@ifnextchar[{\basicssedotar}{\basicssedotar[67]}}%

\def\basicSsedotar[#1]#2{\SSEDOTAR{#2}{}{#100}}%

\newcommand{\Ssedotar}{\@ifnextchar[{\basicSsedotar}{\basicSsedotar[67]}}%

\def\basicssedotaR[#1]#2{\SSEDOTAR{}{#2}{#100}}%

\newcommand{\ssedotaR}{\@ifnextchar[{\basicssedotaR}{\basicssedotaR[67]}}%

\def\basicssemono[#1]{\SSEMONO{}{}{#100}}%

\newcommand{\ssemono}{\@ifnextchar[{\basicssemono}{\basicssemono[67]}}%

\def\basicSsemono[#1]#2{\SSEMONO{#2}{}{#100}}%

\newcommand{\Ssemono}{\@ifnextchar[{\basicSsemono}{\basicSsemono[67]}}%

\def\basicssemonO[#1]#2{\SSEMONO{}{#2}{#100}}%

\newcommand{\ssemonO}{\@ifnextchar[{\basicssemonO}{\basicssemonO[67]}}%

\def\basicsseepi[#1]{\SSEEPI{}{}{#100}}%

\newcommand{\sseepi}{\@ifnextchar[{\basicsseepi}{\basicsseepi[67]}}%

\def\basicSseepi[#1]#2{\SSEEPI{#2}{}{#100}}%

\newcommand{\Sseepi}{\@ifnextchar[{\basicSseepi}{\basicSseepi[67]}}%

\def\basicsseepI[#1]#2{\SSEEPI{}{#2}{#100}}%

\newcommand{\sseepI}{\@ifnextchar[{\basicsseepI}{\basicsseepI[67]}}%

\def\basicssebimo[#1]{\SSEBIMO{}{}{#100}}%

\newcommand{\ssebimo}{\@ifnextchar[{\basicssebimo}{\basicssebimo[67]}}%

\def\basicSsebimo[#1]#2{\SSEBIMO{#2}{}{#100}}%

\newcommand{\Ssebimo}{\@ifnextchar[{\basicSsebimo}{\basicSsebimo[67]}}%

\def\basicssebimO[#1]#2{\SSEBIMO{}{#2}{#100}}%

\newcommand{\ssebimO}{\@ifnextchar[{\basicssebimO}{\basicssebimO[67]}}%

\def\basicsseiso[#1]{\SSEAR{\isosign}{}{#100}}%

\newcommand{\sseiso}{\@ifnextchar[{\basicsseiso}{\basicsseiso[67]}}%

\def\basicSseiso[#1]#2{\SSEAR{#2}{\isosign}{#100}}%

\newcommand{\Sseiso}{\@ifnextchar[{\basicSseiso}{\basicSseiso[67]}}%

\def\basicsseisO[#1]#2{\SSEAR{\isosign}{#2}{#100}}%

\newcommand{\sseisO}{\@ifnextchar[{\basicsseisO}{\basicsseisO[67]}}%

\def\basicsseeql[#1]{\SSEEQL{}{}{#100}}%

\newcommand{\sseeql}{\@ifnextchar[{\basicsseeql}{\basicsseeql[67]}}%

\def\basicSseeql[#1]#2{\SSEEQL{#2}{}{#100}}%

\newcommand{\Sseeql}{\@ifnextchar[{\basicSseeql}{\basicSseeql[67]}}%

\def\basicsseeqL[#1]#2{\SSEEQL{}{#2}{#100}}%

\newcommand{\sseeqL}{\@ifnextchar[{\basicsseeqL}{\basicsseeqL[67]}}%

\def\basicssebiar[#1]{\SSEBIAR{}{}{#100}}%

\newcommand{\ssebiar}{\@ifnextchar[{\basicssebiar}{\basicssebiar[67]}}%

\def\basicSsebiar[#1]#2#3{\SSEBIAR{#2}{#3}{#100}}%

\newcommand{\Ssebiar}{\@ifnextchar[{\basicSsebiar}{\basicSsebiar[67]}}%

\def\basicssebidist[#1]{\SSEBIDIST{}{}{#100}}%

\newcommand{\ssebidist}{\@ifnextchar[{\basicssebidist}{\basicssebidist[67]}}%

\def\basicSsebidist[#1]#2#3{\SSEBIDIST{#2}{#3}{#100}}%

\newcommand{\Ssebidist}{\@ifnextchar[{\basicSsebidist}{\basicSsebidist[67]}}%

\def\basicsseadjar[#1]{\SSEADJAR{}{}{#100}}%

\newcommand{\sseadjar}{\@ifnextchar[{\basicsseadjar}{\basicsseadjar[67]}}%

\def\basicSseadjar[#1]#2#3{\SSEADJAR{#2}{#3}{#100}}%

\newcommand{\Sseadjar}{\@ifnextchar[{\basicSseadjar}{\basicSseadjar[67]}}%

\def\basicsseadjdist[#1]{\SSEADJDIST{}{}{#100}}%

\newcommand{\sseadjdist}{\@ifnextchar[{\basicsseadjdist}{\basicsseadjdist[67]}}%

\def\basicSseadjdist[#1]#2#3{\SSEADJDIST{#2}{#3}{#100}}%

\newcommand{\Sseadjdist}{\@ifnextchar[{\basicSseadjdist}{\basicSseadjdist[67]}}%

\newcommand{\SSWAR}[3]{%
\Z=#3%
\divide\Z by 2%
\begin{picture}(0,0)%
\put(\Z,#3){\line(-1,-2){#3}}%
\put(-\Z,-#3){\sswhead}%
\truex{200}\truey{800}\truez{600}%
\put(-\value{x},\value{x}){\makebox(0,\value{z})[r]{${#1}$}}%
\put(\value{x},-\value{y}){\makebox(0,\value{z})[l]{${#2}$}}%
\end{picture}}%

\newcommand{\SSWDIST}[3]{%
\Z=#3%
\divide\Z by 2%
\begin{picture}(0,0)%
\put(\Z,#3){\line(-1,-2){#3}}%
\put(-\Z,-#3){\sswhead}%
\truex{400}%
\put(0,0){\distsign}%
\truex{200}\truey{800}\truez{600}%
\put(-\value{x},\value{x}){\makebox(0,\value{z})[r]{${#1}$}}%
\put(\value{x},-\value{y}){\makebox(0,\value{z})[l]{${#2}$}}%
\end{picture}}%

\newcommand{\SSWDOTAR}[3]{%
\truex{100}\truey{268}\truez{134}%
\Z=#3%
\divide\Z by 2%
\NUMBEROFDOTS=#3%
\divide\NUMBEROFDOTS by \value{z}%
\advance\NUMBEROFDOTS by 1%
\begin{picture}(0,0)%
\multiput(\Z,#3)(-\value{z},-\value{y}){\NUMBEROFDOTS}%
{\circle*{\value{x}}}%
\put(-\Z,-#3){\sswhead}%
\truex{200}\truey{800}\truez{600}%
\put(-\value{x},\value{x}){\makebox(0,\value{z})[r]{${#1}$}}%
\put(\value{x},-\value{y}){\makebox(0,\value{z})[l]{${#2}$}}%
\end{picture}}%

\newcommand{\SSWMONO}[3]{%
\Z=#3%
\divide\Z by 2%
\truetaiL%
\bimolength=#3%
\advance\bimolength by -\truemonotaiL%
\monolength=\bimolength%
\advance\monolength by -\Z%
\secondmonolength=\monolength%
\multiply\secondmonolength by 2%
\begin{picture}(0,0)%
\put(\monolength,\secondmonolength){\line(-1,-2){\bimolength}}%
\put(\monolength,\secondmonolength){\sswhead}%
\put(-\Z,-#3){\sswhead}%
\truex{200}\truey{800}\truez{600}%
\put(-\value{x},\value{x}){\makebox(0,\value{z})[r]{${#1}$}}%
\put(\value{x},-\value{y}){\makebox(0,\value{z})[l]{${#2}$}}%
\end{picture}}%

\newcommand{\SSWEPI}[3]{%
\Z=#3%
\divide\Z by 2%
\trueheaD%
\bimolength=#3%
\advance\bimolength by -\trueepiheaD%
\epilength=\bimolength%
\advance\epilength by -\Z%
\secondepilength=\epilength%
\multiply\secondepilength by 2%
\begin{picture}(0,0)%
\put(\Z,#3){\line(-1,-2){\bimolength}}%
\put(-\epilength,-\secondepilength){\sswhead}%
\put(-\Z,-#3){\sswhead}%
\truex{200}\truey{800}\truez{600}%
\put(-\value{x},\value{x}){\makebox(0,\value{z})[r]{${#1}$}}%
\put(\value{x},-\value{y}){\makebox(0,\value{z})[l]{${#2}$}}%
\end{picture}}%

\newcommand{\SSWBIMO}[3]{%
\Z=#3%
\divide\Z by 2%
\truetaiL\trueheaD%
\bimolength=#3%
\advance\bimolength by -\truemonotaiL%
\monolength=\bimolength%
\advance\monolength by -\Z%
\advance\bimolength by -\trueepiheaD%
\epilength=\bimolength%
\advance\epilength by -\monolength%
\secondmonolength=\monolength%
\multiply\secondmonolength by 2%
\secondepilength=\epilength%
\multiply\secondepilength by 2%
\begin{picture}(0,0)%
\put(\monolength,\secondmonolength){\line(-1,-2){\bimolength}}%
\put(\monolength,\secondmonolength){\sswhead}%
\put(-\epilength,-\secondepilength){\sswhead}%
\put(-\Z,-#3){\sswhead}%
\truex{200}\truey{800}\truez{600}%
\put(-\value{x},\value{x}){\makebox(0,\value{z})[r]{${#1}$}}%
\put(\value{x},-\value{y}){\makebox(0,\value{z})[l]{${#2}$}}%
\end{picture}}%

\newcommand{\SSWBIAR}[3]{%
\Y=#3%
\divide\Y by 2%
\Z=#3%
\multiply \Z by 2%
\begin{picture}(0,0)%
\put(\Y,#3){\begin{picture}(0,0)%
\truex{313}\truey{156}%
\put(-\value{x},\value{y}){\line(-1,-2){#3}}%
\put(\value{x},-\value{y}){\line(-1,-2){#3}}%
\monolength=#3%
\advance\monolength by -\value{x}%
\epilength=#3%
\advance\epilength by \value{x}%
\secondmonolength=\Z%
\advance\secondmonolength by -\value{y}%
\secondepilength=\Z%
\advance\secondepilength by \value{y}%
\put(-\monolength,-\secondepilength){\sswhead}%
\put(-\epilength,-\secondmonolength){\sswhead}%
\end{picture}}
\truex{300}\truey{1000}\truez{600}%
\put(-\value{x},\value{x}){\makebox(0,\value{z})[r]{${#1}$}}%
\put(\value{x},-\value{y}){\makebox(0,\value{z})[l]{${#2}$}}%
\end{picture}}%

\newcommand{\SSWBIDIST}[3]{%
\Y=#3%
\divide\Y by 2%
\Z=#3%
\multiply \Z by 2%
\begin{picture}(0,0)%
\truex{313}\truey{156}\truez{400}%
\put(\Y,#3){\begin{picture}(0,0)%
\put(-\value{x},\value{y}){\line(-1,-2){#3}}%
\put(\value{x},-\value{y}){\line(-1,-2){#3}}%
\monolength=#3%
\advance\monolength by -\value{x}%
\epilength=#3%
\advance\epilength by \value{x}%
\secondmonolength=\Z%
\advance\secondmonolength by -\value{y}%
\secondepilength=\Z%
\advance\secondepilength by \value{y}%
\put(-\monolength,-\secondepilength){\sswhead}%
\put(-\epilength,-\secondmonolength){\sswhead}%
\end{picture}}
\put(-\value{x},\value{y}){\distsign}%
\put(\value{x},-\value{y}){\distsign}%
\truex{300}\truey{1000}\truez{600}%
\put(-\value{x},\value{x}){\makebox(0,\value{z})[r]{${#1}$}}%
\put(\value{x},-\value{y}){\makebox(0,\value{z})[l]{${#2}$}}%
\end{picture}}%

\newcommand{\SSWADJAR}[3]{%
\Y=#3%
\divide\Y by 2%
\Z=#3%
\multiply \Z by 2%
\begin{picture}(0,0)%
\put(\Y,#3){\begin{picture}(0,0)%
\truex{313}\truey{156}%
\monolength=#3%
\advance\monolength by -\value{x}%
\epilength=#3%
\advance\epilength by \value{x}%
\secondmonolength=\Z%
\advance\secondmonolength by -\value{y}%
\secondepilength=\Z%
\advance\secondepilength by \value{y}%
\put(\value{x},-\value{y}){\line(-1,-2){#3}}%
\put(-\monolength,-\secondepilength){\sswhead}%
\put(-\epilength,-\secondmonolength){\line(1,2){#3}}%
\put(-\value{x},\value{y}){\nnehead}%
\end{picture}}
\truex{300}\truey{1000}\truez{600}%
\put(-\value{x},\value{x}){\makebox(0,\value{z})[r]{${#1}$}}%
\put(\value{x},-\value{y}){\makebox(0,\value{z})[l]{${#2}$}}%
\end{picture}}%

\newcommand{\SSWADJDIST}[3]{%
\Y=#3%
\divide\Y by 2%
\Z=#3%
\multiply \Z by 2%
\begin{picture}(0,0)%
\truex{313}\truey{156}\truez{400}%
\put(\Y,#3){\begin{picture}(0,0)%
\monolength=#3%
\advance\monolength by -\value{x}%
\epilength=#3%
\advance\epilength by \value{x}%
\secondmonolength=\Z%
\advance\secondmonolength by -\value{y}%
\secondepilength=\Z%
\advance\secondepilength by \value{y}%
\put(\value{x},-\value{y}){\line(-1,-2){#3}}%
\put(-\monolength,-\secondepilength){\sswhead}%
\put(-\epilength,-\secondmonolength){\line(1,2){#3}}%
\put(-\value{x},\value{y}){\nnehead}%
\end{picture}}
\put(-\value{x},\value{y}){\distsign}%
\put(\value{x},-\value{y}){\distsign}%
\truex{300}\truey{1000}\truez{600}%
\put(-\value{x},\value{x}){\makebox(0,\value{z})[r]{${#1}$}}%
\put(\value{x},-\value{y}){\makebox(0,\value{z})[l]{${#2}$}}%
\end{picture}}%

\def\basicsswar[#1]{\SSWAR{}{}{#100}}%

\newcommand{\sswar}{\@ifnextchar[{\basicsswar}{\basicsswar[67]}}%

\def\basicSswar[#1]#2{\SSWAR{#2}{}{#100}}%

\newcommand{\Sswar}{\@ifnextchar[{\basicSswar}{\basicSswar[67]}}%

\def\basicsswaR[#1]#2{\SSWAR{}{#2}{#100}}%

\newcommand{\sswaR}{\@ifnextchar[{\basicsswaR}{\basicsswaR[67]}}%

\def\basicsswdist[#1]{\SSWDIST{}{}{#100}}%

\newcommand{\sswdist}{\@ifnextchar[{\basicsswdist}{\basicsswdist[67]}}%

\def\basicSswdist[#1]#2{\SSWDIST{#2}{}{#100}}%

\newcommand{\Sswdist}{\@ifnextchar[{\basicSswdist}{\basicSswdist[67]}}%

\def\basicsswdisT[#1]#2{\SSWDIST{}{#2}{#100}}%

\newcommand{\sswdisT}{\@ifnextchar[{\basicsswdisT}{\basicsswdisT[67]}}%

\def\basicsswdotar[#1]{\SSWDOTAR{}{}{#100}}%

\newcommand{\sswdotar}{\@ifnextchar[{\basicsswdotar}{\basicsswdotar[67]}}%

\def\basicSswdotar[#1]#2{\SSWDOTAR{#2}{}{#100}}%

\newcommand{\Sswdotar}{\@ifnextchar[{\basicSswdotar}{\basicSswdotar[67]}}%

\def\basicsswdotaR[#1]#2{\SSWDOTAR{}{#2}{#100}}%

\newcommand{\sswdotaR}{\@ifnextchar[{\basicsswdotaR}{\basicsswdotaR[67]}}%

\def\basicsswmono[#1]{\SSWMONO{}{}{#100}}%

\newcommand{\sswmono}{\@ifnextchar[{\basicsswmono}{\basicsswmono[67]}}%

\def\basicSswmono[#1]#2{\SSWMONO{#2}{}{#100}}%

\newcommand{\Sswmono}{\@ifnextchar[{\basicSswmono}{\basicSswmono[67]}}%

\def\basicsswmonO[#1]#2{\SSWMONO{}{#2}{#100}}%

\newcommand{\sswmonO}{\@ifnextchar[{\basicsswmonO}{\basicsswmonO[67]}}%

\def\basicsswepi[#1]{\SSWEPI{}{}{#100}}%

\newcommand{\sswepi}{\@ifnextchar[{\basicsswepi}{\basicsswepi[67]}}%

\def\basicSswepi[#1]#2{\SSWEPI{#2}{}{#100}}%

\newcommand{\Sswepi}{\@ifnextchar[{\basicSswepi}{\basicSswepi[67]}}%

\def\basicsswepI[#1]#2{\SSWEPI{}{#2}{#100}}%

\newcommand{\sswepI}{\@ifnextchar[{\basicsswepI}{\basicsswepI[67]}}%

\def\basicsswbimo[#1]{\SSWBIMO{}{}{#100}}%

\newcommand{\sswbimo}{\@ifnextchar[{\basicsswbimo}{\basicsswbimo[67]}}%

\def\basicSswbimo[#1]#2{\SSWBIMO{#2}{}{#100}}%

\newcommand{\Sswbimo}{\@ifnextchar[{\basicSswbimo}{\basicSswbimo[67]}}%

\def\basicsswbimO[#1]#2{\SSWBIMO{}{#2}{#100}}%

\newcommand{\sswbimO}{\@ifnextchar[{\basicsswbimO}{\basicsswbimO[67]}}%

\def\basicsswiso[#1]{\SSWAR{\isosign}{}{#100}}%

\newcommand{\sswiso}{\@ifnextchar[{\basicsswiso}{\basicsswiso[67]}}%

\def\basicSswiso[#1]#2{\SSWAR{#2}{\isosign}{#100}}%

\newcommand{\Sswiso}{\@ifnextchar[{\basicSswiso}{\basicSswiso[67]}}%

\def\basicsswisO[#1]#2{\SSWAR{\isosign}{#2}{#100}}%

\newcommand{\sswisO}{\@ifnextchar[{\basicsswisO}{\basicsswisO[67]}}%

\def\basicsswbiar[#1]{\SSWBIAR{}{}{#100}}%

\newcommand{\sswbiar}{\@ifnextchar[{\basicsswbiar}{\basicsswbiar[67]}}%

\def\basicSswbiar[#1]#2#3{\SSWBIAR{#2}{#3}{#100}}%

\newcommand{\Sswbiar}{\@ifnextchar[{\basicSswbiar}{\basicSswbiar[67]}}%

\def\basicsswbidist[#1]{\SSWBIDIST{}{}{#100}}%

\newcommand{\sswbidist}{\@ifnextchar[{\basicsswbidist}{\basicsswbidist[67]}}%

\def\basicSswbidist[#1]#2#3{\SSWBIDIST{#2}{#3}{#100}}%

\newcommand{\Sswbidist}{\@ifnextchar[{\basicSswbidist}{\basicSswbidist[67]}}%

\def\basicsswadjar[#1]{\SSWADJAR{}{}{#100}}%

\newcommand{\sswadjar}{\@ifnextchar[{\basicsswadjar}{\basicsswadjar[67]}}%

\def\basicSswadjar[#1]#2#3{\SSWADJAR{#2}{#3}{#100}}%

\newcommand{\Sswadjar}{\@ifnextchar[{\basicSswadjar}{\basicSswadjar[67]}}%

\def\basicsswadjdist[#1]{\SSWADJDIST{}{}{#100}}%

\newcommand{\sswadjdist}{\@ifnextchar[{\basicsswadjdist}{\basicsswadjdist[67]}}%

\def\basicSswadjdist[#1]#2#3{\SSWADJDIST{#2}{#3}{#100}}%

\newcommand{\Sswadjdist}{\@ifnextchar[{\basicSswadjdist}{\basicSswadjdist[67]}}%

\newcommand{\NNWAR}[3]{%
\Z=#3%
\divide\Z by 2%
\begin{picture}(0,0)%
\put(\Z,-#3){\line(-1,2){#3}}%
\put(-\Z,#3){\nnwhead}%
\truex{200}\truey{800}\truez{600}%
\put(\value{x},\value{x}){\makebox(0,\value{z})[l]{${#1}$}}%
\put(-\value{x},-\value{y}){\makebox(0,\value{z})[r]{${#2}$}}%
\end{picture}}%

\newcommand{\NNWDIST}[3]{%
\Z=#3%
\divide\Z by 2%
\begin{picture}(0,0)%
\put(\Z,-#3){\line(-1,2){#3}}%
\put(-\Z,#3){\nnwhead}%
\truex{400}%
\put(0,0){\distsign}%
\truex{200}\truey{800}\truez{600}%
\put(\value{x},\value{x}){\makebox(0,\value{z})[l]{${#1}$}}%
\put(-\value{x},-\value{y}){\makebox(0,\value{z})[r]{${#2}$}}%
\end{picture}}%

\newcommand{\NNWDOTAR}[3]{%
\truex{100}\truey{268}\truez{134}%
\Z=#3%
\divide\Z by 2%
\NUMBEROFDOTS=#3%
\divide\NUMBEROFDOTS by \value{z}%
\advance\NUMBEROFDOTS by 1%
\begin{picture}(0,0)%
\multiput(\Z,-#3)(-\value{z},\value{y}){\NUMBEROFDOTS}%
{\circle*{\value{x}}}%
\put(-\Z,#3){\nnwhead}%
\truex{200}\truey{800}\truez{600}%
\put(\value{x},\value{x}){\makebox(0,\value{z})[l]{${#1}$}}%
\put(-\value{x},-\value{y}){\makebox(0,\value{z})[r]{${#2}$}}%
\end{picture}}%

\newcommand{\NNWMONO}[3]{%
\Z=#3%
\divide\Z by 2%
\truetaiL%
\bimolength=#3%
\advance\bimolength by -\truemonotaiL%
\monolength=\bimolength%
\advance\monolength by -\Z%
\secondmonolength=\monolength%
\multiply\secondmonolength by 2%
\begin{picture}(0,0)%
\put(\monolength,-\secondmonolength){\line(-1,2){\bimolength}}%
\put(\monolength,-\secondmonolength){\nnwhead}%
\put(-\Z,#3){\nnwhead}%
\truex{200}\truey{800}\truez{600}%
\put(\value{x},\value{x}){\makebox(0,\value{z})[l]{${#1}$}}%
\put(-\value{x},-\value{y}){\makebox(0,\value{z})[r]{${#2}$}}%
\end{picture}}%

\newcommand{\NNWEPI}[3]{%
\Z=#3%
\divide\Z by 2%
\trueheaD%
\bimolength=#3%
\advance\bimolength by -\trueepiheaD%
\epilength=\bimolength%
\advance\epilength by -\Z%
\secondepilength=\epilength%
\multiply\secondepilength by 2%
\begin{picture}(0,0)%
\put(\Z,-#3){\line(-1,2){\bimolength}}%
\put(-\epilength,\secondepilength){\nnwhead}%
\put(-\Z,#3){\nnwhead}%
\truex{200}\truey{800}\truez{600}%
\put(\value{x},\value{x}){\makebox(0,\value{z})[l]{${#1}$}}%
\put(-\value{x},-\value{y}){\makebox(0,\value{z})[r]{${#2}$}}%
\end{picture}}%

\newcommand{\NNWBIMO}[3]{%
\Z=#3%
\divide\Z by 2%
\truetaiL\trueheaD%
\bimolength=#3%
\advance\bimolength by -\truemonotaiL%
\monolength=\bimolength%
\advance\monolength by -\Z%
\advance\bimolength by -\trueepiheaD%
\epilength=\bimolength%
\advance\epilength by -\monolength%
\secondmonolength=\monolength%
\multiply\secondmonolength by 2%
\secondepilength=\epilength%
\multiply\secondepilength by 2%
\begin{picture}(0,0)%
\put(\monolength,-\secondmonolength){\line(-1,2){\bimolength}}%
\put(\monolength,-\secondmonolength){\nnwhead}%
\put(-\epilength,\secondepilength){\nnwhead}%
\put(-\Z,#3){\nnwhead}%
\truex{200}\truey{800}\truez{600}%
\put(\value{x},\value{x}){\makebox(0,\value{z})[l]{${#1}$}}%
\put(-\value{x},-\value{y}){\makebox(0,\value{z})[r]{${#2}$}}%
\end{picture}}%

\newcommand{\NNWBIAR}[3]{%
\Y=#3%
\divide\Y by 2%
\Z=#3%
\multiply \Z by 2%
\begin{picture}(0,0)%
\put(\Y,-#3){\begin{picture}(0,0)%
\truex{313}\truey{156}%
\put(-\value{x},-\value{y}){\line(-1,2){#3}}%
\put(\value{x},\value{y}){\line(-1,2){#3}}%
\monolength=#3%
\advance\monolength by -\value{x}%
\epilength=#3%
\advance\epilength by \value{x}%
\secondmonolength=\Z%
\advance\secondmonolength by -\value{y}%
\secondepilength=\Z%
\advance\secondepilength by \value{y}%
\put(-\monolength,\secondepilength){\nnwhead}%
\put(-\epilength,\secondmonolength){\nnwhead}%
\end{picture}}
\truex{400}\truey{1000}\truez{600}%
\put(\value{x},\value{x}){\makebox(0,\value{z})[l]{${#1}$}}%
\put(-\value{x},-\value{y}){\makebox(0,\value{z})[r]{${#2}$}}%
\end{picture}}%

\newcommand{\NNWBIDIST}[3]{%
\Y=#3%
\divide\Y by 2%
\Z=#3%
\multiply \Z by 2%
\begin{picture}(0,0)%
\truex{313}\truey{156}\truez{400}%
\put(\Y,-#3){\begin{picture}(0,0)%
\put(-\value{x},-\value{y}){\line(-1,2){#3}}%
\put(\value{x},\value{y}){\line(-1,2){#3}}%
\monolength=#3%
\advance\monolength by -\value{x}%
\epilength=#3%
\advance\epilength by \value{x}%
\secondmonolength=\Z%
\advance\secondmonolength by -\value{y}%
\secondepilength=\Z%
\advance\secondepilength by \value{y}%
\put(-\monolength,\secondepilength){\nnwhead}%
\put(-\epilength,\secondmonolength){\nnwhead}%
\end{picture}}
\put(-\value{x},-\value{y}){\distsign}%
\put(\value{x},\value{y}){\distsign}%
\truex{500}\truey{1000}\truez{600}%
\put(\value{x},\value{x}){\makebox(0,\value{z})[l]{${#1}$}}%
\put(-\value{x},-\value{y}){\makebox(0,\value{z})[r]{${#2}$}}%
\end{picture}}%

\newcommand{\NNWADJAR}[3]{%
\Y=#3%
\divide\Y by 2%
\Z=#3%
\multiply \Z by 2%
\begin{picture}(0,0)%
\put(\Y,-#3){\begin{picture}(0,0)%
\truex{313}\truey{156}%
\monolength=#3%
\advance\monolength by -\value{x}%
\epilength=#3%
\advance\epilength by \value{x}%
\secondmonolength=\Z%
\advance\secondmonolength by -\value{y}%
\secondepilength=\Z%
\advance\secondepilength by \value{y}%
\put(-\value{x},-\value{y}){\line(-1,2){#3}}%
\put(-\epilength,\secondmonolength){\nnwhead}%
\put(-\monolength,\secondepilength){\line(1,-2){#3}}%
\put(\value{x},\value{y}){\ssehead}%
\end{picture}}
\truex{400}\truey{1000}\truez{600}%
\put(\value{x},\value{x}){\makebox(0,\value{z})[l]{${#1}$}}%
\put(-\value{x},-\value{y}){\makebox(0,\value{z})[r]{${#2}$}}%
\end{picture}}%

\newcommand{\NNWADJDIST}[3]{%
\Y=#3%
\divide\Y by 2%
\Z=#3%
\multiply \Z by 2%
\begin{picture}(0,0)%
\truex{313}\truey{156}\truez{400}%
\put(\Y,-#3){\begin{picture}(0,0)%
\monolength=#3%
\advance\monolength by -\value{x}%
\epilength=#3%
\advance\epilength by \value{x}%
\secondmonolength=\Z%
\advance\secondmonolength by -\value{y}%
\secondepilength=\Z%
\advance\secondepilength by \value{y}%
\put(-\value{x},-\value{y}){\line(-1,2){#3}}%
\put(-\epilength,\secondmonolength){\nnwhead}%
\put(-\monolength,\secondepilength){\line(1,-2){#3}}%
\put(\value{x},\value{y}){\ssehead}%
\end{picture}}
\put(\value{x},\value{y}){\distsign}%
\put(-\value{x},-\value{y}){\distsign}%
\truex{500}\truey{1000}\truez{600}%
\put(\value{x},\value{x}){\makebox(0,\value{z})[l]{${#1}$}}%
\put(-\value{x},-\value{y}){\makebox(0,\value{z})[r]{${#2}$}}%
\end{picture}}%

\def\basicnnwar[#1]{\NNWAR{}{}{#100}}%

\newcommand{\nnwar}{\@ifnextchar[{\basicnnwar}{\basicnnwar[67]}}%

\def\basicNnwar[#1]#2{\NNWAR{#2}{}{#100}}%

\newcommand{\Nnwar}{\@ifnextchar[{\basicNnwar}{\basicNnwar[67]}}%

\def\basicnnwaR[#1]#2{\NNWAR{}{#2}{#100}}%

\newcommand{\nnwaR}{\@ifnextchar[{\basicnnwaR}{\basicnnwaR[67]}}%

\def\basicnnwdist[#1]{\NNWDIST{}{}{#100}}%

\newcommand{\nnwdist}{\@ifnextchar[{\basicnnwdist}{\basicnnwdist[67]}}%

\def\basicNnwdist[#1]#2{\NNWDIST{#2}{}{#100}}%

\newcommand{\Nnwdist}{\@ifnextchar[{\basicNnwdist}{\basicNnwdist[67]}}%

\def\basicnnwdisT[#1]#2{\NNWDIST{}{#2}{#100}}%

\newcommand{\nnwdisT}{\@ifnextchar[{\basicnnwdisT}{\basicnnwdisT[67]}}%

\def\basicnnwdotar[#1]{\NNWDOTAR{}{}{#100}}%

\newcommand{\nnwdotar}{\@ifnextchar[{\basicnnwdotar}{\basicnnwdotar[67]}}%

\def\basicNnwdotar[#1]#2{\NNWDOTAR{#2}{}{#100}}%

\newcommand{\Nnwdotar}{\@ifnextchar[{\basicNnwdotar}{\basicNnwdotar[67]}}%

\def\basicnnwdotaR[#1]#2{\NNWDOTAR{}{#2}{#100}}%

\newcommand{\nnwdotaR}{\@ifnextchar[{\basicnnwdotaR}{\basicnnwdotaR[67]}}%

\def\basicnnwmono[#1]{\NNWMONO{}{}{#100}}%

\newcommand{\nnwmono}{\@ifnextchar[{\basicnnwmono}{\basicnnwmono[67]}}%

\def\basicNnwmono[#1]#2{\NNWMONO{#2}{}{#100}}%

\newcommand{\Nnwmono}{\@ifnextchar[{\basicNnwmono}{\basicNnwmono[67]}}%

\def\basicnnwmonO[#1]#2{\NNWMONO{}{#2}{#100}}%

\newcommand{\nnwmonO}{\@ifnextchar[{\basicnnwmonO}{\basicnnwmonO[67]}}%

\def\basicnnwepi[#1]{\NNWEPI{}{}{#100}}%

\newcommand{\nnwepi}{\@ifnextchar[{\basicnnwepi}{\basicnnwepi[67]}}%

\def\basicNnwepi[#1]#2{\NNWEPI{#2}{}{#100}}%

\newcommand{\Nnwepi}{\@ifnextchar[{\basicNnwepi}{\basicNnwepi[67]}}%

\def\basicnnwepI[#1]#2{\NNWEPI{}{#2}{#100}}%

\newcommand{\nnwepI}{\@ifnextchar[{\basicnnwepI}{\basicnnwepI[67]}}%

\def\basicnnwbimo[#1]{\NNWBIMO{}{}{#100}}%

\newcommand{\nnwbimo}{\@ifnextchar[{\basicnnwbimo}{\basicnnwbimo[67]}}%

\def\basicNnwbimo[#1]#2{\NNWBIMO{#2}{}{#100}}%

\newcommand{\Nnwbimo}{\@ifnextchar[{\basicNnwbimo}{\basicNnwbimo[67]}}%

\def\basicnnwbimO[#1]#2{\NNWBIMO{}{#2}{#100}}%

\newcommand{\nnwbimO}{\@ifnextchar[{\basicnnwbimO}{\basicnnwbimO[67]}}%

\def\basicnnwiso[#1]{\NNWAR{\isosign}{}{#100}}%

\newcommand{\nnwiso}{\@ifnextchar[{\basicnnwiso}{\basicnnwiso[67]}}%

\def\basicNnwiso[#1]#2{\NNWAR{#2}{\isosign}{#100}}%

\newcommand{\Nnwiso}{\@ifnextchar[{\basicNnwiso}{\basicNnwiso[67]}}%

\def\basicnnwisO[#1]#2{\NNWAR{\isosign}{#2}{#100}}%

\newcommand{\nnwisO}{\@ifnextchar[{\basicnnwisO}{\basicnnwisO[67]}}%

\def\basicnnwbiar[#1]{\NNWBIAR{}{}{#100}}%

\newcommand{\nnwbiar}{\@ifnextchar[{\basicnnwbiar}{\basicnnwbiar[67]}}%

\def\basicNnwbiar[#1]#2#3{\NNWBIAR{#2}{#3}{#100}}%

\newcommand{\Nnwbiar}{\@ifnextchar[{\basicNnwbiar}{\basicNnwbiar[67]}}%

\def\basicnnwbidist[#1]{\NNWBIDIST{}{}{#100}}%

\newcommand{\nnwbidist}{\@ifnextchar[{\basicnnwbidist}{\basicnnwbidist[67]}}%

\def\basicNnwbidist[#1]#2#3{\NNWBIDIST{#2}{#3}{#100}}%

\newcommand{\Nnwbidist}{\@ifnextchar[{\basicNnwbidist}{\basicNnwbidist[67]}}%

\def\basicnnwadjar[#1]{\NNWADJAR{}{}{#100}}%

\newcommand{\nnwadjar}{\@ifnextchar[{\basicnnwadjar}{\basicnnwadjar[67]}}%

\def\basicNnwadjar[#1]#2#3{\NNWADJAR{#2}{#3}{#100}}%

\newcommand{\Nnwadjar}{\@ifnextchar[{\basicNnwadjar}{\basicNnwadjar[67]}}%

\def\basicnnwadjdist[#1]{\NNWADJDIST{}{}{#100}}%

\newcommand{\nnwadjdist}{\@ifnextchar[{\basicnnwadjdist}{\basicnnwadjdist[67]}}%

\def\basicNnwadjdist[#1]#2#3{\NNWADJDIST{#2}{#3}{#100}}%

\newcommand{\Nnwadjdist}{\@ifnextchar[{\basicNnwadjdist}{\basicNnwadjdist[67]}}%

\newcommand{\EENEAR}[3]{%
\Y=#3%
\divide \Y by 2%
\Z=\Y%
\divide \Z by 3%
\begin{picture}(0,0)%
\put(-\Y,-\Z){\line(3,1){#3}}%
\put(\Y,\Z){\eenehead}%
\truex{200}\truey{800}\truez{600}%
\put(-\value{x},\value{x}){\makebox(0,\value{z})[r]{${#1}$}}%
\put(\value{x},-\value{y}){\makebox(0,\value{z})[l]{${#2}$}}%
\end{picture}}%

\def\basiceenear[#1]{\EENEAR{}{}{#100}}%

\newcommand{\eenear}{\@ifnextchar[{\basiceenear}{\basiceenear[211]}}%

\def\basicEenear[#1]#2{\EENEAR{#2}{}{#100}}%

\newcommand{\Eenear}{\@ifnextchar[{\basicEenear}{\basicEenear[211]}}%

\def\basiceeneaR[#1]#2{\EENEAR{}{#2}{#100}}%

\newcommand{\eeneaR}{\@ifnextchar[{\basiceeneaR}{\basiceeneaR[211]}}%

\newcommand{\EESEAR}[3]{%
\Y=#3%
\divide \Y by 2%
\Z=\Y%
\divide \Z by 3%
\begin{picture}(0,0)%
\put(-\Y,\Z){\line(3,-1){#3}}%
\put(\Y,-\Z){\eesehead}%
\truex{200}\truey{800}\truez{600}%
\put(\value{x},\value{x}){\makebox(0,\value{z})[l]{${#1}$}}%
\put(-\value{x},-\value{y}){\makebox(0,\value{z})[r]{${#2}$}}%
\end{picture}}%

\def\basiceesear[#1]{\EESEAR{}{}{#100}}%

\newcommand{\eesear}{\@ifnextchar[{\basiceesear}{\basiceesear[211]}}%

\def\basicEesear[#1]#2{\EESEAR{#2}{}{#100}}%

\newcommand{\Eesear}{\@ifnextchar[{\basicEesear}{\basicEesear[211]}}%

\def\basiceeseaR[#1]#2{\EESEAR{}{#2}{#100}}%

\newcommand{\eeseaR}{\@ifnextchar[{\basiceeseaR}{\basiceeseaR[211]}}%

\newcommand{\WWNWAR}[3]{%
\Y=#3%
\divide \Y by 2%
\Z=\Y%
\divide \Z by 3%
\begin{picture}(0,0)%
\put(\Y,-\Z){\line(-3,1){#3}}%
\put(-\Y,\Z){\wwnwhead}%
\truex{200}\truey{800}\truez{600}%
\put(\value{x},\value{x}){\makebox(0,\value{z})[l]{${#1}$}}%
\put(-\value{x},-\value{y}){\makebox(0,\value{z})[r]{${#2}$}}%
\end{picture}}%

\def\basicwwnwar[#1]{\WWNWAR{}{}{#100}}%

\newcommand{\wwnwar}{\@ifnextchar[{\basicwwnwar}{\basicwwnwar[211]}}%

\def\basicWwnwar[#1]#2{\WWNWAR{#2}{}{#100}}%

\newcommand{\Wwnwar}{\@ifnextchar[{\basicWwnwar}{\basicWwnwar[211]}}%

\def\basicwwnwaR[#1]#2{\WWNWAR{}{#2}{#100}}%

\newcommand{\wwnwaR}{\@ifnextchar[{\basicwwnwaR}{\basicwwnwaR[211]}}%

\newcommand{\WWSWAR}[3]{%
\Y=#3%
\divide \Y by 2%
\Z=\Y%
\divide \Z by 3%
\begin{picture}(0,0)%
\put(\Y,\Z){\line(-3,-1){#3}}%
\put(-\Y,-\Z){\wwswhead}%
\truex{200}\truey{800}\truez{600}%
\put(-\value{x},\value{x}){\makebox(0,\value{z})[r]{${#1}$}}%
\put(\value{x},-\value{y}){\makebox(0,\value{z})[l]{${#2}$}}%
\end{picture}}%

\def\basicwwswar[#1]{\WWSWAR{}{}{#100}}%

\newcommand{\wwswar}{\@ifnextchar[{\basicwwswar}{\basicwwswar[211]}}%

\def\basicWwswar[#1]#2{\WWSWAR{#2}{}{#100}}%

\newcommand{\Wwswar}{\@ifnextchar[{\basicWwswar}{\basicWwswar[211]}}%

\def\basicwwswaR[#1]#2{\WWSWAR{}{#2}{#100}}%

\newcommand{\wwswaR}{\@ifnextchar[{\basicwwswaR}{\basicwwswaR[211]}}%

\newcommand{\NNNEAR}[3]{%
\Y=#3%
\divide \Y by 2%
\Z=\Y%
\multiply \Z by 3%
\begin{picture}(0,0)%
\put(-\Y,-\Z){\line(1,3){#3}}%
\put(\Y,\Z){\nnnehead}%
\truex{100}\truez{600}%
\put(-\value{x},\value{x}){\makebox(0,\value{z})[r]{${#1}$}}%
\put(\value{x},-\value{z}){\makebox(0,\value{z})[l]{${#2}$}}%
\end{picture}}%

\def\basicnnnear[#1]{\NNNEAR{}{}{#100}}%

\newcommand{\nnnear}{\@ifnextchar[{\basicnnnear}{\basicnnnear[71]}}%

\def\basicNnnear[#1]#2{\NNNEAR{#2}{}{#100}}%

\newcommand{\Nnnear}{\@ifnextchar[{\basicNnnear}{\basicNnnear[71]}}%

\def\basicnnneaR[#1]#2{\NNNEAR{}{#2}{#100}}%

\newcommand{\nnneaR}{\@ifnextchar[{\basicnnneaR}{\basicnnneaR[71]}}%

\newcommand{\SSSWAR}[3]{%
\Y=#3%
\divide \Y by 2%
\Z=\Y%
\multiply \Z by 3%
\begin{picture}(0,0)%
\put(\Y,\Z){\line(-1,-3){#3}}%
\put(-\Y,-\Z){\ssswhead}%
\truex{100}\truez{600}%
\put(-\value{x},\value{x}){\makebox(0,\value{z})[r]{${#1}$}}%
\put(\value{x},-\value{z}){\makebox(0,\value{z})[l]{${#2}$}}%
\end{picture}}%

\def\basicssswar[#1]{\SSSWAR{}{}{#100}}%

\newcommand{\ssswar}{\@ifnextchar[{\basicssswar}{\basicssswar[71]}}%

\def\basicSsswar[#1]#2{\SSSWAR{#2}{}{#100}}%

\newcommand{\Ssswar}{\@ifnextchar[{\basicSsswar}{\basicSsswar[71]}}%

\def\basicssswaR[#1]#2{\SSSWAR{}{#2}{#100}}%

\newcommand{\ssswaR}{\@ifnextchar[{\basicssswaR}{\basicssswaR[71]}}%

\newcommand{\SSSEAR}[3]{%
\Y=#3%
\divide \Y by 2%
\Z=\Y%
\multiply \Z by 3%
\begin{picture}(0,0)%
\put(-\Y,\Z){\line(1,-3){#3}}%
\put(\Y,-\Z){\sssehead}%
\truex{200}\truez{600}%
\put(\value{x},\value{x}){\makebox(0,\value{z})[l]{${#1}$}}%
\put(-\value{x},-\value{z}){\makebox(0,\value{z})[r]{${#2}$}}%
\end{picture}}%

\def\basicsssear[#1]{\SSSEAR{}{}{#100}}%

\newcommand{\sssear}{\@ifnextchar[{\basicsssear}{\basicsssear[71]}}%

\def\basicSssear[#1]#2{\SSSEAR{#2}{}{#100}}%

\newcommand{\Sssear}{\@ifnextchar[{\basicSssear}{\basicSssear[71]}}%

\def\basicssseaR[#1]#2{\SSSEAR{}{#2}{#100}}%

\newcommand{\ssseaR}{\@ifnextchar[{\basicssseaR}{\basicssseaR[71]}}%

\newcommand{\NNNWAR}[3]{%
\Y=#3%
\divide \Y by 2%
\Z=\Y%
\multiply \Z by 3%
\begin{picture}(0,0)%
\put(\Y,-\Z){\line(-1,3){#3}}%
\put(-\Y,\Z){\nnnwhead}%
\truex{200}\truez{600}%
\put(\value{x},\value{x}){\makebox(0,\value{z})[l]{${#1}$}}%
\put(-\value{x},-\value{z}){\makebox(0,\value{z})[r]{${#2}$}}%
\end{picture}}%

\def\basicnnnwar[#1]{\NNNWAR{}{}{#100}}%

\newcommand{\nnnwar}{\@ifnextchar[{\basicnnnwar}{\basicnnnwar[71]}}%

\def\basicNnnwar[#1]#2{\NNNWAR{#2}{}{#100}}%

\newcommand{\Nnnwar}{\@ifnextchar[{\basicNnnwar}{\basicNnnwar[71]}}%

\def\basicnnnwaR[#1]#2{\NNNWAR{}{#2}{#100}}%

\newcommand{\nnnwaR}{\@ifnextchar[{\basicnnnwaR}{\basicnnnwaR[71]}}%

\newcommand{\NEENEAR}[3]{%
\Y=#3%
\divide \Y by 2%
\Z=#3%
\divide \Z by 3%
\begin{picture}(0,0)%
\put(-\Y,-\Z){\line(3,2){#3}}%
\put(\Y,\Z){\neenehead}%
\truex{200}\truey{800}\truez{600}%
\put(-\value{x},\value{x}){\makebox(0,\value{z})[r]{${#1}$}}%
\put(\value{x},-\value{y}){\makebox(0,\value{z})[l]{${#2}$}}%
\end{picture}}%

\def\basicneenear[#1]{\NEENEAR{}{}{#100}}%

\newcommand{\neenear}{\@ifnextchar[{\basicneenear}{\basicneenear[215]}}%

\def\basicNeenear[#1]#2{\NEENEAR{#2}{}{#100}}%

\newcommand{\Neenear}{\@ifnextchar[{\basicNeenear}{\basicNeenear[215]}}%

\def\basicneeneaR[#1]#2{\NEENEAR{}{#2}{#100}}%

\newcommand{\neeneaR}{\@ifnextchar[{\basicneeneaR}{\basicneeneaR[215]}}%

\newcommand{\SEESEAR}[3]{%
\Y=#3%
\divide \Y by 2%
\Z=#3%
\divide \Z by 3%
\begin{picture}(0,0)%
\put(-\Y,\Z){\line(3,-2){#3}}%
\put(\Y,-\Z){\seesehead}%
\truex{200}\truey{800}\truez{600}%
\put(\value{x},\value{x}){\makebox(0,\value{z})[l]{${#1}$}}%
\put(-\value{x},-\value{y}){\makebox(0,\value{z})[r]{${#2}$}}%
\end{picture}}%

\def\basicseesear[#1]{\SEESEAR{}{}{#100}}%

\newcommand{\seesear}{\@ifnextchar[{\basicseesear}{\basicseesear[215]}}%

\def\basicSeesear[#1]#2{\SEESEAR{#2}{}{#100}}%

\newcommand{\Seesear}{\@ifnextchar[{\basicSeesear}{\basicSeesear[215]}}%

\def\basicseeseaR[#1]#2{\SEESEAR{}{#2}{#100}}%

\newcommand{\seeseaR}{\@ifnextchar[{\basicseeseaR}{\basicseeseaR[215]}}%

\newcommand{\NWWNWAR}[3]{%
\Y=#3%
\divide \Y by 2%
\Z=#3%
\divide \Z by 3%
\begin{picture}(0,0)%
\put(\Y,-\Z){\line(-3,2){#3}}%
\put(-\Y,\Z){\nwwnwhead}%
\truex{200}\truey{800}\truez{600}%
\put(\value{x},\value{x}){\makebox(0,\value{z})[l]{${#1}$}}%
\put(-\value{x},-\value{y}){\makebox(0,\value{z})[r]{${#2}$}}%
\end{picture}}%

\def\basicnwwnwar[#1]{\NWWNWAR{}{}{#100}}%

\newcommand{\nwwnwar}{\@ifnextchar[{\basicnwwnwar}{\basicnwwnwar[215]}}%

\def\basicNwwnwar[#1]#2{\NWWNWAR{#2}{}{#100}}%

\newcommand{\Nwwnwar}{\@ifnextchar[{\basicNwwnwar}{\basicNwwnwar[215]}}%

\def\basicnwwnwaR[#1]#2{\NWWNWAR{}{#2}{#100}}%

\newcommand{\nwwnwaR}{\@ifnextchar[{\basicnwwnwaR}{\basicnwwnwaR[215]}}%

\newcommand{\SWWSWAR}[3]{%
\Y=#3%
\divide \Y by 2%
\Z=#3%
\divide \Z by 3%
\begin{picture}(0,0)%
\put(\Y,\Z){\line(-3,-2){#3}}%
\put(-\Y,-\Z){\swwswhead}%
\truex{200}\truey{800}\truez{600}%
\put(-\value{x},\value{x}){\makebox(0,\value{z})[r]{${#1}$}}%
\put(\value{x},-\value{y}){\makebox(0,\value{z})[l]{${#2}$}}%
\end{picture}}%

\def\basicswwswar[#1]{\SWWSWAR{}{}{#100}}%

\newcommand{\swwswar}{\@ifnextchar[{\basicswwswar}{\basicswwswar[215]}}%

\def\basicSwwswar[#1]#2{\SWWSWAR{#2}{}{#100}}%

\newcommand{\Swwswar}{\@ifnextchar[{\basicSwwswar}{\basicSwwswar[215]}}%

\def\basicswwswaR[#1]#2{\SWWSWAR{}{#2}{#100}}%

\newcommand{\swwswaR}{\@ifnextchar[{\basicswwswaR}{\basicswwswaR[215]}}%

\newcommand{\NENNEAR}[3]{%
\Y=#3%
\divide \Y by 2%
\Z=#3%
\multiply \Z by 3%
\divide \Z by 4%
\begin{picture}(0,0)%
\put(-\Y,-\Z){\line(2,3){#3}}%
\put(\Y,\Z){\nennehead}%
\truex{100}\truez{600}%
\put(-\value{x},\value{x}){\makebox(0,\value{z})[r]{${#1}$}}%
\put(\value{x},-\value{z}){\makebox(0,\value{z})[l]{${#2}$}}%
\end{picture}}%

\def\basicnennear[#1]{\NENNEAR{}{}{#100}}%

\newcommand{\nennear}{\@ifnextchar[{\basicnennear}{\basicnennear[143]}}%

\def\basicNennear[#1]#2{\NENNEAR{#2}{}{#100}}%

\newcommand{\Nennear}{\@ifnextchar[{\basicNennear}{\basicNennear[143]}}%

\def\basicnenneaR[#1]#2{\NENNEAR{}{#2}{#100}}%

\newcommand{\nenneaR}{\@ifnextchar[{\basicnenneaR}{\basicnenneaR[143]}}%

\newcommand{\SWSSWAR}[3]{%
\Y=#3%
\divide \Y by 2%
\Z=#3%
\multiply \Z by 3%
\divide \Z by 4%
\begin{picture}(0,0)%
\put(\Y,\Z){\line(-2,-3){#3}}%
\put(-\Y,-\Z){\swsswhead}%
\truex{100}\truez{600}%
\put(-\value{x},\value{x}){\makebox(0,\value{z})[r]{${#1}$}}%
\put(\value{x},-\value{z}){\makebox(0,\value{z})[l]{${#2}$}}%
\end{picture}}%

\def\basicswsswar[#1]{\SWSSWAR{}{}{#100}}%

\newcommand{\swsswar}{\@ifnextchar[{\basicswsswar}{\basicswsswar[143]}}%

\def\basicSwsswar[#1]#2{\SWSSWAR{#2}{}{#100}}%

\newcommand{\Swsswar}{\@ifnextchar[{\basicSwsswar}{\basicSwsswar[143]}}%

\def\basicswsswaR[#1]#2{\SWSSWAR{}{#2}{#100}}%

\newcommand{\swsswaR}{\@ifnextchar[{\basicswsswaR}{\basicswsswaR[143]}}%

\newcommand{\SESSEAR}[3]{%
\Y=#3%
\divide \Y by 2%
\Z=#3%
\multiply \Z by 3%
\divide \Z by 4%
\begin{picture}(0,0)%
\put(-\Y,\Z){\line(2,-3){#3}}%
\put(\Y,-\Z){\sessehead}%
\truex{200}\truez{600}%
\put(\value{x},\value{x}){\makebox(0,\value{z})[l]{${#1}$}}%
\put(-\value{x},-\value{z}){\makebox(0,\value{z})[r]{${#2}$}}%
\end{picture}}%

\def\basicsessear[#1]{\SESSEAR{}{}{#100}}%

\newcommand{\sessear}{\@ifnextchar[{\basicsessear}{\basicsessear[143]}}%

\def\basicSessear[#1]#2{\SESSEAR{#2}{}{#100}}%

\newcommand{\Sessear}{\@ifnextchar[{\basicSessear}{\basicSessear[143]}}%

\def\basicsesseaR[#1]#2{\SESSEAR{}{#2}{#100}}%

\newcommand{\sesseaR}{\@ifnextchar[{\basicsesseaR}{\basicsesseaR[143]}}%

\newcommand{\NWNNWAR}[3]{%
\Y=#3%
\divide \Y by 2%
\Z=#3%
\multiply \Z by 3%
\divide \Z by 4%
\begin{picture}(0,0)%
\put(\Y,-\Z){\line(-2,3){#3}}%
\put(-\Y,\Z){\nwnnwhead}%
\truex{200}\truez{600}%
\put(\value{x},\value{x}){\makebox(0,\value{z})[l]{${#1}$}}%
\put(-\value{x},-\value{z}){\makebox(0,\value{z})[r]{${#2}$}}%
\end{picture}}%

\def\basicnwnnwar[#1]{\NWNNWAR{}{}{#100}}%

\newcommand{\nwnnwar}{\@ifnextchar[{\basicnwnnwar}{\basicnwnnwar[143]}}%

\def\basicNwnnwar[#1]#2{\NWNNWAR{#2}{}{#100}}%

\newcommand{\Nwnnwar}{\@ifnextchar[{\basicNwnnwar}{\basicNwnnwar[143]}}%

\def\basicnwnnwaR[#1]#2{\NWNNWAR{}{#2}{#100}}%

\newcommand{\nwnnwaR}{\@ifnextchar[{\basicnwnnwaR}{\basicnwnnwaR[143]}}%

\newcommand{\Necurve}[2]%
{\begin{picture}(0,0)%
\truex{1300}\truey{2000}\truez{200}%
\put(0,\value{x}){\oval(#200,\value{y})[t]}%
\put(0,\value{x}){\makebox(0,0){\begin{picture}(#200,0)%
\put(#200,0){\line(0,-1){\value{z}}}%
\put(#200,-\value{z}){\shead}%
\put(0,0){\line(0,-1){\value{z}}}\end{picture}}}%
\truex{2500}%
\put(0,\value{x}){\makebox(0,0)[b]{${#1}$}}%
\end{picture}}%

\def\basicnecurvar[#1]{\Necurve{}{#1}}

\newcommand{\necurvar}{\@ifnextchar[{\basicnecurvar}{\basicnecurvar[160]}}%

\def\basicNecurvar[#1]#2{\Necurve{#2}{#1}}%

\newcommand{\Necurvar}{\@ifnextchar[{\basicNecurvar}{\basicNecurvar[160]}}%

\newcommand{\Nwcurve}[2]%
{\begin{picture}(0,0)%
\truex{1300}\truey{2000}\truez{200}%
\put(0,\value{x}){\oval(#200,\value{y})[t]}%
\put(0,\value{x}){\makebox(0,0){\begin{picture}(#200,0)%
\put(#200,0){\line(0,-1){\value{z}}}%
\put(0,0){\line(0,-1){\value{z}}}%
\put(0,-\value{z}){\shead}%
\end{picture}}}%
\truex{2500}%
\put(0,\value{x}){\makebox(0,0)[b]{${#1}$}}%
\end{picture}}%

\def\basicnwcurvar[#1]{\Nwcurve{}{#1}}

\newcommand{\nwcurvar}{\@ifnextchar[{\basicnwcurvar}{\basicnwcurvar[160]}}%

\def\basicNwcurvar[#1]#2{\Nwcurve{#2}{#1}}%

\newcommand{\Nwcurvar}{\@ifnextchar[{\basicNwcurvar}{\basicNwcurvar[160]}}%

\newcommand{\Securve}[2]%
{\begin{picture}(0,0)%
\truex{1300}\truey{2000}\truez{200}%
\put(0,-\value{x}){\oval(#200,\value{y})[b]}%
\put(0,-\value{x}){\makebox(0,0){\begin{picture}(#200,0)%
\put(#200,0){\line(0,1){\value{z}}}%
\put(0,0){\line(0,1){\value{z}}}%
\put(#200,\value{z}){\nhead}%
\end{picture}}}%
\truex{2500}%
\put(0,-\value{x}){\makebox(0,0)[t]{${#1}$}}%
\end{picture}}%

\def\basicsecurvar[#1]{\Securve{}{#1}}

\newcommand{\securvar}{\@ifnextchar[{\basicsecurvar}{\basicsecurvar[160]}}%

\def\basicSecurvar[#1]#2{\Securve{#2}{#1}}%

\newcommand{\Securvar}{\@ifnextchar[{\basicSecurvar}{\basicSecurvar[160]}}%

\newcommand{\Swcurve}[2]%
{\begin{picture}(0,0)%
\truex{1300}\truey{2000}\truez{200}%
\put(0,-\value{x}){\oval(#200,\value{y})[b]}%
\put(0,-\value{x}){\makebox(0,0){\begin{picture}(#200,0)%
\put(#200,0){\line(0,1){\value{z}}}%
\put(0,0){\line(0,1){\value{z}}}%
\put(0,\value{z}){\nhead}%
\end{picture}}}%
\truex{2500}%
\put(0,-\value{x}){\makebox(0,0)[t]{${#1}$}}%
\end{picture}}%

\def\basicswcurvar[#1]{\Swcurve{}{#1}}

\newcommand{\swcurvar}{\@ifnextchar[{\basicswcurvar}{\basicswcurvar[160]}}%

\def\basicSwcurvar[#1]#2{\Swcurve{#2}{#1}}%

\newcommand{\Swcurvar}{\@ifnextchar[{\basicSwcurvar}{\basicSwcurvar[160]}}%

\newcommand{\Escurve}[2]%
{\begin{picture}(0,0)%
\truex{1400}\truey{2000}\truez{200}%
\put(\value{x},0){\oval(\value{y},#200)[r]}%
\put(\value{x},0){\makebox(0,0){\begin{picture}(0,#200)%
\put(0,0){\line(-1,0){\value{z}}}%
\put(0,#200){\line(-1,0){\value{z}}}%
\put(-\value{z},0){\whead}%
\end{picture}}}%
\truex{2500}%
\put(\value{x},0){\makebox(0,0)[l]{${#1}$}}%
\end{picture}}%

\def\basicescurvar[#1]{\Escurve{}{#1}}

\newcommand{\escurvar}{\@ifnextchar[{\basicescurvar}{\basicescurvar[160]}}%

\def\basicEscurvar[#1]#2{\Escurve{#2}{#1}}%

\newcommand{\Escurvar}{\@ifnextchar[{\basicEscurvar}{\basicEscurvar[160]}}%

\newcommand{\Encurve}[2]%
{\begin{picture}(0,0)%
\truex{1400}\truey{2000}\truez{200}%
\put(\value{x},0){\oval(\value{y},#200)[r]}%
\put(\value{x},0){\makebox(0,0){\begin{picture}(0,#200)%
\put(0,0){\line(-1,0){\value{z}}}%
\put(0,#200){\line(-1,0){\value{z}}}%
\put(-\value{z},#200){\whead}%
\end{picture}}}%
\truex{2500}%
\put(\value{x},0){\makebox(0,0)[l]{${#1}$}}%
\end{picture}}%

\def\basicencurvar[#1]{\Encurve{}{#1}}

\newcommand{\encurvar}{\@ifnextchar[{\basicencurvar}{\basicencurvar[160]}}%

\def\basicEncurvar[#1]#2{\Encurve{#2}{#1}}%

\newcommand{\Encurvar}{\@ifnextchar[{\basicEncurvar}{\basicEncurvar[160]}}%

\newcommand{\Wscurve}[2]%
{\begin{picture}(0,0)%
\truex{1300}\truey{2000}\truez{200}%
\put(-\value{x},0){\oval(\value{y},#200)[l]}%
\put(-\value{x},0){\makebox(0,0){\begin{picture}(0,#200)%
\put(0,0){\line(1,0){\value{z}}}%
\put(0,#200){\line(1,0){\value{z}}}%
\put(\value{z},0){\ehead}%
\end{picture}}}%
\truex{2400}%
\put(-\value{x},0){\makebox(0,0)[r]{${#1}$}}%
\end{picture}}%

\def\basicwscurvar[#1]{\Wscurve{}{#1}}

\newcommand{\wscurvar}{\@ifnextchar[{\basicwscurvar}{\basicwscurvar[160]}}%

\def\basicWscurvar[#1]#2{\Wscurve{#2}{#1}}%

\newcommand{\Wscurvar}{\@ifnextchar[{\basicWscurvar}{\basicWscurvar[160]}}%

\newcommand{\Wncurve}[2]%
{\begin{picture}(0,0)%
\truex{1300}\truey{2000}\truez{200}%
\put(-\value{x},0){\oval(\value{y},#200)[l]}%
\put(-\value{x},0){\makebox(0,0){\begin{picture}(0,#200)%
\put(0,0){\line(1,0){\value{z}}}%
\put(\value{z},#200){\ehead}%
\put(0,#200){\line(1,0){\value{z}}}%
\end{picture}}}%
\truex{2400}%
\put(-\value{x},0){\makebox(0,0)[r]{${#1}$}}%
\end{picture}}%

\def\basicwncurvar[#1]{\Wncurve{}{#1}}

\newcommand{\wncurvar}{\@ifnextchar[{\basicwncurvar}{\basicwncurvar[160]}}%

\def\basicWncurvar[#1]#2{\Wncurve{#2}{#1}}%

\newcommand{\Wncurvar}{\@ifnextchar[{\basicWncurvar}{\basicWncurvar[160]}}%

\newcommand{\crosslength}[2]{%
\settowidth{\firstitem}{#1}%
\settowidth{\seconditem}{#2}%
\ifdim\firstitem < \seconditem%
\itemlength=\seconditem%
\else%
\itemlength=\firstitem%
\fi%
\divide\itemlength by 2%
\hspace{\itemlength}}%

\newcommand{\crossarrows}[2]{%
\crosslength{\mbox{$#1$}}{\mbox{$#2$}}%
\begin{picture}(0,0)%
\put(0,0){\makebox(0,0){$#1$}}%
\thinlines%
\put(0,0){\makebox(0,0){$#2$}}%
\thinlines%
\end{picture}%
\crosslength{\mbox{$#1$}}{\mbox{$#2$}}}%

\catcode`\@=12
